\documentclass[a4paper,12pt,amsart,frenchb]{article}

\usepackage{amsmath,amsbsy,amsfonts,amssymb}
\usepackage[french]{babel}
\newcommand{\ass}[2]{\vskip0.3cm\noindent
{\bf {#1}}. { \sl {#2}}\vskip0.3cm\noindent
}

\begin{document}

  \title{ Propri\'et\'es de maximalit\'e concernant une repr\'esentation d\'efinie par Lusztig}
\author{J.-L. Waldspurger}
\date{29 ao\^ut 2017}
 \maketitle
 
 {\bf Introduction}
 
 Soit $N\geq1$ un entier. On note ${\cal P}(2N)$ l'ensemble des partitions  de $2N$, c'est-\`a-dire celui des suites d'entiers $\lambda=(\lambda_{1}\geq \lambda_{2}\geq...\geq\lambda_{r}\geq 0)$  telles que $S(\lambda)=2N$, o\`u  $S(\lambda)=\lambda_{1}+...+\lambda_{r}$. Deux suites sont  identifi\'ees si elles ne diff\`erent que par des termes nuls. Pour une telle partition et pour un entier $i\geq1$, on note $mult_{\lambda}(i)$ le nombre d'indices $j$ tels que $\lambda_{j}=i$. Notons ${\cal P}^{symp}(2N)$ le sous-ensemble des partitions symplectiques, c'est-\`a-dire celles telles que pour tout entier $i\geq1$ impair, $mult_{\lambda}(i)$ est pair. Pour $\lambda\in {\cal P}^{symp}(2N)$, notons $Jord^{bp}(\lambda)$ l'ensemble des entiers $i\geq2$ pairs tels que $mult_{\lambda}(i)\not=0$. On note $\boldsymbol{{\cal P}^{symp}}(2N)$ l'ensemble des couples $(\lambda,\epsilon)$ o\`u $\lambda\in {\cal P}^{symp}(2N)$ et $\epsilon\in \{\pm 1\}^{Jord^{bp}(\lambda)}$. 
 
 Pour tout entier $n\in {\mathbb N}$, notons $W_{n}$ le groupe de Weyl d'un syst\`eme de racines de type $B_{n}$ ou $C_{n}$, avec la convention $W_{0}=\{1\}$. Notons ${\cal W}_{N}$ l'ensemble des couples $(k,\rho)$ o\`u $k$ est un entier tel que $k\geq0$ et $k(k+1)\leq 2N$ et $\rho$ est une repr\'esentation irr\'eductible de $W_{N-k(k+1)/2}$. 
  La correspondance de Springer g\'en\'eralis\'ee, d\'efinie par Lusztig, est une bijection $(\lambda,\epsilon)\mapsto (k(\lambda,\epsilon),\rho(\lambda,\epsilon))$ de $\boldsymbol{{\cal P}^{symp}}(2N)$ sur ${\cal W}_{N}$.

 Il est bien connu que ${\cal P}^{symp}(2N)$ param\`etre les orbites unipotentes dans le groupe symplectique $Sp(2N;{\mathbb C})$. L'ensemble $\boldsymbol{{\cal P}^{symp}}(2N)$ param\`etre les couples $(C,{\cal E})$ form\'es d'une orbite unipotente $C$ et d'un syst\`eme local $Sp(2N;{\mathbb C})$-\'equivariant irr\'eductible sur $C$. On note $(C_{\lambda,\epsilon},{\cal E}_{\lambda,\epsilon})$ le couple param\'etr\'e par $(\lambda,\epsilon)\in \boldsymbol{{\cal P}^{symp}}(2N)$. Pour un tel couple $(C',{\cal E}')$, on d\'efinit le prolongement d'intersection $IC(\bar{C}',{\cal E}')$, qui est support\'e par la fermeture $\bar{C}'$ de $C'$, et ses faisceaux de cohomologie ${\cal H}^mIC(\bar{C}',{\cal E}')$. Lusztig a montr\'e que ces derniers sont nuls si $m$ est impair. Si $C$ est une autre orbite unipotente, la restriction de ${\cal H}^{2m}IC(\bar{C}',{\cal E}')$ \`a $C$ est somme directe de syst\`emes locaux $Sp(2N;{\mathbb C})$-\'equivariants irr\'eductibles. Si ${\cal E}$ est un tel syst\`eme local sur $C$, on note $\mathfrak{mult}(C,{\cal E};C',{\cal E}')$ la multiplicit\'e avec laquelle intervient ${\cal E}$ dans la restriction \`a $C$ de $\oplus_{m\in {\mathbb Z}}{\cal H}^{2m}IC(\bar{C}',{\cal E}')$.  Soient $(\lambda,\epsilon)$ et $(\lambda',\epsilon')$ les \'el\'ements de $\boldsymbol{{\cal P}^{symp}}(2N)$ param\'etrant respectivement $(C,{\cal E})$ et $(C',{\cal E}')$. On pose $\mathfrak{mult}(\lambda,\epsilon;\lambda',\epsilon')=  \mathfrak{mult}(C,{\cal E};C',{\cal E}')$. Supposons que cette multiplicit\'e soit non nulle. Il est clair qu'alors $\lambda\leq \lambda'$ pour l'ordre usuel des partitions. Lusztig a prouv\'e que $k(\lambda,\epsilon)=k(\lambda',\epsilon')$. 

Partons maintenant de $(\lambda,\epsilon)\in \boldsymbol{{\cal P}^{symp}}(2N)$. Posons pour simplifier $k=k(\lambda,\epsilon)$. On d\'efinit la repr\'esentation $\underline{\rho}(\lambda,\epsilon)$ de $W_{N-k(k+1)/2}$ par
$$\underline{\rho}(\lambda,\epsilon)=\oplus_{(\lambda',\epsilon')\in \boldsymbol{{\cal P}^{symp}}(2N)}\mathfrak{mult}( \lambda,\epsilon ;\lambda',\epsilon')\rho(\lambda',\epsilon').$$
Les coefficients $\mathfrak{mult}( \lambda,\epsilon ;\lambda',\epsilon')$ interviennent de fa\c{c}on essentielle dans de nombreux travaux concernant les groupes r\'eductifs  finis, la premi\`ere occurence \'etant peut-\^etre \cite{L1} theorem 24.4. Ce qui nous concerne plus directement est que la repr\'esentation $\underline{\rho}(\lambda,\epsilon)$ contr\^ole les restrictions aux sous-groupes parahoriques des repr\'esentations de r\'eduction unipotente des groupes $p$-adiques (voir \cite{L2} et \cite{W} proposition 5.2). Or on sait peu de choses sur cette repr\'esentation. Il nous semble qu'en g\'en\'eral, outre ce que l'on a d\'ej\`a dit ci-dessus ($\mathfrak{mult}(\lambda,\epsilon;\lambda',\epsilon')\not=0$ entra\^{\i}ne $k(\lambda',\epsilon')=k(\lambda,\epsilon)$), on ne conna\^{\i}t que la propri\'et\'e de minimalit\'e

$\mathfrak{mult}(\lambda,\epsilon;\lambda',\epsilon')\not=0$ entra\^{\i}ne $\lambda'>\lambda$ ou $(\lambda',\epsilon')=(\lambda,\epsilon)$; dans ce dernier cas

\noindent $\mathfrak{mult}(\lambda,\epsilon;\lambda,\epsilon)=1$. 

Le but de l'article est de d\'emontrer que, sous une hypoth\`ese de parit\'e sur $\lambda$, la repr\'esentation $\underline{\rho}(\lambda,\epsilon)$ poss\`ede aussi une propri\'et\'e de maximalit\'e. Pr\'ecis\'ement, nous d\'emontrerons le th\'eor\`eme suivant. 

 \ass{Th\'eor\`eme}{Soit $(\lambda,\epsilon)\in \boldsymbol{{\cal P}^{symp}}(2N)$. Supposons que $\lambda$ n'a que des termes pairs. Alors il existe un unique \'el\'ement $(\lambda^{max},\epsilon^{max})\in \boldsymbol{{\cal P}^{symp}}(2N)$ v\'erifiant les propri\'et\'es suivantes:
 
 (i) $\mathfrak{mult}(\lambda,\epsilon;\lambda^{max},\epsilon^{max})=1$;
 
 (ii) pour tout \'el\'ement $(\lambda',\epsilon')\in \boldsymbol{{\cal P}^{symp}}(2N)$ tel que $\mathfrak{mult}(\lambda,\epsilon;\lambda',\epsilon')\not=0$, on a $\lambda'<\lambda^{max}$ ou $(\lambda',\epsilon')=(\lambda^{max},\epsilon^{max})$. }
 
 En fait, notre m\'ethode permet une certaine latitude sur l'ordre choisi pour comparer les partitions, ou plus exactement les \'el\'ements de $\boldsymbol{{\cal P}^{symp}}(2N)$.    Pour tout entier $n\in {\mathbb N}$, notons simplement $sgn$ le caract\`ere signature de $W_{n}$. Pour $(\lambda,\epsilon)\in \boldsymbol{{\cal P}^{symp}}(2N)$, notons $(^s\lambda,^s\epsilon)$ l'\'el\'ement de $\boldsymbol{{\cal P}^{symp}}(2N)$ tel que $k(^s\lambda,^s\epsilon)=k(\lambda,\epsilon)$ et $\rho(^s\lambda,^s\epsilon)=sgn\otimes \rho(\lambda,\epsilon)$ (la notation ne doit pas pr\^eter \`a confusion: $^s\lambda$, resp. $^s\epsilon$, ne d\'epend pas seulement de $\lambda$, resp. $\epsilon$, mais du couple $(\lambda,\epsilon)$). En g\'en\'eral, cette op\'eration n'est pas d\'ecroissante pour  l'ordre usuel des partitions, c'est-\`a-dire que, pour $(\lambda,\epsilon),(\lambda',\epsilon')\in \boldsymbol{{\cal P}^{symp}}(2N)$, la relation $\lambda\leq \lambda'$ n'entra\^{\i}ne pas ${^s\lambda}'\leq {^s\lambda}$. Toutefois, nous d\'emontrerons le th\'eor\`eme suivant.
 
  \ass{Th\'eor\`eme}{Soit $(\lambda,\epsilon)\in \boldsymbol{{\cal P}^{symp}}(2N)$. Supposons que $\lambda$ n'a que des termes pairs. Alors il existe un unique \'el\'ement $(\lambda^{min},\epsilon^{min})\in \boldsymbol{{\cal P}^{symp}}(2N)$ v\'erifiant les propri\'et\'es suivantes:
 
 (i) $\mathfrak{mult}(\lambda,\epsilon;{^s\lambda}^{min},{^s\epsilon}^{min})=1$;
 
 (ii) pour tout \'el\'ement $(\lambda',\epsilon')\in \boldsymbol{{\cal P}^{symp}}(2N)$ tel que $\mathfrak{mult}(\lambda,\epsilon;{^s\lambda'},{^s\epsilon'})\not=0$, on a $\lambda^{min}<\lambda'$ ou $(\lambda',\epsilon')=(\lambda^{min},\epsilon^{min})$. 
 
 De plus, on a $({^s\lambda}^{min},{^s\epsilon}^{min})=(\lambda^{max},\epsilon^{max})$.}
 
 Dans un article ult\'erieur, nous utiliserons  ce r\'esultat  pour prouver  que certaines repr\'esentations (tr\`es particuli\`eres) de groupes $p$-adiques admettent un front d'onde.
 
 Le couple $(\lambda^{max},\epsilon^{max})$ peut se calculer explicitement. Par contre, nous n'avons pas trouv\'e de formule "simple" pour calculer le couple $(\lambda^{min},\epsilon^{min})$ (il peut \'evidemment se calculer \`a partir de $(\lambda^{max},\epsilon^{max})$ et des formules explicites pour la repr\'esentation de Springer g\'en\'eralis\'ee, mais ce calcul n'est pas "simple").
 
 Les preuves des th\'eor\`emes sont purement combinatoires.  Comme on le sait, les repr\'esentations irr\'eductibles d'un groupe $W_{n}$ sont param\'etr\'es par les couples $(\alpha,\beta)$ de partitions telles que $S(\alpha)+S(\beta)=n$. Une multiplicit\'e $\mathfrak{mult}(\lambda,\epsilon;\lambda',\epsilon')$ peut se r\'ecrire $\mathfrak{mult}(\alpha,\beta;\alpha',\beta')$ o\`u $(\alpha,\beta)$ param\`etre $\rho(\lambda,\epsilon)$ et $(\alpha',\beta')$ param\`etre $\rho(\lambda',\epsilon')$.  Shoji a d\'evelopp\'e une th\'eorie des fonctions sym\'etriques associ\'ees \`a des couples de partitions. En particulier, \`a l'aide d'un couple $(\alpha,\beta)$ et d'un ordre  $<_{I}$ sur l'ensemble d'indices de ce couple, il d\'efinit un 
   polyn\^ome   du type polyn\^ome de Kostka. Pour un autre couple $(\alpha',\beta')$, on peut d\'efinir  la "multiplicit\'e" $mult(<_{I};\alpha,\beta;\alpha',\beta')$ de $(\alpha',\beta')$ dans ce polyn\^ome. Sous certaines hypoth\`eses sur $(\lambda,\epsilon)$ (le couple associ\'e \`a $(\alpha,\beta)$) et pour un ordre $<_{I}$ bien choisi, Shoji prouve que $mult(<_{I};\alpha,\beta;\alpha',\beta')=\mathfrak{mult}(\alpha,\beta;\alpha',\beta')$ pour tout $(\alpha',\beta')$. Nous g\'en\'eraliserons un peu ce r\'esultat en 5.2: il est valable pourvu que $\lambda$ n'ait que des termes pairs. Les assertions des th\'eor\`emes ci-dessus deviennent des assertions concernant les paires $(\alpha',\beta')$ telles que $mult(<_{I};\alpha,\beta;\alpha',\beta')$. Et la d\'efinition de ces multiplicit\'es est purement combinatoire. 
   
   Dans le premier paragraphe, on consid\`ere un couple de partitions $(\alpha,\beta)$ et un ordre $<_{I}$ sur l'ensemble d'indices de ce couple. On d\'efinit des ensembles $P(\alpha,\beta)$ et $P_{A,B;s}(\alpha,\beta)$ form\'es de couples de partitions $(\nu,\mu)$ d\'eduits de $(\alpha,\beta)$ par une combinatoire simple. On montre au paragraphe 4 que ces couples sont pr\'ecis\'ement ceux qui interviennent dans le polyn\^ome  de Shoji associ\'e \`a $(\alpha,\beta)$ avec des propri\'et\'es de maximalit\'e du type de celles intervenant dans les th\'eor\`emes ci-dessus. Entre-temps, les paragraphes 2 et 3 sont consacr\'es \`a prouver divers lemmes techniques concernant ces ensembles $P(\alpha,\beta)$ et $P_{A,B;s}(\alpha,\beta)$. Ils peuvent \^etre n\'eglig\'es en premi\`ere lecture et m\^eme en deuxi\`eme. Au paragraphe 5, on \'etend comme on l'a dit le th\'eor\`eme  de Shoji et on en d\'eduit les th\'eor\`emes, \`a l'aide  des r\'esultats de maximalit\'e prouv\'es au paragraphe 4. La description explicite du couple $(\lambda^{max},\epsilon^{max})$ est donn\'ee au paragraphe 6.

\section{D\'efinitions}

\bigskip

\subsection{Ensembles finis \`a multiplicit\'es}

Soit $n\in {\mathbb N}$. Consid\'erons l'ensemble des suites $\lambda=(\lambda_{1},...,\lambda_{n})$ o\`u les $\lambda_{j}$ sont des nombres r\'eels. Disons que deux telles suites sont \'equivalentes si et seulement si elles diff\`erent par une permutation des indices. Nous appellerons une classe d'\'equivalence un ensemble \`a multiplicit\'es  de nombres r\'eels, \`a $n$ \'el\'ements. Pour une telle suite $\lambda$ et pour $i\in {\mathbb R}$, on note $mult_{\lambda}(i)$ le nombre d'indices $j$ tels que $\lambda_{j}=i$. Notons ${\cal R}_{n}$ l'ensemble des  suites $\lambda=(\lambda_{1},...,\lambda_{n})$ o\`u les $\lambda_{j}$ sont des nombres r\'eels tels que $\lambda_{1}\geq \lambda_{2}...\geq \lambda_{n}$. Pour tout ensemble \`a multiplicit\'es  de nombres r\'eels, \`a $n$ \'el\'ements, qui est donc une classe d'\'equivalence de suites, il existe une unique suite qui appartient \`a la fois \`a cette classe et \`a ${\cal R}_{n}$. Ainsi, l'ensemble des  ensembles \`a multiplicit\'es  de nombres r\'eels, \`a $n$ \'el\'ements, s'identifie \`a ${\cal R}_{n}$. Pour $\lambda\in {\cal R}_{n}$ et pour $k\in \{0,...,n\}$, on pose $S_{k}(\lambda)=\sum_{j=1,...,k}\lambda_{j}$. On pose simplement $S(\lambda)=S_{n}(\lambda)$. On munit ${\cal R}_{n}$ de l'ordre usuel: $\lambda\leq\lambda'$ si et seulement si $S_{k}(\lambda)\leq S_{k}(\lambda')$ pour tout $k=1,...,n$. Pour $\lambda,\lambda'\in {\cal R}_{n}$, on note $\lambda+\lambda'$ la suite $(\lambda_{1}+\lambda'_{1},\lambda_{2}+\lambda'_{2},...,\lambda_{n}+\lambda'_{n})$. Elle appartient \`a ${\cal R}_{n}$. Pour $\lambda\in {\cal R}_{n}$ et $\lambda'\in {\cal R}_{n'}$, on note $\lambda\sqcup \lambda'$ l'\'element de ${\cal R}_{n+n'}$ \'equivalent \`a $(\lambda_{1},...,\lambda_{n},\lambda'_{1},...,\lambda'_{n'})$. On v\'erifie que

 (1) soient $\lambda,\lambda'\in {\cal P}_{n}$ et $\mu,\mu'\in {\cal P}_{m}$; supposons   $\lambda\leq \lambda'$ et $\mu\leq \mu'$; alors $\lambda\sqcup \mu\leq \lambda'\sqcup \mu'$ et on a $\lambda\sqcup \mu=\lambda'\sqcup \mu'$ si et seulement si $\lambda=\lambda'$ et $\mu=\mu'$.
 
  Dans le cas particulier $n=0$, on note $\emptyset$ l'unique \'el\'ement de ${\cal R}_{0}$.

Pour $n\in {\mathbb N}$, on note ${\cal P}_{n}$ l'ensemble des partitions  $\alpha=(\alpha_{1},...,\alpha_{n})$ o\`u $\alpha_{j}\in {\mathbb N}$ pour tout $j$ et o\`u $\alpha_{1}\geq ... \geq \alpha_{n}\geq0$.  C'est un sous-ensemble de ${\cal R}_{n}$. On utilisera parfois la notation $\alpha=(0^{n_{0}},1^{n_{1}},...)$ o\`u $n_{i}=mult_{\alpha}(i)$.  

\bigskip

\subsection{L'ensemble $P(\alpha,\beta)$}
Soient $n,m\in {\mathbb N}$. Posons $I=(\{1,...,n\}\times\{0\}) \cup(\{1,...,m\}\times\{1\})$. On fixe un ordre total $<_{I}$ sur $I$ tel que, pour  $i,j\in \{1,...,n\}$, on ait $(i,0)<_{I}(j,0)\iff i<j$ et, pour $i,j\in \{1,...,m\}$, $(i,1)<_{I}(j,1)\iff i<j$. Pour  $(\alpha,\beta)\in {\cal P}_{n}\times {\cal P}_{m}$, on va d\'efinir un sous-ensemble fini $P(\alpha,\beta)\subset {\cal P}_{n}\times {\cal P}_{m}$. Cette construction se fait par r\'ecurrence sur $n+m$. 

Si $m=n=0$, on a $\alpha=\beta=\emptyset$, $(\emptyset,\emptyset)$ est l'unique \'el\'ement de ${\cal P}_{0}\times {\cal P}_{0}$ et on pose $P(\emptyset,\emptyset)=\{(\emptyset,\emptyset)\}$. 

Supposons $n+m\geq1$. On consid\`ere les deux constructions suivantes.

(a) supposons $n\geq1$. On  pose $a_{1}=1$. On note $b_{1}$ le plus petit \'el\'ement de $\{1,...,m\}$ tel que $(1,0)<_{I}(b_{1},1)$, s'il en existe, puis $a_{2}$ le plus petit \'el\'ement de $\{1,...,n\}$ tel que $(b_{1},1)<_{I}(a_{2},0)$ s'il en existe,  puis $b_{2}$ le plus petit \'el\'ement de $\{1,...,m\}$ tel que $(a_{2},0)<_{I} (b_{2},1)$ s'il en existe etc... On arr\^ete le proc\'ed\'e  quand il n'y a plus d'\'el\'ement v\'erifiant l'in\'egalit\'e requise. On pose $\nu_{1}=\alpha_{a_{1}}+\beta_{b_{1}}+\alpha_{a_{2}}+\beta_{b_{2}}+...$.
On note $i_{1}<...<i_{n'}$ les \'el\'ements de $\{1,...,n\}$ qui ne font pas partie du sous-ensemble $\{a_{1},a_{2},...\}$. On note $j_{1}<... <j_{m'}$ les \'el\'ements de $\{1,...,m\}$ qui ne font pas partie du sous-ensemble $\{b_{1},b_{2},...\}$. On identifie l'ensemble $I'=(\{1,...,n'\}\times\{0\})\cup( \{1,...,m'\}\times\{0\})$ \`a un sous-ensemble de $I$ en identifiant $(k,0)\in I'$ \`a $(i_{k},0)\in I$ et $(k,1)\in I'$ \`a $(j_{k},1)\in I$. On munit $I'$ de l'ordre $<_{I'}$ induit par ce plongement et par  l'ordre $<_{I}
$ sur $I$. On d\'efinit $\alpha'=(\alpha_{i_{1}},\alpha_{i_{2}},...)\in {\cal P}_{n'}$ et $\beta'=(\beta_{j_{1}},\beta_{j_{2}},...)\in {\cal P}_{m'}$. Par r\'ecurrence, l'ensemble $P(\alpha',\beta')$ est d\'efini. On note $P^{a}(\alpha,\beta)$ l'ensemble des couples $(\nu=\{\nu_{1}\}\sqcup \nu'\sqcup\{0^{n-n'-1}\},\mu=\mu'\sqcup\{0^{m-m'}\})$ pour $(\nu',\mu')\in P(\alpha',\beta')$.

(b) Supposons $m\geq1$. On construit un ensemble $P^{b}(\alpha,\beta)$ en inversant les r\^oles de $\alpha$ et $\beta$: on pose $b_{1}=1$. On note $a_{1}$ le plus petit \'el\'ement de $\{1,...,n\}$ tel que $(1,1)<(a_{1},0)$ etc... On pose $\mu_{1}=\beta_{b_{1}}+\alpha_{a_{1}}+\beta_{b_{2}}...$, on d\'efinit les partitions $(\alpha',\beta')$. On note  $P^{b}(\alpha,\beta)$ l'ensemble des couples $(\nu=\nu'\sqcup \{0^{n-n'}\},\mu=\{\mu_{1}\}\sqcup \mu'\sqcup\{0^{m-m'-1}\})$ pour $(\nu',\mu')\in P(\alpha',\beta')$.

On a suppos\'e $n+m\geq1$. Si $m=0$, on pose $P(\alpha,\beta)=P^{a}(\alpha,\beta)$. Si $n=0$, on pose $P(\alpha,\beta)=P^{b}(\alpha,\beta)$. Si $n,m\geq1$, on note  $P(\alpha,\beta)$ la r\'eunion de $P^{a}(\alpha,\beta)$ et $P^{b}(\alpha,\beta)$.  

On a

(1) si $n=0$ ou $m=0$, $P(\alpha,\beta)=\{(\alpha,\beta)\}$; 

(2) $S(\nu)+S(\mu)=S(\alpha)+S(\beta)$ pour tout $(\nu,\mu)\in P(\alpha,\beta)$.

Consid\'erons des ensembles d'indices $1\leq i_{1}<i_{2}<...<i_{n'}\leq n$, $1\leq j_{1}<j_{2}<...<j_{m'}\leq m$, posons $\alpha'=(\alpha_{i_{1}},\alpha_{i_{2}},...)\in {\cal P}_{n'}$ et $\beta'=(\beta_{j_{1}},\beta_{j_{2}},...)\in {\cal P}_{m'}$. Comme dans la construction ci-dessus, identifions l'ensemble d'indices $I'=(\{1,...,n'\}\times\{0\})\cup( \{1,...,m'\}\times\{0\})$ \`a un sous-ensemble de $I$. Supposons $n'\geq1$ et consid\'erons des \'el\'ements $(\nu,\mu)\in P^{a}(\alpha,\beta)$ et $(\nu',\mu')\in P(\alpha',\beta')$. Alors

(3) on a $\nu'_{1}\leq \nu_{1}$.

 Preuve de (3). Comme on l'a dit, on a $\nu_{1}=\alpha_{1}+\beta_{b_{1}}+\alpha_{a_{2}}...$.  La construction entra\^{\i}ne    que $\nu'_{1}$ est de la forme $\alpha_{a'_{1}}+\beta_{b'_{1}}+\alpha_{a'_{2}}...$ o\`u $a'_{1}\geq1$, $b'_{1}$ est un \'el\'ement de $\{1,...,m\}$ tel que $(a'_{1},0)<_{I}(b'_{1},1)$, $a'_{2}$ est un \'el\'ement de $\{1,...,n\}$ tel que $(b'_{1},1)<_{I}(a'_{2},0)$ etc...Par r\'ecurrence, on voit que si $b'_{1}$ existe, $b_{1}$ aussi et $b_{1}\leq b'_{1}$, que si $a'_{2}$ existe, $a_{2}$ aussi et $a_{2}\leq a'_{2}$ etc... L'\'el\'ement $\nu'_{1}=\alpha_{a'_{1}}+\beta_{b'_{1}}+\alpha_{a'_{2}}...$ est alors inf\'erieur ou \'egal \`a $\nu_{1}$. $\square$
 
 En particulier, dans la construction (a), $\nu_{1}$ est bien le plus grand \'el\'ement de la suite $\nu$. Il y a \'evidemment des propri\'et\'es similaires en \'echangeant les r\^oles de $\alpha$ et $\beta$: dans la construction (b), $\mu_{1}$  est bien le plus grand \'el\'ement de la suite $\mu$. 

\bigskip

Supposons $n,m\geq1$. Le proc\'ed\'e (a) cr\'ee  des suites $a_{1}=1$, $a_{2}$,..., $b_{1}$, $b_{2}$, ... un \'el\'ement $\nu_{1}=\alpha_{1}+\beta_{b_{1}}+\alpha_{a_{2}}+...$ et des partitions $(\alpha',\beta')$. Le proc\'ed\'e (b) cr\'ee des objets que nous notons en les soulignant: des suites $\underline{b}_{1}=1$, $\underline{b}_{2}$...,$\underline{a}_{1}$, $\underline{a}_{2}$,... un \'el\'ement $\underline{\mu}_{1}=\beta_{1}+\alpha_{\underline{a}_{1}}+\beta_{\underline{b}_{2}}+...$ et des partitions  $(\underline{\alpha},\underline{\beta})$. Supposons que $(1,0)<_{I}(1,1)$. On voit que $b_{1}=1$, puis, par r\'ecurrence, $\underline{a}_{1}=a_{2}$, $\underline{b}_{2}=b_{2}$, $\underline{a}_{2}=a_{3}$... D'o\`u $\nu_{1}=\alpha_{1}+\underline{\mu}_{1}$, $\beta'=\underline{\beta}'$ et $\underline{\alpha}'=\{\alpha_{1}\}\sqcup\alpha'$. Il y a un r\'esultat similaire si $(1,1)<_{I}(1,0)$. En r\'esum\'e

(4) si $(1,0)<_{I}(1,1)$, on a $\nu_{1}=\alpha_{1}+\underline{\mu}_{1}$, $\underline{\alpha}'=\{\alpha_{1}\}\sqcup\alpha'$  et $\beta'=\underline{\beta}'$; si $(1,1)<_{I}(1,0)$, on a $\underline{\mu}_{1}=\beta_{1}+\nu_{1}$, $\alpha'=\underline{\alpha}'$ et $ \beta'=\{\beta_{1}\}\sqcup\underline{\beta}'$.

Consid\'erons les cas particuliers suivants: 

(5) on suppose que $\beta=\{0^m\}$ et que $\alpha_{i}=0$ pour tout $i\in \{1,...,n\}$ tel que $(1,1)<_{I}(i,0)$; ou bien on suppose que $\alpha=\{0^n\}$ et que $\beta_{i}=0$ pour tout $i\in \{1,...,m\}$ tel que $(1,0)<_{I}(i,1)$; dans les deux cas $P(\alpha,\beta)=\{(\alpha,\beta)\}$.

Cela se prouve par r\'ecurrence sur $n+m$. 

\bigskip

{\bf Remarque.} Un \'el\'ement de $P(\alpha,\beta)$ est obtenu en appliquant  une suite de proc\'ed\'es (a) ou (b). Mais cette suite n'est pas toujours unique. Par exemple, supposons $n=1$, $m=2$, $(1,0)<_{I}(1,1)<_{I}(2,1)$ et supposons $\beta_{1}=\beta_{2}$. Alors le couple $(\{\alpha_{1}+\beta_{1}\},\{\beta_{2},0\})$ appartient \`a $P(\alpha,\beta)$. Il peut \^etre obtenu par application des proc\'ed\'es (a) puis (b) ou par application des proc\'ed\'es (b) puis (a).

 \bigskip

\subsection{L'ensemble $P_{A,B;s}(\alpha,\beta)$}
 On conserve les hypoth\`eses du paragraphe pr\'ec\'edent. Consid\'erons des r\'eels $A,B$ et un r\'eel $s>0$. 
Soit $\alpha\in {\cal P}_{n}$ et $\beta\in {\cal P}_{m}$. On d\'efinit $P_{A,B;s}(\alpha,\beta)\subset P(\alpha,\beta)$ de la fa\c{c}on suivante.  

Si $n=m=0$, ou plus g\'en\'eralement si $n=0$ ou $m=0$, $P_{A,B;s}(\alpha,\beta)=P(\alpha,\beta)=\{(\alpha,\beta)\}$.

Supposons $n\geq1$, $m\geq1$. On reprend la m\^eme construction que dans le paragraphe pr\'ec\'edent. On autorise cette fois le proc\'ed\'e (a)  si et seulement si l'une des conditions suivantes est v\'erifi\'ee

$(1,0)<_{I}(1,1)$ et $\alpha_{1}+A\geq B$ ou $(1,1)<_{I}(1,0)$ et $\beta_{1}+B\leq A$;

Dans ce cas, on d\'efinit $P^{a}_{A,B;s}(\alpha,\beta)$ comme l'ensemble des couples $(\nu=\{\nu_{1}\}\sqcup \nu'\sqcup\{0^{n-n'-1}\},\mu=\mu'\sqcup\{0^{m-m'}\})$ pour $(\nu',\mu')\in P_{A-s,B;s}(\alpha',\beta')$. 

On  autorise le proc\'ed\'e (b) si et seulement si l'une des conditions suivantes est v\'erifi\'ee

$(1,0)<_{I}(1,1)$ et $\alpha_{1}+A\leq B$ ou $(1,1)<_{I}(1,0)$ et $\beta_{1}+B\geq A$.

Dans ce cas, on d\'efinit $P^{b}_{A,B;s}(\alpha,\beta)$ comme l'ensemble des couples $(\nu=\nu'\sqcup\{0^{n-n'}\},\mu=\{\mu_{1}\}\sqcup \mu'\sqcup\{0^{m-m'-1}\})$, pour $(\nu',\mu')\in P_{A,B-s;s}(\alpha',\beta')$.

Si un seul proc\'ed\'e est autoris\'e, disons celui de type (a), on pose $P_{A,B;s}(\alpha,\beta)=P^{a}_{A,B;s}(\alpha,\beta)$. Si les deux proc\'ed\'es sont autoris\'es, on pose $P_{A,B;s}(\alpha,\beta)=P^{a}_{A,B;s}(\alpha,\beta)\cup P^{b}_{A,B;s}(\alpha,\beta)$. 
Au moins un des proc\'ed\'es est loisible donc l'ensemble $P_{A,B;s}(\alpha,\beta)$ est non vide. Il peut avoir plusieurs \'el\'ements puisque les proc\'ed\'es (a) et (b) peuvent \^etre tous deux autoris\'es. On peut toutefois d\'efinir un \'el\'ement canonique "de type (a)"  en modifiant la construction ci-dessus  de la fa\c{c}on suivante: on n'autorise le proc\'ed\'e (b) que si le proc\'ed\'e (a) n'est pas autoris\'e. Dans l'un ou l'autre des proc\'ed\'es, on choisit pour $(\nu',\mu')$ l'\'el\'ement canonique "de type (a)" relatif \`a $(\alpha',\beta')$. On d\'efinit de la m\^eme fa\c{c}on  un \'el\'ement canonique "de type (b)".

Soient $N,M\in {\mathbb N}$. Supposons $N\geq n$, $M\geq m$, $A\geq s(N-1)$, $B\geq s(M-1)$. Appelons $(N,M;s)$-symbole 
 un couple $\Lambda=(\Lambda^{a},\Lambda^{b})\in {\cal R}_{N}\times {\cal R}_{M}$ tel que 
 
 $\Lambda^{a}_{j}\geq0$ et $\Lambda^b_{j}\geq0$ pour tout $j$;

 $\Lambda^{a}_{j}\geq \Lambda^{a}_{j+1}+s$ pour $j=1,...,N-1$, $\Lambda^{b}_{j}\geq \Lambda^{b}_{j+1}+s$ pour $i=1,...M-1$.   
 
 On d\'efinit l'\'el\'ement $p\Lambda\in {\cal R}_{N+M}$ par  $p\Lambda=\Lambda^{a}\sqcup \Lambda^{b}$. 
 
 Pour $R,R'\in {\mathbb R}$ tels que $R\in R'+s{\mathbb N}$, on pose $[R,R']_{s}=\{R,R-s,R-2s,...,R'\}$. Par convention, on pose aussi $[R,R+s]_{s}=\emptyset$. 
  Pour $\nu\in {\cal P}_{n}$, $\mu\in {\cal P}_{m}$, on d\'efinit le $(N,M;s)$-symbole
$$\Lambda^{N,M}_{A,B;s}(\nu,\mu)=( (\nu\sqcup\{0^{N-n}\})+[A,A+s-sN]_{s},(\mu\sqcup\{0^{M-m}\})+[B,B+s-sM]_{s}).$$

Soient $\alpha\in {\cal P}_{n}$ et $\beta\in {\cal P}_{m}$, consid\'erons un \'el\'ement $(\nu,\mu)\in P_{A,B;s}(\alpha,\beta)$. D'apr\`es la d\'efinition de ce terme, on peut trouver un entier $r\geq0$ , des suites de partitions $(\alpha^{i},\beta^{i})_{i=1,...,r+1}$, $(\nu^{i},\mu^{i})_{i=1,...,r+1}$,  des suites de r\'eels $(A_{i})_{i=1,...,r+1}$, $(B_{i})_{i=1,...,r+1}$ et d'entiers $(N_{i})_{i=1,...,r+1}$, $(M_{i})_{i=1,...,r+1}$, v\'erifiant les propri\'et\'es suivantes:

$(\alpha^1,\beta^1)=(\alpha,\beta)$, $(\nu^1,\mu^1)=(\nu,\mu)$, $A_{1}=A$, $B_{1}=B$;

en notant $n_{i}$ et $m_{i}$ les entiers tels que $\alpha^{i}\in {\cal P}_{n_{i}}$ et $\beta^{i}\in {\cal P}_{m_{i}}$, on a $n_{i}+m_{i}\geq1$ si $i\leq r$ et $n_{r+1}=m_{r+1}=0$;

pour $i=1,...,r+1$, $(\nu^{i},\mu^{i})$ appartient \`a $P_{A_{i},B_{i};s}(\alpha^{i},\beta^{i})$;

pour $i=1,...,r$, si $(\nu^{i},\mu^{i})$ est issu du proc\'ed\'e (a), ce proc\'ed\'e cr\'ee l'\'el\'ement $\nu^{i}_{1}$ et les partitions $(\alpha^{i+1},\beta^{i+1})$; on pose $A_{i+1}=A_{i}-s $, $B_{i+1}=B_{i}$, $N_{i+1}=N_{i}-1$, $M_{i+1}=M_{i}$; on a $(\nu^{i},\mu^{i})=(\{\nu^{i}_{1}\}\sqcup \nu^{i+1}\sqcup\{0^{n_{i}-1-n_{i+1}}\},\mu^{i+1}\sqcup\{0^{m_{i}-m_{i+1}}\})$; si $(\nu^{i},\mu^{i})$ est issu du proc\'ed\'e (b), ce proc\'ed\'e cr\'ee l'\'el\'ement $\mu^{i}_{1}$ et les partitions $(\alpha^{i+1},\beta^{i+1})$; on pose $A_{i+1}=A_{i} $, $B_{i+1}=B_{i}-s$, $N_{i+1}=N_{i}$ et $M_{i+1}=M_{i}-1$; on a $(\nu^{i},\mu^{i})=( \nu^{i+1}\sqcup\{0^{n_{i}-n_{i+1}}\},\{\mu^{i}_{1}\}\sqcup\mu^{i+1}\sqcup\{0^{m_{i}-1-m_{i+1}}\})$.

  Comme on l'a remarqu\'e dans le paragraphe pr\'ec\'edent, de telles suites ne sont pas forc\'ement uniques mais fixons-les. Pour $i=1,...,r$, on a
  $$p\Lambda^{N_{i},M_{i}}_{A_{i},B_{i};s}(\nu^{i},\mu^{i})=\{x_{i}\}\sqcup p\Lambda^{N_{i+1},M_{i+1}}_{A_{i+1},B_{i+1};s}(\nu^{i+1},\mu^{i+1}),$$
  o\`u $x_{i}=\nu^{i}_{1}+A_{i}$ si $(\nu^{i},\mu^{i})$ est issu du proc\'ed\'e (a) et $x_{i}=\mu^{i}+B_{i}$ si $(\nu^{i},\mu^{i})$ est issu du proc\'ed\'e (b). On a aussi
   $$p\Lambda^{N_{r+1},M_{r+1}}_{A_{r+1},B_{r+1};s}(\nu^{r+1},\mu^{r+1})=[A_{r+1},A+s-sN]_{s}\sqcup [B_{r+1},B+s-sM]_{s}.$$
   Notons $h$ le plus petit entier $i\in \{1,...,r\}$ tel que $inf(n_{i+1},m_{i+1})=0$. Pour $i=1,...,h$, posons $\lambda_{i}=x_{i}$, o\`u $x_{i}$ est d\'efini comme ci-dessus. Posons
   $$p\Lambda^{N_{h+1},M_{h+1}}_{A_{h+1},B_{h+1};s}(\nu^{h+1},\mu^{h+1})=(\lambda_{h+1},...,\lambda_{N+M}).$$
 Montrons que
 
 (1) on a l'\'egalit\'e
 $$p\Lambda^{N,M}_{A,B;s}(\nu,\mu)= (\lambda_{1},...,\lambda_{N+M}).$$
 
 {\bf Remarque.} Ce qui pr\'ec\`ede montre que les termes  du membre  de gauche sont bien $\lambda_{1},...,\lambda_{N+M}$. L'assertion est que ces termes sont en ordre d\'ecroissant.
 
 \bigskip 
   Par d\'efinition, les termes sont en ordre d\'ecroissant au-del\`a du rang $h+1$. Par r\'ecurrence, il suffit donc de prouver que, si $n,m\geq1$, on a $\lambda_{1}\geq \lambda_{2}$. Par sym\'etrie, on ne perd rien \`a supposer que $(\nu,\mu)$ est issu du proc\'ed\'e (a). Donc $\lambda_{1}=\nu_{1}+A$. On a alors l'une des possibilit\'es suivantes
   
   (2) $r=1$  et $\lambda_{2}= A_{2}=A-s$;
   
   (3) $r=1$  et $\lambda_{2}=B_{2}=B$;

   (4) $r\geq2$, $(\nu^2,\mu^2)$ est construit \`a l'aide du proc\'ed\'e (a) et $\lambda_{2}=\nu^2_{1}+A_{2}=\nu^2_{1}+A-s$;
   
   (5) $r\geq2$, $(\nu^2,\mu^2)$ est construit \`a l'aide du proc\'ed\'e (b) et $\lambda_{2}=\mu^2_{1}+B_{2}=\mu^2_{1}+B$.
   
 Dans le cas (2), l'in\'egalit\'e $\lambda_{1}> \lambda_{2}$ est imm\'ediate. Dans le cas (4), elle r\'esulte de  1.2(3). Notons $\underline{\mu}_{1}$ le terme issu du proc\'ed\'e (b) appliqu\'e \`a $(\alpha,\beta)$. Dans le cas (5), on a $\mu^2_{1}\leq \underline{\mu}_{1}$ d'apr\`es l'assertion sym\'etrique de 1.2(3). Dans les cas (3) et (5), on a donc $\lambda_{2}\leq \underline{\mu}_{1}+B$. Supposons $(1,0)<_{I}(1,1)$. Alors, d'apr\`es 1.2(4), on a $\nu_{1}=\alpha_{1}+\underline{\mu}_{1}$. Puisque le proc\'ed\'e (a) est loisible pour construire $(\nu,\mu)$, on a $\alpha_{1}+A\geq B$. Alors $$\lambda_{2}\leq \underline{\mu}_{1}+B\leq \alpha_{1}+\underline{\mu}_{1}+A=\nu_{1}+A=\lambda_{1}.$$
 Supposons maintenant $(1,1)<_{I}(1,0)$. Alors, d'apr\`es 1.2(4), on a $\nu_{1}=\underline{\mu}_{1}-\beta_{1}$.   Puisque le proc\'ed\'e (a) est loisible pour construire $(\nu,\mu)$, on a  $\beta_{1}+B\leq A$. Alors
   $$\lambda_{2}\leq \underline{\mu}_{1}+B\leq  \underline{\mu}_{1}-\beta_{1}+A=\nu_{1}+A=\lambda_{1}.$$
   Cela prouve (1). 
   
  Supposons que $(\nu,\mu)$ soit l'\'el\'ement canonique de type (b) de $P_{A,B;s}(\alpha,\beta)$. On peut alors choisir de fa\c{c}on canonique les suites $(\alpha^{i},\beta^{i})_{i=1,...,r+1}$ etc..: pour $i=1,...,r$, si $(\nu^{i},\mu^{i})$ peut \^etre construit \`a l'aide de chacun des proc\'ed\'es (a) ou (b), on choisit le proc\'ed\'e (b).    On a alors les in\'egalit\'es plus pr\'ecises suivantes:
  
  (6) soit $i\in \{1,...,h\}$; supposons que $(\nu^{i},\mu^{i})$ soit construit \`a l'aide du proc\'ed\'e (a); alors $\lambda_{i}> \lambda_{i+1}$.
  
  Comme dans la preuve ci-dessus, on peut supposer $i=1$. Par d\'efinition des suites canoniques, puisque $(\nu,\mu)$ est construit \`a l'aide du proc\'ed\'e (a), c'est que le proc\'ed\'e (b) est exclu. On a donc les in\'egalit\'es strictes $\alpha_{1}+A>B$ si $(1,0)<_{I}(1,1)$ et $\beta_{1}+B<A$ si $(1,1)<_{I}(1,0)$. En glissant ces in\'egalit\'es strictes dans la preuve ci-dessus, on obtient (6). 
  
  Dans la suite de l'article, on consid\`ere des donn\'ees $n,m,N,M\in {\mathbb N}$, $A,B\in {\mathbb R}$ et un r\'eel $s>0$ v\'erifiant les hypoth\`eses ci-dessus.  C'est-\`a-dire que l'on 
  suppose fix\'e un ordre $<_{I}$ sur l'ensemble d'indices $I=(\{1,...,n\}\times\{0\})\cup(\{1,...,m\}\times\{1\})$ comme en 1.2; on suppose  
 $N\geq n$, $M\geq m$, 
$A\geq s(N-1)$ et $B\geq s(M-1)$.    
  Le nombre $s$ ne changera pas. Par contre, on raisonnera souvent par r\'ecurrence sur $n$ et $m$ et on s'autorisera \`a faire varier $N$, $M$, $A$ et $B$.
Les r\'eels $A$ et $B$ sont  presque toujours positifs ou nuls. Par exemple, $A$ ne peut \^etre strictement n\'egatif que si $N=0$ donc $n=0$ mais, dans ce cas, $A$ n'interviendra pas r\'eellement dans les constructions. 
     
  \bigskip
  
  \subsection{D\'efinition de $p_{A,B;s}^{N,M}[\alpha,\beta]$}

\ass{Lemme}{Soient $\alpha\in {\cal P}_{n}$ et $\beta\in {\cal P}_{m}$. 

(i) Pour $(\nu,\mu)\in P(\alpha,\beta)$ et $(\boldsymbol{\nu},\boldsymbol{\mu})\in P_{A,B;s}(\alpha,\beta)$, on a $p\Lambda^{N,M}_{A,B;s}(\nu,\mu)\leq p\Lambda_{A,B;s}^{N,M}(\boldsymbol{\nu},\boldsymbol{\mu})$.

 (ii) Pour $(\nu,\mu)\in P(\alpha,\beta)$ et $(\boldsymbol{\nu},\boldsymbol{\mu})\in P_{A,B;s}(\alpha,\beta)$, on a $p\Lambda^{N,M}_{A,B;s}(\nu,\mu)= p\Lambda^{N,M}_{A,B;s}(\boldsymbol{\nu},\boldsymbol{\mu})$ si et seulement si $(\nu,\mu)\in P_{A,B;s}(\alpha,\beta)$.}

Preuve.  L'assertion "si" de (ii) r\'esulte trivialement de (i). On va donc d\'emontrer cette assertion (i) ainsi que l'assertion "seulement si" de (ii).

On raisonne par r\'ecurrence sur $n+m$. Si $n+m=0$, plus g\'en\'eralement si $n=0$ ou $m=0$, $P_{A,B}(\alpha,\beta)=P(\alpha,\beta)=\{(\alpha,\beta)\}$   et les assertions sont  triviales. Supposons $n\geq1$ et $m\geq1$. Les \'el\'ements $(\nu,\mu) $ et $(\boldsymbol{\nu},\boldsymbol{\mu}) $  sont construits selon des proc\'ed\'es (a) ou (b). Supposons d'abord  qu'ils soient construits par le m\^eme proc\'ed\'e, disons par le proc\'ed\'e (a).  Celui-ci construit un terme $\nu_{1}$ et des partitions $\alpha'\in {\cal P}_{n'}$, $\beta'\in {\cal P}_{m'}$. On a $\nu=\{\nu_{1}\}\sqcup \nu'\sqcup \{0\}^{n-n'-1}$  et $\mu=\mu'\sqcup \{0\}^{m-m'}$  pour un \'el\'ement $(\nu',\mu')\in P(\alpha',\beta')$. On en d\'eduit facilement l'\'egalit\'e
$$p\Lambda_{A,B;s}^{N,M}(\nu,\mu)=\{\nu_{1}+A\}\sqcup p\Lambda_{A-s,B;s}^{N-1,M}(\nu',\mu').$$
De m\^eme
$$p\Lambda_{A,B;s}^{N,M}(\boldsymbol{\nu},\boldsymbol{\mu})=\{\nu_{1}+A\}\sqcup p\Lambda_{A-s,B;s}^{N-1,M}(\boldsymbol{\nu}',\boldsymbol{\mu}'),$$
o\`u $(\boldsymbol{\nu}',\boldsymbol{\mu}')$ est un \'el\'ement de 
 $ P_{A-s,B;s}(\alpha',\beta')$. Alors l'in\'egalit\'e $p\Lambda_{A,B;s}^{N,M}(\nu,\mu)\leq p\Lambda^{N,M}_{A,B;s}(\boldsymbol{\nu},\boldsymbol{\mu})$ \'equivaut \`a l'in\'egalit\'e $p\Lambda_{A-s,B;s}^{N-1,M}(\nu',\mu')\leq p\Lambda_{A-s,B;s}^{N-1,M}(\boldsymbol{\nu}',\boldsymbol{\mu}')$ qui est connue par r\'ecurrence. La premi\`ere in\'egalit\'e est une \'egalit\'e si et seulement s'il en est de m\^eme pour la seconde. Pour celle-ci, l'\'egalit\'e entra\^{\i}ne par r\'ecurrence que
  $(\nu',\mu')\in P_{A-s,B;s}(\alpha',\beta')$. Mais alors, par construction, on a $(\nu,\mu)\in P_{A,B;s}(\alpha,\beta)$. 

Supposons maintenant  que $(\nu,\mu) $ et $(\boldsymbol{\nu},\boldsymbol{\mu}) $ soient construits selon des proc\'ed\'es diff\'erents. On ne perd rien \`a supposer $(1,0)< _{I}(1,1)$. Supposons d'abord  que  $(\boldsymbol{\nu},\boldsymbol{\mu}) $  soit construit \`a l'aide du proc\'ed\'e (a) et que $(\nu,\mu)$ le soit selon le proc\'ed\'e (b). Le proc\'ed\'e (a)  d\'efinit   un \'el\'ement que l'on note ici $\boldsymbol{\nu}_{1}$ et un couple que l'on note $(\alpha^1,\beta^1)\in {\cal P}_{n_{1}}\times {\cal P}_{m_{1}}$. Le proc\'ed\'e (b)  d\'efinit   un terme $ \mu_{1}$ et un couple $(\underline{\alpha}^1,\underline{\beta}^1)$. D'apr\`es 1.2(4), l'hypoth\`ese $(1,0)<_{I}(1,1)$ entra\^{\i}ne  que $\boldsymbol{\nu}_{1}=\alpha_{1}+ \mu_{1}$, $\underline{\alpha}^1=\{\alpha_{1}\}\sqcup \alpha^1$, $\underline{\beta}^1=\beta^1$. Puis 
$$(1) \qquad p\Lambda_{A,B;s}^{N,M}(\boldsymbol{\nu},\boldsymbol{\mu})=\{\alpha_{1}+\mu_{1}+A\}\sqcup p\Lambda_{A-s,B;s}^{N-1,M}(\boldsymbol{\nu}^1,\boldsymbol{\mu}^1) ,$$
$$(2)\qquad p\Lambda_{A,B;s}^{N,M}(\nu,\mu)=\{\mu_{1}+B\}\sqcup p\Lambda^{N,M-1}_{A,B-s;s}(\nu^1,\mu^1) ,$$
pour des couples $(\boldsymbol{\nu}^1,\boldsymbol{\mu}^1)\in P_{A-s,B;s}(\alpha^1,\beta^1)$ et $(\nu^1,\mu^1)\in P(\{\alpha_{1}\}\sqcup \alpha^1,\beta^1)$.

Supposons d'abord $m_{1}=0$. Dans ce cas, on a  $(\boldsymbol{\nu}^1,\boldsymbol{\mu}^1)=(\alpha^1,\beta^1=\emptyset)$ et $(\nu^1,\mu^1)=(\{\alpha_{1}\}\sqcup \alpha^1,\beta^1=\emptyset)$. On trouve plus explicitement
$$p\Lambda_{A,B;s}^{N,M}(\boldsymbol{\nu},\boldsymbol{\mu})=\{\alpha_{1}+\mu_{1}+A, \alpha^1_{1}+A-s
,...,\alpha^1_{n_{1}}+A-sn_{1},A-sn_{1}-s,...,A+s-sN\}$$
$$\sqcup\{B, B-s,...,B+s-sM\},$$
$$p\Lambda_{A,B;s}^{N,M}(\nu,\mu)=(\alpha_{1}+A,\alpha^1_{1}+A-s
,...,\alpha^1_{n_{1}}+A-sn_{1},A-sn_{1}-s,...,A+s-sN\}$$
$$\sqcup\{\mu_{1}+B, B-s,...,B+s-sM\}.$$
D'apr\`es la remarque 1.1(1), il suffit de prouver que $\{\alpha_{1}+A,\mu_{1}+B\}\leq  \{\alpha_{1}+\mu_{1}+A,B\}$. Or les hypoth\`eses que $(1,0)<_{I}(1,1)$ et que (a) est autoris\'e pour construire un \'el\'ement de $P_{A,B;s}(\alpha,\beta)$ entra\^{\i}nent que $\alpha_{1}+A\geq B$, ce qui entra\^{\i}ne l'in\'egalit\'e cherch\'ee. On voit aussi que l'on ne peut avoir l'\'egalit\'e  
 $p\Lambda_{A,B;s}^{N,M}(\boldsymbol{\nu},\boldsymbol{\mu})= p\Lambda_{A,B;s}^{N,M}(\nu,\mu)$ que si $\mu_{1}=0$ ou $\alpha_{1}+A=B$. Dans le premier cas, on voit que $(\nu,\mu)=(\boldsymbol{\nu},\boldsymbol{\mu})$ donc $(\nu,\mu)\in P_{A,B;s}(\alpha,\beta)$. Dans le deuxi\`eme cas, le proc\'ed\'e (b) est permis pour construire l'ensemble $P^{b}_{A,B;s}(\alpha,\beta)$ et $(\nu,\mu)$ appartient \`a cet ensemble.
 
 Supposons maintenant $m_{1}\geq1$. Introduisons un couple $(\underline{\nu}^1,\underline{\mu}^1)\in P_{A,B-s;s}(\{\alpha_{1}\}\sqcup \alpha^1,\beta^1)$ et notons $I_{1} =(\{1,...,n_{1}+1\}\times\{0\})\cup (\{1,...,m_{1}\}\times\{1\})$ l'ensemble d'indices intervenant, qui est muni d'un ordre $<I_{1}$ par une injection $I_{1}\to I$.   Puisque $(1,0)\in  I_{1}$ s'identifie \`a $(1,0)\in I$ qui est le plus petit \'el\'ement de $I$, il reste le plus petit \'el\'ement de $I_{1}$.  
 Comme on l'a dit ci-dessus, nos hypoth\`eses entra\^{\i}nent  $\alpha_{1}+A\geq B$, a fortiori $\alpha_{1}+A> B-s$. Donc $(\underline{\nu}^1,\underline{\mu}^1)$ est construit par le proc\'ed\'e (a).   Ce proc\'ed\'e  construit  des \'el\'ements, notons-les $a'_{1}=1$, $a'_{2}$,..., $b'_{1}$, $b'_{2}$... un \'el\'ement $\underline{\nu}_{1}=\alpha_{1}+\beta^1_{b'_{1}}+\alpha^1_{a'_{2}-1}+...$ (il y a un d\'ecalage sur les indices de $\alpha^1$ car le $i$-i\`eme terme de $\alpha^1$ est le $(i+1)$-i\`eme terme de $\{\alpha_{1}\}\sqcup \alpha^1$), puis un couple $(\alpha^2,\beta^2)\in {\cal P}_{n_{2}}\times {\cal P}_{m_{2}}$. Notons en passant que, puisque $(1,0)$ est le plus petit \'el\'ement de $I_{1}$, on a $b'_{1}=1$. On a $(\underline{\nu}^1,\underline{\mu}^1)=(\underline{\nu}_{1}\sqcup \nu^2,\mu^2)$, pour un couple $(\nu^2,\mu^2)\in P_{A-s,B-s;s}(\alpha^2,\beta^2)$. On a
 $$p\Lambda^{N,M-1}_{A,B-s;s}(\underline{\nu}^1,\underline{\mu}^1)=\{\underline{\nu}_{1}+A\}\sqcup p\Lambda^{N-1,M-1}_{A-s,B-s;s}(\nu^2,\mu^2) .$$
 De plus, l'hypoth\`ese de r\'ecurrence entra\^{\i}ne que $p\Lambda^{N,M-1}_{A,B-s;s}(\nu^1,\mu^1)\leq p\Lambda^{N,M-1}_{A,B-s;s}(\underline{\nu}^1,\underline{\mu}^1)$. Avec (2), on obtient
 $$  p\Lambda^{N,M}_{A,B;s}(\nu,\mu)\leq\{\mu_{1}+B\}\sqcup \{\underline{\nu}_{1}+A\}\sqcup p\Lambda^{N-1,M-1}_{A-s,B-s;s}(\nu^2,\mu^2)   .$$
 Posons $ \mu'_{1}=\beta^1_{b'_{1}}+\alpha^1_{a'_{2}-1}+...$. On a donc $\underline{\nu}_{1}=\alpha_{1}+ \mu'_{1}$. D'apr\`es 1.2(3), on a $\mu_{1}\geq  \mu'_{1}$. En utilisant l'in\'egalit\'e $\alpha_{1}+A\geq B$, on voit que $\{\mu_{1}+B\}\sqcup \{\underline{\nu}_{1}+A\}\leq \{\alpha_{1}+\mu_{1}+A\}\sqcup\{ \mu'_{1}+B\}$. Donc
 $$(3)\qquad  p\Lambda^{N,M}_{A,B;s}(\nu,\mu)\leq  \{\alpha_{1}+\mu_{1}+A\}\sqcup\{ \mu'_{1}+B\}\sqcup p\Lambda^{N-1,M-1}_{A-s,B-s;s}(\nu^2,\mu^2) .$$
 
On a vu ci-dessus que $b'_{1}=1$.  Mais alors $ \mu'_{1}$ est le terme issu du proc\'ed\'e (b) appliqu\'e \`a $(\alpha^1,\beta^1)$ et $(\alpha^2,\beta^2)$ est la partition qui se d\'eduit du m\^eme proc\'ed\'e. Ce proc\'ed\'e est loisible pour construire l'ensemble $P^{b}(\alpha^1,\beta^1)$ puisqu'on a suppos\'e $m_{1}\geq1$. Donc $(\nu^2,\{ \mu'_{1}\}\sqcup \mu^2)$ est un \'el\'ement de $ P(\alpha^1,\beta^1)$ et on a
 $$p\Lambda^{N-1,M}_{A-s,B;s}(\nu^2,\{ \mu'_{1}\}\sqcup \mu^2) =\{ \mu'_{1}+B\}\sqcup p\Lambda^{N-1,M-1}_{A-s,B-s;s}(\nu^2,\mu^2) .$$
 Par r\'ecurrence, on a 
 $$p\Lambda^{N-1,M}_{A-s,B;s}(\nu^2,\{ \mu'_{1}\}\sqcup \mu^2)\leq p\Lambda^{N-1,M}_{A-s,B;s}(\boldsymbol{\nu}^1,\boldsymbol{\mu}^1)$$
 et, avec (3), on obtient
 $$p\Lambda^{N,M}_{A,B}(\nu,\mu)\leq \{\alpha_{1}+\mu_{1}+A\}\sqcup p\Lambda^{N-1,M}_{A-s,B;s}(\boldsymbol{\nu}^1,\boldsymbol{\mu}^1) .$$
 En comparant avec (1), on obtient l'in\'egalit\'e cherch\'ee $p\Lambda^{N,M}_{A,B;s}(\nu,\mu)\leq p\Lambda^{N,M}_{A,B;s}(\boldsymbol{\nu},\boldsymbol{\mu})$.
 
 Supposons que $p\Lambda^{N,M}_{A,B;s}(\nu,\mu)= p\Lambda^{N,M}_{A,B;s}(\boldsymbol{\nu},\boldsymbol{\mu})$. Alors toutes les in\'egalit\'es ci-dessus doivent \^etre des \'egalit\'es. D'abord $p\Lambda^{N,M-1}_{A,B-s;s}(\nu^1,\mu^1)= p\Lambda^{N,M-1}_{A,B-s;s}(\underline{\nu}^1,\underline{\mu}^1)$. Par r\'ecurrence, cela entra\^{\i}ne que $(\nu^1,\mu^1)\in P_{A,B-s;s}( \{\alpha_{1}\}\sqcup \alpha^1,\beta^1)$. Puisque le couple $(\underline{\nu}^1,\underline{\mu}^1)$ \'etait un \'el\'ement quelconque de cet ensemble, on peut supposer qu'il est \'egal \`a $(\nu^1,\mu^1)$. Ensuite $\{\mu_{1}+B\}\sqcup \{\underline{\nu}_{1}+A\}= \{\alpha_{1}+\mu_{1}+A\}\sqcup\{ \mu'_{1}+B\}$. Cela entra\^{\i}ne que $\alpha_{1}+A=B$ ou que $ \mu'_{1}=\mu_{1}$. Dans le premier cas, le proc\'ed\'e (b) est permis pour d\'efinir l'ensemble $P^{b}_{A,B;s}(\alpha,\beta)$. Le couple $(\nu,\mu)$ est construit \`a l'aide de ce proc\'ed\'e et de l'\'el\'ement $(\nu^1,\mu^1)\in P_{A,B-s;s}( \{\alpha_{1}\}\sqcup \alpha^1,\beta^1)$. Donc $(\nu,\mu)\in P_{A,B;s}(\alpha,\beta)$. Supposons maintenant $ \mu'_{1}=\mu_{1}$. On doit aussi avoir l'\'egalit\'e $p\Lambda^{N-1,M}_{A-s,B;s}(\nu^2,\{ \mu'_{1}\}\sqcup \mu^2)= p\Lambda^{N-1,M}_{A-s,B;s}(\boldsymbol{\nu}^1,\boldsymbol{\mu}^1)$, ce qui entra\^{\i}ne par r\'ecurrence que $(\nu^2,\{ \mu'_{1}\}\sqcup \mu^2)\in P_{A-s,B;s}(\alpha^1,\beta^1)$. L'\'el\'ement de $P_{A,B;s}(\alpha,\beta)$ construit \`a l'aide du proc\'ed\'e (a) et de ce couple est
 $$(\{\alpha_{1}+\mu_{1}\}\sqcup \nu^2, \{ \mu' _{1}\}\sqcup \mu^2).$$
 Le couple $(\nu,\mu)$ est $ (\{\alpha_{1}+ \mu'_{1}\}\sqcup \nu^2,  \{\mu_{1}\}\sqcup \mu^2)$. Puisque $ \mu'_{1}=\mu_{1}$, ces deux couples sont \'egaux et $(\nu,\mu)$ appartient \`a $P_{A,B;s}(\alpha,\beta)$.

 On suppose maintenant  que $(\nu,\mu)$ est construit \`a l'aide du proc\'ed\'e (a) et que $(\boldsymbol{\nu},\boldsymbol{\mu})$ l'est \`a l'aide du proc\'ed\'e (b) (on suppose toujours  $(1,0)<_{I} (1,1)$). On reprend le raisonnement en \'echangeant les r\^oles de $(\nu,\mu)$ et $(\boldsymbol{\nu}, \boldsymbol{\mu})$. Les \'egalit\'es (1) et (2) sont invers\'ees, c'est-\`a-dire que l'on a
 $$(4) \qquad p\Lambda^{N,M}_{A,B;s}(\boldsymbol{\nu},\boldsymbol{\mu})=\{\mu_{1}+B\}\sqcup p\Lambda_{A,B-s;s}^{N,M-1}(\boldsymbol{\nu}^1,\boldsymbol{\mu}^1) ,$$
$$(5)\qquad p\Lambda^{N,M}_{A,B}(\nu,\mu)=\{\alpha_{1}+\mu_{1}+A\}\sqcup p \Lambda^{N-1,M}_{A-s,B;s}(\nu^1,\mu^1) ,$$
pour des couples $(\boldsymbol{\nu}^1,\boldsymbol{\mu}^1)\in P_{A,B-s;s}(\{\alpha_{1}\}\sqcup\alpha^1,\beta^1)$ et $(\nu^1,\mu^1)\in P( \alpha^1,\beta^1)$.

Supposons d'abord $m_{1}=0$. Dans ce cas , on a  $(\boldsymbol{\nu}^1,\boldsymbol{\mu}^1)=(\{\alpha_{1}\}\sqcup \alpha^1,\beta^1=\emptyset)$ et $(\nu^1,\mu^1)=(  \alpha^1,\beta^1=\emptyset)$. On trouve plus explicitement
$$p\Lambda_{A,B}^{N,M}(\boldsymbol{\nu},\boldsymbol{\mu})=\{\mu_{1}+B\}\sqcup\{ \alpha_{1}+A,\alpha^1_{1}+A-s
,...,\alpha^1_{n_{1}}+A-sn_{1},A-sn_{1}-s,...,A+s-sN\}$$
$$\sqcup\{ B-s,...,B+s-sM\},$$
$$p\Lambda^{N,M}_{A,B;s}(\nu,\mu)=(\alpha_{1}+\mu_{1}+A,\alpha^1_{1}+A-s
,...,\alpha^1_{n_{1}}+A-sn_{1},A-sn_{1}-s,...,A+s-sN\})$$
$$\sqcup\{B, B-s,...,B+s-sM\}.$$
D'apr\`es la remarque 1.1(1), il suffit de prouver que $\{\alpha_{1}+\mu_{1}+A,B\}\leq  \{\alpha_{1}+A,\mu_{1}+B\}$. Or les hypoth\`eses que $(1,0)<_{I}(1,1)$ et que (b) est autoris\'e pour construire des \'el\'ements de $P^{b}_{A,B;s}(\alpha,\beta)$  entra\^{\i}nent que $\alpha_{1}+A\leq B$, d'o\`u  l'in\'egalit\'e cherch\'ee. On voit aussi que l'on ne peut avoir l'\'egalit\'e  
 $p\Lambda^{N,M}_{A,B;s}(\boldsymbol{\nu},\boldsymbol{\mu})= p\Lambda^{N,M}_{A,B;s}(\nu,\mu)$ que si $\mu_{1}=0$ ou $\alpha_{1}+A=B$. Dans le premier cas, on voit que $(\nu,\mu)=(\boldsymbol{\nu},\boldsymbol{\mu})$ donc $(\nu,\mu)\in P_{A,B;s}(\alpha,\beta)$. Dans le deuxi\`eme cas, le proc\'ed\'e (a) est permis pour construire l'ensemble $P^{a}_{A,B;s}(\alpha,\beta)$ et $(\nu,\mu)$ appartient \`a cet ensemble.
 
  Supposons $m_{1}\geq1$. Introduisons un couple $(\underline{\nu}^1,\underline{\mu}^1)\in P_{A-s,B}(  \alpha^1,\beta^1)$.  Montrons que ce couple est obtenu par le proc\'ed\'e (b). C'est \'evident si $n_{1}=0$ puisqu'alors seul ce proc\'ed\'e est autoris\'e. Supposons $n_{1}\geq1$.   On note  $I_{1}=(\{1,...,n_{1}\}\times \{0\})\cup (\{1,...,m_{1}\}\times \{1\})$ l'ensemble d'indices intervenant.  Contrairement \`a ce qui se passait plus haut, l'\'el\'ement $(1,0)$ de $I_{1}$ ne s'identifie plus \`a $(1,0)\in I$. Pour l'ordre sur $I_{1}$, on ne sait pas quel est le plus grand des \'el\'ements $(1,0)$ et $(1,1)$. Supposons que  $(1,0)<_{I_{1}}(1,1)$. 
 Comme on l'a dit ci-dessus, nos hypoth\`eses entra\^{\i}nent  $\alpha_{1}+A\leq B$. Puisque $\alpha^1$ est une partition extraite de $\alpha$, on a $\alpha^1_{1}\leq \alpha_{1}$, a fortiori $\alpha^1_{1}+A-s<B$ et cela impose le proc\'ed\'e (b). Supposons maintenant que $(1,1)<_{I_{1}}(1,0)$. L'in\'egalit\'e $\alpha_{1}+A\leq B$ entra\^{\i}ne $A\leq B$, a fortiori $\beta^1_{1}+B> A-s$. Cela impose encore le proc\'ed\'e (b), ce qui d\'emontre l'assertion. 
   Le proc\'ed\'e (b) construit  des \'el\'ements, notons-les $b'_{1}=1$, $b'_{2}$... $a'_{1}$, $a'_{2}$,...  un \'el\'ement $\underline{\mu}_{1}= \beta^1_{b'_{1}}+\alpha^1_{a'_{1}}+...$, puis un couple $(\alpha^2,\beta^2)\in {\cal P}_{n_{2}}\times {\cal P}_{m_{2}}$. On a $(\underline{\nu}^1,\underline{\mu}^1)=( \nu^2,\{\underline{\mu}_{1}\}\sqcup \mu^2)$, pour un couple $(\nu^2,\mu^2)\in P_{A-s,B-s;s}(\alpha^2,\beta^2)$. On a
 $$p\Lambda^{N-1,M}_{A-s,B;s}(\underline{\nu}^1,\underline{\mu}^1)=\{\underline{\mu}_{1}+B\}\sqcup p\Lambda^{N-1,M-1}_{A-s,B-s;s}(\nu^2,\mu^2) .$$
 De plus, l'hypoth\`ese de r\'ecurrence entra\^{\i}ne que $p\Lambda^{N-1,M}_{A-s,B;s}(\nu^1,\mu^1)\leq p\Lambda^{N-1,M}_{A-s,B;s}(\underline{\nu}^1,\underline{\mu}^1)$. Avec (5), on obtient
 $$  p\Lambda^{N,M}_{A,B;s}(\nu,\mu)\leq\{\alpha_{1}+\mu_{1}+A\}\sqcup \{\underline{\mu}_{1}+B\}\sqcup p\Lambda^{N-1,M-1}_{A-s,B-s;s}(\nu^2,\mu^2)  .$$
  D'apr\`es 1.2(3), on a  $\underline{\mu}_{1}\leq \mu_{1}$. En utilisant l'in\'egalit\'e $\alpha_{1}+A\leq B$, on voit que $\{\alpha_{1}+\mu_{1}+A\}\sqcup \{\underline{\mu}_{1}+B\}\leq \{\alpha_{1}+\underline{\mu}_{1}+A\}\sqcup\{ \mu_{1}+B\}$. 
  Donc
 $$(6)\qquad  p\Lambda^{N,M}_{A,B;s}(\nu,\mu)\leq  \{\mu_{1}+B\}\sqcup\{ \alpha_{1}+\underline{\mu}_{1}+A\}\sqcup p\Lambda^{N-1,M-1}_{A-s,B-s;s}(\nu^2,\mu^2) .$$
 Consid\'erons le couple $(\{\alpha_{1}\}\sqcup \alpha^1,\beta^1)$. Cette fois, l'ordre  sur l'ensemble d'indices  ${\bf I}_{1}=(\{1,...,n_{1}+1\}\times \{0\})\cup (\{1,...,m_{1}\}\times \{1\})$ est tel que $(1,0)<_{{\bf I}_{1}}(1,1)$. On voit que le proc\'ed\'e (a) cr\'ee l'\'el\'ement $\alpha_{1}+\underline{\mu}_{1}$ et les partitions $(\alpha^2,\mu^2)$. Donc le couple $(\{\alpha_{1}+\underline{\mu}_{1}\}\sqcup \nu^2,\mu^2)$ appartient \`a $P(\{\alpha_{1}\}\sqcup \alpha^1,\beta^1)$.
 On a
 $$p\Lambda^{N,M-1}_{A,B-s;s}(\{\alpha_{1}+\underline{\mu}_{1}\}\sqcup \nu^2,\mu^2)=\{\alpha_{1}+\underline{\mu}_{1}+A\}\sqcup p\Lambda^{N-1,M-1}_{A-s,B-s;s}(\nu^2,\mu^2).$$
  D'o\`u, avec (6),
  $$p\Lambda^{N,M}_{A,B;s}(\nu,\mu)\leq  \{\mu_{1}+B\}\sqcup p\Lambda^{N,M-1}_{A,B-s;s}(\{\alpha_{1}+\underline{\mu}_{1}\}\sqcup \nu^2,\mu^2).$$
   Puisque $(\{\alpha_{1}+\underline{\mu}_{1}\}\sqcup \nu^2,\mu^2)$ appartient \`a $P(\{\alpha_{1}\}\sqcup \alpha^1,\beta^1)$, on sait  par r\'ecurrence que $p\Lambda^{N,M-1}_{A,B-s;s}(\{\alpha_{1}+\underline{\mu}_{1}\}\sqcup \nu^2,\mu^2)\leq p\Lambda^{N,M-1}_{A,B-s;s}(\boldsymbol{\nu}^1,\boldsymbol{\mu}^1)$. Alors
  $$p\Lambda_{A,B;s}^{N,M}(\nu,\mu)\leq  \{\mu_{1}+B\}\sqcup  p\Lambda^{N,M-1}_{A,B-s;s}(\boldsymbol{\nu}^1,\boldsymbol{\mu}^1) .$$
  En comparant avec (4), on obtient l'in\'egalit\'e cherch\'ee $p\Lambda^{N,M}_{A,B;s}(\nu,\mu)\leq p\Lambda^{N,M}_{A,B;s}(\boldsymbol{\nu},\boldsymbol{\mu})$.

  Supposons que $p\Lambda_{A,B;s}^{N,M}(\nu,\mu)= p\Lambda_{A,B;s}^{N,M}(\boldsymbol{\nu},\boldsymbol{\mu})$. Alors toutes les in\'egalit\'es ci-dessus doivent \^etre des \'egalit\'es. D'abord $p\Lambda^{N-1,M}_{A-s,B;s}(\nu^1,\mu^1)= p\Lambda^{N-1,M}_{A-s,B;s}(\underline{\nu}^1,\underline{\mu}^1)$. Cela entra\^{\i}ne par r\'ecurrence que $(\nu^1,\mu^1)\in P_{A-s,B;s}(\alpha^1,\beta^1)$. Comme plus haut, on peut alors supposer $(\underline{\nu}^1,\underline{\mu}^1)=(\nu^1,\mu^1)$. Ensuite $\{\alpha_{1}+\mu_{1}+A\}\sqcup \{\underline{\mu}_{1}+B\}= \{\alpha_{1}+\underline{\mu}_{1}+A\}\sqcup\{ \mu_{1}+B\}$. Cela entra\^{\i}ne que $\alpha_{1}+A=B$ ou que $\underline{\mu}_{1}=\mu_{1}$. Dans le premier cas, le proc\'ed\'e (a) est permis pour construire l'ensemble $P^{a}_{A,B;s}(\alpha,\beta)$ et $(\nu,\mu)$ est issu de ce proc\'ed\'e et du couple $(\nu^1,\mu^1)\in P_{A-s,B;s}(\alpha^1,\beta^1)$. Donc $(\nu,\mu)\in P_{A,B;s}(\alpha,\beta)$. Supposons $\underline{\mu}_{1}=\mu_{1}$. On doit encore avoir $p\Lambda^{N,M-1}_{A,B-s;s}(\{\alpha_{1}+\underline{\mu}_{1}\}\sqcup \nu^2,\mu^2)= p\Lambda^{N,M-1}_{A,B-s;s}(\boldsymbol{\nu}^1,\boldsymbol{\mu}^1)$. Par r\'ecurrence, cela entra\^{\i}ne que $(\{\alpha_{1}+\underline{\mu}_{1}\}\sqcup \nu^2,\mu^2)\in P_{A,B-s;s}(\{\alpha_{1}\}\sqcup \alpha^1,\beta^1)$. L'\'el\'ement de $P^{b}_{A,B;s}(\alpha,\beta)$ issu du proc\'ed\'e (b) et de ce couple est
  $$(\{\alpha_{1}+\underline{\mu}_{1}\}\sqcup \nu^2,\{\mu_{1}\}\sqcup\Lambda^2).$$
  Le couple $(\nu,\mu)$ est
  $$(\{\alpha_{1}+\mu_{1}\}\sqcup \nu^2,\{\underline{\mu}_{1}\}\sqcup \mu^2).$$
  L'\'egalit\'e $\underline{\mu}_{1}=\mu_{1}$ entra\^{\i}ne que ces deux couples sont \'egaux, donc $(\nu,\mu)\in P_{A,B;s}(\alpha,\beta)$. Cela ach\`eve la d\'emonstration. $\square$
 
Le lemme autorise la d\'efinition suivante.

\ass{D\'efinition}{Pour $\alpha\in {\cal P}_{n}$ et $\beta\in {\cal P}_{m}$, on note $p_{A,B;s}^{N,M}[\alpha,\beta]$ la partition $p\Lambda_{A,B;s}^{N,M}(\nu,\mu)$ pour un \'el\'ement quelconque $(\nu,\mu)\in P_{A,B;s}(\alpha,\beta)$.}

{\bf Remarque.} Il r\'esulte de 1.2(1) que, si $n=0$ ou $m=0$, on a simplement $p_{A,B;s}^{N,M}[\alpha,\beta]=p\Lambda_{A,B;s}^{N,M}(\alpha,\beta)$. 

\bigskip

\section{Quelques propri\'et\'es des ensembles $P(\alpha,\beta)$}
\bigskip

\subsection{La partition $\alpha^-$}

Pour $n\geq1$ et pour  $\lambda=(\lambda_{1},...,\lambda_{n})\in {\cal R}_{n}$, on d\'efinit l'\'el\'ement $\lambda^-\in {\cal R}_{n-1}$ par $\lambda^-=(\lambda_{2},...,\lambda_{n})$.

 Soient $\alpha\in {\cal P}_{n}$ et $\beta\in {\cal P}_{m}$. On suppose $n\geq1$ et $(1,0)<_{I}(1,1)$ si $m\geq1$.     On plonge l'ensemble d'indices $I^-=(\{1,...,n-1\}\times \{0\})\cup(\{1,...,m\}\times \{1\})$ dans $I$ en envoyant $(i,0)\in I^-$ sur $(i+1,0)\in I$. On munit $I^-$ de l'ordre $<_{I^-}$ induit par ce plongement et l'ordre $<_{I}$ sur $I$. On utilise cet ordre pour construire les objets relatifs \`a $(\alpha^-,\beta)$. 
 
 On utilisera plusieurs fois les propri\'et\'es suivantes, qui se d\'emontrent simplement en appliquant les d\'efinitions. On suppose $m\geq1$. 
  Appliquons \`a $(\alpha,\beta)$ le proc\'ed\'e (a). Il construit un \'el\'ement $\nu_{1}$ et des partitions $(\alpha',\beta')$. Alors
  
  (1) le proc\'ed\'e (b) appliqu\'e \`a $(\alpha^-,\beta)$ construit l'\'el\'ement $\underline{\mu}_{1}=\nu_{1}-\alpha_{1}$ et les m\^emes partitions $(\alpha',\beta')$;
  
  (2) si $(1,0)<_{I^-}(1,1)$, ce qui \'equivaut \`a $(2,0)<_{I}(1,1)$, le proc\'ed\'e (a) appliqu\'e \`a $(\alpha^-,\beta)$ construit l'\'el\'ement $\underline{\nu}_{1}=\nu_{1}-\alpha_{1}+\alpha_{2}$ et les partitions $({\alpha'}^-,\beta')$;
  
  (3) si $(1,1)<_{I^-}(1,0)$, ce qui \'equivaut \`a $(1,1)<_{I}(2,0)$, le proc\'ed\'e (a) appliqu\'e \`a $(\alpha^-,\beta)$ construit l'\'el\'ement $\underline{\nu}_{1}=\nu_{1}-\alpha_{1}-\beta_{1}$ et les partitions $(\alpha',\{\beta_{1}\}\sqcup\beta')$.
  
   Appliquons \`a $(\alpha,\beta)$ le proc\'ed\'e (b). Il construit un \'el\'ement $\mu_{1}$ et des partitions $(\alpha',\beta')$. Alors 
 
 (4) le proc\'ed\'e (b) appliqu\'e \`a $(\alpha^-,\beta)$ construit le m\^eme \'el\'ement $\mu_{1}$ et les partitions $({\alpha'}^-,\beta')$.
 
 Supposons $\alpha_{1}+A\leq B$, auquel cas le proc\'ed\'e (b) est autoris\'e pour construire les \'el\'ements de $P_{A,B;s}(\alpha,\beta)$. Consid\'erons de plus des r\'eels $A',B'$ tels que $A\geq A'$ et $B'\geq B$. Alors
 
 (5) le proc\'ed\'e (b) est autoris\'e pour construire les \'el\'ements de $P_{A',B';s}(\alpha^-,\beta)$.
 
 En effet, si $(1,0)<_{I^-}(1,1)$, on doit voir que $\alpha^-_{1}+A'\leq B'$. Mais $\alpha^-_{1}+A'=\alpha_{2}+A'\leq \alpha_{1}+A\leq B\leq B'$. Si $(1,1)<_{I^-}(1,0)$, on doit voir que $\beta_{1}+B'\geq A'$. Mais $\beta_{1}+B'\geq B'\geq B\geq \alpha_{1}+A\geq A\geq A'$.
 
 \bigskip

 \subsection{L'ensemble $P^{b[c]}(\alpha,\beta)$}
 Soient $\alpha\in {\cal P}_{n}$ et $\beta\in {\cal P}_{m}$ et  soit $c\in {\mathbb N}$. On d\'efinit un sous-ensemble $P^{b[c]}(\alpha,\beta)$ de la fa\c{c}on suivante. Si $m=0$, on pose $P^{b[c]}(\alpha,\beta)=P(\alpha,\beta)$. Supposons $m\geq 1$. Si $c=0$, on pose $P^{b[0]}(\alpha,\beta)= P(\alpha,\beta)$. Supposons $c\geq1$.  Le proc\'ed\'e (b) construit un terme $\mu_{1}$ et des partitions $(\alpha',\beta')\in {\cal P}_{n'}\times {\cal P}_{m'}$. Alors $P^{b[c]}(\alpha,\beta)$ est l'ensemble des couples $(\nu'\sqcup\{0^{n-n'}\};\{\mu_{1}\}\sqcup \mu'\sqcup\{0^{m-1-m'}\})$, pour $(\nu',\mu')\in P^{b[c-1]}(\alpha',\beta')$.  

Montrons que

(1) pour $c,c'\geq m$, on a $P^{b[c]}(\alpha,\beta)=P^{b[c']}(\alpha,\beta)$.

Si $m=0$, les deux ensembles sont \'egaux \`a $P(\alpha,\beta)$. Si $m\geq1$, donc aussi $c,c'\geq1$, on a $m'<m$ dans la construction ci-dessus et la propri\'et\'e r\'esulte par r\'ecurrence de l'\'egalit\'e $P^{b[c-1]}(\alpha',\beta')=P^{b[c'-1]}(\alpha',\beta')$.

 Supposons $m\geq1$, consid\'erons un couple $(\nu,\mu)\in P(\alpha,\beta)$ et un entier $c\in \{1,...,m\}$. Pour tout $d\in \{1,...,c\}$, on suppose donn\'es des r\'eels $A(d)\geq sN-s$, $B(d)\geq sM-s$. On consid\`ere la propri\'et\'e 
  
  (2) $S_{d}((\mu\sqcup\{0^{M-m}\})+[B(d),B(d)+s-sM]_{s})=S_{d}(p^{N,M}_{A(d),B(d);s}[\alpha,\beta])$ pour tout $d\in \{1,...,c\}$. 
   
  \ass{Lemme}{Supposons  (2) v\'erifi\'ee. Alors $(\nu,\mu)$ appartient \`a $P^{b[c]}(\alpha,\beta)$.}
  
Preuve.   Supposons d'abord $c=1$ et posons pour simplifier $A=A(1)$, $B=B(1)$. 
Notons $(\chi,\xi)$ l'\'el\'ement canonique de type (b) de $P_{A,B;s}(\alpha,\beta)$. 
Montrons tout d'abord que
  
(3) cet \'el\'ement est construit \`a l'aide du proc\'ed\'e (b) et on a $\mu_{1}=\xi_{1}$. 

Si $n=0$, le proc\'ed\'e (b) est le seul permis et on a $\mu_{1}=\beta_{1}=\xi_{1}$. Supposons $n\geq1$.  Supposons que  $(1,0)<_{I}(1,1)$  et $\alpha_{1}+A>B$, ou que $(1,1)<_{I}(1,0)$ et $\beta_{1}+B< A$.  Alors $(\chi,\xi)$ est construit \`a l'aide du  proc\'ed\'e (a). Celui-ci  cr\'ee   l'\'el\'ement $\chi_{1} $ et on a
$S_{1}(p_{A,B;s}^{N,M}[\alpha,\beta])=\chi_{1}+A$, cf. 1.3(1). Par l'hypoth\`ese (2), ceci est \'egal \`a $\mu_{1}+B$. Puisque $(\nu,\mu)\in P(\alpha,\beta)$, les relations 1.2(3) et 1.2(4) montrent que 

si $(1,0)<_{I}(1,1)$, $\mu_{1}\leq \chi_{1}-\alpha_{1}$;

si $(1,1)<_{I}(1,0)$, $\mu_{1}\leq \chi_{1}+\beta_{1}$.

Dans le premier cas, on a $\chi_{1}+A=\mu_{1}+B\leq \chi_{1}-\alpha_{1}+B$, c'est-\`a-dire $\alpha_{1}+A\leq B$  contrairement \`a l'hypoth\`ese. Dans le deuxi\`eme, on a $\chi_{1}+A=\mu_{1}+B\leq \chi_{1}+\beta_{1}+B$, c'est-\`a-dire $\beta_{1}+B\geq A$ contrairement \`a l'hypoth\`ese. Ces contradictions prouvent que, si  $(1,0)<_{I}(1,1)$, on a $\alpha_{1}+A\leq B$, et que, si $(1,1)<_{I}(1,0)$, on a $\beta_{1}+B\geq A$. Alors $(\chi,\xi)$ est construit \`a l'aide du proc\'ed\'e (b) d'apr\`es la d\'efinition de cet \'el\'ement canonique. 
 Cela d\'emontre la premi\`ere assertion de (3). On a donc $S_{1}(p_{A,B;s}^{N,M}[\alpha,\beta])=\xi_{1}+B$ d'apr\`es 1.3(1). Par hypoth\`ese, ceci est \'egal \`a $\mu_{1}+B$, d'o\`u $\xi_{1}=\mu_{1}$. Cela prouve (3). 

On va d\'efinir un entier $r\geq 1$ des couples de partitions $(\alpha^{i},\beta^{i})\in {\cal P}_{n_{i}}\times {\cal P}_{m_{i}}$ et $(\nu^{i},\mu^{i}) \in P(\alpha^{i},\beta^{i})$ pour $i=1,...,r$. On pose $\alpha^1=\alpha$, $\beta^1=\beta$, $\nu^1=\nu$, $\mu^1=\mu$.   Supposons construits ces objets jusqu'au rang $i$. Si $m_{i}=0$ ou si $m_{i}>0$ et $(\nu^{i},\mu^{i})$ appartient \`a $ P^{b}(\alpha^{i},\beta^{i})$, on pose $r=i$ et le proc\'ed\'e s'arr\^ete. Si $m_{i}>0$ et $(\nu^{i},\mu^{i})$ n'appartient pas \`a $ P^{b}(\alpha^{i},\beta^{i})$, alors $n_{i}>0$ et $(\nu^{i},\mu^{i})$ appartient \`a $P^{a}(\alpha^{i},\beta^{i})$. Ce proc\'ed\'e cr\'ee le terme $\nu^{i}_{1}$ et un couple de partitions $(\alpha^{i+1},\beta^{i+1})$. On a une \'egalit\'e $\nu^{i}=\{\nu^{i}_{1}\}\sqcup\nu^{i+1}\sqcup\{0^{n_{i}-1-n_{i+1}}\}$, $\mu^{i}=\mu^{i+1}\sqcup\{0^{m_{i}-m_{i+1}}\}$,  pour un couple $(\nu^{i+1},\mu^{i+1})\in P(\alpha^{i+1},\beta^{i+1})$. Cela d\'efinit nos suites. On pose $I_{i}=(\{1,...,n_{i}\}\times\{0\})\cup (\{1,...,m_{i}\}\times\{1\})$. Comme en 1.2, $I_{i+1}$ se plonge dans $I_{i}$. Par r\'ecurrence, $I_{i}$ se plonge  dans $I_{1}=I$ et h\'erite d'un ordre $<_{I_{i}}$.

On va prouver

(4) $r=1$. 

On voit par r\'ecurrence que 

pour $i=1,...,r$, $\nu_{i}=\nu^{i}_{1}$;

si $m_{r}>0$,  $\mu_{1}=\mu_{1}^{i}$ pour $i=1,...,r$;

si $m_{r}=0$,  $\mu_{1}=\mu_{1}^{i}=0$ pour $i=1,...,r-1$;

pour $i=1,...,r$, 
 ${\Lambda_{A,B;s}^{N,M}(\nu,\mu)}^{a}=\{\nu_{1}+A,...,\nu_{i-1}+A+2s-si\}\sqcup {\Lambda_{A+s-si,B}^{N+1-i,M}(\nu^{i},\mu^{i})}^{a}$ et   ${\Lambda_{A,B;s}^{N,M}(\nu,\mu)}^{b}={\Lambda_{A+s-si,B}^{N+1-i,M}(\nu^{i},\mu^{i})}^{b}$, o\`u par exemple ${\Lambda_{A,B;s}^{N,M}(\nu,\mu)}^{a}$ et ${\Lambda_{A,B;s}^{N,M}(\nu,\mu)}^{b}$ sont les deux composantes du $(N,M;s)$-symbole $\Lambda_{A,B;s}^{N,M}(\nu,\mu)$. 
 
 Supposons d'abord $m_{r}=0$. Alors, comme on vient de le dire, $\mu_{1}=0$. D'apr\`es (3), on a $\xi_{1}=0$. Par construction, $\xi_{1}=\beta_{1}+\alpha_{a_{1}}+\beta_{b_{2}}+...$. Donc tous ces termes sont nuls. Rappelons que $a_{2}$ est le plus petit entier $i$ tel que $(1,1)<_{I}(i,0)$ si un tel entier existe.   Cela entra\^{\i}ne $\beta=\{0^m\}$ et, si un terme $\alpha_{i}$ est non nul, on a $(i,0)<_{I}(1,1)$.  D'apr\`es 1.2(5), $P(\alpha,\beta)$ est r\'eduit \`a un seul \'el\'ement, \`a savoir $(\alpha,\beta)$. Donc $(\nu,\mu)=(\chi,\xi)$ et (3) entra\^{\i}ne que $(\nu,\mu)$ est construit \`a l'aide du proc\'ed\'e (b). 
En vertu de notre hypoth\`ese $m>0$, cela entra\^{\i}ne par d\'efinition de $r$  que $r=1$, mais alors les hypoth\`eses $m_{r}=0$ et $m>0$ sont contradictoires. Supposons maintenant $m_{r}>0$. Alors $(\nu^r,\mu^r)$ est construit par le proc\'ed\'e (b) par d\'efinition de $r$. Celui-ci cr\'ee l'\'el\'ement $\mu^{r}_{1}=\mu_{1}$ et des partitions $(\alpha^{r+1},\beta^{r+1})$. On a $\nu^r =\nu^{r+1}\sqcup\{0^{n_{r}-n_{r+1}}\}$ et $\mu^r=\{\mu_{1}\}\sqcup \mu^{r+1}\sqcup\{0^{m_{r}-1-m_{r+1}}\}$, pour un couple $(\nu^{r+1},\mu^{r+1})\in P(\alpha^{r+1},\beta^{r+1})$. Supposons que $r\geq2$ et supposons d'abord $(1,0)<_{I_{r-1}}(1,1)$. Alors $\nu_{r-1}=\nu^{r-1}_{1}=\alpha^{r-1}_{1}+\beta^{r-1}_{1}+\alpha^{r-1}_{a_{2}}+...$.  D'apr\`es  1.2(3), on a $\mu_{1}=\mu^r_{1}\leq \beta^{r-1}_{1}+\alpha^{r-1}_{a_{2}}+...\leq \xi_{1}$. Les deux termes extr\^emes \'etant \'egaux par (3), ces in\'egalit\'es sont des \'egalit\'es. Donc $\nu_{r-1}=\alpha^{r-1}_{1}+\mu_{1}$. Le proc\'ed\'e (b) appliqu\'e \`a $(\alpha^{r-1},\beta^{r-1})$ cr\'ee l'\'el\'ement $\mu_{1}$ et le couple de partitions $(\{\alpha_{1}^{r-1}\}\sqcup \alpha^{r},\beta^r)$. Pour l'ensemble d'indices de ce couple $(1,0)$ reste l'\'el\'ement minimal. Le proc\'ed\'e (a) appliqu\'e \`a ce couple cr\'ee donc l'\'el\'ement $\alpha_{1}+\mu_{1}^r=\alpha_{1}+\mu_{1}$ et le couple de partitions $(\alpha^{r+1},\beta^{r+1})$. Il en r\'esulte que $(\{\alpha_{1}+\mu_{1}\}\sqcup \nu^{r+1}\sqcup\{0^{n_{r}-n_{r+1}}\},\mu^{r+1}\sqcup\{0^{m_{r}-m_{r+1}}\})$ appartient \`a $P(\{\alpha_{1}^{r-1}\}\sqcup \alpha^{r},\beta^r)$. L'\'el\'ement de $P(\alpha^{r-1},\beta^{r-1})$ qui s'en d\'eduit par le proc\'ed\'e (b) est 
$$(\{\alpha_{1}+\mu_{1}\}\sqcup \nu^{r+1}\sqcup\{0^{n_{r-1}-1-n_{r+1}}\},\{\mu_{1}\}\sqcup\mu^{r+1}\sqcup\{0^{m_{r-1}-1-m_{r+1}}\})$$
$$=(\{\alpha_{1}+\mu_{1}\}\sqcup \nu^{r}\sqcup\{0^{n_{r-1}-1-n_{r}}\},\mu^r\sqcup\{0^{m_{r-1}-m_{r}}\})=(\nu^{r-1},\mu^{r-1}).$$
Cela montre que $(\nu^{r-1},\mu^{r-1})$ est obtenu par le proc\'ed\'e (b), ce qui contredit la d\'efinition de $r$. Supposons maintenant que $(1,1)<_{I_{r-1}}(1,0)$. Le proc\'ed\'e (b) appliqu\'e \`a $(\alpha^{r-1},\beta^{r-1})$ cr\'ee l'\'el\'ement  $\beta^{r-1}_{1}+\nu^{r}_{1}$ d'apr\`es 1.2(4). Comme ci-dessus, on a
$$\mu_{1}=\mu^r_{1}\leq \beta^{r-1}_{1}+ \nu^r_{1}\leq \xi_{1},$$
d'o\`u l'\'egalit\'e de ces termes. Donc $\nu^r_{1}=\mu_{1}-\beta^{r-1}_{1}$. N\'ecessairement, le premier terme de $\beta^r$  est $\beta^{r-1}_{1}$ (pr\'ecis\'ement, l'\'el\'ement $(1,1)$ de l'ensemble d'indices $I_{r}$ de $(\alpha^r,\beta^r)$ s'identifie au m\^eme \'el\'ement $(1,1)\in I_{r-1}$) . Le proc\'ed\'e (b) appliqu\'e \`a $(\alpha^{r-1},\beta^{r-1})$ cr\'ee le couple de partitions $(\alpha^r,(\beta^r)^-)$. Supposons $n_{r}>0$. Alors le proc\'ed\'e (a) appliqu\'e \`a $(\alpha^r,(\beta^r)^-)$ cr\'ee l'\'el\'ement $\mu^r_{1}-\beta^r_{1}=\mu_{1}-\beta^{r-1}_{1}$ et le m\^eme couple de partitions $(\alpha^{r+1},\beta^{r+1})$. L'\'el\'ement $(\{\mu_{1}-\beta^{r-1}_{1}\}\sqcup\nu^{r+1}\sqcup\{0^{n_{r}-1-n_{r+1}}\},\mu^{r+1}\sqcup\{0^{m_{r}-1-m_{r+1}}\})$ appartient \`a $P(\alpha^r,(\beta^r)^-)$ et l'\'el\'ement de $P(\alpha^{r-1},\beta^{r-1})$ qui s'en d\'eduit par le proc\'ed\'e (b) est 
$$(\{\mu_{1}-\beta^{r-1}_{1}\}\sqcup\nu^{r+1}\sqcup\{0^{n_{r-1}-1-n_{r+1}}\},\{\mu_{1}\}\sqcup \mu^{r+1}\sqcup\{0^{m_{r}-1-m_{r+1}}\})$$
$$=(\{\mu_{1}-\beta^{r-1}_{1}\}\sqcup\nu^{r}\sqcup\{0^{n_{r-1}-1-n_{r}}\},  \mu^{r}\sqcup\{0^{m_{r-1}-m_{r}}\})=(\nu^{r-1},\mu^{r-1}).$$
Cela montre que $(\nu^{r-1},\mu^{r-1})$ est obtenu par le proc\'ed\'e (b), ce qui contredit la d\'efinition de $r$. Il reste le cas o\`u $n_{r}=0$. On a alors $\mu^{r}_{1}=\beta^{r}_{1}$. L'\'egalit\'e
$$\mu^r_{1}= \beta^{r-1}_{1}+\alpha^{r-1}_{1}+\beta^{r-1}_{b_{2}}+...$$
entra\^{\i}ne  que $\beta^r_{1}=\beta^{r-1}_{1}$, $\alpha^{r-1}_{1}=\beta^{r-1}_{b_{2}}=...=0$. Autrement dit $\alpha^{r-1}=\{0^{n_{r-1}}\}$ et les termes $\beta^{r-1}_{i}$ ne sont non nuls que si $(i,1)<_{I_{r-1}}(1,0)$.  D'apr\`es 1.2(5), il n'y a qu'un seul \'el\'ement dans $P(\alpha^{r-1},\beta^{r-1})$, qui peut \^etre obtenu par le proc\'ed\'e (b) puisque $m_{r-1}>0$. Cela contredit encore la d\'efinition de $r$. Cela prouve (4).

D'apr\`es l'hypoth\`ese $m>0$, l'\'egalit\'e $r=1$ entra\^{\i}ne   
 que $(\nu,\mu)$ est construit \`a l'aide du proc\'ed\'e (b).  Cela d\'emontre le lemme dans le cas $c=1$. 
 
 Supposons maintenant $c>1$. L'hypoth\`ese (2) pour $c$ est  plus forte que celle pour $c-1$. En raisonnant par r\'ecurrence, on sait que $(\nu,\mu)$ appartient \`a $P^{b[c-1]}(\alpha,\beta)$. On pose $(\alpha^1,\beta^1)=(\alpha,\beta)$, $(\nu^1,\mu^1)=(\nu,\mu)$ et on d\'efinit maintenant un entier $t$ et, pour $i=1,...,t$, des couples $(\alpha^{i},\beta^{i})$ et $(\nu^{i},\beta^{i})\in P(\alpha^{i},\beta^{i})$ de la fa\c{c}on suivante. Si $(\nu^{i},\mu^{i})$ ne peut  pas \^etre construit par le proc\'ed\'e (b), on pose $t=i$ et la construction s'arr\^ete. Si $(\nu^{i},\mu^{i})$ est construit par le proc\'ed\'e (b), ce proc\'ed\'e cr\'ee un terme $\mu^{i}_{1}$ et des partitions $(\alpha^{i+1},\beta^{i+1})$. On a $\nu^{i}=\nu^{i+1}\sqcup\{0^{n_{i}-n_{i+1}}\}$ et $\mu^{i}=\{\mu^{i}_{1}\}\sqcup\mu^{i+1}\sqcup\{0^{m_{i}-1-m_{i+1}}\}$ pour un couple $(\nu^{i+1},\mu^{i+1})\in P(\alpha^{i+1},\beta^{i+1})$. On voit que $\mu^{i}_{1}=\mu_{i}$ (sauf \'eventuellement si $m_{t}=0$ auquel cas $\mu^{t}_{1}$ n'existe pas; dans ce cas, on a $\mu_{t}=0$). La condition $(\nu,\mu) \in P^{b[c-1]}(\alpha,\beta)$ signifie que $t\geq c$ ou que $t<c$ et $m_{t}=0$. Supposons d'abord $t<c$ et $m_{t}=0$. Dans ce cas, on a $(\nu^t,\mu^t)\in P^{b[c']}(\alpha^{t},\beta^t)$ pour tout $c'\in {\mathbb N}$. En particulier $(\nu^t,\mu^t)\in P^{b[c+1-t]}(\alpha^{t},\beta^t)$. Alors par r\'ecurrence, $(\nu^{i},\mu^{i})\in P^{b[c+1-i]}(\alpha^{i},\beta^{i})$, d'o\`u $(\nu,\mu)\in P^{b[c]}(\alpha,\beta)$. Remarquons que ce raisonnement vaut aussi si $t=c$ et $m_{t}=0$. Supposons maintenant $t\geq c$. Si $t\geq c+1$, alors $(\nu,\mu)\in P^{b[c]}(\alpha,\beta)$ comme on le voulait. On suppose $t=c$. Comme on vient de le dire, on peut aussi supposer $m_{t}>0$.   
 On pose pour simplifier $A=A(c)$, $ B=B(c)$. Notons $(\chi,\xi)$ l'\'el\'ement canonique de type (b) de $P_{A,B}(\alpha,\beta)$. On applique la m\^eme construction en rempla\c{c}ant $(\nu,\mu)$ par $(\chi,\xi)$. On note $T$ l'entier similaire \`a $t$ et  $(\underline{\alpha}^{i},\underline{\beta}^{i})$ et $(\chi^{i},\xi^{i})$ les analogues de $(\alpha^{i},\beta^{i})$ et $(\nu^{i},\mu^{i})$. Evidemment, pour $i\leq inf(t,T)$, les couples de partitions $(\alpha^{i},\beta^{i})$ et $(\underline{\alpha}^{i},\underline{\beta}^{i})$  sont les m\^emes. Les termes $\xi_{i}=\xi^{i}_{1}$ et $\mu_{i}=\mu^{i}_{1}$ sont aussi les m\^emes pour $i=1,...inf(t,T)-1$. Montrons que
 
 (5) on a $T\geq t$. 
 
 Supposons $T<t$. On a $m_{T}>0$ puisqu'on a suppos\'e $m_{t}>0$. Par d\'efinition de $T$, $(\chi^{T},\xi^T)$ est construit \`a l'aide du proc\'ed\'e (a). En particulier $n_{T}>0$. Supposons $(1,0)<_{I_{T}}(1,1)$. Puisque $(\chi,\xi)$ est l'\'el\'ement canonique de type (b), et que ce proc\'ed\'e n'est pas utilis\'e, on a $\alpha^T_{1}+A> B+s-sT$. Puisque les partitions $(\alpha^{i},\beta^{i})$ pour $i=T+1,...,t$ se d\'eduisent de $(\alpha^T,\beta^T)$ par applications successives de proc\'ed\'es (b), le terme $\alpha_{1}^T$ se conserve, c'est-\`a-dire $\alpha_{1}^t=\alpha_{1}^T$ et $(1,0)$ est encore le plus grand terme de l'ensemble d'indices du couple $(\alpha^t,\beta^t)$.
 Par d\'efinition de $t$, $(\nu^t,\mu^t)$ est construit \`a l'aide du proc\'ed\'e (a). Ce proc\'ed\'e cr\'ee le terme $\nu_{1} $ et des partitions $(\alpha^{t+1},\beta^{t+1})$. On a $\nu^t=\{\nu_{1}\}\sqcup\nu^{t+1}\sqcup\{0^{n_{t}-1-n_{t+1}}\}$ et $\mu^t=\mu^{t+1}$ pour un couple $(\nu^{t+1},\mu^{t+1})\in P(\alpha^{t+1},\beta^{t+1})$. En particulier, le terme $\mu_{t}$ de notre partition $\mu$ de d\'epart est \'egal \`a $\mu^{t+1}_{1}$ ou est nul si $m_{t+1}=0$. Toujours en appliquant 1.2(3) et 1.2(4), on a $\mu_{t}\leq \nu_{1}-\alpha_{1}^t$. Avec l'in\'egalit\'e $\alpha^T_{1}+A> B+s-sT$, on en d\'eduit $\nu_{1}+A>\mu_{t}+B+s-st$. Mais alors
 $$(6) \qquad \nu_{1}+A+S_{t-1}(\{\mu_{1}+B,...,\mu_{t-1}+B+2s-st\})> S_{t}(\{\mu_{1}+B,...,\mu_{t}+B+s-st\}).$$
 Le membre de gauche est inf\'erieur ou \'egal \`a $S_{t}(p\Lambda_{A,B;s}^{N,M}(\nu,\mu))$. D'apr\`es le lemme 1.4, c'est inf\'erieur ou \'egal \`a $S_{t}(p_{A,B;s}^{N,M}[\alpha,\beta])$. Or, puisque $c=t$, l'hypoth\`ese (2) nous dit que ce dernier terme est \'egal au  membre de droite de (6). Cela contredit cette in\'egalit\'e (6) qui est stricte. Supposons maintenant $(1,1)<_{I_{T}}(1,0)$. Puisque $(\chi,\xi)$ est l'\'el\'ement canonique de type (b), et que ce proc\'ed\'e n'est pas utilis\'e pour construire $(\chi^T,\xi^T)$, on a $\beta^T_{1}+B+s-sT< A$, a fortiori $\beta^t_{1}+B+s-st<A$. Par d\'efinition de $t$, $(\nu^t,\mu^t)$ est construit \`a l'aide du proc\'ed\'e (a). Si $(1,0)<_{I_{t}}(1,1)$, on voit  comme ci-dessus que $\mu_{t}\leq \nu_{1}-\alpha^t_{1}$. Si au contraire $(1,1)<_{I_{t}}(1,0)$, un calcul similaire montre que $\mu_{t}\leq \nu_{1}+\beta^t_{1}$. En tout cas, $\nu_{1}+A> \mu_{t}+B+s-st$. On obtient une contradiction comme ci-dessus. Cela prouve (5).
 
 Puisqu'on a suppos\'e $t=c$, (5) et 1.3(1) entra\^{\i}nent que les $c-1$ plus grands termes de $p\Lambda_{A,B;s}^{N,M}[\alpha,\beta]$ sont $\mu_{1}+B$,...,$\mu_{c-1}+B+2s-cs$ et
 $$S_{c}(p\Lambda_{A,B;s}^{N,M}[\alpha,\beta])=S_{c-1}(\{\mu_{1}+B,...,\mu_{c-1}+B+2s-cs\})+S_{1}(p\Lambda_{A,B+s-cs;s}^{N,M+1-c}[\alpha^c,\beta^c]).$$
 On a aussi
 $$S_{c}(\{\mu_{1}+B,...,\mu_{c}+B+s-cs\})=S_{c-1}(\{\mu_{1}+B,...,\mu_{c-1}+B+2s-cs\})+S_{1}(\{\mu^c_{1}+B+s-cs\}).$$
 D'apr\`es (2), on a donc
$$S_{1}(\{\mu^c_{1}+B+s-cs\})=  S_{1}(p\Lambda_{A,B+s-cs;s}^{N,M+1-c}[\alpha^c,\beta^c]).$$
Mais alors les couples $(\alpha^c,\beta^c)$ et $(\nu^c,\mu^c)$ satisfont les hypoth\`eses de l'\'enonc\'e pour le nouvel entier $c=1$. Ce cas a d\'ej\`a \'et\'e trait\'e: cela entra\^{\i}ne que $(\nu^c,\mu^c)\in P^{b}(\alpha^{c},\beta^c)$. Cela contredit la d\'efinition de $t$ et l'\'egalit\'e $t=c$. Cette contradiction ach\`eve la d\'emonstration. $\square$

\bigskip

\subsection{Relation entre $P(\alpha,\beta)$ et $P(\alpha^-,\beta)$}
  
Supposons $n\geq1$, soit $x\in {\mathbb N}$ et 
 soit $c\in {\mathbb N}$ avec  $c\geq1$.   Soit $(\underline{\nu},\underline{\mu})\in{\mathbb Z}^{n-1}\times {\mathbb Z}^m$.  On d\'efinit $(\nu,\mu)\in  {\mathbb Z}^n\times {\mathbb Z}^m$ par:
 
  si $c\leq m$, $\nu_{1}=x+\underline{\mu}_{c}$, $\nu_{i}=\underline{\nu}_{i-1}$ pour $i=2,...,n$, $\mu_{j}=\underline{\mu}_{j}$ pour $j=1,...,c-1$, $\mu_{j}=\underline{\mu}_{j+1}$ pour $j=c,...,m-1$, $\mu_{m}=0$;
  
  si $c\geq m+1$,  $\nu_{1}=x$, $\nu_{i}=\underline{\nu}_{i-1}$ pour $i=2,...,n$, $\mu_{j}=\underline{\mu}_{j}$ pour $j=1,...,m$.
  
  On d\'efinit note $\iota_{c,x}:{\mathbb Z}^{n-1}\times {\mathbb Z}^m\to {\mathbb Z}^n\times {\mathbb Z}^m$ l'application qui, \`a $(\underline{\nu},\underline{\mu})$, associe $(\nu,\mu)$. Remarquons que $\iota_{c,x}=\iota_{c',x}$ pour $c,c'\geq m+1$. 
  
  Soit maintenant $(\alpha,\beta)\in {\cal P}_{n}\times {\cal P}_{m}$ et supposons $(1,0)<_{I}(1,1)$. 

\ass{Lemme}{(i) Pour tout $(\nu,\mu)\in P(\alpha,\beta)$, il existe $(\underline{\nu},\underline{\mu})\in P(\alpha^-,\beta)$ et $c\in \{1,...,m+1\}$ de sorte que $(\nu,\mu)=\iota_{c,\alpha_{1}}(\underline{\nu},\underline{\mu})$.

(ii) Soit $c\in {\mathbb N}$ avec $c\geq1$ et soit  $(\underline{\nu},\underline{\mu})\in P^{b[c]}(\alpha^-,\beta)$.  Alors $\iota_{c,\alpha_{1}}(\underline{\nu},\underline{\mu})$ appartient \`a $P(\alpha,\beta)$.}

Preuve de (i).  Si  $m=0$, les ensembles $P(\alpha,\beta)$ et $P(\alpha^-,\beta)$ sont r\'eduits \`a l'\'el\'ement $(\alpha,\emptyset)$, resp. $(\alpha^-,\emptyset)$ et on a $\iota_{1,\alpha_{1}}(\alpha^-,\emptyset)=(\alpha,\emptyset)$.

On suppose d\'esormais  $m>0$ et on raisonne par r\'ecurrence sur $n+m$. Soit $(\nu,\mu)\in P(\alpha,\beta)$. Supposons que cet \'el\'ement soit construit \`a l'aide du proc\'ed\'e (a). Ce proc\'ed\'e cr\'ee   un \'el\'ement $\nu_{1 }$ et des partitions $(\alpha',\beta')$. On a $\nu=\{\nu_{1}\}\sqcup \nu'\sqcup\{0^{n'-n-1}\}$, $\mu=\mu'\sqcup\{0^{m'-m}\}$ pour un couple $(\nu',\mu')\in P(\alpha',\beta')$.    Puisque $m\geq1$, on peut appliquer le proc\'ed\'e (b) pour construire les \'el\'ements de $P^b(\alpha^-,\beta)$. D'apr\`es 2.1(1), ce proc\'ed\'e cr\'ee  un \'el\'ement $\underline{\mu}_{1}=\nu_{1}-\alpha_{1}$ et les m\^emes partitions $(\alpha',\beta')$. En cons\'equence, le couple $(\underline{\nu},\underline{\mu})=(\nu'\sqcup\{0^{n-1-n'}\}, \{\mu_{1}\}\sqcup \mu'\sqcup\{0^{m-1-m'}\})$ appartient \`a $P(\alpha^-,\beta)$. On v\'erifie que $(\nu,\mu)=\iota_{1,\alpha_{1}}(\underline{\nu},\underline{\mu})$, ce qui conclut.

Supposons au contraire que $(\nu,\mu)$ soit construit \`a l'aide du proc\'ed\'e (b).   Ce proc\'ed\'e cr\'ee   un \'el\'ement $\mu_{1} $ et des partitions $(\alpha',\beta')$. L'hypoth\`ese $(1,0)<_{I}(1,1)$ implique que   le plus grand terme de $\alpha'$ est $\alpha_{1}$ (ou plus pr\'ecis\'ement, que le terme $(1,0)$ de l'ensemble d'indices $I'$ associ\'e \`a $(\alpha',\beta')$  s'identifie \`a $(1,0)\in I$). On a $\nu=\nu'\sqcup\{0^{n-n'}\}$, $\mu=\{\mu_{1}\}\sqcup\mu'\sqcup\{0^{m-m'-1}\}$, pour un couple $(\nu',\mu')\in P(\alpha',\beta')$. Par hypoth\`ese de r\'ecurrence, il existe $c\in \{1,...,m'+1\}$ et $(\underline{\nu}',\underline{\mu}')\in P(\alpha^{'-},\beta')$ de sorte que $(\nu',\mu')=\iota_{c,\alpha_{1}}(\underline{\nu}',\underline{\mu}')$. On peut appliquer le proc\'ed\'e (b) pour construire les \'el\'ements de $P^b(\alpha^-,\beta)$. D'apr\`es 2.1(4), ce proc\'ed\'e cr\'ee le m\^eme \'el\'ement $\mu_{1}$ et les partitions $(\alpha^{'-},\beta')$. En cons\'equence, le couple $(\underline{\nu},\underline{\mu})=((\underline{\nu}'\sqcup\{0^{n-n'}\},\{\mu_{1}\}\sqcup\underline{\mu}'\sqcup\{0^{m-m'-1}\})$ appartient \`a $P(\alpha^-,\beta)$. Mais on voit que $(\alpha,\beta)=\iota_{c+1,\alpha_{1}}(\underline{\nu},\underline{\mu})$ et on a $c+1\leq m+1$ puisque $m'<m$. Cela ach\`eve la preuve de (i). 

Preuve de (ii). Le cas $m=0$ se traite comme ci-dessus. On suppose $m>0$.  Le proc\'ed\'e (b) est loisible pour d\'efinir l'ensemble $P^b(\alpha^-,\beta)$. Ce proc\'ed\'e cr\'ee un \'el\'ement $\mu_{1} $ et des partitions $(\alpha',\beta')\in {\cal P}_{n'}\times {\cal P}_{m'}$. 
Par d\'efinition de l'ensemble $P^{b[c]}(\alpha^-,\beta)$, on a $\underline{\nu}=\nu'\sqcup\{0^{n-1-n'}\}$, $\underline{\mu}=\{\mu_{1}\}\sqcup \mu'\sqcup\{0^{m-1-m'}\}$ pour un \'el\'ement $(\nu',\mu')\in P^{b[c-1]}(\alpha',\beta')$.   Supposons d'abord $c=1$. Parce que $(1,0)<_{I}(1,1)$ et d'apr\`es 2.1(1),  le proc\'ed\'e (a) appliqu\'e \`a $(\alpha,\beta)$ cr\'ee l'\'el\'ement $\nu_{1}=\alpha_{1}+\mu_{1}$ et les m\^emes partitions $(\alpha',\beta')$. Alors le couple $(\nu_{1}\sqcup\nu'\sqcup\{0^{n-1-n'}\},\mu')$ appartient \`a $P(\alpha,\beta)$. Or ce couple est \'egal \`a $\iota_{1,\alpha_{1}}(\underline{\nu},\underline{\mu})$. Supposons maintenant $c\geq2$. Toujours parce que $(1,0)<_{I}(1,1)$ et d'apr\`es 2.1(4), le proc\'ed\'e (b) appliqu\'e \`a $(\alpha,\beta)$ cr\'ee le m\^eme \'el\'ement $\mu_{1}$ et le couple de partitions $(\alpha'',\beta'')=(\{\alpha_{1}\}\sqcup \alpha',\beta')$. On a $\alpha'=\alpha^{''-}$. 
 0n applique le lemme par r\'ecurrence: l'\'el\'ement $(\underline{\nu}'',\underline{\mu}'')=\iota_{c-1,\alpha_{1}}(\nu',\mu')$ appartient \`a $P(\alpha'',\beta'')$. Alors, par le proc\'ed\'e (b), le couple $(\nu,\mu)=(\underline{\nu}''\sqcup\{0^{n-n''}\},\{\mu_{1}\}\sqcup\underline{\mu}''\sqcup\{0^{m-1-m''}\})$ appartient \`a $P(\alpha,\beta)$. On v\'erifie que les couples $(\nu,\mu)$ et
  $\iota_{c,\alpha_{1}}(\underline{\nu},\underline{\mu})$ sont tous deux \'egaux \`a
  
 $(\{\alpha_{1}+\mu'_{c-1},\nu'_{1} ,...,\nu'_{n'},0,...,0\}, \{\mu_{1},\mu'_{1},...,\mu'_{c-2},\mu'_{c},...,\mu'_{m'},0...,0\})$ si $c\leq m'+1$,
 
 $(\{\alpha_{1},\nu'_{1} ,...,\nu'_{n'},0,...,0\}, \{\mu_{1},\mu'_{1},...,\mu'_{m'},0...,0\})$ si $c> m'+1$.
  
  Donc $\iota_{c,\alpha_{1}}(\underline{\nu},\underline{\mu})=(\nu,\mu)$ appartient \`a $P(\alpha,\beta)$.
 Cela ach\`eve la d\'emonstration. $\square$
 
 \bigskip

\subsection{Nombres d'\'el\'ements "de type $b$"}

Soient $\alpha\in {\cal P}_{n}$, $\beta\in {\cal P}_{m}$. Introduisons l'\'el\'ement canonique de type (b) $(\nu,\mu)\in P_{A,B;s}(\alpha,\beta)$. Soit $k\in \{1,...,N+M\}$. Il existe des entiers $a,b\in{\mathbb N}$ tels que $a+b=k$, de sorte que l'ensemble des  $k$ plus grands \'el\'ements de $p_{A,B;s}^{N,M}[\alpha,\beta]=p\Lambda_{A,B;s}^{N,M}(\nu,\mu)$ soient la r\'eunion de celui des  $a$ plus grands termes de $(\nu\sqcup\{0^{N-n}\})+[A,A+s-sN]_{s}$ et de celui des $b$ plus grands termes de $(\mu\sqcup\{0^{M-m}\})+[B,B+s-sM]_{s}$. Les entiers $a$ et $b$ ne sont pas compl\`etement d\'etermin\'es car les deux partitions intervenant peuvent avoir des termes communs: si le $a$-i\`eme terme de la premi\`ere est \'egal au $(b+1)$-i\`eme terme de la seconde, le couple $(a-1,b+1)$ convient lui-aussi. On choisit $a$ et $b$ de sorte que $b$ soit le plus grand possible. On note alors $b_{A,B;s}^{N,M}(\alpha,\beta;k)=b$. 

\ass{Lemme}{(i) Il existe un entier $e\in {\mathbb N}$ tel que $b_{A,B+es/2;s}^{N,M}(\alpha,\beta;k)=inf(k,M)$ pour tout $k=1,...N+M$.

(ii) On a 
$$b_{A,B;s}^{N,M}(\alpha,\beta;k)\leq b_{A,B+s/2;s}^{N,M}(\alpha,\beta;k)\leq b_{A,B;s}^{N,M}(\alpha,\beta;k)+1$$
pour tout $k=1,...,N+M$.}

Preuve. Les termes des \'el\'ements de $P_{A,B;s}(\alpha,\beta)$  sont  des sommes de termes de $\alpha$ et de termes de $\beta$. Ils restent donc dans un ensemble fini d'entiers ind\'ependant de $A$ et $B$. Notons $c$ le plus grand \'el\'ement de cet ensemble. Choisissons  $e$ tel que $B+es/2+s-sM> c+A$.  Pour  l'\'el\'ement canonique de type (b) $(\nu,\mu)\in P_{A,B+es/2;s}(\alpha,\beta)$, il est clair que les termes de $(\mu\sqcup\{0^{M-m}\})+[B+es/2,B+es/2+s-sM]_{2}$ sont tous strictement sup\'erieurs aux termes de $(\nu\sqcup\{0^{N-n}\})+[A,A+s-sN]_{2}$. Alors $e$ v\'erifie la propri\'et\'e (i).

Fixons $k\in \{1,...,N+M\}$. On pose $b=b_{A,B;s}^{N,M}(\alpha,\beta;k)$, $\underline{b}=b_{A,B+s/2;s}^{N,M}(\alpha,\beta;k)$, $a=k-b$, $\underline{a}=k-\underline{b}$. L'assertion (ii) est \'evidente si $N=0$ (on a $b=\underline{b}=k$) ou $M=0$ (on a $b=\underline{b}=0$). On suppose $N,M\geq1$. 

Supposons d'abord $n=0$ ou $m=0$. On a $(\nu,\mu)=(\alpha,\beta)$.  Pour simplifier, posons $\alpha_{i}=0$ pour $i=n+1,...,N$ et $\beta_{i}=0$ pour $i=m+1,...,M$.   Les entiers $a$ et $b$ sont caract\'eris\'es par les propri\'et\'es

 $0\leq a\leq N$, $0\leq b\leq M$, $a+b=k$; 

(1)  les termes $\beta_{1}+B, \beta_{2}+B-s,...,\beta_{b}+B+s-sb$ sont sup\'erieurs ou \'egaux \`a $\alpha_{a+1}+A-sa,...,\alpha_{N}+A+s-sN$ ainsi qu'\`a $\beta_{b+1}+B-sb,...,\beta_{M}+B+s-sM$;
 
 (2) les termes $\alpha_{1}+A,\alpha_{2}+A-s,...,\alpha_{a}+A+s-sa$ sont strictement sup\'erieurs  \`a $\alpha_{a+1}+A-sa,...,\alpha_{N}+A+s-sN$ ainsi qu'\`a $\beta_{b+1}+B-sb,...,\beta_{M}+B+s-sM$.
 
 La propri\'et\'e (1) se conserve quand on remplace $B$ par $B+s/2$. Si la propri\'et\'e (2) se conserve \'egalement, on a $\underline{b}=b$ et $\underline{a}=a$. Supposons que la propri\'et\'e (2) ne se conserve pas. Cela entra\^{\i}ne $a\geq1$, $b<M$ et $\alpha_{a}+A+s-sa\leq\beta_{b+1}+B+s/2-sb$. Mais alors
 
  les termes $\beta_{1}+B+s/2, \beta_{2}+B-s/2,...,\beta_{b+1}+B+s/2-sb$ sont sup\'erieurs ou \'egaux \`a $\alpha_{a}+A-sa,...,\alpha_{N}+A+s-sN$ ainsi qu'\`a $\beta_{b+2}+B-s/2-sb,...,\beta_{M}+B+3s/2-sM$;
 
  les termes $\alpha_{1}+A,\alpha_{2}+A-s,...,\alpha_{a-1}+A+2s-sa$ sont strictement sup\'erieurs  \`a $\alpha_{a}+A+s-sa,...,\alpha_{N}+A+s-sN$ ainsi qu'\`a $\beta_{b+2}+B-s/2-sb,...,\beta_{M}+B+3s/2-sM$.
  
  D'o\`u $\underline{b}=b+1$ et $\underline{a}=a-1$.

On suppose maintenant $n,m\geq1$. Introduisons les \'el\'ements canoniques de type (b) $(\nu,\mu)\in P_{A,B;s}(\alpha,\beta)$ et $(\underline{\nu},\underline{\mu})\in P_{A,B+s/2;s}(\alpha,\beta)$. On utilise les notations habituelles pour $(\nu,\mu)$ et les m\^emes notations soulign\'ees pour $(\underline{\nu},\underline{\mu})$. Supposons d'abord que les \'el\'ements $(\nu,\mu)$ et $(\underline{\nu},\underline{\mu})$ soient construits \`a l'aide du m\^eme pr\'ec\'ed\'e, par  exemple le proc\'ed\'e (b). Ce proc\'ed\'e cr\'ee un terme $\mu_{1}$ et des partitions $(\alpha',\beta')$. On a $(\nu,\mu)=(\nu',\{\mu_{1}\}\sqcup \mu')$ o\`u $(\nu',\mu')$ est l'\'el\'ement canonique de type (b) de $P_{A,B-s;s}(\alpha',\beta')$ et $(\underline{\nu},\underline{\mu})=(\underline{\nu}',\{\mu_{1}\}\sqcup\underline{\mu}')$ o\`u $(\underline{\nu}',\underline{\mu}')$  est l'\'el\'ement canonique de type (b) de $P_{A,B-s/2;s}(\alpha',\beta')$. D'o\`u $p\Lambda_{A,B;s}^{N,M}(\nu,\mu)=\{\mu_{1}+B\}\sqcup p\Lambda_{A,B-s;s}^{N,M-1}(\nu',\mu')$ et $p\Lambda_{A,B+s/2;s}^{N,M}(\underline{\nu},\underline{\mu})=
\{\mu_{1}+B+s/2\}\sqcup 
p\Lambda_{A,B-s/2;s}^{N,M-1}(\underline{\nu}',\underline{\mu}')$. Si $k=1$, on voit que $b=\underline{b}=1$. Si $k\geq2$, on voit que $b=1+b_{A,B-s;s}^{N,M}(\alpha',\beta';k-1)$, $\underline{b}=1+b_{A,B-s/2;s}^{N,M-1}(\alpha',\beta';k-1)$. L'assertion \`a prouver s'ensuit par r\'ecurrence.

On suppose maintenant que les \'el\'ements $(\nu,\mu)$ et $(\underline{\nu},\underline{\mu})$ soient construits \`a l'aide de proc\'ed\'es diff\'erents. Supposons d'abord $(1,0)<_{I}(1,1)$. Si $(\nu,\mu)$ est construit \`a l'aide du proc\'ed\'e (b), on a  $\alpha_{1}+A\leq B$.  Mais alors $\alpha_{1}+A< B+s/2$ et $(\underline{\nu},\underline{\mu})$ est aussi construit \`a l'aide du proc\'ed\'e (b), contrairement \`a l'hypoth\`ese.  Donc $(\nu,\mu)$ est construit \`a l'aide du proc\'ed\'e (a). On a donc $\alpha_{1}+A>B$ (si on a \'egalit\'e, la construction de l'\'el\'ement canonique de type (b) privil\'egie le proc\'ed\'e (b)). L'\'el\'ement $(\underline{\nu},\underline{\mu})$ est construit \`a l'aide du proc\'ed\'e (b) donc $\alpha_{1}+A\leq B+s/2$. Le proc\'ed\'e (a) construit un \'el\'ement $\nu_{1} $ et des partitions $(\alpha',\beta')$. Puisque $(1,0)<_{I}(1,1)$, le proc\'ed\'e (b) construit d'apr\`es 1.2(4) l'\'el\'ement $\underline{\mu}_{1}= \nu_{1}-\alpha_{1}$ et les partitions $(\{\alpha_{1}\}\sqcup \alpha',\beta')$. On obtient $p\Lambda_{A,B;s}^{N,M}(\nu,\mu)=\{\nu_{1}+A\}\sqcup p\Lambda_{A-s,B;s}^{N-1,M}(\nu',\mu')$   o\`u $(\nu',\mu')$ est l'\'el\'ement canonique de type (b) de $P_{A-s,B;s}(\alpha',\beta')$et $p\Lambda_{A,B+s/2;s}^{N,M}(\underline{\nu},\underline{\mu})=\{\underline{\mu}_{1}+B+s/2\}\sqcup p\Lambda_{A,B-s/2;s}^{N,M-1}( \underline{\nu}',\underline{\mu}')$ o\`u $(\underline{\nu}',\underline{\mu}')$  est l'\'el\'ement canonique de type (b) de $P_{A,B-s/2;s}( \{\alpha_{1}\}\sqcup\alpha',\beta')$. Si $k=1$, on voit que $b=0$ et $\underline{b}=1$. Supposons $k\geq2$. Supposons d'abord $n', m'\geq1$. L'\'el\'ement $(1,0)$ reste le plus petit \'el\'ement de l'ensemble d'indices de $(\{\alpha_{1}\}\sqcup\alpha',\beta')$. L'in\'egalit\'e $\alpha_{1}+A>B$ entra\^{\i}ne $\alpha_{1}+A> B-s/2$. Donc l'\'el\'ement $(\underline{\nu}',\underline{\mu}')$ est construit \`a l'aide du proc\'ed\'e (a).   Ce proc\'ed\'e cr\'ee un terme 
$\underline{\nu}_{1} $ et des partitions $(\alpha'',\beta'')$. On obtient
$$(3) \qquad p\Lambda_{A,B+s/2;s}^{N,M}(\underline{\nu},\underline{\mu})=\{\underline{\mu}_{1}+B+s/2\}\sqcup \{\underline{\nu}_{1}+A\}\sqcup p\Lambda_{A-s,B-s/2;s}^{N-1,M-1}(\underline{\nu}'',\underline{\mu}''),$$
o\`u $(\underline{\nu}'',\underline{\mu}'')$ est l'\'el\'ement canonique de type (b) de $P_{A-s,B-s/2;s}(\alpha'',\beta'')$. Consid\'erons l'ensemble d'indices $I'$ relatif au couple $(\alpha',\beta')$. Si $(1,0)<_{I'}(1,1)$, l'\'egalit\'e $\alpha_{1}+A\leq B+s/2$ entra\^{\i}ne $\alpha'_{1}+A-s<B$ et $(\nu',\mu')$ est construit \`a l'aide du proc\'ed\'e de type (b). Si $(1,1)<_{I'}(1,0)$, la m\^eme in\'egalit\'e entra\^{\i}ne $\beta'_{1}+B> A-s$ et la m\^eme conclusion. D'apr\`es 2.1(1), ce proc\'ed\'e (b) cr\'ee l'\'el\'ement $\mu_{1}= \underline{\nu}_{1}-\alpha_{1}$ et les m\^emes partitions $(\alpha'',\beta'')$. On obtient
$$(4)\qquad p\Lambda_{A,B;s}^{N,M}(\nu,\mu)=\{\nu_{1}+A\}\sqcup \{\mu_{1}+B\}\sqcup p\Lambda_{A-s,B-s;s}^{N-1,M-1}(\nu'',\mu''),$$
o\`u $(\nu'',\mu'')$ est l'\'el\'ement canonique de type (b) de $P_{A-s,B-s}(\alpha'',\beta'')$. Remarquons que les termes de (3) et (4) sont en ordre d\'ecroissant, au sens renforc\'e qui intervient dans la d\'efinition de $b$ et $\underline{b}$: par exemple, dans (3), on a $\underline{\mu}_{1}+B+s/2\geq \underline{\nu}_{1}+A$ et $\underline{\nu}_{1}+A$ est strictement sup\'erieur aux termes suivants, ainsi qu'il r\'esulte de 1.3(6).  

On a suppos\'e $n',m'\geq1$.  Supposons maintenant $n'=0$ ou $m'=0$ et montrons que l'on a encore l'\'egalit\'e (3). 
 L'hypoth\`ese $n'\geq1$ n'a pas servi pour obtenir (3).  Si $m'=0$, seul le proc\'ed\'e (a) est autoris\'e pour construire les \'el\'ements de $P_{A,B-s/2}(\{\alpha_{1}\}\sqcup \alpha',\beta')$. Cela conduit \`a l'\'egalit\'e (3) o\`u on a simplement $\underline{\nu}_{1}=\alpha_{1}$ et $\underline{\mu}''=\emptyset$. Le plus grand terme de $p\Lambda_{A-s,B-s/2;s}^{N-1,M-1}(\underline{\nu}'',\underline{\mu}'')$ est soit $\underline{\nu}''_{1}+A-s$ (avec la convention $\underline{\nu}''_{1}=0$ si $n''=0$) soit $B-s/2$. Ces termes sont strictement inf\'erieurs \`a $\alpha_{1}+A$ en vertu de l'in\'egalit\'e  $\alpha_{1}+A>B$ et parce que $\underline{\nu}''_{1}$ est ici soit nul, soit un  terme de la partition $\alpha'$.  Notons que, si $n'=0$ ou $m'=0$, on a $\alpha''=\alpha'$ et $\beta''$ est \'egal soit \`a $\beta'=\emptyset$ si $m'=0$, soit \`a $\beta'$ priv\'e de son plus grand terme $\beta'_{1}$ si $m'\geq1$. 
Tournons-nous vers l'\'egalit\'e (4). La partition $p\Lambda_{A-s,B;s}^{N-1,M}(\alpha',\beta')$ est alors $\{\alpha'_{1}+A-s,...,\alpha'_{N-1}+A+s-sN\}\sqcup \{\beta'_{1}+B,...,\beta'_{M}+B+s-sM\}$ avec la convention $\alpha'_{i}=0$ si $i\geq n'+1$ et $\beta'_{i}=0$ si $i\geq m'+1$. Le plus grand terme de cette partition est $\beta'_{1}+B$ en vertu de l'in\'egalit\'e  $\alpha_{1}+A\leq B+s/2$ qui implique $\alpha'_{1}+A-s<B$. On obtient l'\'egalit\'e (4) o\`u $\mu_{1}=\beta'_{1}$ et o\`u  $(\alpha'',\beta'')$ est le couple d\'ecrit ci-dessus.

 On suppose maintenant $(1,1)<_{I}(1,0)$. Si $(\nu,\mu)$ est construit \`a l'aide du proc\'ed\'e (b), on a $\beta_{1}+B\geq A$ mais alors $\beta_{1}+B+s/2> A$ et $(\underline{\nu},\underline{\mu})$ est construit \`a l'aide du proc\'ed\'e (b) contrairement \`a l'hypoth\`ese. Donc $(\nu,\mu)$ est construit \`a l'aide du proc\'ed\'e (a). On a donc $\beta_{1}+B<A$ (si on a \'egalit\'e, la construction de l'\'el\'ement canonique de type (b) privil\'egie le proc\'ed\'e (b)). L'\'el\'ement $(\underline{\nu},\underline{\mu})$ est construit \`a l'aide du proc\'ed\'e (b) donc $\beta_{1}+B+s/2\geq A$.  Le proc\'ed\'e (b) construit un \'el\'ement $\underline{\mu}_{1} $ et des partitions $(\alpha',\beta')$. D'apr\`es 1.2(4), le proc\'ed\'e (a) construit un \'el\'ement $\nu_{1} =\underline{\mu}_{1}-\beta_{1}$ et les partitions $(\alpha',\{\beta_{1}\}\sqcup \beta')$. On obtient $p\Lambda_{A,B;s}^{N,M}(\nu,\mu)=\{\nu_{1}+A\}\sqcup p\Lambda_{A-s,B;s}^{N-1,M}(\nu',\mu')$   o\`u $(\nu',\mu')$ est l'\'el\'ement canonique de type (b) de $P_{A-s,B;s}(\alpha',\{\beta_{1}\}\sqcup \beta')$et $p\Lambda_{A,B+s/2;s}^{N,M}(\underline{\nu},\underline{\mu})=\{\underline{\mu}_{1}+B+s/2\}\sqcup p\Lambda_{A,B-s/2;s}^{N,M-1}( \underline{\nu}',\underline{\mu}')$ o\`u $(\underline{\nu}',\underline{\mu}')$  est l'\'el\'ement canonique de type (b) de $P_{A,B-s/2}( \alpha',\beta')$. Si $k=1$, on voit que $b=0$ et $\underline{b}=1$. Supposons $k\geq2$. Supposons d'abord $n', m'\geq1$. L'\'el\'ement $(1,1)$ reste le plus petit \'el\'ement de l'ensemble d'indices de $(\alpha',\{\beta_{1}\}\sqcup\beta')$. L'in\'egalit\'e $\beta_{1}+B+s/2\geq A$ entra\^{\i}ne $\beta_{1}+B> A-s$. Donc l'\'el\'ement $( \nu',\mu')$ est construit \`a l'aide du proc\'ed\'e (b).   Ce proc\'ed\'e cr\'ee un terme 
$ \mu_{1} $ et des partitions $(\alpha'',\beta'')$. On obtient l'\'egalit\'e (4) 
 o\`u $(\nu'',\mu'')$ est encore l'\'el\'ement canonique de type (b) de $P_{A-s,B-s;s}(\alpha'',\beta'')$. 
Consid\'erons l'ensemble d'indices $I'$ relatif au couple $(\alpha',\beta')$. Si $(1,0)<_{I'}(1,1)$, l'\'egalit\'e $\beta_{1}+B<A$ entra\^{\i}ne $\alpha'_{1}+A>B-s$ et $(\nu',\mu')$ est construit \`a l'aide du proc\'ed\'e de type (a). Si $(1,1)<_{I'}(1,0)$, la m\^eme \'egalit\'e entra\^{\i}ne $\beta'_{1}+B-s<A$ et la m\^eme conclusion. On applique 2.1(1) en \'echangeant les r\^oles des deux partitions: le proc\'ed\'e (a) cr\'ee l'\'el\'ement $\underline{\nu}_{1} =\mu_{1}-\beta_{1}$ et les m\^emes partitions $(\alpha'',\beta'')$. On obtient l'\'egalit\'e (3) 
 o\`u $(\underline{\nu}'',\underline{\mu}'')$ est encore  l'\'el\'ement canonique de type (b) de $P_{A-s,B-s/2}(\alpha'',\beta'')$.  Les termes de (3) et (4) sont encore en ordre d\'ecroissant, au sens renforc\'e qui intervient dans la d\'efinition de $b$ et $\underline{b}$. On a suppos\'e $n',m'\geq1$ mais, comme plus haut, on v\'erifie que les \'egalit\'es (3) et (4) et leurs propri\'et\'es de d\'ecroissance  restent vraies si on l\`eve cette hypoth\`ese.

Les \'egalit\'es (3) et (4) sont donc \'etablies en tout cas (pour $k\geq2$). Si $k=2$, on voit que $b=\underline{b}=1$. Si $k\geq3$, on a
$b=1+b_{A-s,B-s;s}^{N-1,M-1}(\alpha'',\beta'';k-2)$ et $\underline{b}=1+b_{A-s,B-s/2}^{N-1,M-1}(\alpha'',\beta'';k-2)$. Le r\'esultat cherch\'e s'ensuit par r\'ecurrence. $\square$

 \bigskip
 
 \subsection{Modification de l'ensemble d'indices}
   Soient $n',m'\in {\mathbb N}$ tels que $n'\leq n$ et $m'\leq m$. On suppose $(n'+1,0)>_{I}(m',1)$ si $n'<n$ et $m'\geq1$ et $(m'+1,1)>_{I}(n',0)$ si $m'<m$ et $n'\geq1$. Soient $\alpha'\in {\cal P}_{n'}$, $\beta'\in {\cal P}_{m'}$, posons $\alpha=\alpha'\sqcup\{0^{n-n'}\}$, $\beta=\beta'\sqcup\{0^{m-m'}\}$. L'ensemble d'indices $I'=(\{1,...,n'\}\times\{0\})\times (\{1,...,m'\}\times\{1\})$ est un sous-ensemble de $I$, on le munit de l'ordre induit. Cela permet de d\'efinir les ensembles $P(\alpha',\beta')$ et $P_{A,B;s}(\alpha',\beta')$. 
   
   \ass{Lemme}{Sous ces hypoth\`eses, on a
   $$P(\alpha,\beta)=\{(\nu'\sqcup\{0^{n-n'}\},\mu'\sqcup\{0^{m-m'}\}; (\nu',\mu')\in P(\alpha',\beta')\},$$ 
$$P_{A,B;s}(\alpha,\beta)=\{(\nu'\sqcup\{0^{n-n'}\},\mu'\sqcup\{0^{m-m'}\}); (\nu',\mu')\in P_{A,B;s}(\alpha',\beta')\}.$$}

Preuve. Si $n=0$ ou $m=0$, on a $P(\alpha,\beta)=P_{A,B;s}(\alpha,\beta)=\{(\alpha,\beta)\}$ et $P(\alpha',\beta')=P_{A,B;s}(\alpha',\beta')=\{(\alpha',\beta')\}$ et le r\'esultat s'ensuit. On suppose $n,m\geq1$. On  raisonne par r\'ecurrence sur  $n+m-n'-m'$. Supposons que l'on ait trait\'e le cas o\`u cet entier vaut $1$, et traitons le cas o\`u il est strictement sup\'erieur \`a $1$.  Supposons par exemple $(n,0)<_{I}(m,1)$. Avec les hypoth\`eses pos\'ees avant l'\'enonc\'e, cela entra\^{\i}ne $m'<m$. On pose $\beta''=\beta\sqcup\{0^{m-m'-1}\}$. On voit alors que les paires $(\alpha',\beta')$ et $(\alpha,\beta'')$ v\'erifient les m\^emes conditions que ci-dessus, $m$ \'etant remplac\'e par $m-1$. Et que les paires $(\alpha,\beta'')$ et $(\alpha,\beta)$ v\'erifient aussi ces conditions, $(n',m')$ \'etant remplac\'e par $(n,m-1)$. Par r\'ecurrence, on peut appliquer le lemme \`a chacun des deux couples de  paires et le r\'esultat s'ensuit. Il nous reste \`a traiter le cas o\`u $n+m-n'-m'=1$. On ne perd rien \`a supposer $n'=n$ et $m'=m-1$.  Avec les hypoth\`eses pos\'ees avant l'\'enonc\'e, cela entra\^{\i}ne $(n,0)<_{I}(m,1)$. Supposons d'abord $m=1$, donc $m'=0$. On a encore $P(\alpha,\beta')=P_{A,B;s}(\alpha,\beta')=\{(\alpha,\beta')\}$. D'apr\`es 1.2(5), on a aussi $P(\alpha,\beta)=P_{A,B;s}(\alpha,\beta)=\{(\alpha,\beta)\}$ et le r\'esultat s'ensuit. Supposons maintenant $m\geq2$. Le proc\'ed\'e (a) appliqu\'e \`a $(\alpha,\beta')$ cr\'ee un \'el\'ement $\nu_{1}=\alpha_{1}+\beta'_{b_{1}}+\alpha_{a_{2}}+...$ et des partitions, notons-les ici $(\underline{\alpha},\underline{\beta}')\in {\cal P}_{\underline{n}}\times {\cal P}_{\underline{m}'}$.   Le terme $(m,1)$ est le plus grand terme de l'ensemble d'indices $I$. On voit alors que

(1)   si le dernier terme de $\nu_{1}$ est $\alpha_{a_{t}}$, le proc\'ed\'e (a) appliqu\'e \`a $(\alpha,\beta)$ cr\'ee le terme $\alpha_{1}+\beta'_{b_{1}}+\alpha_{a_{2}}+...+\alpha_{a_{t}}+\beta_{m}=\nu_{1}$ (puisque $\beta_{m}=0$) et les m\^emes partitions $(\underline{\alpha},\underline{\beta}')$;

(2) si le dernier terme de $\nu_{1}$ est $\beta'_{b_{t}}$, le proc\'ed\'e (a) appliqu\'e \`a $(\alpha,\beta)$ cr\'ee le m\^eme terme $\nu_{1}$ et les partitions $(\underline{\alpha},\underline{\beta}'\sqcup\{0\})$.

Dans le cas (1), on a
$$(3) \qquad P^{a}(\alpha,\beta')=\{(\nu_{1}\sqcup\underline{\nu}\sqcup\{0^{n-\underline{n}-1}\},\underline{\mu}'\sqcup\{0^{m'-\underline{m}'}\}); (\underline{\nu},\underline{\mu}')\in P(\underline{\alpha},\underline{\beta}')\};$$
$$P^{a}(\alpha,\beta)=\{(\nu_{1}\sqcup\underline{\nu}\sqcup\{0^{n-\underline{n}-1}\},\underline{\mu}'\sqcup\{0^{m-\underline{m}'}\}); (\underline{\nu},\underline{\mu}')\in P(\underline{\alpha},\underline{\beta}')\}.$$
On en d\'eduit l'analogue de la premi\`ere \'egalit\'e de l'\'enonc\'e o\`u $P(\alpha,\beta)$ et $P(\alpha,\beta')$ sont remplac\'es par $P^{a}(\alpha,\beta)$ et $P^{a}(\alpha,\beta')$. Dans le cas (2), on a encore (3) et 
$$P^{a}(\alpha,\beta)=\{(\nu_{1}\sqcup\underline{\nu}\sqcup\{0^{n-\underline{n}-1}\},\underline{\mu}\sqcup\{0^{m-\underline{m}'-1}\}); (\underline{\nu},\underline{\mu})\in P(\underline{\alpha},\underline{\beta})\},$$
o\`u $\underline{\beta}=\underline{\beta}'\sqcup\{0\}$.
Mais $(\underline{\alpha},\underline{\beta})$ et $(\underline{\alpha},\underline{\beta}')$ sont dans la m\^eme situation que $(\alpha,\beta)$ et $(\alpha,\beta')$. En raisonnant par r\'ecurrence sur $n+m$, on peut appliquer l'\'enonc\'e pour d\'ecrire $P(\underline{\alpha},\underline{\beta})$ \`a l'aide de $P(\underline{\alpha},\underline{\beta}')$. Les \'egalit\'es ci-dessus conduisent alors \`a l'analogue de la premi\`ere \'egalit\'e de l'\'enonc\'e o\`u $P(\alpha,\beta)$ et $P(\alpha,\beta')$ sont remplac\'es par $P^{a}(\alpha,\beta)$ et $P^{a}(\alpha,\beta')$. Cette analogue est donc d\'emontr\'ee dans les deux cas (1) et (2). On a raisonn\'e avec le proc\'ed\'e (a). Le m\^eme raisonnement s'applique au pr\'ec\'ed\'e (b), avec une conclusion similaire. La r\'eunion de ces deux conclusions entra\^{\i}ne la premi\`ere \'egalit\'e de l'\'enonc\'e.  Pour la seconde, le fait que le proc\'ed\'e (a) ou (b) soit autoris\'e pour construire les \'el\'ements de $P_{A,B;s}(\alpha,\beta')$ ou $P_{A,B;s}(\alpha,\beta)$ ne d\'epend que de $\alpha_{1}$ et $\beta_{1}$. Ces termes sont les m\^emes pous nos deux couples de partitions (rappelons que $n\geq1$ et $m'=m-1\geq1$) donc les m\^emes proc\'ed\'es sont autoris\'es pour ces deux couples. Le raisonnement se poursuit alors comme ci-dessus et conduit \`a la conclusion. $\square$

\bigskip

\section{Trois lemmes de majorations}

  \subsection{Premier  lemme } 
  
  On suppose $n\geq1$ et $(1,0)<_{I}(1,1)$ si $m\geq1$. Soient $\alpha\in {\cal P}_{n}$ et $\beta\in {\cal P}_{m}$.

 \ass{Lemme}{ Pour tout $k=1,...,N+M$, 
 on a l'in\'egalit\'e 
 $$S_{k}(p_{A,B;s}^{N,M}[\alpha^-,\beta])\leq S_{k}(p_{A,B;s}^{N,M}[\alpha,\beta])\leq S_{k}(p_{A,B;s}^{N,M}[\alpha^-,\beta])+\alpha_{1}+\beta_{1},$$
 avec la convention $\beta_{1}=0$ si $m=0$.}
 
 Il y a un lemme sym\'etrique en \'echangeant les r\^oles de $\alpha$ et $\beta$. On doit alors supposer $m\geq1$ et $(1,1)<_{I}(1,0)$. On d\'efinit la partition $\beta^-$ et, dans l'\'enonc\'e ci-dessus, on remplace le couple $(\alpha^-,\beta)$ par $(\alpha,\beta^-)$. La d\'emonstration de ces lemmes sera donn\'ee en 3.3.
 
 \bigskip
 
 \subsection{Deuxi\`eme lemme }
 
   \ass{Lemme}{On suppose $N>n$ et $M>m$. Soient $\alpha\in {\cal P}_{n}$ et $\beta\in {\cal P}_{m}$. Pour tout $k\in\{1,...,N+M-1\}$, on a l'in\'egalit\'e
  $$S_{k}(p_{A,B-s;s}^{N,M-1}[\alpha,\beta])\leq S_{k}(p_{A-s,B;s}^{N-1,M}[\alpha,\beta])+ sup(\alpha_{1}+A-B,0),$$
  avec la convention $\alpha_{1}=0$ si $n=0$.} 
  
 L\`a encore, il y a un sym\'etrique de ce lemme obtenu en \'echangeant les r\^oles de $\alpha$ et $\beta$ et, simultan\'ement, ceux de $A$ et $B$ ainsi que de $N$ et $M$. La d\'emonstration sera donn\'ee en 3.4.
 
 Pour d\'emontrer ce lemme, le lemme pr\'ec\'edent et leurs sym\'etriques, on raisonnne  par r\'ecurrence sur $n+m$. C'est-\`a-dire que l'on v\'erifiera les quatre lemmes dans le cas $n=m=0$ (en fait dans les cas $n=0$ ou $m=0$). Ensuite, on supposera que  $n+m\geq1$ et que les quatre lemmes sont v\'erifi\'es pour des couples $(n',m')$ tels que $n'+m'<n+m$.  On ne fera les d\'emonstrations que des lemmes \'enonc\'es, celles des  sym\'etriques \'etant \'evidemment analogues. Pour simplifier l'\'ecriture, on adopte  dans ces d\'emonstrations la convention suivante: pour une partition, par exemple $\alpha\in {\cal P}_{n}$, on pose $\alpha_{i}=0$ pour $i>n$. 
  
  \bigskip
  
  \subsection{Preuve du lemme 3.1}

 Si $m=0$, on a
 $$p_{A,B;s}^{N,M}[\alpha,\beta])=\{\alpha_{1}+A,...,\alpha_{N}+A+s-sN\}\sqcup \{B,...,B+s-sM\}$$
 et
 $$p_{A,B;s}^{N,M}[\alpha^-,\beta])=\{\alpha_{2}+A,...,\alpha_{N+1}+A+s-sN\}\sqcup \{B,...,B+s-sM\}.$$
 La premi\`ere in\'egalit\'e de l'\'enonc\'e r\'esulte des in\'egalit\'es $\alpha_{i+1} \leq \alpha_{i}$ pour $i=1,...,N$. Pour la seconde, il existe $a\in \{0,...,N\}$ et $b\in \{0,...,M\}$ tels que $a+b=k$ et
 $$S_{k}(p_{A,B;s}^{N,M}[\alpha,\beta])=S_{a}(\{\alpha_{1}+A,...,\alpha_{N}+A+s-sN\})+S_{b}(\{B,...,B+s-sM\}).$$
 On a 
 $$S_{a}(\{\alpha_{1}+A,...,\alpha_{N}+A+s-sN\})=\alpha_{1}-\alpha_{a+1}+S_{a}(\{\alpha_{2}+A,...,\alpha_{N+1}+A+s-sN\}).$$
 D'o\`u
 $$S_{k}(p_{A,B;s}^{N,M}[\alpha,\beta])\leq \alpha_{1}+S_{a}(\{\alpha_{2}+A,...,\alpha_{N+1}+A+s-sN\})+S_{b}(\{B,...,B+s-sM\})$$
 $$\leq \alpha_{1}+S_{k}(p_{A,B;s}^{N,M}[\alpha^-,\beta]).$$
 
   On suppose $m\geq1$. 
  Supposons d'abord $\alpha_{1}+A\leq B$. Le proc\'ed\'e (b) est autoris\'e pour construire les \'el\'ements de $P_{A,B;s}(\alpha,\beta)$.  D'apr\`es 2.1(5), il en est de m\^eme pour l'ensemble $P_{A,B;s}(\alpha^-,\beta)$. Le proc\'ed\'e (b) appliqu\'e \`a $\alpha,\beta$ construit  un terme $\mu_{1}$ et des partitions $(\alpha',\beta')\in {\cal P}_{n'}\times {\cal P}_{m'}$.  D'apr\`es 2.1(4), le m\^eme proc\'ed\'e appliqu\'e \`a $(\alpha^-,\beta)$ construit    le m\^eme terme $\mu_{1}$ et le couple de partitions $({\alpha'}^-,\beta')\in {\cal P}_{n'-1}\times {\cal P}_{m'}$. Il y a donc un \'el\'ement de $P_{A,B;s}(\alpha,\beta)$ de la forme $(\nu',\{\mu_{1}\}\sqcup \mu')$ pour un \'el\'ement $(\nu',\mu')\in P_{A,B-s;s}(\alpha',\beta')$, et un  \'el\'ement de $P_{A,B;s}(\alpha^-,b)$  de la forme $(\nu'',\{\mu_{1}\}\sqcup \mu'')$ pour un \'el\'ement $(\nu'',\mu'')\in P_{A,B-s;s}(\alpha^{'-},\beta')$. 
 On obtient $p_{A,B;s}^{N,M}[\alpha,\beta]=\{\mu_{1}+B\}\sqcup p_{A,B-s;s}^{N,M-1}[\alpha',\beta'] $ et $p_{A,B;s}^{N,M}[\alpha^-,\beta]=\{\mu_{1}+B\}\sqcup p_{A,B-s;s}^{N,M-1}[\alpha^{'-},\beta'] $. Les in\'egalit\'es \`a prouver sont alors claires si $k=1$. Si $k\geq2$, on a par r\'ecurrence
 $$S_{k-1} (p_{A,B-s;s}^{N,M-1}[\alpha^{'-},\beta'])\leq S_{k-1}(  p_{A,B-s;s}^{N,M-1}[\alpha',\beta'] )\leq 
 S_{k-1} (p_{A,B-s;s}^{N,M-1}[\alpha^{'-},\beta'])+\alpha'_{1}+\beta'_{1}.$$
 Les in\'egalit\'es de l'\'enonc\'e s'ensuivent puisque $\alpha'_{1}\leq \alpha_{1}$ et $\beta'_{1}\leq \beta_{1}$.

 Supposons maintenant $\alpha_{1}+A> B$. On construit un \'el\'ement de $P_{A,B;s}(\alpha,\beta)$ par le proc\'ed\'e (a). Celui-ci nous fournit   un terme $\nu_{1} $ et des partitions $(\alpha',\beta')\in {\cal P}_{n'}\times {\cal P}_{m'}$. Un \'el\'ement de $P_{A,B;s}(\alpha,\beta)$ est de la forme $(\{\nu_{1}\}\sqcup \nu',\mu')$ pour un \'el\'ement $(\nu',\mu')\in P_{A-s,B;s}(\alpha',\beta')$ d'o\`u
 $$(1) \qquad p^{N,M}_{A,B;s}[\alpha,\beta]=\{\nu_{1}+A\}\sqcup p^{N-1,M}_{A-s,B;s}[\alpha',\beta'] .$$
 
 Supposons $n=1$.  Alors $\nu_{1}=\alpha_{1}+\beta_{1}$, $\alpha'=\alpha^-=\emptyset$ et $\beta'=\beta^-$. On a $n'=0$ d'o\`u 
 $$p^{N-1,M}_{A-s,B;s}[\alpha',\beta']=[A-s,A+s-sN]_{s}\sqcup((\beta^-\sqcup\{0^{M-m+1}\}+[B,B+s-sM]_{s}).$$
 D'o\`u
 $$p_{A,B;s}^{N,M}[\alpha,\beta]=\{\alpha_{1}+\beta_{1}+A,A-s,...,A+s-sN\}\sqcup \{\beta_{2}+B,\beta_{3}+B-s,...,\beta_{M+1}+B+s-sM\}.$$
   On a aussi
  $$p_{A,B;s}^{N,M}[\alpha^-,\beta]= \{A,A-s,...,A+s-sN\}\sqcup\{\beta_{1}+B,\beta_{2}+B-s,...,\beta_{M}+B+s-sM\}.$$ 
 On peut fixer des entiers  $a\in \{0,...,N\}$ et $b\in \{0,...,M\}$ tels que $a+b=k$ et
 $$S_{k}(p_{A,B;s}^{N,M}[\alpha^-,\beta])=S_{a}( \{A,A-s,...,A+s-sN\})+S_{b}(\{\beta_{1}+B,\beta_{2}+B-s,...,\beta_{M}+B+s-sM\}).$$
 Supposons d'abord $a\geq1$. Alors 
 $$S_{b}(\{\beta_{1}+B,\beta_{2}+B-s,...,\beta_{M}+B+s-sM\})=$$
 $$S_{b}(\{\beta_{2}+B,\beta_{3}+B-s,...,\beta_{M+1}+B+s-sM\}) +\beta_{1}-\beta_{b+1}$$
 et
 $$S_{a}( \{A,A-s,...,A+s-sN\})+\beta_{1}\leq S_{a}( \{\alpha_{1}+\beta_{1}+A,A-s,...,A+s-sN\}).$$
 D'o\`u
 $$S_{k}(p_{A,B;s}^{N,M}[\alpha^-,\beta])\leq S_{a}( \{\alpha_{1}+\beta_{1}+A,A-s,...,A+s-sN\})$$
 $$+S_{b}(\{\beta_{2}+B,\beta_{3}+B-s,...,\beta_{M+1}+B+s-sM\})\leq S_{k}[p_{A,B;s}^{N,M}[\alpha,\beta]).$$
 Supposons maintenant $a=0$, donc $b=k\geq1$.Alors
 $$S_{k}(\{\beta_{1}+B,\beta_{2}+B-s,...,\beta_{M}+B+s-sM\})=$$
 $$S_{k-1}(\{\beta_{2}+B,\beta_{3}+B-s,...,\beta_{M+1}+B+s-sM\})+\beta_{1}+B+s-sk.$$
 D'apr\`es l'hypoth\`ese $\alpha_{1}+A>B$, on a 
 $$\beta_{1}+B+s-sk<\beta_{1}+\alpha_{1}+A=S_{1}(\{\alpha_{1}+\beta_{1}+A,A-s,...,A+s-sN\}).$$
 D'o\`u
 $$S_{k}(p_{A,B;s}^{N,M}[\alpha^-,\beta])\leq S_{1}(\{\alpha_{1}+\beta_{1}+A,A-s,...,A+s-sN\})$$
 $$+S_{k-1}(\{\beta_{2}+B,\beta_{3}+B-s,...,\beta_{M+1}+B+s-sM\})\leq S_{k}[p_{A,B;s}^{N,M}[\alpha,\beta]).$$
 C'est la premi\`ere in\'egalit\'e de l'\'enonc\'e. 
On peut fixer de nouveaux  entiers  $a\in \{0,...,N\}$ et $b\in \{0,...,M\}$ tels que $a+b=k$ et
$$S_{k}(p_{A,B;s}^{N,M}[\alpha,\beta])=S_{a}( \{\alpha_{1}+\beta_{1}+A,A-s,...,A+s-sN\})$$
$$+S_{b}(\{\beta_{2}+B,\beta_{3}+B-s,...,\beta_{M+1}+B+s-sM\}).$$
Il est clair que  
$$S_{a}( \{\alpha_{1}+\beta_{1}+A,A-s,...,A+s-sN\})\leq \alpha_{1}+\beta_{1}+S_{a}( \{A,A-s,...,A+s-sN\}),$$
tandis que
$$S_{b}(\{\beta_{2}+B,\beta_{3}+B-s,...,\beta_{M+1}+B+s-sM\})\leq S_{b}(\{\beta_{1}+B,\beta_{2}+B-s,...,\beta_{M}+B+s-sM\}).$$
D'o\`u
$$S_{k}(p_{A,B;s}^{N,M}[\alpha,\beta])\leq \alpha_{1}+\beta_{1}+S_{a}( \{A,A-s,...,A+s-sN\})$$
$$+S_{b}(\{\beta_{1}+B,\beta_{2}+B-s,...,\beta_{M}+B+s-sM\}\leq \alpha_{1}+\beta_{1}+S_{k}(p_{A,B;s}^{N,M}[\alpha^-,\beta]).$$
C'est la seconde in\'egalit\'e de l'\'enonc\'e.

 On suppose maintenant $n\geq2$. Supposons d'abord $(1,0)<_{I^-}(1,1)$, c'est-\`a-dire $(2,0)<_{I}(1,1)$. Cela entra\^{\i}ne que le premier terme de $\alpha'$ est $\alpha_{2}=\alpha_{1}^-$. Supposons  $\alpha^-_{1}+A\geq B$. On construit un \'el\'ement de $P_{A,B;s}(\alpha^-,\beta)$ par le proc\'ed\'e (a).  D'apr\`es 2.1(2), ce proc\'ed\'e construit le terme $\nu_{1}-\alpha_{1}+\alpha_{2}$, et les partitions   sont $(\alpha^{'-},\beta')$. Un \'el\'ement de $P_{A,B;s}(\alpha^-,\beta)$ est de la forme $(\{\nu_{1}-\alpha_{1}+\alpha_{2}\}\sqcup \nu'',\mu'')$ pour un \'el\'ement $(\nu'',\mu'')\in P_{A-s,B;s}(\alpha^{'-},\beta')$, d'o\`u
 $$p_{A,B;s}^{N,M}[\alpha^-,\beta]=\{\nu_{1}-\alpha_{1}+\alpha_{2}+A\}\sqcup p_{A-s,B;s}^{N-1,M}[\alpha^{'-},\beta'] .$$
  En comparant avec (1), les in\'egalit\'es de l'\'enonc\'e s'ensuivent si $k=1$. Supposons $k\geq2$.   
  Par r\'ecurrence, on a
  $$S_{k-1}(p_{A-s,B;s}^{N-1,M}[\alpha^{'-},\beta'])\leq S_{k-1}(p_{A-s,B;s}^{N-1,M}[\alpha^{'},\beta'])\leq S_{k-1}(p_{A-s,B;s}^{N-1,M}[\alpha^{'-},\beta'])+\alpha'_{1}+\beta'_{1}.$$
   D'o\`u
   $$S_{k}(p_{A,B;s}^{N,M}[\alpha^-,\beta])=\nu_{1}-\alpha_{1}+\alpha_{2}+A+ S_{k-1}(p_{A-s,B;s}^{N-1,M}[\alpha^{'-},\beta'])$$
   $$\leq \nu_{1}+A+S_{k-1}(p_{A-s,B;s}^{N-1,M}[\alpha^{'},\beta'])=S_{k}(p_{A,B;s}^{N,M}[\alpha,\beta]).$$
   Et
   $$S_{k}(p_{A,B;s}^{N,M}[\alpha,\beta])=\nu_{1}+A+S_{k-1}(p_{A-s,B;s}^{N-1,M}[\alpha^{'},\beta'])\leq \nu_{1}+A+\alpha'_{1}+\beta'_{1}+S_{k-1}(p_{A-s,B;s}^{N-1,M}[\alpha^{'-},\beta']).$$
   On se rappelle que $\alpha'_{1}=\alpha_{2}$ et on a $\beta'_{1}\leq \beta_{1}$. D'o\`u
  $$S_{k}(p_{A,B;s}^{N,M}[\alpha,\beta])\leq \nu_{1}+A+\alpha_{2}+\beta_{1}+  S_{k-1}(p_{A-s,B;s}^{N-1,M}[\alpha^{'-},\beta'])$$
  $$=\alpha_{1}+\beta_{1}+(\nu_{1}-\alpha_{1}+\alpha_{2}+A)+S_{k-1}(p_{A-s,B;s}^{N-1,M}[\alpha^{'-},\beta'])=\alpha_{1}+\beta_{1}+
  S_{k}(p_{A,B;s}^{N,M}[\alpha^-,\beta]).$$

 On suppose maintenant $(1,0)<_{I^-}(1,1)$  mais $\alpha^-_{1}+A<B$. On construit un \'el\'ement de $P_{A,B;s}(\alpha^-,\beta)$ par le proc\'ed\'e (b).  D'apr\`es 2.1(1), ce proc\'ed\'e construit l'\'el\'ement $\mu_{1}=\nu_{1}-\alpha_{1}$ et les m\^emes partitions $(\alpha',\beta')$. D'o\`u
$$(2) \qquad p_{A,B;s}^{N,M}[\alpha^-,\beta]=\{\nu_{1}-\alpha_{1}+B\}\sqcup p_{A,B-s;s}^{N,M-1}[\alpha',\beta'].$$
  
Interrompons notre raisonnement pour consid\'erer le cas o\`u $(1,0)>_{I^-}(1,1)$, c'est-\`a-dire $(1,1)<_{I}(2,0)$  et o\`u $\beta_{1}+B\geq A$.   On construit un \'el\'ement de $P_{A,B}(\alpha^-,\beta)$ par le proc\'ed\'e (b), qui construit encore   l'\'el\'ement $\mu_{1} =\nu_{1}-\alpha_{1}$ et les partitions $(\alpha',\beta')$. D'o\`u encore (2). On traite ensemble les deux situations. 

Si $k=1$, on doit montrer que $\nu_{1}-\alpha_{1}+B\leq \nu_{1}+A\leq \nu_{1}-\alpha_{1}+B+\alpha_{1}+\beta_{1}=\nu_{1}+\beta_{1}+B$. La premi\`ere in\'egalit\'e r\'esulte de 
$\alpha_{1}+A> B$. La seconde n'est autre que $A\leq \beta_{1}+B$. Or dans les deux cas consid\'er\'es ci-dessus, on a soit  directement cette in\'egalit\'e, soit $\alpha^-_{1}+A<B$ mais alors $A\leq \alpha^-_{1}+A<B\leq \beta_{1}+B$. Cela conclut.

Si $k\geq2$, le lemme 3.2 appliqu\'e par r\'ecurrence nous dit que
$$S_{k-1}(p_{A,B-s;s}^{N,M-1}[\alpha',\beta'])\leq S_{k-1}(p_{A-s,B;s}^{N-1,M}[\alpha',\beta']) +sup(\alpha'_{1}+A-B,0).$$
 Toujours d'apr\`es l'hypoth\`ese $\alpha_{1}+A> B$ et parce que $\alpha'_{1}\leq \alpha_{1}$, on a
$sup(\alpha'_{1}+A-B,0)\leq \alpha_{1}+A-B$. Alors
$$S_{k}(p_{A,B;s}^{N,M}[\alpha^-,\beta])=\nu_{1}-\alpha_{1}+B+S_{k-1}(p_{A,B-s;s}^{N,M-1}[\alpha',\beta'])\leq \nu_{1}+A+ S_{k-1}(p_{A-s,B;s}^{N-1,M}[\alpha',\beta'])$$
$$=S_{k}(p_{A,B;s}^{N,M}[\alpha,\beta]).$$
Le sym\'etrique du m\^eme lemme 3.2 nous dit que
$$S_{k-1}(p_{A-s,B;s}^{N-1,M}[\alpha',\beta'])\leq S_{k-1}(p_{A,B-s;s}^{N,M-1}[\alpha',\beta'])+sup(\beta'_{1}+B-A,0).$$
Comme on l'a vu ci-dessus, on a ici $\beta_{1}+B-A\geq0$. On a aussi $\beta'_{1}\leq \beta_{1}$, donc $sup(\beta'_{1}+B-A,0)\leq \beta_{1}+B-A$. D'o\`u
$$S_{k}(p_{A,B;s}^{N,M}[\alpha,\beta])=\nu_{1}+A+ S_{k-1}(p_{A-s,B;s}^{N-1,M}[\alpha',\beta'])\leq \nu_{1}+\beta_{1}+B+S_{k-1}(p_{A,B-s;s}^{N,M-1}[\alpha',\beta'])$$
$$=\alpha_{1}+\beta_{1}+(\nu_{1}-\alpha_{1}+B)+S_{k-1}(p_{A,B-s;s}^{N,M-1}[\alpha',\beta'])=\alpha_{1}+\beta_{1}+
S_{k}(p_{A,B;s}^{N,M}[\alpha^-,\beta]).$$
Cela conclut.

Il reste le cas o\`u  $(1,1)<_{I}(2,0)$  mais $\beta_{1}+B< A$.   On construit un \'el\'ement de $P_{A,B}(\alpha^-,\beta)$ par le proc\'ed\'e (a). D'apr\`es 2.1(3), ce proc\'ed\'e cr\'ee l'\'el\'ement  $\nu_{1}-\alpha_{1}-\beta_{1}$, les partitions  $(\alpha',\{\beta_{1}\}\sqcup \beta')$. On obtient
$$(3) \qquad p_{A,B;s}^{N,M}[\alpha^-,\beta]=\{\nu_{1}-\alpha_{1}-\beta_{1}+A\}\sqcup p_{A-s,B;s}^{N-1,M}[\alpha',\{\beta_{1}\}\sqcup \beta'].$$
Si $k=1$, les in\'egalit\'es de l'\'enonc\'e r\'esultent de cette formule et de (1). Supposons $k\geq2$. 
 Posons $\alpha^1=\alpha'$, $\beta^1=,\{\beta_{1}\}\sqcup \beta'$,d'o\`u $\beta'=\beta^{1-}$. L'ensemble d'indices $I_{1}=(\{1,...,n_{1}\}\times \{0\})\cup(\{1,...,m_{1}\}\times \{1\})$ est par construction plong\'e dans $I$ et h\'erite d'un ordre $<_{I_{1}}$ pour lequel $(1,1)<_{I_{1}}(1,0)$. On peut donc appliquer le sym\'etrique de notre lemme 3.1 au couple de  partitions $(\alpha^1,\beta^1)$: on a
 $$S_{k-1}(p_{A-s,B;s}^{N-1,M}[\alpha', \beta'])\leq S_{k-1}(p_{A-s,B;s}^{N-1,M}[\alpha',\{\beta_{1}\}\sqcup \beta'])\leq S_{k-1}(p_{A-s,B;s}^{N-1,M}[\alpha', \beta'])+\alpha'_{1}+\beta_{1}.$$
 Alors
 $$S_{k}(p_{A,B;s}^{N,M}[\alpha^-,\beta])=\nu_{1}-\alpha_{1}-\beta_{1}+A+S_{k-1}(p_{A-s,B;s}^{N-1,M}[\alpha',\{\beta_{1}\}\sqcup \beta'])$$
 $$\leq \nu_{1}-\alpha_{1}+\alpha'_{1}+A+S_{k-1}(p_{A-s,B;s}^{N-1,M}[\alpha', \beta'])\leq -\alpha_{1}+\alpha'_{1}+S_{k}(p_{A,B;s}^{N,M}[\alpha,\beta]).$$
 Puisque $\alpha'_{1}\leq \alpha_{1}$, la premi\`ere in\'egalit\'e de l'\'enonc\'e s'ensuit. On a aussi
 $$S_{k}(p_{A,B;s}^{N,M}[\alpha,\beta])=\nu_{1}+A+ S_{k-1}(p_{A-s,B;s}^{N-1,M}[\alpha',\beta'])\leq  \nu_{1}+A+S_{k-1}(p_{A-s,B;s}^{N-1,M}[\alpha',\{\beta_{1}\}\sqcup \beta'])$$
 $$=\alpha_{1}+\beta_{1}+S_{k}(p_{A,B;s}^{N,M}[\alpha^-,\beta]).$$
 C'est la deuxi\`eme in\'egalit\'e de l'\'enonc\'e, ce qui ach\`eve la d\'emonstration. $\square$
 
 \bigskip

   \subsection{Preuve du lemme 3.2}
   Supposons $n=0$ ou $m=0$. On a
   $$p_{A,B-s;s}^{N,M-1}[\alpha,\beta]=\{\alpha_{1}+A,...,\alpha_{N}+A+s-sN\}\sqcup \{\beta_{1}+B-s,...,\beta_{M-1}+B+s-sM\},$$
  $$p_{A-s,B;s}^{N-1,M}[\alpha,\beta]=\{\alpha_{1}+A-s,...,\alpha_{N-1}+A+s-sN\}\sqcup \{\beta_{1}+B,...,\beta_{M}+B+s-sM\}.$$ 
  Introduisons des entiers $a\in \{0,...,N\}$, $b\in \{0,...,M-1\}$ tels que
  $$S_{k}(p_{A,B-s}^{N,M-1}[\alpha,\beta])=S_{a}(\{\alpha_{1}+A,...,\alpha_{N}+A+s-sN\})+S_{b}(\{\beta_{1}+B-s,...,\beta_{M-1}+B+s-sM\}).$$
  Si $a=0$, on a $b=k$ et
  $$S_{k}(p_{A,B-s}^{N,M-1}[\alpha,\beta])=S_{k}(\{\beta_{1}+B-s,...,\beta_{M-1}+B+s-sM\})$$
  $$=S_{k}(\{\beta_{1}+B,...,\beta_{M}+B+s-sM\})-sk\leq S_{k}(p_{A-s,B;s}^{N-1,M}[\alpha,\beta]).$$
  Supposons $a\geq1$. Alors 
  $$S_{a}(\{\alpha_{1}+A,...,\alpha_{N}+A+s-sN\})=S_{a}(\{\alpha_{1}+A-s,...,\alpha_{N-1}+A+s-sN\})+sa$$
  $$=S_{a-1}(\{\alpha_{1}+A-s,...,\alpha_{N-1}+A+s-sN\})+\alpha_{a}+A,$$
  $$S_{b}(\{\beta_{1}+B-s,...,\beta_{M-1}+B+s-sM\})=S_{b}(\{\beta_{1}+B,...,\beta_{M}+B+s-sM\})-sb$$
  $$=S_{b+1}(\{\beta_{1}+B,...,\beta_{M}+B+s-sM\})-\beta_{b+1}-B.$$
  D'o\`u 
  $$ S_{k}(p_{A,B-s;s}^{N,M-1}[\alpha,\beta])= \alpha_{a}-\beta_{b+1}+A-B+S_{a-1}(\{\alpha_{1}+A-s,...,\alpha_{N-1}+A+s-sN\})$$
  $$+S_{b+1}(\{\beta_{1}+B,...,\beta_{M}+B+s-sM\})\leq \alpha_{1}+A-B+S_{k}(p_{A-2,B}^{N-1,M}[\alpha,\beta]).$$
  
   On suppose maintenant $n\geq1$, $m\geq1$. Supposons que l'une des conditions suivantes soient v\'erifi\'ees:
   
 (1)  $(1,0)<_{I}(1,1)$ et $\alpha_{1}+A-s\geq B$;
   
  (2)  $(1,1)<_{I}(1,0)$ et $\beta_{1}+B\leq A-s$. 
  
  Remarquons que (1) implique $\alpha_{1}+A> B-s$ et que (2) implique $\beta_{1}+B-s<A$. Alors les \'el\'ements de $P_{A,B-s;s}(\alpha,\beta)$ comme de $P_{A-s,B;s}(\alpha,\beta)$ peuvent se construire par le proc\'ed\'e (a). Celui-ci cr\'ee un \'el\'ement $\nu_{1}$ et des partitions $(\alpha',\beta')$. On a
  $$p_{A,B-s;s}^{N,M-1}[\alpha,\beta]=\{\nu_{1}+A\}\sqcup p_{A-s,B-s;s}^{N-1,M-1}[\alpha',\beta'],$$
  $$p_{A-s,B;s}^{N-1,M}[\alpha,\beta]=\{\nu_{1}+A-s\}\sqcup p_{A-2s,B;s}^{N-2,M}[\alpha',\beta'].$$
 Si $k=1$, on a
 $$S_{1}(p_{A,B-s;s}^{N,M-1}[\alpha,\beta])=S_{1}(p_{A-s,B;s}^{N-1,M}[\alpha,\beta])+s$$
 Mais les hypoth\`eses (1) ou (2) impliquent que $s\leq \alpha_{1}+A-B$ et on obtient l'in\'egalit\'e cherch\'ee. Si $k\geq2$, on a par
   r\'ecurrence
  $$S_{k-1}(p_{A-s,B-s;s}^{N-1,M-1}[\alpha',\beta'])\leq S_{k-1}(p_{A-2s,B;s}^{N-2,M}[\alpha',\beta'])+sup(\alpha'_{1}+A-s-B,0).$$
 D'o\`u
 $$S_{k}(p_{A,B-s;s}^{N,M-1}[\alpha,\beta])\leq S_{k}(p_{A-s,B;s}^{N-1,M}[\alpha,\beta])+s+sup(\alpha'_{1}+A-s-B,0).$$
 Comme ci-dessus, nos hypoth\`eses impliquent $s\leq \alpha_{1}+A-B$. On a aussi $s+\alpha'_{1}+A-s-B\leq \alpha_{1}+A-B$. Donc $s+sup(\alpha'_{1}+A-s-B,0)\leq sup(\alpha_{1}+A-B,0)$ et, avec l'in\'egalit\'e ci-dessus, on obtient celle cherch\'ee.
 
 Supposons que l'une des conditions suivantes soient v\'erifi\'ees:
   
 (3)  $(1,0)<_{I}(1,1)$ et $\alpha_{1}+A\leq B-s$;
   
  (4)  $(1,1)<_{I}(1,0)$ et $\beta_{1}+B-s\geq A$. 
   
   Remarquons que (3) implique $\alpha_{1}+A-s<B$ et que (4) implique $\beta_{1}+B>A-s$. Alors les \'el\'ements de $P_{A,B-s;s}(\alpha,\beta)$ comme de $P_{A-s,B;s}(\alpha,\beta)$ peuvent se construire par le proc\'ed\'e (b). Celui-ci cr\'ee un \'el\'ement $\mu_{1}$ et des partitions $(\alpha',\beta')$. On a 
  $$p_{A,B-s;s}^{N,M-1}[\alpha,\beta]=\{\mu_{1}+B-s\}\sqcup p_{A,B-2s;s}^{N,M-2}[\alpha',\beta'],$$
  $$p_{A-s,B;s}^{N-1,M}[\alpha,\beta]=\{\mu_{1}+B\}\sqcup p_{A-s,B-s;s}^{N-1,M-1}[\alpha',\beta'].$$
 Si $k=1$, on a
 $$S_{1}(p_{A,B-s;s}^{N,M-1}[\alpha,\beta])=S_{1}(p_{A-s,B;s}^{N-1,M}[\alpha,\beta])-s,$$
 et, \'evidemment, $-s<sup(\alpha_{1}+A-B,0)$. 
  Si $k\geq2$, on a par
   r\'ecurrence
  $$S_{k-1}(p_{A,B-2s;s}^{N,M-2}[\alpha',\beta'])\leq S_{k-1}(p_{A-s,B-s;s}^{N-1,M-1}[\alpha',\beta'])+sup(\alpha'_{1}+A-B+s,0).$$
 D'o\`u
 $$S_{k}(p_{A,B-s;s}^{N,M-1}[\alpha,\beta])\leq S_{k}(p_{A-s,B;s}^{N-1,M}[\alpha,\beta])-s+sup(\alpha'_{1}+A+s-B,0).$$
 Comme ci-dessus,  $-s<sup(\alpha_{1}+A-B,0)$. On a aussi $-s+\alpha'_{1}+A+s-B\leq \alpha_{1}+A-B$. Donc $-s+sup(\alpha'_{1}+A+s-B,0)\leq sup(\alpha_{1}+A-B,0)$ et, avec l'in\'egalit\'e ci-dessus, on obtient celle cherch\'ee. 
 
 Supposons enfin qu'aucune des relations (1) \`a (4) ne soit satisfaite. On a alors l'une des possibilit\'es suivantes:
 
 (5) $(1,0)<_{I}(1,1)$ et $B-s<\alpha_{1}+A< B+s$;
 
 (6) $(1,1)<_{I}(1,0)$ et $A-s<\beta_{1}+B<A+s$.
 
 On voit  que les \'el\'ements de $P_{A,B-s;s}(\alpha,\beta)$ se construisent \`a l'aide du proc\'ed\'e (a) tandis que ceux de $P_{A-s,B;s}(\alpha,\beta)$ se construisent \`a l'aide du proc\'ed\'e (b). Le proc\'ed\'e (a) cr\'ee un terme $\nu_{1}$ et des partitions $(\alpha',\beta')$. Le proc\'ed\'e (b) cr\'ee un terme $\mu_{1}$ et des partitions $(\alpha'',\beta'')$. D'apr\`es 1.2(4), sous l'hypoth\`ese (5), on a $\nu_{1}=\alpha_{1}+\mu_{1}$, $\alpha''$ contient le terme $\alpha_{1}$ et on a $\alpha'=\alpha^{''-}$ et $\beta'=\beta''$. Alors
 $$p_{A,B-s;s}^{N,M-1}[\alpha,\beta]=\{\alpha_{1}+\mu_{1}+A\}\sqcup p_{A-s,B-s;}^{N-1,M-1}[\alpha^{''-},\beta''],$$
 $$p_{A-s,B;s}^{N-1,M}[\alpha,\beta]=\{\mu_{1}+B\}\sqcup p_{A-s,B-s;s}^{N-1,M-1}[\alpha'',\beta''].$$
 D'apr\`es 1.2(4), sous l'hypoth\`ese (6), on a $\mu_{1}=\beta_{1}+\nu_{1}$, $\beta'$ contient le terme $\beta_{1}$ et on a $\beta''=\beta^{'-}$ et $\alpha''=\alpha'$. Alors, 
$$p_{A,B-s;s}^{N,M-1}[\alpha,\beta]=\{\nu_{1}+A\}\sqcup p_{A-s,B-s;s}^{N-1,M-1}[\alpha',\beta'],$$
 $$p_{A-s,B;s}^{N-1,M}[\alpha,\beta]=\{\beta_{1}+\nu_{1}+B\}\sqcup p_{A-s,B-s;s}^{N-1,M-1}[\alpha',\beta^{'-}].$$
 
 Si $k=1$, on a
$$S_{1}(p_{A,B-s;s}^{N,M-1}[\alpha,\beta])=S_{1}(p_{A-s,B;s}^{N-1,M}[\alpha,\beta])+\alpha_{1}+A-B$$
sous l'hypoth\`ese (5) et
$$S_{1}(p_{A,B-s;s}^{N,M-1}[\alpha,\beta])=S_{1}(p_{A-s,B;s}^{N-1,M}[\alpha,\beta])+A-B-\beta_{1}$$
sous l'hypoth\`ese (6). Les termes $\alpha_{1}+A-B$, resp. $A-B-\beta_{1}$ sont \'evidemment inf\'erieurs ou \'egaux \`a $sup(\alpha_{1}+A-B,0)$ et on conclut.

Si $k\geq2$, on voit que, dans le cas (5),  les partitions $(\alpha'',\beta'')$ v\'erifient les hypoth\`eses du lemme 3.1 tandis que, dans le cas (6), les partitions $(\alpha',\beta')$ v\'erifient celles du lemme sym\'etrique. Dans le cas (5), on a donc
$$S_{k-1}(p_{A-s,B-s;s}^{N-1,M-1}[\alpha^{''-},\beta''])\leq S_{k-1}(p_{A-s,B-s;s}^{N-1,M-1}[\alpha'',\beta'']),$$
d'o\`u
$$S_{k}(p_{A,B-s;s}^{N,M-1}[\alpha,\beta])=\alpha_{1}+\mu_{1}+A+S_{k-1}(p_{A-s,B-s;s}^{N-1,M-1}[\alpha^{''-},\beta''])$$
$$\leq  \alpha_{1}+A-B+(\mu_{1}+B)+S_{k-1}(p_{A-s,B-s;s}^{N-1,M-1}[\alpha'',\beta''])=\alpha_{1}+A-B+S_{k}(p_{A-s,B;s}^{N-1,M}[\alpha,\beta]),$$
ce qui entra\^{\i}ne l'in\'egalit\'e cherch\'ee. Dans le cas (6), on a 
$$S_{k-1}(p_{A-s,B-s;s}^{N-1,M-1}[\alpha',\beta'])\leq S_{k-1}(p_{A-s,B-s;s}^{N-1,M-1}[\alpha',\beta^{'-}])+\alpha'_{1}+\beta_{1}.$$
D'o\`u
$$S_{k}(p_{A,B-s;s}^{N,M-1}[\alpha,\beta])=\nu_{1}+A+S_{k-1}(p_{A-s,B-s;s}^{N-1,M-1}[\alpha',\beta'])\leq \nu_{1}+A+\alpha'_{1}+\beta_{1}+S_{k-1}(p_{A-s,B-s;s}^{N-1,M-1}[\alpha',\beta^{'-}])$$
$$=\alpha'_{1}+A-B+(\beta_{1}+\nu_{1}+B)+S_{k-1}(p_{A-s,B-s;s}^{N-1,M-1}[\alpha',\beta^{'-}])=\alpha'_{1}+A-B+S_{k}(p_{A-s,B;s}^{N-1,M}[\alpha,\beta]).$$
Puisque $\alpha'_{1}\leq \alpha_{1}$, on obtient l'in\'egalit\'e cherch\'ee. $\square$
 
 \bigskip

\subsection{Troisi\`eme  lemme }
On suppose $n\geq1$ et $(1,0)<_{I}(1,1)$ si $m\geq1$.

\ass{Lemme}{ 
(i) Soit $k\in {1,...,N+M}$.
On a l'in\'egalit\'e
$$S_{k-1}(p^{N-1,M}_{A-s,B;s}[\alpha^-,\beta])+\alpha_{1}+A\leq S_{k}(p^{N,M}_{A,B;s}[\alpha,\beta]).$$

(ii) Soit $e\in {\mathbb N}$ et $k\in {1,...,N+M-1}$.
Posons $b=b_{A-s,B+es/2;s}^{N-1,M}(\alpha^-,\beta;k)$. Si $b\geq1$, on a l'in\'egalit\'e
$$S_{k}(p^{N-1,M}_{A-s,B+es/2;s}[\alpha^-,\beta])+\alpha_{1}+A-B-s+bs(1-e/2)\leq S_{k}(p^{N,M}_{A,B;s}[\alpha,\beta]).$$
Si $b=0$, on a l'in\'egalit\'e
$$S_{k}(p^{N-1,M}_{A-s,B+es/2;s}[\alpha^-,\beta]) \leq S_{k}(p^{N,M}_{A,B;s}[\alpha,\beta]).$$}

 Les deux assertions se traitent par r\'ecurrence en reprenant la construction de l'\'el\'ement canonique de type (b) de $P_{A-s,B;s}(\alpha^-,\beta)$ pour l'assertion (i), de $P_{A-s,B+es/2;s}(\alpha^-,\beta)$ pour l'assertion (ii). La premi\`ere construction est le cas particulier $e=0$ de la seconde. On  unifie les notations en posant $e=0$ quand on traite l'assertion (i). On pose aussi dans ce cas $b^-=b_{A-s,B;s}^{N-1,M}(\alpha^-,\beta;k-1)$. La preuve est s\'epar\'ee en cinq cas trait\'es dans les paragraphes suivants.
 
 \bigskip
 
 \subsection{Preuve du lemme 3.5, cas $m=0$}

Supposons  $m=0$. Alors 
$$S_{k-1}(p^{N-1,M}_{A-s,B;s}[\alpha^-,\beta])=S_{k-1-b^-}((\alpha^-\sqcup\{0^{N-n}\})+[A-s,A+s-sN]_{s})+S_{b^-}([B,B+s-sM]_{s}).$$
On a 
$$S_{k-1-b^-}((\alpha^-\sqcup\{0^{N-n}\})+[A-s,A+s-sN]_{s})=S_{k-b^-}((\alpha\sqcup\{0^{N-n}\})+[A,A+s-sN]_{s})-\alpha_{1}-A,$$
d'o\`u
$$S_{k-1}(p^{N-1,M}_{A-s,B;s}[\alpha^-,\beta])=-\alpha_{1}-A+S_{k-b^-}((\alpha\sqcup\{0^{N-n}\}))+[A,A+s-sN]_{s})+S_{b^-}([B,B+s-sM]_{s})$$
$$\leq -\alpha_{1}-A+S_{k}(p^{N,M}_{A,B;s}[\alpha,\beta]).$$
Cela prouve (i).

On a aussi
$$S_{k}(p^{N-1,M}_{A-s,B+es/2;s}[\alpha^-,\beta])=S_{k-b}((\alpha^-\sqcup\{0^{N-n}\})+[A-s,A+s-sN]_{s})$$
$$+S_{b}([B+es/2,B+es/2+s-sM]_{s}).$$
Si $b=0$, on a \'evidemment $S_{k}((\alpha^-\sqcup\{0^{N-n}\})+[A-s,A+s-sN]_{s})<S_{k}((\alpha\sqcup\{0^{N-n}\})+[A,A+s-sN]_{s})$, d'o\`u 
$$S_{k}(p^{N-1,M}_{A-s,B+es/2;s}[\alpha^-,\beta])<S_{k}((\alpha\sqcup\{0^{N-n}\})+[A,A+s-sN]_{s})\leq S_{k}(p^{N,M}_{A,B;s}[\alpha,\beta])$$
et (ii).

Si $b\geq1$, on a 
$$S_{b}([B+es/2,B+es/2+s-sM]_{s})=S_{b-1}([B+es/2,B+es/2+s-sM]_{s})+B+es/2+s-sb$$
$$=S_{b-1}([B,B+s-sM]_{2})+B+s+bs(e/2-1),$$ 
$$S_{k-b}((\alpha^-\sqcup\{0^{N-n}\})+[A-s,A+s-sN]_{s})=S_{k-b+1}((\alpha\sqcup\{0^{N-n}\})+[A,A+s-sN]_{s})-\alpha_{1}-A.$$
 D'o\`u
$$S_{k}(p^{N-1,M}_{A-s,B+es/2;s}[\alpha^-,\beta])=-\alpha_{1}-A+B+s+bs(e/2-1)+S_{k-b+1}((\alpha\sqcup\{0^{N-n}\})+[A,A+s-sN]_{s})$$
$$+S_{b-1}([B,B+s-sM]_{s})\leq -\alpha_{1}-A+B+s+bs(e/2-1)+S_{k}(p^{N,M}_{A,B;s}[\alpha,\beta]).$$
Cela prouve (ii) dans le cas $m=0$.

 \bigskip
 
 \subsection{Preuve du lemme 3.5, cas $m\geq1$, $\alpha_{1}+A-B\leq 0$}

Supposons $m\geq1$ et $\alpha_{1}+A-B\leq 0$. Dans ce cas, le proc\'ed\'e (b) est autoris\'e pour construire les \'el\'ements de $P_{A,B;s}(\alpha,\beta)$. Ce proc\'ed\'e construit un \'el\'ement $\mu_{1} $ et des partitions $(\alpha',\beta')$. Le premier terme de $\alpha'$ est $\alpha_{1}$.
On a $p^{N,M}_{A,B;s}[\alpha,\beta]=\{\mu_{1}+B\}\sqcup p^{N,M-1}_{A,B-s;s}[\alpha',\beta']$. D'apr\`es 2.1(5),  le proc\'ed\'e (b) est autoris\'e pour construire des \'el\'ements de $P_{A-s,B+es/2}(\alpha^-,\beta)$.  D'apr\`es 2.1(4), ce proc\'ed\'e construit le m\^eme \'el\'ement $\mu_{1}$ et les partitions $(\alpha^{'-},\beta)$. D'o\`u 
$$(1) \qquad p^{N-1,M}_{A-s,B+es/2;s}[\alpha^-,\beta]=\{\mu_{1}+B+es/2\}\sqcup p^{N-1,M-1}_{A-s,B+se/2-s;s}[\alpha^{'-},\beta'].$$

Pour l'assertion (i), on a $e=0$ comme on l'a dit. Si $k=1$, on a $S_{k-1}(p^{N-1,M}_{A-s,B;s}[\alpha^-,\beta])=0$ et $S_{k}(p^{N,M}_{A,B;s}[\alpha,\beta])=\mu_{1}+B\geq \alpha_{1}+A$ d'apr\`es l'hypoth\`ese $\alpha_{1}+A-B\leq 0$. Si $k\geq2$, on a $S_{k-1}(p^{N-1,M}_{A-s,B;s}[\alpha^-,\beta])=\mu_{1}+B+S_{k-2}(p^{N-1,M-1}_{A-s,B-s;s}[\alpha^{'-},\beta'])$ et $S_{k}(p^{N,M}_{A,B;s}[\alpha,\beta])=\mu_{1}+B+S_{k-1}(p^{N,M-1}_{A,B-s;s}[\alpha^-,\beta])$. Puisque le premier terme de $\alpha'$ est $\alpha_{1}$, on a par r\'ecurrence 
$$S_{k-2}(p^{N-1,M-1}_{A-s,B-s;s}[\alpha^{'-},\beta'])+\alpha_{1}+A\leq S_{k-1}(p^{N,M-1}_{A,B-s;s}[\alpha^-,\beta]),$$
d'o\`u (i).

Tournons-nous vers l'assertion (ii). Si $k=1$, la formule (1) montre que  $b=1$ et on doit voir que $\mu_{1}+B+es/2+\alpha_{1}+A-B-es/2\leq \mu_{1}+B$, ce qui r\'esulte de l'hypoth\`ese $\alpha_{1}+A-B\leq 0$. Si $k\geq2$,  la formule (1) montre que $b\geq1$ et que $b_{A-s,B+se/2-s;s}^{N-1,M-1}(\alpha^{'-},\beta')=b-1$.    Si $b=1$, on a  par r\'ecurrence $S_{k-1}(p^{N-1,M-1}_{A-s,B+se/2-s;s}[\alpha^{'-},\beta'])\leq S_{k-1}(p^{N,M-1}_{A,B-s;s}[\alpha',\beta'])$ et le calcul est identique au cas $k=1$. Supposons $b\geq2$. Alors, par r\'ecurrence, on a
$$S_{k-1}(p^{N-1,M-1}_{A-s,B+se/2-s;s}[\alpha^{'-},\beta'])+\alpha_{1}+A-B+s(b-1)(1-e/2)\leq S_{k-1}(p^{N,M-1}_{A,B-s;s}[\alpha',\beta']),$$
d'o\`u
$$S_{k}(p^{N-1,M}_{A-s,B+es/2;s}[\alpha^-,\beta]+\alpha_{1}+A-B-s+bs(1-e/2)$$
$$=\mu_{1}+B+es/2+S_{k-1}(p^{N-1,M-1}_{A-s,B+s(e/2-1);s}[\alpha^{'-},\beta'])+\alpha_{1}+A-B-s+bs(1-e/2)$$
$$\leq S_{k-1}(p^{N,M-1}_{A,B-s;s}[\alpha',\beta'])+\mu_{1}+B=S_{k}(p^{N,M}_{A,B;s}[\alpha,\beta]),$$
ce qui est l'in\'egalit\'e cherch\'ee.

 \bigskip
 
 \subsection{Preuve du lemme 3.5, cas $m\geq1$, $\alpha_{1}+A-B>0$, utilisation du proc\'ed\'e (b) pour $(\alpha^-,\beta)$}

Supposons $m\geq1$ et  $\alpha_{1}+A-B>0$. Le proc\'ed\'e (a) appliqu\'e \`a $(\alpha,\beta)$ construit un \'el\'ement $\nu_{1}=\alpha_{1}+\beta_{1}+\alpha_{a_{2}}+... $ et des partitions $(\alpha',\beta')$. On a
$$(1) \qquad p^{N,M}_{A,B;s}[\alpha,\beta]=\{\nu_{1}+A\}\sqcup p^{N-1,M}_{A-s,B;s}[\alpha',\beta'].$$
Remarquons que cette formule suffit \`a d\'emontrer l'assertion (i) dans le cas $k=1$ puisqu'alors $S_{k-1}(p^{N-1,M}_{A-s,B;s}[\alpha^-,\beta])=0$ et $S_{k}(p^{N,M}_{A,B;s}[\alpha,\beta])=\nu_{1}+A\geq \alpha_{1}+A$. 

Supposons que le proc\'ed\'e (b) soit loisible pour construire des \'el\'ements de $P_{A-s,B+es/2;s}(\alpha^-,\beta)$. D'apr\`es 2.1(1), ce proc\'ed\'e construit l'\'el\'ement $\mu_{1}=\nu_{1}-\alpha_{1}$ et les m\^emes partitions $\alpha',\beta'$. On a
$$(2) \qquad p^{N-1,M}_{A-s,B+es/2;s}[\alpha^-,\beta]=\{\mu_{1}+B+es/2\}\sqcup p^{N-1,M-1}_{A-s,B+es/2-s;s}[\alpha',\beta'].$$

Pour l'assertion (i), on fait $e=0$. On a d\'ej\`a trait\'e le cas $k=1$. Si $k\geq2$, on introduit  l'\'el\'ement canonique de type (b) $(\nu',\mu')\in P_{A-s,B-s;s}(\alpha',\beta')$. La formule (2) montre que $b^-\geq1$ et $b_{A-s,B-s;s}^{N-1,M-1}(\alpha',\beta';k-2)=b^--1$. D'o\`u
$$S_{k-2}(p^{N-1,M-1}_{A-s,B-s;s}[\alpha',\beta'])=S_{k-1-b^-}((\nu'\sqcup\{0^{N-n}\})+[A-s,A+s-sN]_{s})$$
$$+S_{b^--1}((\mu'\sqcup\{0^{M-m}\})+[B-s,B+s-sM]_{s}).$$
On a 
$$S_{b^--1}((\mu'\sqcup\{0^{M-m}\})+[B-s,B+s-sM]_{s})=S_{b^-}((\mu'\sqcup\{0^{M-m}\})+[B,B+s-sM]_{s})-\mu'_{b^-}-B,$$
avec la convention $\mu'_{b^-}=0$ si $b^->m-1$. D'o\`u
$$S_{k-2}(p^{N-1,M-1}_{A-s,B-s;s}[\alpha',\beta'])\leq S_{k-1-b^-}((\nu'\sqcup\{0^{N-n}\})+[A-s,A+s-sN]_{s})$$
$$+S_{b^-}((\mu'\sqcup\{0^{M-m}\})+[B,B+s-sM]_{s})-\mu'_{b^-}-B\leq S_{k-1}(p\Lambda_{A-s,B;s}^{N-1,M}(\nu',\mu'))-B.$$
Puisque $(\nu',\mu')\in P (\alpha',\beta')$, on a  $S_{k-1}(p\Lambda_{A-s,B;s}^{N-1,M}(\nu',\mu'))\leq S_{k-1}(p_{A-s,B;s}^{N-1,M}[\alpha',\beta'])$ d'apr\`es le lemme 1.4. D'o\`u
$$S_{k-2}(p^{N-1,M-1}_{A-s,B-s;s}[\alpha',\beta'])\leq S_{k-1}(p_{A-s,B;s}^{N-1,M}[\alpha',\beta'])-B,$$
puis, gr\^ace \`a (1) et (2), 
$$S_{k-1}(p^{N-1,M}_{A-s,B;s}[\alpha^-,\beta])=\mu_{1}+B+ S_{k-2}(p^{N-1,M-1}_{A-s,B-s;s}[\alpha',\beta'])\leq \nu_{1}-\alpha_{1}+S_{k-1}(p_{A-s,B;s}^{N-1,M}[\alpha',\beta'])$$
$$-\alpha_{1}-A+(\nu_{1}+A+S_{k-1}(p_{A-s,B;s}^{N-1,M}[\alpha',\beta']))=-\alpha_{1}-A+S_{k}(p_{A,B;s}^{N,M}[\alpha,\beta]).$$
C'est l'in\'egalit\'e cherch\'ee.

Traitons maintenant (ii). D'apr\`es (2), on a $b\geq1$ et $b_{A-s,B+se/2-s;s}^{N-1,M-1}(\alpha',\beta';k-1)=b-1$.    Introduisons l'\'el\'ement canonique de type (b) $(\nu',\mu')\in P_{A-s,B+se/2-s;s}[\alpha',\beta']$. On a 
$$S_{k-1}(p^{N-1,M-1}_{A-s,B+se/2-s;s}[\alpha',\beta'])$$
$$=S_{k-b}(\nu'+[A-s,A+s-sN]_{s})+S_{b-1}(\mu'+[B+se/2-s,B+s+es/2-sM]_{s})$$
$$=S_{k-b}(\nu'+[A-s,A+s-sN]_{s})+S_{b-1}(\mu'+[B,B+2s-sM]_{s})+s(e/2-1)(b-1)$$
$$\leq S_{k-1}(p\Lambda^{N-1,M-1}_{A-s,B;s}(\nu',\mu'))+s(e/2-1)(b-1).$$
Puisque $(\nu',\mu')\in P(\alpha',\beta')$, on a $p\Lambda^{N-1,M-1}_{A-s,B;s}(\nu',\mu')\leq p^{N-1,M-1}_{A-s,B;s}[\alpha',\beta']$,  d'o\`u
$$S_{k-1}(p^{N-1,M-1}_{A-s,B+es/2-s;s}[\alpha',\beta'])\leq S_{k-1}(p^{N,M}_{A-s,B;s}[\alpha',b'])+s(e/2-1)(b-1),$$
puis
$$S_{k}(p^{N-1,M}_{A-s,B+es/2}[\alpha^-,\beta])=\mu_{1}+B+es/2+S_{k-1}(p^{N-1,M-1}_{A-s,B+es/2-s;s}[\alpha',\beta'])$$
$$=\nu_{1}-\alpha_{1}+B+es/2+S_{k-1}(p^{N-1,M-1}_{A-s,B+es/2-s;s}[\alpha',\beta'])$$
$$\leq B+s+s(e/2-1)b-\alpha_{1}-A+(\nu_{1}+A+S_{k-1}(p^{N-1,M}_{A-s,B;s}[\alpha',\beta']))$$
$$=B+s+s(e/2-1)b-\alpha_{1}-A+S_{k}(p^{N,M}_{A,B;s}[\alpha,\beta]).$$
C'est  l'in\'egalit\'e cherch\'ee.

 \bigskip
 
 \subsection{Preuve du lemme 3.5, cas $m\geq1$, $\alpha_{1}+1-B>0$, utilisation du proc\'ed\'e (a) pour $(\alpha^-,\beta)$, $(1,0)<_{I^-}(1,1)$}

On suppose que $m\geq1$, $\alpha_{1}+A-B>0$ mais que le proc\'ed\'e (b) n'est pas autoris\'e pour construire des \'el\'ements de $P_{A-s,B+es/2}(\alpha^-,\beta)$.   On a $n\geq2$ sinon seul ce proc\'ed\'e (b) est autoris\'e. On a encore l'\'egalit\'e (1) de 3.8. Supposons $(1,0)<_{I^-}(1,1)$, c'est-\`a-dire $(2,0)<_{I}(1,1)$. Le premier terme de $\alpha'$ est $\alpha_{2}=\alpha^-_{1}$. Puisque le proc\'ed\'e (b) n'est pas autoris\'e, on a $\alpha_{2}+A-s> B+es/2$. D'apr\`es 2.1(2), le proc\'ed\'e (a) construit l'\'el\'ement $\underline{\nu}_{1}=\nu_{1}-\alpha_{1}+\alpha_{2}$ et les partitions $(\alpha^{'-},\beta')$. On a
$$(1) \qquad p^{N-1,M}_{A-s,B+es/2;s}[\alpha^-,\beta]=\{\nu_{1}-\alpha_{1}+\alpha_{2}+A-s\}\sqcup p^{N-2,M}_{A-2s,B+es/2;s}[\alpha^{'-},\beta'].$$

Pour l'assertion (i), on fait $e=0$.  On  a d\'ej\`a trait\'e le cas $k=1$. Si $k\geq2$, parce que $\alpha'_{1}=\alpha_{2}$,  on a par r\'ecurrence 
$$S_{k-2}(p^{N-2,M}_{A-2s,B;s}[\alpha^{'-},\beta'])\leq -\alpha_{2}-A+s+S_{k-1}(p^{N-1,M}_{A-s,B;s}[\alpha',\beta']) $$
d'o\`u
$$S_{k-1}(p^{N-1,M}_{A-s,B;s}[\alpha^-,\beta])=\nu_{1}-\alpha_{1}+\alpha_{2}+A-s+S_{k-2}(p^{N-2,M}_{A-2s,B;s}[\alpha^{'-},\beta'])\leq \nu_{1}-\alpha_{1}+S_{k-1}(p^{N-1,M}_{A-s,B;s}[\alpha',\beta'])$$
$$=-\alpha_{1}-A+(\nu_{1}+A+S_{k-1}(p^{N-1,M}_{A-s,B;s}[\alpha',\beta']))=-\alpha_{1}-A+S_{k}(p^{N,M}_{A,B;s}[\alpha,\beta]),$$
ce qui est l'in\'egalit\'e cherch\'ee. 

Traitons (ii). Si $k=1$,  la relation (1) montre que $b=0$. On a
$$S_{1}(p^{N-1,M}_{A-s,B+es/2}[\alpha^-,\beta])=\nu_{1}-\alpha_{1}+\alpha_{2}+A-s<\nu_{1}+A=S_{1}(p^{N,M}_{A,B;s}[\alpha,\beta])$$
ce qui est l'in\'egalit\'e cherch\'ee. Supposons $k\geq2$. La relation (4) montre que
$b=b_{A-2s,B+es/2;s}^{N-2,B}(\alpha^{'-},\beta';k-1)$. Si $b=0$, on a par r\'ecurrence
$$S_{k-1}(p^{N-2,M}_{A-2s,B+es/2;s}[\alpha^{'-},\beta'])\leq S_{k-1}(p^{N-1,M}_{A-s,B;s}[\alpha',\beta']),$$
d'o\`u
 $$S_{k}(p^{N-1,M}_{A-s,B+es/2}[\alpha^-,\beta])=\nu_{1}-\alpha_{1}+\alpha_{2}+A-s+S_{k-1}(p^{N-2,M}_{A-2s,B+es/2;s}[\alpha^{'-},\beta'])$$
 $$<\nu_{1}+A+S_{k-1}(p^{N-1,M}_{A-s,B;s}[\alpha',\beta'])=S_{k}(p^{N,M}_{A,B;s}[\alpha,\beta])$$
 ce qui est l'in\'egalit\'e cherch\'ee. 
Si $b>0$, parce que $\alpha'_{1}=\alpha_{2}$, on a par r\'ecurrence  
$$S_{k-1}(p^{N-2,M}_{A-2s,B;s}[\alpha^{'-},\beta'])\leq -\alpha_{2}-A+s+B+s+sb(e/2-1)+S_{k-1}(p^{N-1,M}_{A-s,B;s}[\alpha',\beta']),$$
d'o\`u
$$S_{k}(p^{N-1,M}_{A-s,B;s}[\alpha^-,\beta])=\nu_{1}-\alpha_{1}+\alpha_{2}+A-s+S_{k-1}(p^{N-2,M}_{A-2s,B;s}[\alpha^{'-},\beta'])$$
$$\leq -\alpha_{1} -A+B+s+sb(e/2-1)+(\nu_{1}+A+S_{k-1}(p^{N-1,M}_{A-s,B;s}[\alpha',\beta']))$$
$$=-\alpha_{1}-A+B+s+sb(e/2-1)+S_{k}(p^{N,M}_{A,B;s}[\alpha,\beta]),$$
ce qui est l'in\'egalit\'e cherch\'ee. 

\bigskip

 \subsection{Preuve du lemme 3.5, cas $m\geq1$, $\alpha_{1}+1-B>0$, utilisation du proc\'ed\'e (a) pour $(\alpha^-,\beta)$, $(1,1)<_{I^-}(1,0)$}

Supposons  que $m\geq1$,  que $\alpha_{1}+1-B>0$ et que le proc\'ed\'e (b) n'est pas autoris\'e pour construire les \'el\'ements de $P_{A-s,B+es/2;s}(\alpha^-,\beta)$. Supposons de plus que   $(1,1)<_{I^-}(1,0)$, c'est-\`a-dire $(1,1)<_{I}(2,0)$. On a toujours l'\'egalit\'e (1) de 3.8.  Puisque le proc\'ed\'e (b) n'est pas autoris\'e, on a $\beta_{1}+B+es/2<A-s$. D'apr\`es 2.1(3), le proc\'ed\'e (a) construit  l'\'el\'ement $\underline{\nu}_{1}=\nu_{1}-\alpha_{1}-\beta_{1}$ et les partitions $(\alpha',\{\beta_{1}\}\sqcup \beta')$. On a
$$(1) \qquad p^{N-1,M}_{A-s,B+es/2;s}[\alpha^-,\beta]=\{\nu_{1}-\alpha_{1}-\beta_{1}+A-s\} \sqcup p^{N,M}_{A-2s,B+es/2;s}[\alpha',\{\beta_{1}\}\sqcup \beta'].$$
Posons $\alpha^1=\alpha'$, $\beta^1=\{\beta_{1}\}\sqcup \beta'$ d'o\`u $\beta'=(\beta^{1})^-$. On va d\'efinir un entier $r$, des termes $\nu_{i+1}$ pour $i=1,...,r$ et des partitions $(\alpha^{i},\beta^{i})\in {\cal P}_{n_{i}}\times {\cal P}_{m_{i}}$ pour $i=1,...,r+1$, extraites de $(\alpha',\beta')$.  Les propri\'et\'es suivantes seront v\'erifi\'ees. L'ensemble d'indices $I_{i}=(\{1,...,n_{i}\}\times \{0\})\cup(\{1,...,m_{i}\}\times \{1\})$ se plonge dans $I^-$ et est donc muni d'un ordre $<_{I_{i}}$.  On a $m_{i}\geq1$ et le plongement identifie $(1,1)\in I_{i}$ \`a $(1,1)\in I^-$. En particulier, $(1,1)$ est le plus petit terme de $I_{i}$. Le premier terme de $\beta^{i}$ est $\beta_{1}$.   Supposons d\'efinie $(\alpha^{i},\beta^{i})$, ce qui est le cas pour $i=1$. Si $n_{i}=0$, on pose $r=i-1$ et on a termin\'e. Supposons $n_{i}\geq1$. Si $\beta_{1}+B+es/2\geq A-s-si$, on pose $r=i-1$ et on a termin\'e. Si $\beta_{1}+B+es/2< A-s-si$, on applique le proc\'ed\'e (a) pour construire l'\'el\'ement canonique de type (b) de $P_{A-s-si,B+es/2}(\alpha^{i},\beta^{i})$. On pose $I_{i}^-=(\{1,...,n_{i}\}\times \{0\})\cup(\{1,...,m_{i}-1\}\times \{1\})$. On plonge cet ensemble dans $I_{i}$ en identifiant $(j,1)\in I_{i}^-$ \`a $(j+1,1)\in I_{i}$, d'o\`u un ordre $<_{I_{i}^-}$ sur $I_{i}^-$. Montrons que l'on peut aussi appliquer le proc\'ed\'e (a) pour construire un \'el\'ement de $P_{A-si,B;s}(\alpha^{i},\beta^{i-})$.   Si $(1,0)<_{I^-_{i}}(1,1)$, on doit montrer que $\alpha^{i}+A-si\geq B$; si $(1,1)<_{I^-_{i}}(1,0)$, on doit montrer que  $\beta^{i-}_{1}+B\leq A-si$; dans les deux cas, cela r\'esulte de l'in\'egalit\'e $\beta_{1}+B+es/2< A-s-si$. D'apr\`es l'assertion sym\'etrique de 2.1(4), le proc\'ed\'e (a) appliqu\'e \`a $(\alpha^{i},\beta^{i})$, resp. $(\alpha^{i},\beta^{i-})$, construit le m\^eme \'el\'ement  que l'on note $\nu_{i+1}$   et des partitions $(\alpha^{i+1},\beta^{i+1})$, resp. $(\alpha^{i+1},\beta^{i+1,-})$. A ce point, on obtient les \'egalit\'es
$$p_{A,B;s}^{N,M}[\alpha,\beta]=\{\nu_{1}+A,\nu_{2}+A-s,...,\nu_{r+1}+A-sr\}\sqcup p^{N-r-1,M}_{A-sr-s,B;s}[\alpha^{r+1},\beta^{r+1,-}],$$
$$(2) \qquad p_{A-s,B+es/2;s}^{N-1,M}[\alpha^-,\beta]$$
$$=\{\nu_{1}-\alpha_{1}-\beta_{1}+A-s,\nu_{2}+A-2s,...,\nu_{r+1}+A-sr-s\}\sqcup p^{N-r-2,M}_{A-sr-2s,B+es/2;s}[\alpha^{r+1},\beta^{r+1}].$$
Par d\'efinition de $r$, on a soit $n_{r+1}=0$, soit $\beta_{1}+B+es/2\geq A-2s-sr$. Dans ce dernier cas, on a aussi $\beta_{1}+B+es/2<A-s-sr$: si $r=0$, cela r\'esulte de l'hypoth\`ese $\beta_{1}+B+es/2<A-s$; si $r>0$ cela r\'esulte de ce que la construction ne s'est pas arr\^et\'ee en $r-1$. On pose $\nu=(\nu_{1},...,\nu_{r+1})$. 

Supposons d'abord $k\leq r+1$.  Pour l'assertion (i),  on a d\'ej\`a trait\'e le cas $k=1$, cf. 3.8. Si $k\geq2$, on a
$$S_{k-1}(p_{A-s,B;s}^{N-1,M}[\alpha^-,\beta])=S_{k-1}(\nu)-\alpha_{1}-\beta_{1}+S_{k-1}([A-s,A+s-sk]_{s})$$
$$\leq S_{k}(\nu)-\alpha_{1}-A+S_{k}([A,A+s-sk]_{s})=S_{k}(p_{A,B;s}^{N,M}[\alpha,\beta])-\alpha_{1}-A.$$
Pour (ii), (2) entra\^{\i}ne que $b=0$ et on a
$$S_{k}(p_{A-s,B+es/2;s}^{N-1,M}[\alpha^-,\beta])=S_{k}(\nu)-\alpha_{1}-\beta_{1}+S_{k}([A-s,A-sk]_{s}$$
$$= S_{k}(\nu)-\alpha_{1}-\beta_{1}-sk+S_{k}([A,A+s-sk]_{s})<S_{k}(p_{A,B;s}^{N,M}[\alpha,\beta]).$$

Supposons $k>r+1$ et posons $h=k-r-1$. On voit que, pour prouver (i), il suffit de prouver que
$$(3)\qquad S_{h-1}(p^{N-r-2,M}_{A-sr-2s,B;s}[\alpha^{r+1},\beta^{r+1}])+A-\beta_{1}-sr-s\leq S_{h}(p^{N-r-1,M}_{A-sr-s,B;s}[\alpha^{r+1},\beta^{r+1,-}]).$$
Dans le cas de l'assertion (ii), on voit que $b=b_{A-sr-2s,B+es/2;s}^{N-r-2,M}(\alpha^{r+1},\beta^{r+1};h)$. Il suffit de prouver:

dans le cas $b=0$,
$$(4) \qquad S_{h}(p^{N-r-2,M}_{A-sr-2s,B+es/2;s}[\alpha^{r+1},\beta^{r+1}])-\alpha_{1}-\beta_{1}-sr-s< S_{h}(p^{N-r-1,M}_{A-sr-s,B;s}[\alpha^{r+1},\beta^{r+1,-}]);$$

dans le cas $b\geq1$, 
$$(5)\qquad S_{h}(p^{N-r-2,M}_{A-sr-2s,B+es/2;s}[\alpha^{r+1},\beta^{r+1}])-\beta_{1}+A-B-sr-2s+sb(1-e/2)$$
$$\leq S_{h}(p^{N-r-1,M}_{A-sr-s,B;s}[\alpha^{r+1},\beta^{r+1,-}]).$$

Supposons d'abord $n_{r+1}=0$. Posons par convention  $\beta^{r+1}_{i}=0$ si $i\geq m_{r+1}+1$. Rappelons que $\beta^{r+1}_{1}=\beta_{1}$. Alors 
$$p^{N-r-2,M}_{A-sr-2s,B+es/2;s}[\alpha^{r+1},\beta^{r+1}]=[A-sr-2s,A+s-sN]_{s}$$
$$\sqcup \{\beta_{1}+B+es/2,\beta^{r+1}_{2}+B+es/2-s,...,\beta^{r+1}_{M}+B+es/2+s-sM\},$$
et
$$p^{N-r-1,M}_{A-sr-s,B;s}[\alpha^{r+1},\beta^{r+1,-}]=[A-sr-s,A+s-sN]_{s}\sqcup \{\beta^{r+1}_{2}+B,\beta^{r+1}_{2}+B-s,...,\beta^{r+1}_{M+1}+B+s-sM\}.$$

Traitons (3). On a $b^-=b_{A-sr-2s,B;s}^{N-r-2,M}(\alpha^{r+1},\beta^{r+1};h-1)$ d'apr\`es (2). D'o\`u 
$$S_{h-1}(p^{N-r-2,M}_{A-sr-2s,B;s}[\alpha^{r+1},\beta^{r+1}])=S_{h-1-b^-}([A-sr-2s,A+s-sN]_{s})$$
$$+S_{b^-}(\{\beta_{1}+B,\beta^{r+1}_{2}+B-s,...,\beta^{r+1}_{M}+B+s-sM\}).$$
On a 
$$S_{h-1-b^-}([A-sr-2s,A+s-sN]_{s})=S_{h-b^-}([A-sr-s,A+s-sN]_{s})-A+sr+s,$$
$$S_{b^-}(\{\beta_{1}+B,\beta^{r+1}_{2}+B-s,...,\beta^{r+1}_{M}+B+s-sM\})=$$
$$S_{b^-}(\{\beta^{r+1}_{2}+B,\beta^{r+1}_{3}+B-s,...,\beta^{r+1}_{M+1}+B+s-sM\})+\beta_{1}-\beta^{r+1}_{b^-+1}.$$
Enfin
$$S_{h-b^-}([A-sr-s,A+s-sN]_{s})+S_{b^-}(\{\beta^{r+1}_{2}+B,\beta^{r+1}_{3}+B-s,...,\beta^{r+1}_{M+1}+B+s-sM\})\leq $$
$$S_{h}(p^{N-r-1,M}_{A-sr-s,B;s}[\alpha^{r+1},\beta^{r+1,-}]).$$
En mettant bout-\`a-bout ces relations, on obtient (3). 

Traitons (4). Puisque $b=0$ pour cette assertion, on a
$$S_{h}(p^{N-r-2,M}_{A-sr-2s,B+es/2;s}[\alpha^{r+1},\beta^{r+1}])=S_{h}([A-sr-2s,A+s-sN]_{s})$$
$$=S_{h}([A-sr-s,A+s-sN]_{s})-sh\leq S_{h}(p^{N-r-1,M}_{A-sr-s,B;s}[\alpha^{r+1},\beta^{r+1,-}])-sh,$$
et (4) s'en d\'eduit. 

Traitons (5). On a
$$S_{h}(p^{N-r-2,M}_{A-sr-2s,B+es/2;s}[\alpha^{r+1},\beta^{r+1}])=S_{h-b}([A-sr-2s,A+s-sN]_{s})$$
$$+S_{b}( \{\beta_{1}+B+es/2,\beta^{r+1}_{2}+B+es/2-s,...,\beta^{r+1}_{M}+B+es/2+s-sM\}).$$
Puisque $b\geq1$, on a
 $$S_{b}( \{\beta_{1}+B+es/2,\beta^{r+1}_{2}+B+es/2-s,...,\beta^{r+1}_{M}+B+es/2+s-sM\})=$$
 $$S_{b-1}( \{\beta^{r+1}_{2}+B,\beta^{r+1}_{3}+B-s,...,\beta^{r+1}_{M+1}+B+s-sM\})+\beta_{1}+B+s-sb+ebs/2.$$
et
$$S_{h-b}([A-sr-2s,A+s-sN]_{s})=S_{h-b+1}([A-sr-s,A+s-sN]_{s})-A+sr+s.$$
On a aussi
$$S_{b-1}( \{\beta^{r+1}_{2}+B,\beta^{r+1}_{3}+B-s,...,\beta^{r+1}_{M+1}+B+s-sM\})+S_{h-b+1}([A-sr-s,A+s-sN]_{s})\leq $$
$$S_{h}(p^{N-r-1,M}_{A-sr-s,B;s}[\alpha^{r+1},\beta^{r+1,-}]),$$
et (5) se d\'eduit de ces relations.

Supposons enfin $n_{r+1}\geq1$ donc $A-2s-sr\leq \beta_{1}+B+es/2< A-s-sr$ ainsi qu'on l'a dit plus haut. La premi\`ere in\'egalit\'e entra\^{\i}ne que le proc\'ed\'e (b)  s'applique pour construire l'\'el\'ement canonique  de type (b) de $P_{A-2s-sr,B+es/2;s}(\alpha^{r+1},\beta^{r+1})$.  On obtient un terme $\mu_{1} $ et des partitions $(\alpha'',\beta'')$. Montrons que    l'\'el\'ement canonique de type (b) de $P_{A-sr-s,B;s}(\alpha^{r+1},\beta^{r+1,-})$ se construit \`a l'aide du proc\'ed\'e (a). Si $m_{r+1}=1$, c'est le seul proc\'ed\'e autoris\'e. Si $m_{r+1}\geq2$ et $(1,0)<_{I^-_{r+1}}(1,1)$, on a $\alpha_{1}^{r+1}+A-s-sr> B$ parce que $B\leq \beta_{1}+B+es/2<A-s-sr$. Si $(1,1)< _{I_{r+1}^-}(1,0)$, on a $\beta^{r+1}_{2}+B\leq \beta_{1}+B\leq \beta_{1}+B+es/2< A-s-sr$. D'apr\`es l'assertion sym\'etrique de 2.1(1), ce proc\'ed\'e (a) construit l'\'el\'ement $\nu_{r+2}=\mu_{1}-\beta_{1}$ et les partitions $(\alpha'',\beta'')$.On obtient
$$(6) \qquad p^{N-r-2,M}_{A-sr-2s,B+es/2;s}[\alpha^{r+1},\beta^{r+1}]=\{\mu_{1}+B+es/2\}\sqcup p^{N-r-2,M-1}_{A-sr-2s,B+es/2-s;s}[\alpha'',\beta'']$$
$$(7) \qquad p^{N-r-1,M}_{A-sr-s,B;s}[\alpha^{r+1},\beta^{r+1,-}]=\{\mu_{1}-\beta_{1}+A-sr-s\}\sqcup p^{N-r-2,M}_{A-sr-2s,B;s}[\alpha'',\beta''].$$

Traitons (3). Si $h=1$, on a
$$S_{h-1}(p^{N-r-2,M}_{A-sr-2s,B;s}[\alpha^{r+1},\beta^{r+1}])=0$$
et 
$$S_{h}(p^{N-r-1,M}_{A-sr-s,B;s}[\alpha^{r+1},\beta^{r+1,-}])=\mu_{1}-\beta_{1}+A-sr-s\geq -\beta_{1}+A-sr-s,$$
d'o\`u l'in\'egalit\'e cherch\'ee. Si $h\geq2$, on introduit l'\'el\'ement canonique de type (b)  $(\nu'',\mu'')\in P_{A-sr-2s,B-s;s}(\alpha'',\beta'')$. On a
$$S_{h-1}(p^{N-r-2,M}_{A-sr-2s,B;s}[\alpha^{r+1},\beta^{r+1}])=\mu_{1}+B+S_{h-1-b^-}((\nu''\sqcup\{0^{N-r-2-n''}\})+[A-sr-2s,A+s-sN]_{s})$$
$$+S_{b^--1}((\mu''\sqcup\{0^{M-1-m''}\})+[B-s,B+s-sM]_{s}).$$
On a
$$S_{b^--1}((\mu''\sqcup\{0^{M-1-m''}\})+[B-s,B+s-sM]_{s})=S_{b^-}((\mu''\sqcup\{0^{M-m''}\})+[B,B+s-sM]_{s})-B-\mu''_{b^-},$$
avec toujours la convention $\mu''_{i}=0$ si $i\geq m''+1$. 
On a aussi
$$S_{h-1-b^-}((\nu''\sqcup\{0^{N-r-2-n''}\})+[A-sr-2s,A+s-sN]_{s})+S_{b^-}((\mu''\sqcup\{0^{M-m''}\})+[B,B+s-sM]_{s})$$
$$\leq S_{h-1}(p\Lambda_{A-sr-2s,B;s}^{N-r-2,M}(\nu'',\mu''))\leq S_{h-1}(p_{A-sr-2s,B;s}^{N-r-2,M}[\alpha'',\beta'']),$$
puisque $(\nu'',\mu'')\in P(\alpha'',\beta'')$. D'o\`u
$$S_{h-1}(p^{N-r-2,M}_{A-sr-2s,B;s}[\alpha^{r+1},\beta^{r+1}])\leq \mu_{1}-\mu''_{b^-}+S_{h-1}(p_{A-sr-2s,B;s}^{N-r-2,M}[\alpha'',\beta''])\leq $$
$$-A+\beta_{1}+sr+s+(\mu_{1}-\beta_{1}+A-sr-s)
+S_{h-1}(p_{A-sr-2s,B;s}^{N-r-2,M}[\alpha'',\beta''])=$$
$$-A+\beta_{1}+sr+s+S_{h}(p^{N-r-1,M}_{A-sr-s,B;s}[\alpha^{r+1},\beta^{r+1,-}]).$$
La relation (3) en r\'esulte.

La relation (6) entra\^{\i}ne que  $b\geq1$ et il reste seulement \`a prouver (5). On voit aussi que $b_{A-sr-2s,B-s;s}^{N-r-2,M-1}(\alpha'',\beta'';h-1)=b-1$.  Si $h=1$, on a $b=1$. D'apr\`es (6) et (7), on a
$$S_{h}(p^{N-r-2,M}_{A-sr-2s,B+es/2;s}[\alpha^{r+1},\beta^{r+1}])=\mu_{1}+B+es/2,$$
$$S_{h}(p^{N-r-1,M-1}_{A-sr-s,B;s}[\alpha^{r+1},\beta^{r+1,-}])=\mu_{1}-\beta_{1}+A-sr-s,$$
et un calcul simple conduit \`a (5). Supposons $h\geq2$. Introduisons l'\'el\'ement canonique de type (b)   $(\nu'',\mu'')$ dans $P_{A-sr-2s,B+es/2-s;s}(\alpha'',\beta'')$. Puisque $b_{A-sr-2s,B-s;s}^{N-r-2,M-1}(\alpha'',\beta'';h-1)=b-1$, on a
$$S_{h-1} (p^{N-r-2,M-1}_{A-sr-2s,B+es/2-s}[\alpha'',\beta''])=S_{h-b}((\nu''\sqcup\{0^{N-r-2-n''}\})+[A-sr-2s,A+s-sN]_{s})$$
$$+S_{b-1}((\mu''\sqcup\{0^{M-1-m''}\})+[B+es/2-s,B+es/2+s-sM]_{s}),$$
et 
 $$S_{b-1}((\mu''\sqcup\{0^{M-1-m''}\})+[B+es/2-s,B+es/2+s-sM]_{s})=$$
$$S_{b-1}((\mu''\sqcup\{0^{M-1-m''}\})+[B,B+s-sM]_{s})+s(b-1)(e/2-1).$$
D'o\`u
$$S_{h-1} (p^{N-r-2,M-1}_{A-sr-2s,B+es/2-s}[\alpha'',\beta''])\leq S_{h-1}(p\Lambda_{A-sr-2s,B;s}^{N-r-2,M}(\nu'',\mu''))+s(b-1)(e/2-1).$$
Parce que $(\nu'',\mu'')$ appartient \`a $P(\alpha'',\beta'')$, on a
$$S_{h-1}(p\Lambda_{A-sr-2s,B;s}^{N-r-2,M}(\nu'',\mu''))\leq S_{h-1}(
p^{N-r-2,M}_{A-sr-2s,B;s}[\alpha'',\beta'']).$$
D'o\`u
$$S_{h}(p^{N-r-2,M}_{A-sr-2s,B+es/2;s}[\alpha^{r+1},\beta^{r+1}])=\mu_{1}+B+es/2+S_{h-1} (p^{N-r-2,M-1}_{A-sr-2s,B+es/2-s;s}[\alpha'',\beta''])$$
$$\leq \mu_{1}+B+es/2+s(b-1)(e/2-1)+S_{h-1}(
p^{N-r-2,M}_{A-sr-2s,B;s}[\alpha'',\beta''])$$
$$=\beta_{1}-A+B+sr+2s+sb(e/2-1)+ (\mu_{1}-\beta_{1}+A-sr-s)+S_{h-1}(
p^{N-r-2,M}_{A-sr-2s,B;s}[\alpha'',\beta''])$$
$$=\beta_{1}-A+B+sr+2s+sb(e/2-1)+S_{h}(p_{A-sr-s,B;s}^{N-r-1,M}[\alpha^{r+1},\beta^{r+1,-}]).$$
C'est la relation (5). Cela ach\`eve la d\'emonstration. $\square$

\section{Multiplicit\'es}

\bigskip

\subsection{D\'efinition et \'enonc\'e de la proposition principale}

 Soient $\alpha\in {\cal P}_{n}$ et $\beta\in {\cal P}_{m}$.  Rappelons que $I=(\{1,...,n\}\times\{0\}) \cup(\{1,...,m\}\times\{1\})$ est muni d'un ordre $<_{I}$.  Soit $J$ l'ensemble des couples $((i,e),(j,f))\in I$ tels que $(i,e)<_{I} (j,f)$ et $e\not=f$. Posons $X={\mathbb N}^J$. Un \'el\'ement de $X$ s'\'ecrit $x=(x_{(i,e),(j,f)})_{((i,e),(j,f))\in J}$, avec des $x_{(i,e),(j,f)}\in {\mathbb N}$. Pour $x\in X$, on d\'efinit un couple $(\alpha[x],\beta[x])\in {\mathbb Z}^n\times {\mathbb Z}^m$ par
 $$\alpha[x]_{i}=\alpha_{i}+(\sum_{j\in \{1,...,m\}, (i,0)<_{I} (j,1)}x_{(i,0),(j,1)})-(\sum_{j\in \{1,...,m\}, (j,1)<_{I} (i,0)}x_{(j,1),(i,0)}),$$
 $$\beta[x]_{j}=\beta_{j}+(\sum_{i\in \{1,...n\},(j,1)<_{I}(i,0)}x_{(j,1),(i,0)})-(\sum_{i\in \{1,...n\},(i,0)<_{I}(j,1)}x_{(i,0),(j,1)}).$$
 Pour $(\nu,\mu)\in {\mathbb Z}^n\times {\mathbb Z}^m$, notons $X(\alpha,\beta;\nu,\mu)$ l'ensemble des \'el\'ements $x\in X$ tels que $\alpha[x]=\nu$ et $\beta[x]=\mu$.
 
 {\bf Remarques.} (1) On a $S(\alpha[x])+S(\beta[x])=S(\alpha)+S(\beta)$ pour tout $x\in X$. On peut donc se limiter aux $(\nu,\mu)$ tels que $S(\nu)+S(\mu)=S(\alpha)+S(\beta)$. 
 
 (2) L'ensemble $X(\alpha,\beta;\nu,\mu)$ est fini. En effet, identifions l'ensemble d'indices $I$ \`a $\{1,...,n+m\}$ de fa\c{c}on croissante. Nos couples $(\alpha,\beta)$ et $(\alpha[x],\beta[x])$ deviennent des \'el\'ements de ${\mathbb Z}^{n+m}$. L'ensemble $X$ devient un ensemble d'\'el\'ements $(x_{k,l})_{(k,l)}$ o\`u $(k,l)$ parcourt un sous-ensemble $\Sigma$ de celui des couples tels que $k,l\in \{1,...,n+m\}$ et $k<l$. Pour $\gamma=(\gamma_{1},...,\gamma_{n+m})\in {\mathbb Z}^{n+m}$, posons $<\gamma,\delta>=\sum_{i=1,...,n+m}(n+m+1-i)\gamma_{i}$. En particulier, pour $x\in X$, $<x,\delta>=\sum_{(k,l)\in \Sigma}(l-k)x_{k,l}$. Tous les coefficients sont positifs donc, pour un entier $c\in {\mathbb Z}$ quelconque, il n'y a qu'un nombre fini de $x\in X$ tels que $<x,\delta>=c$. Or un \'el\'ement de $X(\alpha,\beta;\nu,\mu)$
doit v\'erifier $<x,\delta>=<(\nu,\mu),\delta>-<(\alpha,\beta),\delta>$. 

\bigskip

On note $\mathfrak{S}_{n}$ et $\mathfrak{S}_{m}$ les groupes de permutations de $\{1,...,n\}$ et $\{1,...,m\}$. On note $sgn$ leur caract\`ere signature usuel. Pour $\nu\in {\mathbb Z}^n$ et $w\in \mathfrak{S}_{n}$, on d\'efinit $\nu^w$ et $\nu[w]$ par $\nu^w_{i}=\nu_{wi}$ et  $\nu[w]_{i}=\nu_{wi}+i-wi$.  On d\'efinit de m\^eme $\mu^v$ et $\mu[v]$ pour $\mu\in {\mathbb Z}^m$ et $v\in \mathfrak{S}_{m}$.

Soient $(\nu,\mu)\in {\cal P}_{n}\times {\cal P}_{m}$. On suppose $S(\nu)+S(\mu)=S(\alpha)+S(\beta)$. On pose
$$mult(\alpha,\beta;\nu,\mu)=\sum_{w\in \mathfrak{S}_{n},v\in \mathfrak{S}_{m}}sgn(w)sgn(v)\vert X(\alpha,\beta;\nu[w],\mu[v])\vert.$$
On a $mult(\alpha,\beta;\nu,\mu)\in {\mathbb Z}$.

\ass{Proposition}{Soient $(\nu,\mu)\in {\cal P}_{n}\times {\cal P}_{m}$. On suppose $S(\nu)+S(\mu)=S(\alpha)+S(\beta)$. 

(i) Supposons $mult(\alpha,\beta;\nu,\mu)\not=0$. Alors $p\Lambda^{N,M}_{A,B;s}(\nu,\mu)\leq p^{N,M}_{A,B;s}[\alpha,
\beta]$.

(ii) Supposons $p\Lambda^{N,M}_{A,B;s}(\nu,\mu)= p_{A,B;s}^{N,M}[\alpha,\beta]$. Alors $mult(\alpha,\beta;\nu,\mu)\not=0$ si et seulement si $(\nu,\mu)\in P_{A,B;s}(\alpha,\beta)$. 

(iii) Supposons $(\nu,\mu)\in P_{A,B;s}(\alpha,\beta)$. Alors $mult(\alpha,\beta;\nu,\mu)=1$.}

\bigskip

\subsection{D\'ebut de la preuve de la proposition}
Preuve. On raisonne par r\'ecurrence sur $n+m$, tout \'etant trivial si $n+m=0$. Plus g\'en\'eralement, tout est trivial si $n=0$ ou $m=0$ (car $X=\emptyset$). On suppose donc $n\geq1$ et $m\geq1$ et on fixe $(\nu,\mu)\in {\cal P}_{n}\times {\cal P}_{m}$ avec  $S(\nu)+S(\mu)=S(\alpha)+S(\beta)$. On ne perd rien \`a supposer $(1,0)< _{I}(1,1)$. On d\'ecompose 
$$(1) \qquad  mult(\alpha,\beta;\nu,\mu)=\sum_{k=1,...,n}mult_{k}(\alpha,\beta;\nu,\mu),$$
o\`u 
$$mult_{k}(\alpha,\beta;\nu,\mu)=\sum_{w\in \mathfrak{S}_{n}, w1=k,v\in \mathfrak{S}_{m}}sgn(w)sgn(v)\vert X(\alpha,\beta;\nu[w],\mu[v])\vert.$$

On fixe $k\in \{1,...,n\}$.  
 On identifie $\mathfrak{S}_{n-1}$ \`a l'ensemble des $w\in \mathfrak{S}_{n}$ tels que $w1=k$: un \'el\'ement $w'\in \mathfrak{S}_{n-1}$ s'identifie \`a $w$ tel que $w1=k$, $wi=w'(i-1)$ si $i=2...,n$ et $w'(i-1)<k$, $wi=w'(i-1)+1$ si $i=2,...,n$ et $w'(i-1)\geq k$. Ainsi, on a $sgn(w)=(-1)^{k+1}sgn(w')$.  On note $J'$ le sous-ensemble des $((i,e),(j,f))\in J$ tels que $(i,e)\not=(1,0)$ et $J''$ le compl\'ementaire, c'est-\`a-dire l'ensemble des \'el\'ements $((1,0),(j,1))$, pour $j\in \{1,...,m\}$. L'ensemble $X$ se d\'ecompose conform\'ement en $X_{J'}\times X_{J''}$ et on peut identifier $X_{J''}$ \`a ${\mathbb N}^m$. On \'ecrit tout \'el\'ement  de $X$ sous la forme $(x',(x_{j})_{j=1,...m})$ avec $x'\in X_{J'}$. Comme en 2.1, introduisons la partition $\alpha^-=(\alpha_{2},...,\alpha_{n})$ et l'ordre $<_{I^-}$ sur l'ensemble d'indices  $I^-=(\{1,...,n-1\}\times \{0\})\cup (\{1,...,m\}\times \{1\})$. 
Pour $x'\in X_{J'}$, on d\'efinit $(\alpha^-[x'],\beta[x'])$ de m\^eme que l'on a d\'efini $(\alpha[x],\beta[x])$. Pour $x=(x',(x_{j})_{j=1,...,m})$, on calcule
$$\alpha[x]_{1}=\alpha_{1}+\sum_{j=1,...,m}x_{j};$$
pour $i=2,...,n$ et $j=1,...,m$,
$$\alpha[x]_{i}=\alpha'^-[x']_{i-1},$$
$$\beta[x]_{j}=\beta[x']_{j}-x_{j}.$$

Fixons $w'\in \mathfrak{S}_{n-1}$, soit $w\in \mathfrak{S}_{n}$ l'\'el\'ement associ\'e et d\'ecrivons l'ensemble $X(\alpha,\beta;\nu[w],\mu[v])$. Un \'el\'ement $(x',(x_{j})_{j=1,...m})\in X$ appartient \`a cet ensemble si et seulement si on a les \'egalit\'es suivantes

(2) $ \alpha_{1}+\sum_{j=1,...,m}x_{j}=\nu_{k}+1-k$;

(3) $\alpha^-[x']_{i-1}=\nu_{wi}-wi+i$ pour 
 $i=2,...,n$;
 
(4) $ \beta[x']-x_{j}=\mu_{vj}-vj+j$ pour $j=1,...,m$.

 La  relation (3) se d\'ecompose en
$\alpha^-[x']_{i-1}=\nu_{w'(i-1)}-w'(i-1)+i$ si $w'(i-1)<k$ et
$\alpha^-[x']_{i-1}=\nu_{w'(i-1)+1}-w'(i-1)+i-1$ si $w'(i-1)\geq k$. D\'efinissons $\nu'\in {\cal P}_{n-1}$ par
$\nu'_{i}=\nu_{i}+1$ si $i< k$ et $\nu'_{i}=\nu_{i+1}$ si $i\geq k$. La relation (3) \'equivaut \`a

(5) $\alpha^-[x']_{i}=\nu'_{w'i}-w'i+i$
pour tout $i\in \{1,...,n-1\}$. 

Pour tout \'el\'ement $y=(y_{1},...,y_{m})\in {\mathbb Z}^m$, notons $stab(y)$ le groupe des $u\in \mathfrak{S}_{m}$ tels que $y^{u}=y$. Notons $({\mathbb Z}^m)^+$ l'ensemble des $y$ tels que $y_{1}\geq y_{2}...\geq y_{m}$. Tout \'el\'ement $y\in {\mathbb Z}^m$ s'\'ecrit $y=z^{u}$ pour un unique $z\in ({\mathbb Z}^m)^+$ et un $u\in \mathfrak{S}_{m}$ dont la classe dans $stab(z)\backslash \mathfrak{S}_{m}$ est uniquement d\'etermin\'ee. On d\'ecompose $X(\alpha,\beta;\nu[w],\mu[v])$ en r\'eunion
$$\sqcup_{z,u}X(\alpha,\beta;\nu[w],\mu[v];z,u),$$
o\`u $z$ parcourt $({\mathbb Z}^m)^+$ et $u$ parcourt $stab(z)\backslash \mathfrak{S}_{m}$ et o\`u $X(\alpha,\beta;\nu[w],\mu[v];z,u)$ est le sous-ensemble des $x\in X(\alpha,\beta;\nu[w],\mu[v])$ tels que 

(6) $\beta[x']_{j}-j=z_{uj}$ pour tout $j=1,...,m$.  

L'ensemble $X(\alpha,\beta;\nu[w],\mu[v];z,u)$ se d\'ecompose en produit  d'un sous-ensemble de $X_{J'} $ et d'un sous-ensemble de $X_{J''} $.
Le premier est l'ensemble des \'el\'ements $x'\in X_{J'}$ v\'erifiant les relations (5) et (6).  Celles-ci ne d\'ependent que de $\alpha^-,\beta,\nu'[w'], z, u$, notons l'ensemble en question $X_{J'}(\alpha^-,\beta;\nu'[w'];z,u)$. Le deuxi\`eme ensemble est celui des $(x_{j})_{j=1,...,m}\in X_{J''}$ v\'erifiant la relation (2) ainsi que

(7) $z_{uj}-x_{j}=\mu_{vj}-vj$ pour tout $j=1,...,m$.

Remarquons que, si  cette derni\`ere relation est v\'erifi\'ee, on peut remplacer (2) par l'\'egalit\'e

(8) $\alpha_{1}+S(z)-S(\mu)+\sum_{j=1,...,m}j=\nu_{k}+1-k$.

Celle-ci est ind\'ependante des $x_{j}$. En cons\'equence, $X(\alpha,\beta;\nu[w],\mu[v];z,u)$ est vide si (8) n'est pas v\'erifi\'ee.  Si elle l'est, le sous-ensemble de $X_{J''}$ est d\'efini par la seule relation (7) et il est loisible de noter $X_{J''}(\mu[v];z,u)$ cet ensemble. Notons $Z(\alpha_{1},\nu_{k},\mu)$ l'ensemble des $z\in ({\mathbb Z}^m)^+$ v\'erifiant (8). On obtient
$$ X(\alpha,\beta;\nu[w],\mu[v])=\sqcup_{z,u; z\in Z(\alpha_{1},\nu_{k},\mu)}X_{J'}(\alpha',\beta;\nu'[w'];z,u)\times X_{J''}(\mu[v];z,u).$$
Puis
$$(9) \qquad \sum_{v\in \mathfrak{S}_{m}}sgn(v)\vert X(\alpha,\beta;\nu[w],\mu[v])\vert =\sum_{ z\in Z(\alpha_{1},\nu_{k},\mu), u\in stab(z)\backslash \mathfrak{S}_{m}}\vert X_{J'}(\alpha^-,\beta;\nu'[w'];z,u)\vert$$
$$ \sum_{v\in \mathfrak{S}_{m}}sgn(v)\vert X_{J''}(\mu[v];z,u)\vert .$$

Notons $Q(\alpha_{1},\nu_{k}+1-k,\mu)$ l'ensemble des partitions $\mu'=(\mu'_{1},...,\mu'_{m})\in {\cal P}_{m}$ v\'erifiant
$$(10) \qquad \mu'_{1}\geq \mu_{1}\geq \mu'_{2}\geq \mu_{2}\geq...\geq\mu'_{m}\geq\mu_{m}$$
et
$$(11) \qquad \alpha_{1}+S(\mu')-S(\mu)=\nu_{k}+1-k.$$
Pour $\mu'\in Q(\alpha_{1},\nu_{k}+1-k,\mu)$, notons $z(\mu')\in {\mathbb Z}^m$ l'\'el\'ement tel que $z(\mu')_{i}=\mu'_{i}-i$. On a $z(\mu')\in Z(\alpha_{1},\nu_{k},\mu)$ et  $stab(z(\mu'))=\{1\}$. Montrons que l'on a

(12) soit $z\in Z(\alpha_{1},\nu_{k},\mu)$; si $\sum_{v\in \mathfrak{S}_{m}}sgn(v)\vert X_{J''}(\mu[v];z,u)\vert\not=0$,  il existe $\mu'\in Q(\alpha_{1},\nu_{k}+1-k,\mu)$ tel que $z=z(\mu')$; inversement, pour $\mu'\in Q(\alpha_{1},\nu_{k}+1-k,\mu)$, on a  
$$\sum_{v\in \mathfrak{S}_{m}}sgn(v)\vert X_{J''}(\mu[v];z(\mu'),u)\vert=sgn(u).$$

Preuve. La suite $\mu_{1}-1,\mu_{2}-2,...$ est strictement d\'ecroissante. Puisque $z\in({\mathbb Z}^m)^+$, il existe d'unique entiers $0\leq l_{1}\leq l_{2}\leq ...\leq l_{m}\leq m$ tels que
$$z_{1}\geq...\geq z_{l_{1}}\geq \mu_{1}-1> z_{l_{1}+1}\geq...\geq z_{l_{2}}\geq \mu_{2}-2>z_{l_{2}+1}\geq...$$
$$\geq z_{l_{m}}\geq \mu_{m}-m>z_{l_{m}+1}\geq...\geq z_{m}.$$ 
Soit $v\in \mathfrak{S}_{m}$. La relation (6) se r\'ecrit
$x_{uv^{-1} h}=z_{uv^{-1} h}-\mu_{h}+h$ pour tout $h=1,...,m$. Puisque les $x_{ h}$ sont des entiers positifs ou nuls quelconques, l'ensemble $ X_{J''}(\mu[v];z,u)$ a au plus un \'el\'ement et il en a un si et seulement si
$z_{uv^{-1} h}\geq \mu_{h}-h$
pour tout $h=1,...,m$, ce qui \'equivaut \`a 

(13) $uv^{-1}h\leq l_{h}$ pour tout $h=1,...,m$. 

Notons $V$ l'ensemble des $v\in \mathfrak{S}_{m}$ qui v\'erifient cette condition. On obtient
$$\sum_{v\in \mathfrak{S}_{m}}sgn(v)\vert X_{J''}(\mu[v];z,u)\vert=\sum_{v\in V}sgn(v).$$

Supposons la somme de gauche non nulle. Alors $V$ n'est pas vide. La condition (13) pour $v\in V$  entra\^{\i}ne $uv^{-1}1,...uv^{-1}h\leq l_{h}$, en particulier $l_{m}=m$.  Introduisons le sous-groupe  $\mathfrak{T}$ des \'el\'ements de $\mathfrak{ S}_{m}$ qui conservent chacun des intervalles $\{1,...,l_{1}\}$, $\{l_{1}+1,...,l_{2}\}$ .... On voit alors que l'action du groupe $u^{-1}\mathfrak{T}u$ par multiplication \`a droite sur $\mathfrak{S}_{m}$ conserve l'ensemble $V$. 
  En fixant un ensemble $\underline{V}$ de repr\'esentants de cette action dans cet ensemble, on  obtient
 $$\sum_{v\in  V}sgn(v) =\sum_{v\in \underline{V}}sgn(v)\sum_{v'\in u^{-1}\mathfrak{T}u}sgn(v').$$
 La non-nullit\'e de cette expression entra\^{\i}ne que $\mathfrak{T}=\{1\}$.  On en d\'eduit $l_{h}=h$ pour tout $h$. On a alors
 $$z_{1}\geq \mu_{1}-1>z_{2}\geq \mu_{2}-2>...>z_{m}\geq \mu_{m}-m.$$
 En posant $\mu'_{i}=z_{i}+i$, on a $\mu'\in Q(\alpha_{1},\nu_{k}+1-k,\mu)$ et $z=z(\mu')$. Inversement, si $z=z(\mu')$ pour $\mu'\in Q(\alpha_{1},\nu_{k}+1-k,\mu)$, on a $l_{h}=h$ pour tout $h$. On   voit par r\'ecurrence sur $h$ que la condition (13) est v\'erifi\'ee  si et seulement si $uv^{-1}h=h$ pour tout $h$, autrement dit $u=v$.  On a alors $ V=\{u\}$ d'o\`u 
$\sum_{v\in  V}sgn(v) =sgn(u)$. Cela d\'emontre (12).  
  
 Pour $\mu'\in Q(\alpha_{1},\nu_{k}+1-k,\mu)$ l'ensemble $X_{J'}(\alpha^-,\beta;\nu'[w'];z(\mu'),u)$ n'est autre que $X(\alpha^-,\beta;\nu'[w'],\mu'[u])$. L'\'egalit\'e (9) et la propri\'et\'e (12) entra\^{\i}nent
 $$\sum_{w\in \mathfrak{S}_{n}, w1=k}\sum_{v\in \mathfrak{S}_{m}}sgn(w)sgn(v)\vert X(\alpha,\beta;\nu[w],\mu[v])\vert =(-1)^{k+1}\sum_{\mu'\in Q(\alpha_{1},\nu_{k}+1-k,\mu)}\sum_{w'\in \mathfrak{S}_{n-1}}$$
 $$\sum_{v\in \mathfrak{S}_{m}}sgn(w')sgn(v)\vert X(\alpha^-,\beta;\nu'[w'],\mu'[v])\vert.$$
 Autrement dit
 $$(14) \qquad mult_{k}(\alpha,\beta;\nu,\mu)=(-1)^{k+1}\sum_{\mu'\in Q(\alpha_{1},\nu_{k}+1-k,\mu)}mult(\alpha^-,\beta;\nu',\mu').$$
 
 \bigskip
 
 \subsection{Deuxi\`eme \'etape de la preuve}
  On suppose donn\'e un entier $k\in \{1,...,n\}$ et un \'el\'ement $\mu'\in Q(\alpha_{1},\nu_{k}+1-k,\mu)$. On suppose que 
 
 (1) $mult(\alpha^-,\beta;\nu',\mu')\not=0$ ou $(\nu',\mu')\in P(\alpha^-,\beta)$. 
 
 Cette hypoth\`ese entra\^{\i}ne
 
 (2) pour tout $c\in \{1,...,N+M-1\}$, tous r\'eels $A'\geq sN-2s$ et  $B'\geq sM-s$, on a
 $$S_{c}(p_{A',B';s}^{N-1,M}(\nu',\mu'))\leq S_{c}(p_{A',B';s}^{N-1,M}[\alpha^-,\beta]).$$
 
 En effet, si  $mult(\alpha^-,\beta;\nu',\mu')\not=0$ , c'est le (i) de la proposition que l'on applique par r\'ecurrence. Si $(\nu',\mu')\in P(\alpha^-,\beta)$, c'est le lemme 1.4.
 
    Pour tout entier $c\in \{1,...,N+M\}$, on fixe des r\'eels $A(c)\geq sN-s$, $B(c)\geq sM-s$ . On va prouver que 
 
 (3) pour tout $c\in \{1,...,N+M\}$, on a 
 $$S_{c}(p\Lambda_{A(c),B(c);s}^{N,M}(\nu,\mu))\leq S_{c}(p_{A(c),B(c);s}^{N,M}[\alpha,\beta]).$$
 
 {\bf Remarque.} L'entier $c$ \'etant fix\'e, on pourrait aussi bien simplifier la notation en posant $A=A(c)$ et $B=B(c)$. Mais on reprendra la d\'emonstration qui suit en 4.5  et il nous semble que les preuves de ce paragraphe seront plus claires si nous introduisons d\`es maintenant les notations qui seront alors n\'ecessaires.
 
 \bigskip

  Soit $c\in \{1,...,N+M\}$. On fixe 
 des entiers $a(c)\in \{0,...,N\}$, $b(c)\in \{0,...,M\}$ tels que $a(c)+b(c)=c$ et
 $$S_{c}(p\Lambda_{A(c),B(c);s}^{N,M}(\nu,\mu))=S_{a(c)}((\nu\sqcup\{ 0^{N-n}\})+[A(c),A(c)+s-sN]_{s})$$
 $$+S_{b(c)}((\mu\sqcup \{0^{M-m}\})+[B(c),B(c)+s-sM]_{s}).$$

 Supposons d'abord que $a(c)<k$. On a $\nu_{i}=\nu'_{i}-1$ pour $i=1,...,a(c)$, donc $S_{a(c)}((\nu\sqcup\{0^{N-n}\})+[A(c),A(c)+s-sN]_{s})=S_{a}((\nu'\sqcup\{0^{N+1-n}\})+[A(c),A(c)+2s-sN]_{s})-a(c)$ (la premi\`ere partition a $N$ termes et la seconde en a $N-1$ mais on ne consid\`ere que les $a$ premiers termes et on $a<k\leq n\leq N$). D'apr\`es 4.2(10), on a $S_{b(c)}((\mu\sqcup\{0^{M-m}\}+[B(c),B(c)+s-sM]_{s})\leq S_{b(c)}((\mu'\sqcup\{0^{M-m}\}+[B(c),B(c)+s-sM]_{s})$ d'o\`u
 $$ S_{c}(p\Lambda^{N,M}_{A(c),B(c);s}(\nu,\mu))\leq S_{a(c)}((\nu'\sqcup\{0^{N-n}\})+[A(c),A(c)+2s-sN]_{s})$$
 $$+S_{b(c)}((\mu'\sqcup\{0^{M-m}\})+[B(c),B(c)+s-sM]_{s})-a(c)$$
 $$\leq S _{c}(p\Lambda^{N-1,M}_{A(c),B(c);s}(\nu',\mu'))-a(c)\leq 
 S_{c}( p_{A(c),B(c);s}^{N-1,M}[\alpha^-,\beta])-a(c),$$
 la derni\`ere in\'egalit\'e r\'esultant de (2).  On a $S_{c}( p_{A(c),B(c);s}^{N-1,M}[\alpha^-,\beta])\leq S_{c}( p_{A(c),B(c);s}^{N,M}[\alpha^-,\beta])$: cette in\'egalit\'e est \'equivalente \`a $S_{c}(p\Lambda_{A(c),B(c);s}^{N-1,M}(\chi,\xi))\leq S_{c}(p\Lambda_{A(c),B(c);s}^{N,M}(\chi,\xi))$ o\`u 
  $(\chi,\xi)$ est un \'el\'ement de $P_{A(c),B(c);s}(\alpha^-,\beta)$,  et cette derni\`ere in\'egalit\'e est \'evidente puisque $p\Lambda_{A(c),B(c);s}^{N-1,M}(\chi,\xi))$ est extraite de $p\Lambda_{A(c),B(c);s}^{N,M}(\chi,\xi)$. 
 On applique le lemme 3.1:
 
 \noindent $S_{c}( p_{A(c),B(c);s}^{N,M}[\alpha^-,\beta])\leq S_{c}(p_{A(c),B(c);s}^{N,M}[\alpha,\beta])$ d'o\`u l'in\'egalit\'e  (3).
 
 On suppose maintenant $a(c)\geq k$. On a $\nu_{i}=\nu'_{i}-1$ pour $i=1,...,k-1$, $\nu_{k}=k-1+\alpha_{1}+S(\mu')-S(\mu)$ (d'apr\`es 4.2(11)) et $\nu_{i}=\nu'_{i-1}$ pour $i=k+1,...,n$. C'est-\`a-dire
 $$\nu\sqcup\{0^{N-n}\}=(\nu'_{1}-1,...,\nu'_{k-1}-1,k-1+\alpha_{1}+S(\mu')-S(\mu),\nu'_{k},...,\nu'_{n-1},0,...,0).$$
 D'o\`u
 $$S_{a(c)}(\nu\sqcup\{0^{N-n}\})=\alpha_{1}+S(\mu')-S(\mu)+S_{a(c)-1}(\nu'\sqcup\{0^{N-n}\}).$$
 On a aussi
 $$S_{a(c)}([A(c),A(c)+s-sN]_{s})=A(c)+S_{a(c)-1}([A(c)-s,A(c)+s-sN]_{s}).$$
 D'o\`u 
 $$(4) \qquad S_{a(c)}((\nu\sqcup\{0^{N-n}\})+[A(c),A(c)+s-sN]_{s})=\alpha_{1}+S(\mu')-S(\mu)+A(c)$$
 $$+S_{a(c)-1}((\nu'\sqcup\{0^{N-n}\})+[A(c)-s,A(c)+s-sN]_{s}).$$
 
 Supposons $b(c)\geq m$. Dans ce cas, 
 $$S_{b(c)}((\mu\sqcup\{0^{M-m}\})+[B(c),B(c)+s-sM]_{s})=S(\mu)+S_{b(c)}([B(c),B(c)+s-sM]_{s})$$
 $$=S(\mu)-S(\mu')+
 S_{b(c)}((\mu'\sqcup\{0^{M-m}\})+[B(c),B(c)+s-sM]_{s}).$$
 D'o\`u
 $$S_{c}(p\Lambda^{N,M}_{A(c),B(c);s}(\nu,\mu))=\alpha_{1}+A(c)+S_{a(c)-1}((\nu'\sqcup\{0^{N-n}\})+[A-s,A+s-sN]_{s})$$
 $$+S_{b(c)}((\mu'\sqcup\{0^{M-m}\})+[B(c),B(c)+s-sM]_{s})$$
 $$\leq \alpha_{1}+A(c)+S_{c-1}(p\Lambda^{N-1,M}_{A(c)-s,B(c);s}(\nu',\mu'))\leq   \alpha_{1}+A(c)+S_{c-1}(p^{N-1,M}_{A(c)-s,B(c);s}[\alpha^-,\beta]),$$
 la derni\`ere in\'egalit\'e r\'esultant de (2) si $c\geq2$ et \'etant triviale si $c=1$. 
 On applique le lemme 3.5(i):  
 $$ \alpha_{1}+A(c)+S_{c-1}(p^{N-1,M}_{A(c)-s,B(c);s}[\alpha^-,\beta])\leq S_{c}(p^{N,M}_{A(c),B(c);s}[\alpha,\beta]),$$
 d'o\`u l'in\'egalit\'e (3). 
 
 Supposons maintenant $b(c)<m$. On a donc $b(c)<M$. On a aussi  $b(c)<c$ car $c=b(c)+a(c)$ et $a(c)\geq k\geq1$. Donc $b(c)<inf(c,M)$. 
 On a d\'efini en 2.4 l'entier $b_{A(c)-s,B(c);s}^{N-1,M}(\alpha^-,\beta;c)$. S'il est sup\'erieur ou \'egal \`a $b(c)+1$, on pose $e(c)=0$ et $\underline{b}(c)=b_{A(c)-s,B(c);s}^{N-1,M}(\alpha^-,\beta;c)$. Sinon, le lemme 2.4 entra\^{\i}ne l'existence d'un entier $e\geq1$ tel que $b_{A(c)-s,B(c)+es/2;s}^{N-1,M}(\alpha^-,\beta;c)=b(c)+1$. On fixe cet entier  que l'on note $e(c)$ et on pose $\underline{b}(c)=b_{A(c)-s,B(c)+e(c)s/2;s}^{N-1,M}(\alpha^-,\beta;c)=b(c)+1$. On a
$$S_{b(c)}((\mu\sqcup\{0^{M-m}\})+[B(c),B(c)+s-sM]_{s})=S_{b(c)}(\mu)+S_{b(c)}([B(c),B(c)+s-sM]_{s})$$
$$=S_{b(c)}(\mu)+S_{b(c)+1}([B(c),B(c)+s-sM]_{s})-B(c)+sb(c)$$
$$=S_{b(c)}(\mu)+S_{b(c)+1}([B(c)+e(c)s/2,B(c)+e(c)s/2+s-sM]_{2})-B(c)+sb(c)-(b(c)+1)e(c)s/2$$
$$=S_{b(c)}(\mu)-S_{b(c)+1}(\mu')+S_{b(c)+1}((\mu'\sqcup\{0^{M-m\}})+[B(c)+e(c)s/2,B(c)+e(c)s/2+s-sM]_{s})$$
$$-B(c)+sb(c)-(b(c)+1)e(c)s/2.$$
Avec (4), on obtient
$$S_{c}(p\Lambda^{N,M}_{A(c),B(c);s}(\nu,\mu))=\alpha_{1}+S(\mu')-S(\mu)+S_{b(c)}(\mu)-S_{b(c)+1}(\mu')+A(c) -B(c)+sb(c)$$
$$-(b(c)+1)e(c)s/2+S_{a(c)-1}((\nu'\sqcup\{0^{N-n}\})+[A(c)-s,A(c)+s-sN]_{s})$$
$$+S_{b(c)+1}((\mu'\sqcup\{0^{M-m\}})+[B(c)+e(c)s/2,B(c)+e(c)s/2+s-sM]_{s})$$
$$\leq \alpha_{1}+S(\mu')-S(\mu)+S_{b(c)}(\mu)-S_{b(c)+1}(\mu')+A -B+sb(c)-(b(c)+1)e(c)s/2$$
$$+S_{c}(p\Lambda_{A(c)-s,B(c)+e(c)s/2;s}^{N-1,M}(\nu',\mu')).$$
D'apr\`es (2), on a
$$S_{c}(p\Lambda_{A(c)-s,B(c)+e(c)s/2;s}^{N-1,M}(\nu',\mu'))\leq S_{c}(p_{A(c)-s,B(c)+e(c)s/2;s}^{N-1,M}[\alpha^-,\beta]).$$
On a aussi 
$$S(\mu')-S(\mu)+S_{b(c)}(\mu)-S_{b(c)+1}(\mu')=-\mu_{b(c)+1}+\mu'_{b(c)+2}-\mu_{b(c)+2}+\mu'_{b(c)+3}...-\mu_{m-1}+\mu'_{m}-\mu_{m}.$$
Les sommes $-\mu_{b(c)+1}+\mu'_{b(c)+2}$, $-\mu_{b(c)+2}+\mu'_{b(c)+3}$ ... sont n\'egatives ou nulles d'apr\`es 4.2(10) et le dernier terme $-\mu_{m}$ l'est aussi. D'o\`u
$$S_{c}(p\Lambda^{N,M}_{A(c),B(c);s}(\nu,\mu))\leq \alpha_{1}+A(c) -B(c)+sb(c)-(b(c)+1)e(c)s/2+S_{c}(p_{A(c)-s,B(c)+e(c)s/2;s}^{N-1,M}[\alpha^-,\beta]).$$
On a
$$\alpha_{1}+A(c)-B(c)+sb(c)-(b(c)+1)e(c)s/2\leq \alpha_{1}+A(c)-B(c)-s+s\underline{b}(c)(1-e(c)/2).$$
En effet, cela \'equivaut \`a $s(\underline{b}(c)-1-b(c))+s(b(c)+1-\underline{b}(c))e(c)/2\geq0$. Le premier terme est positif ou nul par construction et le second est nul (si $e(c)\not=0$, $\underline{b}(c)=b(c)+1$ par construction). D'o\`u
$$S_{c}(p\Lambda^{N,M}_{A(c),B(c)}(\nu,\mu))\leq \alpha_{1}+A(c) -B(c)-s+s\underline{b}(c)(1-e(c)/2)+S_{c}(p_{A(c)-s,B(c)+e(c)s/2;s}^{N-1,M}[\alpha^-,\beta]).$$
On a $\underline{b}(c)\geq b(c)+1\geq 1$. Le lemme 3.5(ii) entra\^{\i}ne alors l'in\'egalit\'e (3). Cela ach\`eve la preuve de cette in\'egalit\'e. 

\bigskip

\subsection{Preuve du (i) de la proposition 4.1}

Supposons $mult(\alpha,\beta;\nu,\mu)\not=0$. D'apr\`es 4.2(1), on peut fixer $k\in \{1,...,n\}$ tel que 
$mult_{k}(\alpha,\beta;\nu,\mu)\not=0$. D'apr\`es 4.2(14), on peut fixer $\mu'\in Q(\alpha_{1},\nu_{k}+1-k,\mu)$ tel que $mult(\alpha^-,\beta;\nu',\mu')\not=0$. Dans les constructions de 4.3, on prend simplement  $A(c)=A$ et $B(c)=B$ pour tout $c$. Les hypoth\`eses de ce paragraphe \'etant satisfaites, la relation (3) de ce paragraphe nous dit que
$$S_{c}(p\Lambda_{A,B;s}^{N,M}(\nu,\mu))\leq S_{c}(p_{A,B;s}^{N,M}[\alpha,\beta])$$
pour tout $c\in \{1,...,N+M\}$. Autrement dit $p\Lambda_{A,B;s}^{N,M}(\nu,\mu)\leq p_{A,B;s}^{N,M}[\alpha,\beta]$. C'est l'assertion (i) de la proposition.

\bigskip

\subsection{Une nouvelle \'etape}
Comme en 4.3, fixons pour tout $c\in \{1,...,N+M\}$ des  r\'eels $A(c)\geq sN-s$, $B(c)\geq sM-s$, $a(c)\in \{0,...,N\}$, $b(c)\in \{0,...M\}$. Disons que ces donn\'ees satisfont l'hypoth\`ese $(Hyp)(\alpha,\beta;\nu,\mu)$ si les conditions suivantes sont satisfaites 

(1) pour tout $c\in \{1,...,N+M\}$, $a(c)+b(c)=c$ et 
$$S_{c}(p\Lambda_{A(c),B(c);s}^{N,M}(\nu,\mu))=S_{a(c)}((\nu\sqcup\{ 0^{N-n}\})+[A(c),A(c)+s-sN]_{s})$$
$$+S_{b(c)}((\mu\sqcup \{0^{M-m}\})+[B(c),B(c)+s-sM]_{s});$$

(2) pour tout $c\in \{1,...,N+M-1\}$,  $a(c)\leq a(c+1)\leq a(c)+1$ et $b(c)\leq b(c+1)\leq b(c)+1$;

(3)  pour tout $c\in \{1,...,N+M\}$,  
$$S_{c}(p\Lambda_{A(c),B(c);s}^{N,M}(\nu,\mu))=S_{c}(p_{A(c),B(c);s}^{N,M}[\alpha,\beta]).$$

Supposons ces hypoth\`eses satisfaites. On note $c^0$ le plus petit entier $c$ tel que $a(c)=1$ et $c^1$ le plus petit entier $c$ tel que $b(c)= m$  (de tels entiers existent d'apr\`es (2)).  On a forc\'ement $c^1\geq m$. 

{\bf Remarque} (4)  $c^0\not=c^1$.

En effet, si $c^1=1$, on a forc\'ement $m=1$ et $b(1)=1$. Alors $a(1)=1-b(1)=0$ et $1\not=c^0$ puisque $a(c^0)=1$. Si $c^1\geq2$, on a
  $b(c^1-1)=m-1$ donc $a(c^1-1)=c^1-1-b(c^1-1)=c^1-b(c^1)=a(c^1)$ ce qui interdit \`a $c^1$ d'\^etre  le plus petit entier $c$ tel que $a(c)=1$. 

\bigskip
 
On d\'efinit un \'el\'ement $\underline{\mu}\in {\mathbb Z}^m$ par

si $c^0\leq c^1$, 
$\underline{\mu}_{i}=\mu_{i}$ pour $i=1,..., c^0-1$, $\underline{\mu}_{c^0}=\nu_{1}-\alpha_{1}$,  $\underline{\mu}_{i}=\mu_{i-1}$ pour $i=c^0+1,...,m$;

si $c^0>c^1$, $\underline{\mu}=\mu$. 

On d\'efinit des r\'eels $\underline{A}(c)$, $\underline{B}(c)$ et des entiers  $\underline{a}(c)$, $\underline{b}(c)$ pour $c\in \{1,...,N+M-1\}$ de la fa\c{c}on suivante. 
Pour $c=1,...,c^0-1$, on pose $\underline{A}(c)=A(c)$, $\underline{B}(c)=B(c)$, $\underline{a}(c)=a(c)=0$, $\underline{b}(c)=b(c)=c$.
Si $c^1< c^0$, on pose $\underline{A}(c)=A(c+1)-s$, $\underline{B}(c)=B(c+1)$, $\underline{a}(c)=a(c+1)-1$ et $\underline{b}(c)=b(c+1)$
pour $c=c^0,...,N+M-1$. Supposons $c^0<c^1$.  Soit $c\in \{c^0,...,c^1-1\}$. On pose $\underline{A}(c)=A(c)-s$, $\underline{a}(c)=a(c)-1$ et $\underline{b}(c)=b(c)+1$. Si $b_{A(c)-s,B(c);s}^{N-1,M}(\alpha^-,\beta;c)\geq b(c)+1$, on pose $\underline{B}(c)=B(c)$. Si $b_{A(c)-s,B(c);s}^{N-1,M}(\alpha^-,\beta;c)\leq b(c)$, on choisit un entier $e(c)$ tel que $b_{A(c)-s,B(c)+e(c)s/2;s}^{N-1,M}(\alpha^-,\beta;c)= b(c)+1$ et on pose $\underline{B}(c)=B(c)+e(c)s/2$. Soit $c\in \{c^1,...,N+M-1\}$. 
On pose $\underline{A}(c)=A(c+1)-s$, $\underline{B}(c)=B(c+1)$, $\underline{a}(c)=a(c+1)-1$ et $\underline{b}(c)=b(c+1)$. 

Il est clair que, pour tout $c\in \{1,...,N+M-1\}$, on a $\underline{a}(c)\in \{0,...,N-1\}$, $\underline{b}(c)\in \{0,...,M\}$ et $\underline{a}(c)+\underline{b}(c)=c$. Montrons que l'on a l'analogue de (2), c'est-\`a-dire que

(5)  si $c\not=N+M-1$, on a $\underline{a}(c)\leq\underline{a}(c+1)\leq\underline{a}(c)+1$, $\underline{b}(c)\leq\underline{b}(c+1)\leq\underline{b}(c)+1$.

Il suffit de prouver la premi\`ere s\'erie d'in\'egalit\'es  car la seconde se d\'eduit de la premi\`ere et de l'\'egalit\'e $\underline{a}(c)+\underline{b}(c)=c$. Supposons $c^1<c^0$. En vertu de (2), (5) est clair sauf si $c=c^0-1$.  Ce cas n'intervient que si $c^0\not=N+M$. Dans ce cas, on a $\underline{a}(c^0-1)=0$ et $\underline{a}(c^0)=a(c^0+1)-1$. Or
$$0=a(c^0)-1\leq a(c^0+1)-1\leq a(c^0)=1$$
d'o\`u l'in\'egalit\'e cherch\'ee. Supposons $c^0<c^1$. En vertu de (2), (5) est clair sauf si $c=c^0-1$ (ce  cas n'intervient que si $c^0\not=1$)  ou $c=c^1-1$ (ce cas n'intervient que si $c^1\not=N+M$). Si $c^0\not=1$, on a $\underline{a}(c^0-1)=0$ et $\underline{a}(c^0)=a(c^0)-1=0$, d'o\`u  l'in\'egalit\'e cherch\'ee. Si $c^1\not=N+M$, on a $\underline{a}(c^1-1)=a(c^1-1)-1$ et $\underline{a}(c^1)=a(c^1+1)-1$. On a $b(c^1-1)< b(c^1)$, d'o\`u $a(c^1-1)=a(c^1)$ et l'in\'egalit\'e cherch\'ee se d\'eduit des in\'egalit\'es $a(c^1)\leq a(c^1+1)\leq a(c^1)+1$. Cela prouve (5).  

\ass{Lemme}{On suppose satisfaite l'hypoth\`ese $(Hyp)(\alpha,\beta;\nu,\mu)$. Soient $k\in \{1,...,n\}$ et $\mu'\in {\cal P}_{m}$. Supposons que $\mu'$ appartienne \`a $Q(\alpha_{1},\nu_{k}+1-k,\mu)$ et que $mult_{k}(\alpha^-,\beta;\nu',\mu')\not=0$ ou que $(\nu',\mu')\in P(\alpha^-,\beta)$. Alors

(i) $k=1$;

(ii) $\mu'=\underline{\mu}$ et, si $c^0<c^1$,  $\mu_{m}=0$;

(iii) les donn\'ees $\underline{A}(c)$, $\underline{B}(c)$, $\underline{a}(c)$ et $\underline{b}(c)$ pour $c\in \{1,...,N+M-1\}$ satisfont l'hypoth\`ese $(Hyp)(\alpha^-,\beta;\nu',\mu')$. }

{\bf Remarque.} Les assertions (i) et (ii) impliquent que, sous les hypoth\`eses de l'\'enonc\'e, $\underline{\mu}$ est une partition, qui appartient \`a $Q(\alpha_{1},\nu_{1},\mu)$ (en particulier, si $c^1<c^0$, on
a alors $\nu_{1}-\alpha_{1}=S(\underline{\mu})-S(\mu)=0$). On voit alors que $(\nu,\mu)=\iota_{c^0,\alpha_{1}}(\nu',\underline{\mu})$ avec la notation de 2.3.

\bigskip

Preuve. On reprend la preuve de 4.3 dont les hypoth\`eses sont satisfaites. Pour tout $c\in \{1,...,N+M\}$, on y a \'ecrit diff\'erentes suites  d'in\'egalit\'es dont le premier terme \'etait 
$S_{c}(p\Lambda_{A(c),B(c);s}^{N,M}(\nu,\mu))$ et le dernier $ S_{c}(p_{A(c),B(c);s}^{N,M}[\alpha,\beta])$.
 Maintenant, on sait par (3) que les termes extr\^emes sont \'egaux. Donc toutes les in\'egalit\'es de ces suites sont des \'egalit\'es. Supposons que $k>1$. On prend $c=c^0$. On a alors $a(c^0)=1<k$. Mais on a alors prouv\'e que 
 $$S_{c}(p\Lambda_{A(c),B(c);s}^{N,M}(\nu,\mu))\leq  S_{c}(p_{A(c),B(c);s}^{N,M}[\alpha,\beta])-a(c^0),$$
 ce qui contredit (3). Cela prouve que $k=1$. Supposons donc $k=1$. On veut prouver (ii) et que, pour tout $c\in \{1,...,N+M-1\}$, on a
 
$$ (6)\qquad S_{c}(p\Lambda_{\underline{A}(c),\underline{B}(c);s}^{N-1,M}(\nu',\mu'))=S_{\underline{a}(c)}((\nu'\sqcup\{ 0^{N-n}\})+[\underline{A}(c),\underline{A}(c)+2s-sN]_{s})$$
$$+S_{\underline{b}(c)}((\mu'\sqcup \{0^{M-m}\})+[\underline{B}(c),\underline{B}(c)+s-sM]_{s});$$
  et
  $$(7) \qquad S_{c}(p\Lambda_{\underline{A}(c),\underline{B}(c);s}^{N-1,M}(\nu',\mu'))=S_{c}(p_{\underline{A}(c),\underline{B}(c);s}^{N-1,M}[\alpha^-,\beta]).$$
 
 Soit $c\in \{1,...,c^0-1\}$. Alors $a(c)=0<k=1$. On a alors utilis\'e l'in\'egalit\'e
 $S_{b(c)}(\mu\sqcup\{0^{M-m}\})\leq S_{b(c)}(\mu'\sqcup\{0^{M-m}\})$. Elle devient une \'egalit\'e. De plus, on a $b(c)=c$. Ces \'egalit\'es sont vraies pour tout $c<c^0$ donc les termes de $ \mu\sqcup\{0^{M-m}\}$ et de $\mu'\sqcup\{0^{M-m}\}$ co\^{\i}ncident jusqu'\`a l'ordre $c^0-1$. Si $c^1<c^0$, on a $c^1=b(c^1)=m$, donc $\mu'=\mu$. On a aussi $\underline{\mu}=\mu$ par d\'efinition d'o\`u $\mu'=\underline{\mu}$. Si $c^0<c^1$, on a $c^0-1=b(c^0-1)<m$ et les termes de $\mu$ et $\mu'$ co\^{\i}ncident jusqu'\`a l'ordre $c^0-1$. On a aussi utilis\'e les in\'egalit\'es
 $$S_{a(c)}((\nu'\sqcup\{0^{N-n}\})+[A(c),A(c)+2s-sN]_{s})+S_{b(c)}((\mu'\sqcup\{0^{M-m}\})+[B(c),B(c)+s-sM]_{s})$$
 $$\leq S_{c}(p\Lambda^{N-1,M}_{A(c),B(c);s}(\nu',\mu'))\leq S_{c}(p_{A(c),B(c);s}^{N-1,M}[\alpha^-,\beta]).$$
 Ces in\'egalit\'es deviennent des \'egalit\'es, qui ne sont autres que (6) et (7). 
 
 Supposons $c^0<c^1$ et prenons maintenant $c\in \{c^0,...,c^1-1\}$. On a $a(c)\geq 1=k$ et $b(c)<m$. On a utilis\'e les in\'egalit\'es $\mu_{b(c)+1}\geq \mu'_{b(c)+2}$, $\mu_{b(c)+2}\geq \mu'_{b(c)+3}$,... et $\mu_{m}\geq 0$. Ces in\'egalit\'es deviennent des \'egalit\'es. Donc $\mu_{m}=0$.  En prenant $c=c^0$, on a $b(c)=b(c^0-1)=c^0-1$ et on  voit que  les termes $\mu'_{c^0+1},... \mu'_{m}$ co\^{\i}ncident avec $\mu_{c^0},...,\mu_{m-1}$. Puisqu'on a d\'ej\`a vu que $\mu'$ et $\mu$ co\"{\i}ncidaient jusqu'\`a l'ordre $c^0-1$, seul le terme $\mu'_{c^0}$ reste ind\'etermin\'e. Mais on a $\mu'_{c^0}=S(\mu')-S(\mu)$ et, d'apr\`es l'hypoth\`ese $\mu'\in Q(\alpha_{1},\nu_{1},\mu)$, cette expression vaut $\nu_{1}-\alpha_{1}$. Cela prouve que $\mu'=\underline{\mu}$, d'o\`u (ii). On a aussi utilis\'e les in\'egalit\'es 
  $$S_{a(c)-1}((\nu'\sqcup\{0^{N-n}\})+[A(c)-s,A(c)+s-sN]_{s})$$
  $$+S_{b(c)+1}((\mu'\sqcup\{0^{M-m}\})+[B(c)+e(c)s/2,B(c)+e(c)s/2+s-sM]_{s})$$
  $$\leq S_{c}(p\Lambda^{N-1,M}_{A(c)-s,B(c)+e(c)s/2;s}(\nu',\mu'))\leq S_{c}(p_{A(c)-s,B(c)+e(c)s/2;s}^{N-1,M}[\alpha^-,\beta]).$$
  Elles deviennent des \'egalit\'es qui ne sont autres que (6) et (7). 
  
  Prenons maintenant $c\in \{sup(c^0,c^1),...,N+M-1\}$. On  utilise notre calcul  de 4.3 en $c'=c+1$. On a $a(c')\geq 1=k$ et $b(c')\geq m$.  On a utilis\'e les in\'egalit\'es
  $$S_{a(c')-1}((\nu'\sqcup\{0^{N-n}\})+[A(c')-s,A(c')+s-sN]_{s})+$$
  $$S_{b(c')}((\mu'\sqcup\{0^{M-m}\})+[B(c'),B(c')+s-sM]_{s})\leq S_{c'}(p\Lambda^{N-1,M}_{A(c')-s,B(c');s}(\nu',\mu'))$$
  $$\leq S_{c'}(p_{A(c')-s,B(c');s}^{N-1,M}[\alpha^-,\beta]).$$
  Elles deviennnent des \'egalit\'es qui ne sont autres que (6) et (7). Cela ach\`eve la d\'emonstration. $\square$

\bigskip

\subsection{Avant-derni\`ere \'etape}
On  l\`eve pour cette \'etape les hypoth\`eses  $n\geq1$, $m\geq1$ et $(1,0)<_{I}(1,1)$ (c'est n\'ecessaire pour raisonner par r\'ecurrence dans la preuve du lemme ci-dessous). Soit $(\nu,\mu)\in {\cal P}_{n}\times {\cal P}_{m}$. Supposons qu'il existe des entiers $A(c)$, $B(c)$, $a(c)$, $b(c)$ pour $c=1,...,N+M$, satisfaisant l'hypoth\`ese $(Hyp)(\alpha,\beta;\nu,\mu)$. 

\ass{Lemme}{Sous ces hypoth\`eses, les deux conditions $mult(\alpha,\beta;\nu,\mu)\not=0$ et $(\nu,\mu)\in P(\alpha,\beta)$ sont \'equivalentes. Si elles sont v\'erifi\'ees, on a $mult(\alpha,\beta;\nu,\mu)=1$.}

Preuve. Si $n=0$ ou $m=0$, les deux conditions \'equivalent \`a $(\nu,\mu)=(\alpha,\beta)$ et on a bien $mult(\alpha,\beta;\alpha,\beta)=1$. Supposons $n\geq1$ et $m\geq1$. De nouveau, on ne perd rien \`a supposer $(1,0)<_{I}(1,1)$.   D'apr\`es 4.2(1) et (14), on a
$$mult(\alpha,\beta;\nu,\mu)=\sum_{k=1,...,n}\sum_{\mu'\in Q(\alpha_{1},\nu_{k}+1-k,\mu)}(-1)^{k-1}mult(\alpha^-,\beta;\nu',\mu').$$
Le lemme 4.5 entra\^{\i}ne qu'il y a au plus un terme non nul dans cette somme: on a forc\'ement $k=1$ et $\mu'=\underline{\mu}$ (ce terme ne contribue pas forc\'ement et d'ailleurs, $\underline{\mu}$ n'est pas forc\'ement une partition).  

Supposons $mult(\alpha,\beta;\nu,\mu)\not=0$. Alors il y a au moins un terme qui contribue et on a
$$(1) \qquad mult(\alpha,\beta;\nu,\mu)=mult(\alpha^-,\beta;\nu',\underline{\mu}).$$
Le lemme 4.5 entra\^{\i}ne aussi  que les donn\'ees $\underline{A}(c)$, $\underline{B}(c)$, $\underline{a}(c)$, $\underline{b}(c)$ pour $c=1,...,N+M-1$ v\'erifient l'hypoth\`eses $(Hyp)(\alpha^-,\beta;\nu',\underline{\mu})$. On applique l'\'enonc\'e par r\'ecurrence: on a donc $(\nu',\underline{\mu})\in P(\alpha^-,\beta)$ et $mult(\alpha^-,\beta;\nu',\underline{\mu})=1$. Avec (1), on en d\'eduit $mult(\alpha,\beta;\nu,\mu)=1$. Par d\'efinition, on a 

pour $c=1,...,c^0-1$, $\underline{b}(c)=c$;

si $c^0<c^1$, 
$\underline{b}(c^0)=b(c^0)+1=c^0$;

  si $c^1< c^0$, auquel cas $c^1=m$, $\underline{b}(c^1)=   b(c^1)=m$.
  
  Alors   $(Hyp)(\alpha^-,\beta;\nu',\underline{\mu})$ entra\^{\i}ne que $(\nu',\underline{\mu})$ v\'erifie les hypoth\`eses de 2.2 o\`u $c=c^0$ si $c^0<c^1$ et $c=m$ si $c^1<c^0$. 
  En se rappelant que, d'apr\`es 2.2(1), on a $P^{b[c]}(\alpha^-,\beta)=P^{b[m]}(\alpha^-,\beta)$ si $c\geq m$, le lemme 2.2 implique  en tout cas que $(\nu',\underline{\mu})\in P^{b[c^0]}(\alpha^-,\beta)$. Le lemme  2.3 implique que $\iota_{c^0,\alpha_{1}}(\nu',\underline{\mu})$ appartient \`a $P(\alpha,\beta)$. Or la remarque qui suit le lemme 4.5 dit que $\iota_{c^0,\alpha_{1}}(\nu',\underline{\mu})=(\nu,\mu)$. Donc $(\nu,\mu)\in P(\alpha,\beta)$, ce qui d\'emontre un sens de l'\'equivalence de l'\'enonc\'e. 
  
Au lieu de supposer que  $mult(\alpha,\beta;\nu,\mu)\not=0$, on suppose maintenant que $(\nu,\mu)\in P(\alpha,\beta)$. D'apr\`es le lemme 2.3, on peut fixer $c\in \{1,...,m+1\}$ et $(\nu'',\mu'')\in P(\alpha^-,\beta)$ tel que $(\nu,\mu)=\iota_{c,\alpha_{1}}(\nu'',\mu'')$. Il r\'esulte de la d\'efinition de 2.3 que $\nu''=\nu'$. On voit aussi que $(\nu',\mu'')$ appartient \`a $Q(\alpha_{1},\nu_{1},\mu)$. Alors le lemme 4.5 entra\^{\i}ne de nouveau   que les donn\'ees $\underline{A}(c)$, $\underline{B}(c)$, $\underline{a}(c)$, $\underline{b}(c)$ pour $c=1,...,N+M-1$ v\'erifient l'hypoth\`ese $(Hyp)(\alpha^-,\beta;\nu',\mu'')$. On applique l'\'enonc\'e par r\'ecurrence: on a donc $mult(\alpha^-,\beta;\nu',\mu'')=1$. Comme on l'a vu au d\'ebut de la preuve, le terme $(\nu',\mu'')$ est alors l'unique terme contribuant \`a $mult(\alpha,\beta;\nu,\mu)$. D'o\`u $mult(\alpha,\beta;\nu,\mu)=1$, ce qui ach\`eve la d\'emonstration. $\square$

\bigskip
 \subsection{Fin de la preuve de la proposition}
On revient aux hypoth\`eses de 4.2.  Supposons $p\Lambda_{A,B;s}^{N,M}(\nu,\mu)=p_{A,B;s}^{N,M}[\alpha,\beta]$. Pour tout $c=1,...,N+M$, on pose $A(c)=A$ et $B(c)=B$. Alors l'\'egalit\'e (3) de 4.5 est satisfaite. On voit que l'on peut trouver des entiers $a(c)$ et $b(c)$ v\'erifiant (1) et (2) de ce paragraphe (par exemple, on choisit $a(c) $ le plus grand possible pour que (1) soit v\'erifi\'e; alors (2) l'est aussi). L'hypoth\`ese $(Hyp)(\alpha,\beta;\nu,\mu)$ est donc v\'erifi\'ee. Le lemme 4.6 implique que les conditions $mult(\alpha,\beta;\nu,\mu)\not=0$ et $(\nu,\mu)\in P(\alpha,\beta)$ sont \'equivalentes. Mais, compte tenu de l'hypoth\`ese $p\Lambda_{A,B;s}^{N,M}(\nu,\mu)=p_{A,B;s}^{N,M}[\alpha,\beta]$, la condition $(\nu,\mu)\in P(\alpha,\beta)$ est \'equivalente \`a $(\nu,\mu)\in P_{A,B;s}(\alpha,\beta)$. Cela prouve le (ii) de la proposition. Pour le (iii), l'hypoth\`ese $(\nu,\mu)\in P_{A,B;s}(\alpha,\beta)$ est plus forte que celle du (ii). Comme on vient de le voir, on peut appliquer le lemme 4.6, qui entra\^{\i}ne $mult(\alpha,\beta;\nu,\mu)=1$. Cela ach\`eve la preuve de la proposition. $\square$

 \bigskip

\section{Composantes extr\'emales d'une repr\'esentation construite par Lusztig}
\bigskip

\subsection{La correspondance de Springer g\'en\'eralis\'ee}

Pour $n\in {\mathbb N}$, on note ${\cal P}(n)$ l'ensemble des partitions de $n$, c'est-\`a-dire l'ensemble des partitions $\lambda=(\lambda_{1},...,\lambda_{t})$ telles que $S(\lambda)=n$.

{\bf Remarque.} On prendra garde de ne pas confondre  l'ensemble ${\cal P}(n)$ avec l'ensemble ${\cal P}_{n}$ des paragraphes pr\'ec\'edents:  pour le premier ensemble, $n$ est la somme des termes de la partition tandis que, pour le second,  $n$ est son nombre de termes. D'autre part, il convient ici d'identifier deux partitions ne diff\'erant que par des termes nuls: $(\lambda_{1},...,\lambda_{t})$ est la m\^eme partition que $(\lambda_{1},...,\lambda_{t},0,...,0)$. 

\bigskip

On note ${\cal P}_{2}(n)$ l'ensemble des couples $(\alpha,\beta)$ de partitions telles que $S(\alpha)+S(\beta)=n$. On note $W_{n}$ le groupe de Weyl d'un syst\`eme de racines de type $C_{n}$, avec la convention $W_{0}=\{1\}$. On note $W_{n}^{\vee}$ l'ensemble des repr\'esentations irr\'eductibles de $W_{n}$. Il est bien connu que $W_{n}^{\vee}$ est param\'etr\'e par ${\cal P}_{2}(n)$. Pour $(\alpha,\beta)\in {\cal P}_{2}(n)$, on note $\rho_{\alpha,\beta}$ la repr\'esentation de $W_{n}$ param\'etr\'ee par ce couple. En particulier, $\rho_{(n),\emptyset}$ est la repr\'esentation triviale  et $\rho_{\emptyset,(1^n)}$ est le caract\`ere signature usuel, que l'on note $sgn$. 

Soit $N\in {\mathbb N}$.  Notons ${\cal P}^{symp}(2N)$ l'ensemble des partitions symplectiques de $2N$, c'est-\`a-dire l'ensemble des partitions $\lambda \in {\cal P}(2N)$ telles  que  $mult_{\lambda}(i)$ est pair pour tout entier $i$ impair.
Pour $\lambda=(\lambda_{1},...,\lambda_{t})\in {\cal P}^{symp}(2N)$, on note $Jord(\lambda)$ l'ensemble des entiers $i\geq1$ tels que $mult_{\lambda}(i)>0$ et $Jord^{bp}(\lambda)$ le sous-ensemble des entiers pairs. On note $\boldsymbol{{\cal P}^{symp}}(2N)$ l'ensemble des couples $(\lambda,\epsilon)$ o\`u $\lambda\in {\cal P}^{symp}(2N)$ et $\epsilon\in \{\pm 1\}^{Jord^{bp}(\lambda)}$.

La correspondance de Springer g\'en\'eralis\'ee, d\'efinie par Springer et Lusztig, \'etablit une bijection entre $\boldsymbol{{\cal P}^{symp}}(2N)$ et l'ensemble
$\bigcup_{k} W_{N-k(k+1)/2}^{\vee}$, o\`u $k$ parcourt les \'el\'ements de ${\mathbb N}$ tels que $k(k+1)\leq 2N$. Ou encore avec l'ensemble $\bigcup_{k}{\cal P}_{2}(N-k(k+1)/2)$. Rappelons la d\'efinition de cette bijection. On fixe un entier $r$ assez grand 

{\bf Remarque.} On peut prendre $r= N$  mais il peut \^etre utile de laisser \`a $r$ une certaine libert\'e. La construction qui suit ne d\'epend pas du choix de $r$.
\bigskip

Soit $k\in {\mathbb N}$ tel que $k(k+1)\leq 2N$ et soit $(\alpha,\beta)\in {\cal P}_{2}(N-k(k+1)/2)$. On d\'efinit un couple de partitions $(A_{\alpha,\beta},B_{\alpha,\beta})$ par

si $k$ est pair, $A_{\alpha,\beta}=\alpha+[2r+k,0]_{2}$, $B_{\alpha,\beta}=\beta+[2r-1-k,1]_{2}$;

si $k$ est impair, $A_{\alpha,\beta}=\beta+[2r-1-k,0]_{2}$, $B_{\alpha,\beta}=\alpha+[2r+k,1]_{2}$.

{\bf Remarque.} Pour d\'efinir ces sommes, on  compl\`ete $\alpha$ et $\beta$ 
 par suffisamment de $0$ pour que $\alpha$ ait  $r+[k/2]+1$ termes et que $\beta$ ait $r-[(k+1)/2]$ termes, o\`u, pour  tout r\'eel $x$, $[x]$ d\'esigne sa partie enti\`ere.  
  
 \bigskip
Soit $(\lambda,\epsilon)\in \boldsymbol{{\cal P}^{symp}}(2N)$. Dans la partition $\lambda+[2r-1,0]_{1}$, il y a autant d'entiers pairs que d'impairs. On note les pairs $2z_{1}>...>2z_{r}$ et les
impairs $2z'_{1}+1>...>2z'_{r}+1$.  Cela d\'efinit des partitions $z$ et $z'$. On pose
$A^{\sharp}_{\lambda}=(z'+[r+1,2]_{1})\sqcup\{0\}$, $B^{\sharp}_{\lambda}=z+[r,1]_{1}$. 
On calcule facilement
$$(1)(a)\qquad A^{\sharp}_{\lambda,j}=\left\lbrace\begin{array}{cc}\lambda_{2j-1}/2+2r+2-2j,&\rm{ si\,\, }\lambda_{2j-1}\rm{\,\, est\,\,pair,}\\ (\lambda_{2j}-1)/2+2r+2-2j,&\rm{ si \,\,}\lambda_{2j}\rm{\,\, est \,\,impair\,\,et}\,\,S_{2j}(\lambda)\rm{\,\, est\,\,pair, }\\ (\lambda_{2j-2}+1)/2+2r+2-2j,&\rm{ si \,\,}j\geq2,\,\,\lambda_{2j-2}\rm{ \,\,est \,\,impair\,\,et\,\,}\\ & S_{2j-2}(\lambda)\rm{ \,\,est\,\,impair,}\\ \end{array}\right.$$
pour $j=1,...,r+1$ et
$$(1)(b)\qquad B^{\sharp}_{\lambda,j}=\left\lbrace\begin{array}{cc}\lambda_{2j}/2+2r+1-2j,&\rm{ si\,\, }\lambda_{2j}\rm{\,\, est\,\,pair,}\\ (\lambda_{2j-1}+1)/2+2r+1-2j,&\rm{ si \,\,}\lambda_{2j-1}\rm{\,\, est \,\,impair\,\,et}\,\,S_{2j-1}(\lambda)\rm{\,\, est\,\,impair, }\\ (\lambda_{2j+1}-1)/2+2r+1-2j,&\rm{ si \,\,} \lambda_{2j+1}\rm{ \,\,est \,\,impair\,\,et\,\,} S_{2j+1}(\lambda)\rm{ \,\,est\,\,pair,}\\ \end{array}\right.$$
pour $j=1,...,r$. On v\'erifie que les divers cas ci-dessus sont exclusifs l'un de l'autre.

Consid\'erons l'ensemble $(A^{\sharp}_{\lambda}\cup B^{\sharp}_{\lambda})-(A^{\sharp}_{\lambda}\cap B^{\sharp}_{\lambda})$, c'est-\`a-dire l'ensemble des entiers qui interviennent dans $A^{\sharp}_{\lambda}$ ou $B^{\sharp}_{\lambda}$ mais pas dans les deux. Notons $D_{\lambda}$ l'ensemble des sous-ensembles de $(A^{\sharp}_{\lambda}\cup B^{\sharp}_{\lambda})-(A^{\sharp}_{\lambda}\cap B^{\sharp}_{\lambda})$ form\'es d'entiers cons\'ecutifs, maximaux pour cette propri\'et\'e, et ne contenant pas $0$. On s'aper\c{c}oit qu'il y a une bijection croissante $i\mapsto \Delta_{i}$ de $Jord^{bp}(\lambda)$ sur $D_{\lambda}$. On pose
$$A_{\lambda,\epsilon}=(A^{\sharp}_{\lambda}-\bigcup_{i\in Jord^{bp}(\lambda); \epsilon_{i}=-1}(A^{\sharp}_{\lambda}\cap \Delta_{i}))\cup(\bigcup_{i\in Jord^{bp}(\lambda); \epsilon_{i}=-1}(B^{\sharp}_{\lambda}\cap \Delta_{i})),$$
$$B_{\lambda,\epsilon}=(B^{\sharp}_{\lambda}-\bigcup_{i\in Jord^{bp}(\lambda); \epsilon_{i}=-1}(B^{\sharp}_{\lambda}\cap \Delta_{i}))\cup(\bigcup_{i\in Jord^{bp}(\lambda); \epsilon_{i}=-1}(A^{\sharp}_{\lambda}\cap \Delta_{i})).$$

On montre qu'il existe un unique entier $k\in {\mathbb N}$ tel que $k(k+1)\leq 2N$ et un unique couple $(\alpha,\beta)\in {\cal P}_{2}(N-k(k+1)/2)$ tels que $(A_{\alpha,\beta},B_{\alpha,\beta})=(A_{\lambda,\epsilon},B_{\lambda,\epsilon})$. La correspondance de Springer g\'en\'eralis\'ee associe \`a $(\lambda,\epsilon)$ le couple $(k,\rho_{\alpha,\beta})$.  On note  ce couple $(k(\lambda,\epsilon),\rho(\lambda,\epsilon))$. 

L'entier $k(\lambda,\epsilon)$ se calcule facilement. Notons $i_{1}>...>i_{m}$ les \'el\'ements $i\in Jord^{bp}(\lambda)$ tels que $mult_{\lambda}(i)$ soit impair. Posons
$$M(\lambda,\epsilon)=\sum_{l=1,...,m; \epsilon_{i_{l}}=-1}(-1)^l.$$
Alors $k(\lambda,\epsilon)=2M(\lambda,\epsilon)$ si $M(\lambda,\epsilon)\geq0$ et $k(\lambda,\epsilon)=-2M(\lambda,\epsilon)-1$ si $M(\lambda,\epsilon)<0$.

Consid\'erons le cas particulier o\`u $Jord(\lambda)=Jord^{bp}(\lambda)$, c'est-\`a-dire que $\lambda$ n'a que des termes pairs. Ecrivons $\lambda=(\lambda_{1},...,\lambda_{2r+1})$, avec $\lambda_{2r+1}=0$. On peut consid\'erer que $\epsilon$ est une fonction d\'efinie sur l'ensemble d'indices $\{1,...,2r+1\}$: si $j\in \{1,...,2r+1\}$ est tel que $\lambda_{j}\not=0$, on pose $\epsilon(j)=\epsilon_{\lambda_{j}}$; si $\lambda_{j}=0$, on pose $\epsilon(j)=1$. On prendra garde de distinguer $\epsilon(j)$ o\`u $j$ est un indice et $\epsilon_{i}$ o\`u $i$ est un terme de $\lambda$. 
Avec la recette ci-dessus, on calcule
$$A_{\lambda,\epsilon}=\{\lambda_{j}/2+2r+1-j; j=1,...,2r+1, \epsilon(j)=(-1)^{j+1}\},$$
$$B_{\lambda,\epsilon}=\{\lambda_{j}/2+2r+1-j; j=1,...,2r+1, \epsilon(j)=(-1)^j\}.$$
 On constate que $A_{\lambda,\epsilon}\sqcup B_{\lambda,\epsilon}$ est "sans multiplicit\'es", c'est-\`a-dire que tout entier y intervient avec multiplicit\'e au plus $1$.
Remarquons qu'avec les m\^emes hypoth\`eses, on a

(2) $M(\lambda,\epsilon)=\sum_{j=1,...,2r+1}(-1)^{j+1}\frac{\epsilon(j)-1}{2}$.

Preuve. Soit $i\in Jord^{bp}(\lambda)\sqcup \{0\}$, notons $\{j_{i}^-,...,j_{i}^+\}$ le sous-ensemble des $j\in \{1,...,2r+1\}$ tels que $\lambda_{j}=i$. Pour tous ces $j$, on a $\epsilon(j)=\epsilon_{i}$ et la contribution de ces $j$  au membre de droite de (2) est
$$\frac{\epsilon_{i}-1}{2}\sum_{j=j_{i}^-,...,j_{i}^+}(-1)^{j+1}=\left\lbrace\begin{array}{cc}0,&\text{ si }j_{i}^+-j_{i}^- \text{ est impair},\\ \frac{\epsilon_{i}-1}{2}(-1)^{j_{i}^-+1},&\text{ si }j_{i}^+-j_{i}^-\text{ est pair}.\\ \end{array}\right.$$
Si $\epsilon_{i}=1$,  on voit que  l'intervalle $\{j_{i}^-,...,j_{i}^+\}$ ne contribue pas au membre de droite de (2). L'entier $i$ ne contribue pas non plus \`a $M(\lambda,\epsilon)$. Remarquons que cela traite le cas $i=0$. Supposons $i>0$. Alors $j^+-j^-=mult_{\lambda}(i)-1$. Si $mult_{\lambda}(i)$ est pair, ou si $mult_{\lambda}(i)$ est pair et $\epsilon_{i}=1$, l'intervalle $\{j_{i}^-,...,j_{i}^+\}$ ne contribue pas au membre de droite de (1). L'entier $i$ ne contribue pas non plus \`a $M(\lambda,\epsilon)$. Si $mult_{\lambda}(i)$ est pair et $\epsilon_{i}=-1$, on a $i=i_{l}$ pour un $l\in \{1,...,m\}$ et on v\'erifie par r\'ecurrence sur $l$ que $j_{i}^-$ est de la m\^eme parit\'e que $l$. Alors l'intervalle $\{j_{i}^-,...,j_{i}^+\}$ contribue  au membre de droite de (2) par $(-1)^l$ et $i$ contribue de la m\^eme fa\c{c}on \`a $M(\lambda,\epsilon)$. Cela prouve (2).

\bigskip

\subsection{Un r\'esultat de Shoji}
Dans les quatre premiers paragraphes, nous avons d\'efini divers termes d\'ependant d'un ordre fix\'e $<_{I}$ sur un ensemble d'indices $I$. Pour plus de pr\'ecision, nous ins\'ererons maintenant cet ordre dans la notation.

Soit $(\lambda,\epsilon)\in \boldsymbol{{\cal P}^{symp}}(2N)$. On suppose que $\lambda$ n'a que des termes pairs. Notons $k=k(\lambda,\epsilon)$ et $(\alpha,\beta)$ l'\'el\'ement de ${\cal P}_{2}(N-k(k+1)/2)$ qui correspond \`a $(\lambda,\epsilon)$ par la correspondance de Springer g\'en\'eralis\'ee. Choisissons un entier $r$ assez grand, posons   $n=r+[k/2]+1$, $m=r-[k/2]$ et consid\'erons que $\alpha$ a $n$ termes et que $\beta$ a $m$ termes. Pour $(\nu,\mu)\in {\cal P}_{n}\times {\cal P}_{m}$, on a d\'efini une multiplicit\'e $mult(\alpha,\beta;\nu,\mu)$ en 4.1, que l'on note maintenant $mult(<_{I};\alpha,\beta;\nu,\mu)$ puisqu'elle  d\'epend du choix  de l'ordre $<_{I}$  sur l'ensemble d'indices $I=(\{1,...,n\}\times\{0\})\cup(\{1,...,m\}\times\{1\})$. Nous allons d\'efinir un  ordre, d\'ependant lui-m\^eme de $(\alpha,\beta)$, que nous noterons $<_{I,\alpha,\beta}$. Posons $\Lambda=A_{\alpha,\beta}\sqcup B_{\alpha,\beta}$. 
D'apr\`es les formules du paragraphe 5.1, on a l'\'egalit\'e
$$ \Lambda=p\Lambda_{A,B;2}^{n,m}(\alpha,\beta),$$
o\`u $A=2r+k$ et $B=2r-k-1$. 

{\bf Remarque.} Il y a ici une difficult\'e de notation. Si $k$ est pair, on a $A_{\alpha,\beta}=\alpha+[A,A+2-2n]_{2}$ et $B_{\alpha,\beta}=\beta+[B,B+2-2m]_{2}$ mais, si $k$ est impair, on doit permuter ces termes: $A_{\alpha,\beta}=\beta+[B,B+2-2m]_{2}$ et $B_{\alpha,\beta}=\alpha+[A,A+2-2n]_{2}$.

\bigskip
On a aussi $ \Lambda=A_{\lambda,\epsilon}\sqcup B_{\lambda,\epsilon}$ par d\'efinition et on a remarqu\'e dans le paragraphe 5.1 que cette partition \'etait sans multiplicit\'es (parce que $\lambda$ est \`a termes pairs). Ecrivons cette partition $\Lambda=(\Lambda_{1},...,\Lambda_{2r+1})$. Pour $j=1,...,n$, il y a un unique $a_{j}\in \{1,...,2r+1\}$ tel que $\alpha_{j}+A+2-2j
=\Lambda_{a_{j}}$. Pour $l=1,...,m$, il y a un unique $b_{l}\in \{1,...,2r+1\}$ tel que $\beta_{l}+B+2-2l
=\Lambda_{b_{l}}$. On d\'efinit l'ordre sur $I$ de sorte que $(j,0)<_{I,\alpha,\beta}(l,1)$ si et seulement si $a_{j}< b_{l}$. 

Soit $(\lambda',\epsilon')\in \boldsymbol{{\cal P}^{symp}}(2N)$, supposons $k(\lambda',\epsilon')=k$. Dans l'introduction, on a d\'efini la multiplicit\'e  $\mathfrak{mult}(\lambda,\epsilon;\lambda',\epsilon')$. Notons $(\alpha',\beta')\in {\cal P}_{2}(N-k(k+1)/2)$ le couple correspondant \`a $(\lambda',\epsilon')$. On consid\`ere comme ci-dessus que $\alpha'$ a $n$ termes et que $\beta'$ en a $m$.

\ass{Proposition (Shoji)}{On rappelle que $\lambda$ est \`a termes pairs. Pour $(\lambda',\epsilon')$ comme ci-dessus, on a l'\'egalit\'e
$$\mathfrak{mult}(\lambda,\epsilon;\lambda',\epsilon')=mult(<_{I,\alpha,\beta};\alpha,\beta;\alpha',\beta').$$}

 Shoji d\'emontre ce r\'esultat en \cite{Shoji} p.685. En fait, ses hypoth\`eses sont un peu plus contraignantes que les n\^otres aussi allons-nous reprendre sa d\'emonstration dans les deux paragraphes suivants. Notre situation est plus simple que celle consid\'er\'ee  par Shoji car il consid\`ere des familles de partitions $\boldsymbol{\alpha}$ \`a $e\geq 1$ termes tandis que nous nous limitons aux paires de partitions, c'est-\`a-dire au cas $e=2$. Nous continuerons donc \`a noter nos familles comme des paires et n'utiliserons pas la notation $\boldsymbol{\alpha}$ de Shoji. De plus, les exposants $\pm$ qui interviennent  parfois dans \cite{Shoji} sont inutiles dans le cas $e=2$ et nous les faisons dispara\^{\i}tre.  

 \bigskip
 
 \subsection{D\'emonstration de la proposition 5.2, premi\`ere \'etape}

La d\'emonstration de Shoji s'applique pourvu que sa proposition 6.1  soit v\'erifi\'ee. C'est-\`a-dire que l'on doit prouver l'\'egalit\'e $Q_{\Lambda}=R_{\Lambda}$. Ainsi que Shoji l'explique, il suffit de prouver l'assertion suivante. Ecrivons $R_{\Lambda}$ comme combinaison lin\'eaire de fonctions de Schur $s_{\alpha',\beta'}$, \`a coefficients dans ${\mathbb Q}(t)$:
$$R_{\Lambda}=\sum_{(\alpha',\beta')\in {\cal P}_{n}\times {\cal P}_{m}}r_{\Lambda}(\alpha',\beta')s_{\alpha',\beta'}.$$
On doit prouver que, pour tout couple $(\alpha',\beta')\in {\cal P}_{n}\times {\cal P}_{m} $  tel que $r_{\Lambda}(\alpha',\beta')\not=0$, on a 
$$(1) \qquad p\Lambda_{A,B;2}^{n,m}(\alpha',\beta')\leq \Lambda.$$ 

On commence par calculer les coefficients $r_{\Lambda}(\alpha',\beta')$.
On note $(v^{0}_{i})_{i=1,...,n}$ la base canonique de ${\mathbb Z}^n$ et $(v^{1}_{i})_{i=1,...,m}$ celle de ${\mathbb Z}^m$. On identifie naturellement les \'el\'ements de ces bases \`a des \'el\'ements de ${\mathbb Z}^n\times {\mathbb Z}^m$.
 Pour simplifier, nous notons simplement $<_{I}$ l'ordre $<_{I,\alpha,\beta}$.  Notons $i^{0}_{max}$ le plus grand indice $i\in \{1,...,n\}$ tel que $\alpha_{i}\not=0$ et $i^{1}_{max}$ le plus grand indice $i\in \{1,...,m\}$ tel  que $\beta_{i}\not=0$. Notons $\nu_{0}$ le plus grand des deux \'el\'ements $(i^{0}_{max},0)$ et $(i^{1}_{max},1)$ pour l'ordre $<_{I}$. Notons $J$ le sous-ensemble de ${\mathbb Z}^n\times {\mathbb Z}^m$ form\'e des  \'el\'ements $v_{i}^{e}-v_{j}^{f}$, pour les couples $((i,e),(j,f))\in I^2$ tels que $e\not=f$, $(i,e)<_{I}(j,f)$ et $(i,e)\leq_{I} \nu_{0}$. Notons $K$ le sous-ensemble de ${\mathbb Z}^n\times {\mathbb Z}^m$ form\'e des  \'el\'ements $v_{i}^{e}-v_{j}^{f}$, pour les couples
  $((i,e),(j,f))\in I^2$ tels que  $e=f$, $(i,e)<_{I}(j,f)$ et $\nu_{0}<_{I}(i,e)$. On pose $\Delta_{n}=(n-1,n-2,...,0)=[n-1,0]_{1}$ et on d\'efinit de m\^eme $\Delta_{m}$. Le couple $(\alpha+\Delta_{n},\beta+\Delta_{m})$ est un \'el\'ement de ${\mathbb Z}^n\times {\mathbb Z}^m$. Le groupe $\mathfrak{ S}_{n}$ agit sur ${\mathbb Z}^n$: pour $x=(x_{1},...,x_{n})\in {\mathbb Z}^n$ et $\sigma\in \mathfrak{ S}_{n}$, on pose $x^{\sigma}=(x_{\sigma 1},... ,x_{\sigma n})$. De m\^eme, $\mathfrak{S}_{m}$ agit sur ${\mathbb Z}^m$. Soit $(\alpha',\beta')\in {\cal P}_{n}\times {\cal P}_{m}$. Pour $d\in {\mathbb N}$, posons
  $$(2) \qquad r^{d}_{\Lambda}(\alpha',\beta')=\sum_{\sigma\in \mathfrak{S}_{n},\tau\in \mathfrak{S}_{m}}sgn(\sigma)sgn(\tau)$$
  $$\vert \{X\subset J\cup K; \vert X\vert =d, (\alpha+\Delta_{n},\beta+\Delta_{m})-\sum_{x\in X}x=((\alpha'+\Delta_{n})^{\sigma},(\beta'+\Delta_{m})^{\tau})\}\vert .$$
D'apr\`es \cite{Shoji} 3.13.1, on a l'\'egalit\'e
$$(3) \qquad r_{\Lambda}(\alpha',\beta')=v_{\alpha,\beta}(t)^{-1}\sum_{d\in {\mathbb N}}(-1)^dr^{d}_{\Lambda}(\alpha',\beta')t^d,$$
o\`u $v_{\alpha,\beta}(t)$ est un certain polyn\^ome non nul. 
On peut r\'ecrire l'expression (2) sous la forme
$$(4) \qquad  r^{d}_{\Lambda}(\alpha',\beta')=\sum_{\sigma\in \mathfrak{S}_{n},\tau\in \mathfrak{S}_{m}}sgn(\sigma)sgn(\tau)\sum_{X_{J}\subset J; \vert X_{J}\vert \leq d}\vert {\cal X}_{K}(X_{J},\sigma,\tau)\vert ,$$
o\`u
$${\cal X}_{K}(X_{J},\sigma,\tau)=\{X_{K}\subset K; \vert X_{K}\vert =d-\vert X_{J}\vert ,$$
$$ (\alpha+\Delta_{n},\beta+\Delta_{m})-\sum_{x\in X_{K}}x=((\alpha'+\Delta_{n})^{\sigma},(\beta'+\Delta_{m})^{\tau})+\sum_{x\in X_{J}}x\}.$$
 Introduisons le plus petit nombre $k\in \{0,...,n\}$ tel que $(n-k,0)\leq_{I}\nu_{0}$  et le plus petit nombre $h\in \{0,...,m\}$ tel que $(m-h,1)\leq \nu_{0}$. Pour $y=(y_{1},...,y_{n})\in {\mathbb Z}^n$, on pose $y_{\leq n-k}=(y_{1},...,y_{n-k})$ et $y_{>n-k}=(y_{n-k+1},...,y_{n})$. Pour $y\in {\mathbb Z}^m$, on d\'efinit de m\^eme $y_{\leq m-h}$ et $y_{>m-h}$. Pour $y=(y^0,y^1)\in {\mathbb Z}^n\times {\mathbb Z}^m$, on pose $p_{\leq}(y)=(y^0_{\leq n-k},y^1_{\leq m-h})$ et $p_{>}(y)=(y^0_{>n-k},y^1_{>m-h})$. Remarquons que, pour $y\in K$, on a $p_{\leq}(y)=0$. Donc $K$  s'identifie par $p_{>}$ \`a un sous-ensemble de ${\mathbb Z}^{k}\times {\mathbb Z}^{h}$. De plus, on a $p_{>}(\alpha+\Delta_{n},\beta+\Delta_{m})=(\Delta_{k},\Delta_{h})$.  Pour $J$, $\sigma$ et $\tau$ comme ci-dessus, l'ensemble ${\cal X}_{K}(X_{J},\sigma,\tau)$ est donc vide si la condition
 
 (5) $p_{\leq}(\alpha+\Delta_{n},\beta+\Delta_{m})=p_{\leq}(((\alpha'+\Delta_{n})^{\sigma},(\beta'+\Delta_{m})^{\tau})+\sum_{x\in X_{J}}x)$
 
 \noindent n'est pas v\'erifi\'ee. Si elle l'est, ${\cal X}_{K}(X_{J},\sigma,\tau)$ s'identifie \`a l'ensemble des $X_{K}\subset K$ tels que $\vert X_{K}\vert =d-\vert X_{J}\vert $ et
 $$(\Delta_{k},\Delta_{h})-\sum_{x\in X_{K}}x=p_{>}(((\alpha'+\Delta_{n})^{\sigma},(\beta'+\Delta_{m})^{\tau})+\sum_{x\in X_{J}}x).$$
 Notons $ \mathfrak{S}[X_{J}]$ l'ensemble des $ (\sigma,\tau)\in \mathfrak{S}_{n}\times \mathfrak{ S}_{m}$ tels que (5) soit v\'erifi\'ee. Supposons que $ (\sigma,\tau)$ appartienne \`a cet ensemble.  
  Par coh\'erence avec ce qui pr\'ec\`ede, un  \'el\'ement de ${\mathbb Z}^{k}$ sera not\'e  $y=(y_{n-k+1},...,y_{n})$. On note $({\mathbb Z}^k)^+$ l'ensemble des $y=(y_{n-k+1},...,y_{n})\in {\mathbb Z}^{k}$ tels que $y_{n-k+1}\geq...\geq y_{n}$. Le groupe $\mathfrak{S}_{k}$ agit sur ${\mathbb Z}^{k}$. Comme pr\'ec\'edemment, on note $stab(y)$ le fixateur dans ce groupe d'un \'el\'ement $y\in {\mathbb Z}^{k}$. Tout \'el\'ement $y\in {\mathbb Z}^{k}$ s'\'ecrit $y=z^{\sigma'}$ pour un unique \'el\'ement $z\in ({\mathbb Z}^k)^+$ et un \'el\'ement $\sigma'\in \mathfrak{S}_{k}$ dont l'image dans $stab(z)\backslash \mathfrak{S}_{k}$ est uniquement d\'etermin\'ee. De m\^emes d\'efinitions s'appliquent \`a ${\mathbb Z}^h$. Alors 
       ${\cal X}_{K}(X_{J},\sigma,\tau)$ se d\'ecompose en  la r\'eunion sur les $z\in ({\mathbb Z}^k)^+$, $\sigma'\in stab(z)\backslash \mathfrak{S}_{k}$, $y\in ({\mathbb Z}^h)^+$, $\tau'\in stab(y)\backslash \mathfrak{S}_{h}$ du sous-ensemble des $X_{K}\in {\cal X}_{K}(X_{J},\sigma,\tau)$ tels que 
$$(6) \qquad (\Delta_{h},\Delta_{k})-\sum_{x\in X_{K}}x=(z^{\sigma'},y^{\tau'}).$$
Ce sous-ensemble est non vide seulement si 
$$(7) \qquad (z^{\sigma'},y^{\tau'})=p_{>}(((\alpha'+\Delta_{n})^{\sigma},(\beta'+\Delta_{m})^{\tau})+\sum_{x\in X_{J}}x).$$
Cette condition ne d\'epend pas de l'ensemble $X_{K}$. Pour $\sigma'\in \mathfrak{S}_{k}$ et $\tau'\in \mathfrak{S}_{h}$, notons ${\cal Z}(\sigma',\tau',\sigma,\tau,X_{J})$ l'ensemble des $(z,y)\in ({\mathbb Z}^k)^+\times ({\mathbb Z}^h)^+$ qui v\'erifient (7). D'autre part, pour $(z,y)\in ({\mathbb Z}^k)^+\times ({\mathbb Z}^h)^+$, $\sigma'\in \mathfrak{S}_{k}$,  $\tau'\in \mathfrak{S}_{h}$ et $f\in {\mathbb N}$, notons ${\cal X}_{K}(f,z,y,\sigma',\tau')$ l'ensemble des $X_{K}\subset K$ qui v\'erifient (6) et ont $f$ \'el\'ements.  Si $(z,y)\not\in {\cal Z}(\sigma',\tau',\sigma,\tau,X_{J})$, l'ensemble des $X_{K}\in {\cal X}_{K}(X_{J},\sigma,\tau)$ qui v\'erifient (6) est vide. Si $(z,y)\in {\cal Z}(\sigma',\tau',\sigma,\tau,X_{J})$, cet ensemble est \'egal \`a ${\cal X}_{K}(d-\vert X_{J}\vert,z,y,\sigma',\tau')$. En rassemblant ces calculs, on obtient la formule suivante
$$(8) \qquad r^{d}_{\Lambda}(\alpha',\beta')=\sum_{X_{J}\subset J, \vert X_{J}\vert \leq d}\sum_{(\sigma,\tau)\in \mathfrak{S}[X_{J}]}sgn(\sigma)sgn(\tau)\sum_{\sigma'\in \mathfrak{S}_{k},\tau'\in \mathfrak{S}_{h}}\sum_{(z,y)\in {\cal Z}(\sigma',\tau',\sigma,\tau,X_{J})}$$
$$\vert stab(z)\vert ^{-1}\vert stab(y)\vert ^{-1}\vert {\cal X}_{K}(d-\vert X_{J}\vert,z,y,\sigma',\tau')\vert .$$
On plonge naturellement $\mathfrak{S}_{k}$ dans $\mathfrak{S}_{n}$: le groupe $\mathfrak{S}_{k}$ agit sur les $k$ derni\`eres coordonn\'ees. De m\^eme, on plonge $\mathfrak{S}_{h}$ dans $\mathfrak{S}_{m}$. On s'aper\c{c}oit que l'action de $\mathfrak{S}_{k}\times \mathfrak{S}_{h}$ conserve l'ensemble $J$. On en d\'eduit une action $X_{J}\mapsto X_{J}^{\sigma',\tau'}$ sur l'ensemble des  sous-ensembles de $X_{J}$.   
La relation (7) \'equivaut \`a 
$$(z,y)=p_{>}(((\alpha'+\Delta_{n})^{\sigma{\sigma'}^{-1}},(\beta'+\Delta_{m})^{\tau{\tau'}^{-1}})+\sum_{x\in X_{J}^{{\sigma'}^{-1},{\tau'}^{-1}}}x).$$
Autrement dit ${\cal Z}(\sigma',\tau',\sigma,\tau,X_{J})={\cal Z}(1,1,\sigma{\sigma'}^{-1},\tau{\tau'}^{-1},X_{J}^{{\sigma'}^{-1},{\tau'}^{-1}})$. D'autre part,  l'action de $\mathfrak{S}_{k}\times \mathfrak{ S}_{h}$ ne change pas l'image d'un \'el\'ement de ${\mathbb Z}^n\times {\mathbb Z}^m$ par la projection $p_{\leq}$. Donc $(\sigma,\tau)\in \mathfrak{S}[X_{J}]$ si et seulement si $(\sigma{\sigma'}^{-1},\tau{\tau'}^{-1})\in \mathfrak{S}[X_{J}^{{\sigma'}^{-1},{\tau'}^{-1}}]$. Dans la formule (8), on intervertit les sommes en commen\c{c}ant par sommer en $\sigma',\tau'$. On remplace ensuite $X_{J}$ par $X_{J}^{\sigma',\tau'}$ et $(\sigma,\tau)$ par $(\sigma\sigma',\tau\tau')$. Ceci fait appara\^{\i}tre des termes $sgn(\sigma')sgn(\tau')$. On utilise les propri\'et\'es ci-dessus et on intervertit de nouveau les sommes. On obtient alors
$$(9) \qquad r^{d}_{\Lambda}(\alpha',\beta')=\sum_{X_{J}\subset J, \vert X_{J}\vert \leq d}\sum_{(\sigma,\tau)\in \mathfrak{S}[X_{J}]}sgn(\sigma)sgn(\tau)\sum_{(z,y)\in {\cal Z}(1,1,\sigma,\tau,X_{J})}\vert stab(z)\vert ^{-1}\vert stab(y)\vert ^{-1}$$
$$\sum_{\sigma'\in \mathfrak{S}_{k},\tau'\in \mathfrak{S}_{h}}sgn(\sigma')sgn(\tau')\vert {\cal X}_{K}(d-\vert X_{J}\vert,z,y,\sigma',\tau')\vert.$$
Posons
$$v_{k}(t)=\prod_{i=1,...,k}\frac{1-t^{i}}{1-t},$$
et d\'efinissons de m\^eme $v_{h}(t)$. Notons $c_{f}$ le coefficient de $t^f$ dans $v_{h}(t)v_{k}(t)$. Montrons que, pour $f\in {\mathbb N}$ et $(z,y)\in ({\mathbb Z}^k)^+\times ({\mathbb Z}^h)^+$, on a
$$(10)\qquad  \sum_{\sigma'\in \mathfrak{S}_{k},\tau'\in \mathfrak{S}_{h}}sgn(\sigma')sgn(\tau')\vert {\cal X}_{K}(f,z,y,\sigma',\tau')\vert=\left\lbrace\begin{array}{cc}0,&\text{ si }(z,y)\not=(\Delta_{k},\Delta_{h}),\\ (-1)^fc_{f},&\text{ si }(z,y)=(\Delta_{k},\Delta_{h}).\\ \end{array}\right.$$
En reprenant les d\'efinitions, on voit que la somme \`a calculer est \'egale au produit de 
$(-1)^f$ et  du coefficient de $X_{n-k+1}^{z_{n-k+1}}...X_{n}^{z_{n}} Y_{m-h+1}^{y_{m-h+1}}...Y_{m}^{y_{m}}t^f$ dans le polyn\^ome 

\noindent $P(X_{n-k+1},...,X_{n},t)Q(Y_{m-h+1},...,Y_{m},t)$, o\`u
$$P(X_{n-k+1},...,X_{n},t)=\sum_{\sigma'\in \mathfrak{S}_{k}}sgn(\sigma')(\prod_{i=n-k+1,...,n}X_{\sigma'(n-k+1)}^{k-1}...X_{\sigma'(n-1)})$$
$$(\prod_{n-k+1\leq i<j\leq n}(1-tX_{\sigma'(j)}X_{\sigma'(i)}^{-1})),$$
$$Q(Y_{m-h+1},...,Y_{m},t)=\sum_{\tau'\in \mathfrak{S}_{h}}sgn(\tau')(\prod_{i=m-h+1,...,m}Y_{\tau'(m-h+1)}^{h-1}...Y_{\tau'(m-1)})$$
$$(\prod_{m-h+1\leq i<j\leq m}(1-tY_{\tau'(j)}Y_{\tau'(i)}^{-1})).$$
On r\'ecrit
$$P(X_{n-k+1},...,X_{n},t)=\sum_{\sigma'\in \mathfrak{S}_{k}}sgn(\sigma')\prod_{n-k+1\leq i<j\leq n}(X_{\sigma'(i)}-tX_{\sigma'(j)}).$$
On sait que
$$\prod_{n-k+1\leq i<j\leq n}(X_{\sigma'(i)}-X_{\sigma'(j)})=sgn(\sigma')\prod_{n-k+1\leq i<j\leq n}(X_{i}-X_{j}).$$
On a donc aussi
$$P(X_{n-k+1},...,X_{n},t)=(\prod_{n-k+1\leq i<j\leq n}(X_{i}-X_{j}))\sum_{\sigma'\in \mathfrak{S}_{k}}\prod_{n-k+1\leq i<j\leq n}\frac{X_{\sigma'(i)}-tX_{\sigma'(j)}}{X_{\sigma'(i)}-X_{\sigma'(j)}}.$$
 D'apr\`es \cite{Mac} III.(1.4), on a
 $$P(X_{n-k+1},...,X_{n},t)=v_{k}(t)\prod_{n-k+1\leq i<j\leq n}(X_{i}-X_{j}).$$
 Puisque $z$ appartient \`a $({\mathbb Z}^k)^+$, il est bien connu que le coefficient de $X_{n-k+1}^{z_{n-k+1}}...X_{n}^{z_{n}} $ dans ce polyn\^ome est nul sauf si $z=\Delta_{k}$, auquel cas ce coefficient vaut $v_{k}(t)$. Un m\^eme calcul s'applique \`a $Q(Y_{m-h+1},...,Y_{m},t)$. La relation (10) s'ensuit.

On utilise (10) pour simplifier (9). Seul le couple $(z,y)=(\Delta_{k},\Delta_{h})$ peut contribuer. Pour ce couple, $\vert stab(z)\vert=\vert stab(y) \vert =1$.  Ce couple contribue si et seulement s'il appartient \`a ${\cal Z}(1,1,\sigma,\tau,X_{J})$, autrement dit si 
$$(\Delta_{k},\Delta_{h})=p_{>}(((\alpha'+\Delta_{n})^{\sigma},(\beta'+\Delta_{m})^{\tau})+\sum_{x\in X_{J}}x).$$
Mais la conjonction de cette relation et de la condition que $(\sigma,\tau)$ appartienne \`a $\mathfrak{S}[X_{J}]$, c'est-\`a-dire que la relation (5) soit v\'erifi\'ee, \'equivaut \`a l'\'egalit\'e
$$(11) \qquad (\alpha+\Delta_{n},\beta+\Delta_{m})=((\alpha'+\Delta_{n})^{\sigma},(\beta'+\Delta_{m})^{\tau})+\sum_{x\in X_{J}}x.$$
Pour tous $\sigma\in \mathfrak{S}_{n}$, $\tau\in \mathfrak{S}_{m}$ et $f\in {\mathbb N}$, notons ${\cal X}_{J}(f,\sigma,\tau)$ l'ensemble des sous-ensembles $X_{J}\subset J$ tels que $\vert X_{J}\vert =f$ et (11) soit v\'erifi\'ee. On obtient
$$r^{d}_{\Lambda}(\alpha',\beta')=\sum_{\sigma\in \mathfrak{S}_{n},\tau\in \mathfrak{S}_{m}}sgn(\sigma)sgn(\tau)\sum_{f=0,...,d}(-1)^{d-f}c_{d-f}\vert {\cal X}_{J}(f,\sigma,\tau)\vert .$$
En reportant cette expression dans (3), on obtient
$$(12) \qquad r_{\Lambda}(\alpha',\beta')=v_{\alpha,\beta}(t)^{-1}v_{k}(t)v_{h}(t)\sum_{f\in {\mathbb N}}(-1)^fs_{\Lambda}^f(\alpha',\beta')t^f,$$
o\`u
$$s_{\Lambda}^f(\alpha',\beta')=\sum_{\sigma\in \mathfrak{S}_{n},\tau\in \mathfrak{S}_{m}}sgn(\sigma)sgn(\tau)\vert {\cal X}_{J}(f,\sigma,\tau)\vert.$$
La formule (12) est similaire \`a (3), on s'est simplement d\'ebarrass\'e de l'ensemble $K$. 

\bigskip

\subsection{Fin de la d\'emonstration de la proposition 5.2}

A ce point, nous allons utiliser l'identification de $I$ \`a $\{1,...,2r+1\}$ d\'efinie en 5.2: on identifie $(j,0)$ \`a $a_{j}$ et $(l,1)$ \`a $b_{l}$. Cette identification induit un isomorphisme $\iota:{\mathbb R}^n\times {\mathbb R}^m\simeq {\mathbb R}^{2r+1}$. Les \'el\'ements de ${\cal P}_{n}\times {\cal P}_{m}$ s'identifient ainsi \`a des \'el\'ements de ${\mathbb R}^{2r+1}$. On note $I_{a}=\{a_{1},...,a_{n}\}$ et $I_{b}=\{b_{1},...,b_{m}\}$. Introduisons le groupe $G=GL(2r+1,{\mathbb C})$, son sous-tore diagonal $T$ et le groupe $X^*(T)$ de ses caract\`eres alg\'ebriques. On a naturellement $X^*(T)\otimes_{{\mathbb Z}}{\mathbb R}={\mathbb R}^{2r+1}$. On note $V={\mathbb R}^{2r+1}$ et  $(\epsilon_{i})_{i=1,...,2r+1}$ la base standard de cet espace. La d\'ecomposition $\{1,...,2r+1\}=I_{a}\sqcup I_{b}$ d\'efinit un sous-groupe de Levi semi-standard $M$ de $G$: les \'el\'ements de $M$ conservent les deux  sous-espaces de ${\mathbb C}^{2r+1}$ engendr\'es  par les  $\epsilon_{i}$ pour  $i\in I_{a}$, resp. $i\in I_{b}$. On introduit la chambre positive ferm\'ee $ {\bar C}^+$ associ\'ee au sous-groupe de Borel triangulaire sup\'erieur de $G$. Autrement dit 
  $ {\bar C}^+$ est l'ensemble des $x=(x_{1},...,x_{2r+1})\in  V$ tels que 
$x_{1}\geq... \geq x_{2r+1}$. Notons $W=\mathfrak{S}_{2r+1}$ le groupe de Weyl de $G$ relatif \`a $T$ et $W^M$ le groupe de Weyl de $M$. Via le plongement $\iota$, le groupe $\mathfrak{S}_{n}\times \mathfrak{S}_{m}$ du paragraphe pr\'ec\'edent s'identifie \`a $W^M$. Comme on le sait, pour tout $x\in  V$, il existe un unique \'el\'ement de ${\bar C}^+$, que l'on note $x^+$, de sorte qu'il existe $w\in W$ tel que $w(x)=x^+$. On introduit l'ensemble $\Sigma$ des racines $\epsilon_{i}-\epsilon_{j}$ pour $i\not=j$ et le sous-ensemble $\Sigma^+$ des racines positives, c'est-\`a-dire des $\epsilon_{i}-\epsilon_{j}$ pour $i<j$. On note aussi $\Sigma^M$ l'ensemble des racines dans $M$, c'est-\`a-dire les $\epsilon_{i}-\epsilon_{j}$ pour $i\not= j$, $i,j\in I_{a}$ ou $i,j\in I_{b}$. On pose $\Sigma^{M,+}=\Sigma^M\cap \Sigma^+$. Pour tout entier $e\in {\mathbb N}$, posons $\delta_{e}=\{(e-1)/2,(e-3)/2,...,(1-e)/2)=[(e-1)/2,(1-e)/2]_{1}$. L'\'el\'ement $\delta_{2r+1}$ est \'egal \`a la demi-somme des racines positives et  $\iota(\delta_{n},\delta_{m})$ est \'egal \`a la demi-somme des racines dans $\Sigma^{M,+}$. Nous noterons ces \'el\'ements $\delta$ et $\delta^M$. Posons ${\cal J}=\iota(J)$. On voit que ${\cal J}$ est un sous-ensemble de $\Sigma^+-\Sigma^{M,+}$.  
On note $z$ l'\'el\'ement de ${\mathbb R}^{2r+1}$ d\'efini par

si $k$ est pair, $z_{i}=n-1$ pour $i\in I_{a}$ et $z_{i}=m$ pour $i\in I_{b}$;

si $k$ est impair, $z_{i}=n$ pour $i\in I_{a}$ et $z_{i}=m-1$ pour $i\in I_{b}$.

Avec ces d\'efinitions, on voit que, pour $(\alpha',\beta')\in {\cal P}_{n}\times {\cal P}_{m}$, $\iota(\Lambda_{A,B;2}^{n,m}(\alpha',\beta'))=\iota(\alpha',\beta')+2\delta^M+z$ et $p\Lambda_{A,B;2}^{n,m}(\alpha',\beta') =\iota(\Lambda_{A,B;2}^{n,m}(\alpha',\beta'))^+$. Remarquons que l'ordre $<_{I,\alpha,\beta}$ a \'et\'e d\'efini de sorte que $\iota(\Lambda_{A,B;2}^{n,m}(\alpha,\beta))$ appartienne \`a ${\bar C}^+$, autrement dit $\Lambda=p\Lambda_{A,B;2}^{n,m}(\alpha,\beta)=\iota(\Lambda_{A,B;2}^{n,m}(\alpha,\beta))$. A ce point, on se rappelle l'hypoth\`ese que $\Lambda$ est sans multiplicit\'e. Puisque les termes de $\Lambda$ sont des entiers, on a $\Lambda_{i}\geq \Lambda_{i+1}+1$ pour tout $i=1,...,2r$, ce qui entra\^{\i}ne

(1) $\Lambda-\rho\in {\bar C}^+$. 

Soit $(\alpha',\beta')\in {\cal P}_{n}\times {\cal P}_{m}$, supposons $r_{\Lambda}(\alpha',\beta')\not=0$. D'apr\`es 5.3(12), on peut fixer $\sigma\in \mathfrak{S}_{n}$, $\tau\in \mathfrak{S}_{m}$ et $X_{J}\subset J$ de sorte que 5.3(11) soit v\'erifi\'ee. On voit que, dans cette relation, on peut aussi bien remplacer $\Delta_{n}$ et $\Delta_{m}$ par $\delta_{n}$ et $\delta_{m}$ car $(\Delta_{n}-\delta_{n},\Delta_{m}-\delta_{m})$ est fixe par $\sigma\times \tau$.   En posant  $X=\iota(X_{J})$, $v=\iota(\alpha,\beta)$, $v'=\iota(\alpha',\beta')$, et en notant $w^M$ l'\'el\'ement de $W^M$ auquel s'identifie $\sigma\times \tau$, cette relation s'\'ecrit
$$v'+\delta^M=w^M (v+\delta^M-\sum_{x\in X}x).$$
Alors 
$$\iota(\Lambda_{A,B;2}^{n,m}(\alpha',\beta'))=v'+2\delta^M+z=w^M(v+\delta^M-\sum_{x\in X}x)+\delta^M+z=w^M(\Lambda-\delta^M-z-\sum_{x\in X}x)+\delta^M+z.$$
L'\'el\'ement $z$ dispara\^{\i}t de cette relation: il est central dans $M$ et fixe par $w^M$. Pour d\'emontrer que $p\Lambda_{A,B;2}^{n,m}(\alpha',\beta') \leq \Lambda$, cf. 5.3(1), on est ramen\'e au probl\`eme suivant. L'\'el\'ement $\Lambda$ v\'erifie (1). On se donne un sous-ensemble $X\subset \Sigma^{+}-\Sigma^{M,+}$ et un \'el\'ement $w^M\in W^M$. On pose
$$\underline{\Lambda}=w^M(\Lambda-\delta^M-\sum_{x\in X}x)+\delta^{M}$$
et on veut prouver
$$(2) \qquad \underline{\Lambda}^+\leq \Lambda.$$
On a 
$$\delta-\delta^M=\frac{1}{2}\sum_{x\in \Sigma^+-\Sigma^{M,+}}x=\frac{1}{2}(\sum_{x\in X}x+\sum_{x\in \Sigma^+-(\Sigma^{M,+}\cup X)}x).$$
  D'o\`u
$$\delta^M+\sum_{x\in X}x=\delta+\mu,$$
o\`u 
$$(3)\qquad \mu=\frac{1}{2}(\sum_{x\in X}x-\sum_{x\in \Sigma^+-(\Sigma^{M,+}\cup X)}x).$$
  On a
$$\underline{\Lambda}=w^M(\Lambda-\delta-\mu)+\delta^M=w^M(\Lambda-\delta+\mu'),$$
o\`u $\mu'=-\mu+(w^M)^{-1}(\delta^M)$.
Par d\'efinition, l'\'el\'ement $\underline{\Lambda}^+$ est insensible \`a la transformation de $\underline{\Lambda}$ par n'importe quel \'el\'ement de $W$. Cela entra\^{\i}ne
$$\underline{\Lambda}^+=(\Lambda-\delta+\mu' )^+.$$
Notons $^+{\bar C}$ le c\^one ferm\'e  engendr\'e par les \'el\'ements de $\Sigma^+$. Pour prouver (2), il suffit de prouver que $(\Lambda-\delta+\mu')^+\in \Lambda-{^+{\bar C}}$. Par d\'efinition, il existe $w\in W$ tel que $(\Lambda-\delta+\mu')^+=w(\Lambda-\delta+\mu')$. D'o\`u
$$(\Lambda-\delta+\mu')^+= w(\Lambda-\delta)+(\Lambda-\delta)-(\Lambda-\delta)+w(\mu')=\Lambda+w(\Lambda-\delta)-(\Lambda-\delta)+w(\mu')-\delta.$$
Puisque $\Lambda-\delta\in {\bar C}^+$, il est bien connu que $w(\Lambda-\delta)-(\Lambda-\delta)$ appartient \`a $-{^+{\bar C}}$. Il reste \`a prouver que $w(\mu')-\delta$ appartient aussi \`a ce c\^one. 
 Notons $E$ l'ensemble des \'el\'ements de $V$ qui s'\'ecrivent $\frac{1}{2}(\sum_{y\in Y}y-\sum_{y\in \Sigma^+-Y}y)$ pour un sous-ensemble $Y\subset \Sigma^+$ et notons $E^M$ l'ensemble des \'el\'ements de $V$ qui s'\'ecrivent $\frac{1}{2}(\sum_{y\in Y}y-\sum_{y\in \Sigma^{M,+}-Y}y)$ pour un sous-ensemble $Y\subset \Sigma^{M,+}$. En utilisant (3), on voit que $-\mu+H\in E$ pour tout $H\in E^M$. D'autre part, $E$ est conserv\'e par l'action de $W$ et $E^M$ l'est par celle de $W^M$. Puisque $\delta^M$ appartient \`a $E^M$, l'\'el\'ement $\mu'=-\mu+(w^M)^{-1}(\delta^M)$ appartient \`a $E$ et aussi $w(\mu')\in E$. Mais on voit que $H-\delta$ appartient \`a $-{^+{\bar C}}$ pour tout $H\in E$. En particulier $w(\mu')-\delta$ appartient \`a $-{^+{\bar C}}$, ce qui ach\`eve la d\'emonstration. $\square$

\bigskip
 \subsection{Un th\'eor\`eme de maximalit\'e}
 \ass{Th\'eor\`eme}{Soit $(\lambda,\epsilon)\in \boldsymbol{{\cal P}^{symp}}(2N)$. Supposons que $\lambda$ n'a que des termes pairs. Alors il existe un unique \'el\'ement $(\lambda^{max},\epsilon^{max})\in \boldsymbol{{\cal P}^{symp}}(2N)$ v\'erifiant les propri\'et\'es suivantes:
 
 (i) $\mathfrak{mult}(\lambda,\epsilon;\lambda^{max},\epsilon^{max})=1$;
 
 (ii) pour tout \'el\'ement $(\lambda',\epsilon')\in \boldsymbol{{\cal P}^{symp}}(2N)$ tel que $\mathfrak{mult}(\lambda,\epsilon;\lambda',\epsilon')\not=0$, on a $\lambda'<\lambda^{max}$ ou $(\lambda',\epsilon')=(\lambda^{max},\epsilon^{max})$. }
 
 Preuve. Posons $k=k(\lambda,\epsilon)$. Puisqu'on se pr\'eoccupe de paires $(\lambda',\epsilon')$ telles que $\mathfrak{mult}(\lambda,\epsilon;\lambda',\epsilon')\not=0$, on peut se limiter au sous-ensemble $\boldsymbol{{\cal P}^{symp}}(2N)_{k}$ des \'el\'ements $(\lambda',\epsilon')$ tels que $k(\lambda',\epsilon')=k$. Montrons tout d'abord que, pour $(\lambda',\epsilon'),(\lambda'',\epsilon'')\in \boldsymbol{{\cal P}^{symp}}(2N)_{k}$, on a l'\'equivalence
 
 (1) $\lambda'\leq \lambda''$ si et seulement si $A_{\lambda',\epsilon'}\sqcup B_{\lambda',\epsilon'}\leq A_{\lambda'',\epsilon''}\sqcup B_{\lambda'',\epsilon''}$.
 
 On applique  \`a $\lambda'$ la construction de $A_{\lambda',\epsilon'}$ et $B_{\lambda',\epsilon'}$ rappel\'ee en 5.1. Pour $c\in \{1,...,2r\}$, on a
 $$S_{c}(\lambda')=S_{c}(\lambda'+[2r-1,0]_{1})-S_{c}([2r-1,0]_{1}).$$
 On a d\'ecompos\'e $\lambda'+[2r-1]_{1}$ en la r\'eunion des $2z_{i}$ et des $2z'_{i}+1$ pour $i=1,...,r$. On v\'erifie par r\'ecurrence les propri\'et\'es suivantes:
 
 pour $c$ pair, les $c$ plus grands termes de $\lambda'+[2r-1]_{1}$ sont les $2z_{i}$ et $2z'_{i}+1$ pour $i=1,...,c/2$, \`a l'exception du cas o\`u   $S_{c}(\lambda)$ est impair; dans ce cas, $\lambda_{c}$ est forc\'ement impair, les $c$ plus grands termes de $\lambda'+[2r-1]_{1}$ sont les $2z_{i}$ pour $i=1,...,c/2-1$ et les $2z'_{i}+1$ pour $i=1,...,c/2+1$ et on a $z'_{c/2+1}=z_{c/2}$;
 
 pour $c$ impair, les $c$ plus grands termes de $\lambda'+[2r-1]_{1}$ sont les $2z_{i}$ pour $i=1,...,(c-1)/2$ et les $2z'_{i}+1$ pour $i=1,...,(c+1)/2$, \`a l'exception du cas o\`u  $S_{c}(\lambda)$ est impair; dans ce cas $\lambda_{c}$ est forc\'ement impair,  les $c$ plus grands termes de $\lambda'+[2r-1]_{1}$ sont les $2z_{i}$ pour $i=1,...,(c+1)/2$ et les $2z'_{i}+1$ pour $i=1,...,(c-1)/2$ et on a $z_{(c+1)/2}=z'_{(c+1)/2}+1$.
 
 Posons $c^+=c^-=c/2$ si $c$ est pair, $c^+=(c+1)/2$, $c^-=(c-1)/2$ si $c$ est impair. Posons aussi $\delta_{c}(\lambda')=1$ si $S_{c}(\lambda)$ est impair. Les propri\'et\'es ci-dessus entra\^{\i}nent  que
 $$S_{c}(\lambda'+[2r-1,0]_{1})=2S_{c^+}(z')+c^++2S_{c^-}(z)+\delta_{c}(\lambda').$$
 On a
 $$S_{c^+}(z')=S_{c^+}(A_{\lambda'}^{\sharp})-S_{c^+}([r+1,2]_{1}),$$
 $$S_{c^-}(z)=S_{c^-}(B_{\lambda'}^{\sharp})-S_{c^-}([r,1]_{1}).$$
 D'o\`u
 $$S_{c}(\lambda')=2S_{c^+}(A_{\lambda'}^{\sharp})+2S_{c^-}(B_{\lambda'}^{\sharp})+\delta_{c}(\lambda')+C_{c},$$
 o\`u $C_{c}$ est un nombre ind\'ependant de $\lambda'$. On v\'erifie qu'en notant $a^{\sharp}_{j}$ et $b^{\sharp}_{j}$ les termes de $A_{\lambda'}^{\sharp}$ et $B^{\sharp}_{\lambda'}$, on a
 $$a^{\sharp}_{1}\geq b^{\sharp}_{1}\geq a^{\sharp}_{2}\geq...\geq a^{\sharp}_{r}\geq b^{\sharp}_{r}\geq a^{\sharp}_{r+1}.$$
 Donc 
 $$S_{c^+}(A_{\lambda'}^{\sharp})+S_{c^-}(B_{\lambda'}^{\sharp})=S_{c}(A_{\lambda'}^{\sharp}\sqcup B_{\lambda'}^{\sharp}).$$
 D'apr\`es les d\'efinitions, on a 
  $$A_{\lambda'}^{\sharp}\sqcup B_{\lambda'}^{\sharp}= A_{\lambda',\epsilon'}\sqcup B_{\lambda',\epsilon'}.$$
  D'o\`u
$$(2) \qquad S_{c}(\lambda')=2S_{c}(A_{\lambda',\epsilon'}\sqcup B_{\lambda',\epsilon'})+\delta_{c}(\lambda')+C_{c}.$$ 
Supposons $\lambda'\leq \lambda''$. Alors $S_{c}(\lambda')\leq S_{c}(\lambda'')$, d'o\`u 
 $$S_{c}(A_{\lambda',\epsilon'}\sqcup B_{\lambda',\epsilon'})\leq S_{c}(A_{\lambda'',\epsilon''}\sqcup B_{\lambda'',\epsilon''})+\delta_{c}(\lambda'')/2-\delta_{c}(\lambda')/2\leq S_{c}(A_{\lambda'',\epsilon''}\sqcup B_{\lambda'',\epsilon''})+1/2.$$
 Puisque les sommes intervenant ici sont enti\`eres, on en d\'eduit
  $$S_{c}(A_{\lambda',\epsilon'}\sqcup B_{\lambda',\epsilon'})\leq S_{c}(A_{\lambda'',\epsilon''}\sqcup B_{\lambda'',\epsilon''}).$$
  Cette in\'egalit\'e \'etant v\'erifi\'ee pour tout $c=1,...,2r$, on conclut que $A_{\lambda',\epsilon'}\sqcup B_{\lambda',\epsilon'}\leq A_{\lambda'',\epsilon''}\sqcup B_{\lambda'',\epsilon''}$. Inversement supposons v\'erifi\'ee cette derni\`ere in\'egalit\'e. Alors 
   $$S_{c}(A_{\lambda',\epsilon'}\sqcup B_{\lambda',\epsilon'})\leq S_{c}(A_{\lambda'',\epsilon''}\sqcup B_{\lambda'',\epsilon''}).$$
   Avec (2), on en d\'eduit
   $$(3) \qquad S_{c}(\lambda')\leq S_{c}(\lambda'')+\delta_{c}(\lambda')-\delta_{c}(\lambda'')\leq S_{c}(\lambda'')+\delta_{c}(\lambda').$$
   Si $S_{c}(\lambda')$ est pair, $\delta_{c}(\lambda')=0$, d'o\`u  l'in\'egalit\'e
   $$(4) \qquad S_{c}(\lambda')\leq S_{c}(\lambda'').$$
 Supposons que $S_{c}(\lambda')$ soit impair et que (4) ne soit pas v\'erifi\'ee. Alors (3) force l'\'egalit\'e $S_{c}(\lambda')=S_{c}(\lambda'')+1$. Puisque $\lambda'$ est symplectique, le fait que $S_{c}(\lambda')$ soit impair entra\^{\i}ne que $\lambda_{c}$ est impair, que $\lambda'_{c}=\lambda'_{c+1}$ et que $S_{c-1}(\lambda')$ et $S_{c+1}(\lambda')$ sont pairs. On peut appliquer (4) pour $c-1$ (dans le cas particulier $c=1$, (4) est triviale pour $c-1=0$). Puisque $S_{c}(\lambda')=S_{c-1}(\lambda')+\lambda'_{c}$ et $S_{c}(\lambda'')=S_{c-1}(\lambda'')+\lambda''_{c}$, cette in\'egalit\'e (4) pour $c-1$ et l'\'egalit\'e $S_{c}(\lambda')=S_{c}(\lambda'')+1$ entra\^{\i}nent $\lambda'_{c}\geq \lambda''_{c}+1$. Puisque $\lambda'_{c+1}=\lambda'_{c}$ et $\lambda''_{c+1}\leq \lambda''_{c}$, on a aussi $\lambda'_{c+1}\geq \lambda''_{c+1}+1$. Alors l'\'egalit\'e $S_{c}(\lambda')=S_{c}(\lambda'')+1$ entra\^{\i}ne 
 $$S_{c+1}(\lambda')=S_{c}(\lambda')+\lambda'_{c+1}\geq S_{c}(\lambda'')+\lambda''_{c+1}+2> S_{c+1}(\lambda'').$$
 Puisque $S_{c+1}(\lambda')$ est pair, les deux termes extr\^emes sont \'egaux d'apr\`es (4) pour $c+1$ (dans le cas particulier $c=2r$,cette \'egalit\'e (4) pour $c+1$ est triviale). Cette contradiction
invalide notre hypoth\`ese, donc (4) est forc\'ement v\'erifi\'ee. Cette in\'egalit\'e \'etant maintenant prouv\'ee pour tout $c$, on en d\'eduit $\lambda'\leq \lambda''$, ce qui prouve (1). 

Notons $(\alpha,\beta)\in {\cal P}_{2}(N-k(k+1)/2)$ le couple correspondant \`a $(\lambda,\epsilon)$. On consid\`ere que $\alpha$ et $\beta$ ont respectivement $n$  et $m$ termes comme en 5.2. La proposition 5.2 et la relation (1) permettent de r\'ecrire l'\'enonc\'e sous la forme

il existe un unique couple $(\alpha^{max},\beta^{max})\in {\cal P}_{n}\times {\cal P}_{m}$ v\'erifiant les conditions

(i) $mult(<_{I,\alpha,\beta};\alpha,\beta;\alpha^{max},\beta^{max})=1$;

(ii) pour tout couple $(\alpha',\beta')\in {\cal P}_{n}\times {\cal P}_{m}$ tel que $mult(<_{I,\alpha,\beta};\alpha,\beta;\alpha',\beta')\not=0$, on a $p\Lambda_{A,B;2}^{n,m}(\alpha',\beta')< p\Lambda_{A,B;2}^{n,m}(\alpha^{max},\beta^{max})$ ou $(\alpha',\beta')=(\alpha^{max},\beta^{max})$.

La proposition 4.1 entra\^{\i}ne que ceci est v\'erifi\'e si et seulement si l'ensemble $P_{A,B;2}(<_{I,\alpha,\beta};\alpha,\beta)$ a un unique \'el\'ement et, dans ce cas, $(\alpha^{max},\beta^{max})$ est cet unique \'el\'ement (comme on l'a dit plus haut, on a gliss\'e l'ordre $<_{I,\alpha,\beta}$ dans la notation). Pour prouver cette assertion, g\'en\'eralisons un peu les hypoth\`eses. On consid\`ere des entiers $n,m\in {\mathbb N}$ et des partitions $\alpha\in {\cal P}_{n}$, $\beta\in {\cal P}_{m}$. On fixe des entiers $A,B$ tels que $A\geq 2(n-1)$ et $B\geq 2(m-1)$, on d\'efinit la  partition $p\Lambda_{A,B;2}^{n,m}(\alpha,\beta)$ et on suppose que cette suite est sans multiplicit\'es. On d\'efinit l'ordre $<_{I,A,B,\alpha,\beta}$ sur l'ensemble d'indices $I=(\{1,...,n\}\times\{0\})\cup (\{1,...,m\}\times \{1\})$ de la fa\c{c}on suivante: pour $i\in \{1,...,n\}$ et $j\in \{1,...,m\}$, on a $(i,0)<_{I,A,B,\alpha,\beta}(j,1)$ si et seulement si $\alpha_{i}+A+2-2i> \beta_{j}+B+2-2j$. On veut alors prouver que $P_{A,B;2}(<_{I,A,B,\alpha,\beta};\alpha,\beta)$ a un unique \'el\'ement. 

{\bf Remarque.} Les hypoth\`eses que l'on vient de poser sont bien v\'erifi\'ees dans la situation pr\'ec\'edente, $A$ et $B$ \'etant d\'efinis comme en 5.2. En particulier, l'ordre $<_{I,\alpha,\beta}$ co\"{\i}ncide par d\'efinition avec l'ordre $<_{I,A,B,\alpha,\beta}$. 

\bigskip
Si $n=0$ ou $m=0$, la conclusion est triviale puisque $P(\alpha,\beta)$ n'a qu'un \'el\'ement. On suppose $n\geq1$, $m\geq1$. Pour simplifier, on note $<_{I}$ l'ordre $<_{I,A,B,\alpha,\beta}$. On ne perd rien \`a supposer que $(1,0)<_{I}(1,1)$. On a alors $\alpha_{1}+A>\beta_{1}+B\geq B$ et le proc\'ed\'e (a) est l'unique proc\'ed\'e autoris\'e pour construire les \'el\'ements de $P_{A,B;2}(<_{I};\alpha,\beta)$.   Ce proc\'ed\'e construit un \'el\'ement $\nu_{1}$ et des partitions $(\alpha',\beta')\in {\cal P}_{n'}\times {\cal P}_{m'}$. De l'ordre $<_{I}$ sur $I$ se d\'eduit un ordre sur l'ensemble d'indices $I'=(\{1,...,n'\}\times\{0\})\cup (\{1,...,m'\}\times \{1\})$, cf. 1.2, notons-le $<_{I'}$. Alors le proc\'ed\'e (a)  d\'efinit une bijection entre $P_{A,B;2}(<_{I};\alpha,\beta)$ et $P_{A-2,B;2}(<_{I'};\alpha',\beta')$ et il faut d\'emontrer que ce dernier ensemble n'a qu'un \'el\'ement. De nouveau, c'est trivial si $n'=0$ ou $m'=0$. On suppose $n'\geq1$, $m'\geq1$. En raisonnant par r\'ecurrence, il suffit
de prouver que les donn\'ees $(\alpha',\beta')$, $A-2$, $B$ v\'erifient les m\^emes hypoth\`eses que $(\alpha,\beta)$, $A$ et $B$ et que l'ordre $<_{I'}$ sur $I'$ co\"{\i}ncide avec l'ordre $<_{I',A-2,B,\alpha',\beta'}$.  Rappelons que l'\'el\'ement $\nu_{1}$ est de la forme $\alpha_{a_{1}}+\beta_{b_{1}}+\alpha_{a_{2}}+\beta_{b_{2}}+...$ o\`u $a_{1}=b_{1}=1$, $a_{2}$ est le plus petit $i\in \{1,...,n\}$ tel que $(b_{1},1)<_{I}(i,0)$ si un tel \'el\'ement existe, $b_{2}$ est le plus petit $j\in \{1,...,m\}$ tel que $(a_{2},0)<_{I}(j,1)$ si un tel \'el\'ement existe etc... Les derniers \'el\'ements de $\nu_{1}$ peuvent \^etre $...+\alpha_{a_{t}}+\beta_{b_{t}}$ ou $...+\alpha_{a_{t}}+\beta_{b_{t}}+\alpha_{a_{t+1}}$. Dans le premier cas, on pose $ T=t$; dans le second, on pose $ T=t+1$. La partition $\alpha'$ est form\'ee des termes $\alpha_{a_{1}+1},...,\alpha_{a_{2}-1}$, $\alpha_{a_{2}+1},...,\alpha_{a_{3}-1}$,..., que l'on r\'eindexe en $\alpha'_{1},...,\alpha'_{n'}$. On en d\'eduit la formule suivante. Pour $i\in \{1,...,n'\}$, notons $k(i)$ le plus grand entier $k\in \{1,...,T\}$ tel que $a_{k}-k<i$. Alors $\alpha'_{i}=\alpha_{i+k(i)}$. Remarquons que l'on a $a_{k(i)}< i+k(i)$. Si $k(i)<T$, on a aussi $i+k(i)<a_{k(i)+1}$ (sinon, on aurait $a_{k(i)+1}-k(i)-1<i$, contredisant la maximalit\'e de $k(i)$). De la m\^eme fa\c{c}on, pour $j\in \{1,...,m'\}$, notons $h(j)$ le plus grand entier $h\in \{1,...,t\}$ tel que $b_{h}-h<j$. Alors $\beta'_{j}=\beta_{j+h(j)}$. On a $b_{h(j)}<j+h(j)$. Si $h(j)<t$, on a $j+h(j)<b_{h(j)+1}$. 
Soient $i\in \{1,...,n'\}$ et $j\in \{1,...,m'\}$. Montrons que

(5) les conditions $(i+k(i),0)<_{I}(j+h(j),1)$ et $k(i)\leq h(j)$ sont \'equivalentes; si elles sont v\'erifi\'ees, on a $k(i)\leq t$ et $(i+k(i),0)<_{I}(b_{k(i)},1)$;

(6) les conditions $(j+h(j),1)<_{I}(i+k(i),0)$ et $h(j)< k(i)$ sont \'equivalentes; si elles sont v\'erifi\'ees, on a  $h(j)+1\leq T$ et $(j+h(j),1)<_{I}(a_{h(j)+1},0)$.

Supposons $(i+k(i),0)<_{I}(j+h(j),1)$. Puisque $a_{k(i)}<i+k(i)$, on a aussi $(a_{k(i)},0)<_{I}(j+h(j),1)$.   Ainsi $j+h(j)$ est un \'el\'ement $l$ de $\{1,...,m\}$ tel que $(a_{k(i)},0)<_{I}(l,1)$. L'existence d'un tel \'el\'ement implique que $b_{k(i)}$ existe et $b_{k(i)}$ est le plus petit $l$ v\'erifiant cette condition. Cela entra\^{\i}ne que $k(i)\leq t$ et que $b_{k(i)}\leq j+h(j)$. Si $h(j)=t$, on a  $k(i)\leq t=h(j)$. Si $h(j)<t$, on a vu que $j+h(j)< b_{h(j)+1}$. Alors $b_{k(i)}<b_{h(j)+1}$, ce qui entra\^{\i}ne encore $k(i)\leq h(j)$. Enfin, si on avait $(b_{k(i)},1)<_{I}(i+k(i),0)$,  cela entra\^{\i}nerait pour la m\^eme raison que ci-dessus l'existence de l'entier $a_{k(i)+1}$ et on aurait $a_{k(i)+1}\leq i+k(i)$. On a vu que l'on avait au contraire $i+k(i)< a_{k(i)+1}$. Cela d\'emontre que $(i+k(i),0)<_{I}(b_{k(i)},1)$. On a ainsi d\'emontr\'e la deuxi\`eme assertion de (5) et un sens de l'\'equivalence de la premi\`ere assertion. Remarquons que les premi\`eres assertions de (5) et (6) sont les  n\'egatives l'une de l'autre.  Pour d\'emontrer compl\`etement ces premi\`eres assertions, il suffit de prouver le m\^eme sens de celle de (6). Supposons donc $(j+h(j),1)<_{I}(i+k(i),0)$. Comme ci-dessus, cela implique $(b_{h(j)},1)<_{I}(i+k(i),0)$, puis que $a_{h(j)+1}$ existe, donc $h(j)+1\leq T$, et que $a_{h(j)+1}\leq i+k(i)$. Si $k(i)=T$, on a $h(j)<h(j)+1\leq T=k(i)$. Si $k(i)<T$, on a vu que $i+k(i)<a_{k(i)+1}$. D'o\`u aussi $a_{h(j)+1}<a_{k(i)+1}$ puis $h(j)<k(i)$. Cela   d\'emontre le sens voulu de la premi\`ere assertion de  (6). Enfin, si on avait $(a_{h(j)+1},0)<_{I}(j+h(j),1)$, cela entra\^{\i}nerait comme ci-dessus l'existence de $b_{h(j)+1}$ et on aurait $b_{h(j)+1}\leq j+h(j)$, alors que l'on au contraire $j+h(j)< b_{h(j)+1}$. Cela ach\`eve la d\'emonstration de (5) et (6). 

Montrons maintenant que

(7) les conditions $(i+k(i),0)<_{I}(j+h(j),1)$ et $\alpha'_{i}+A-2i>\beta'_{j}+B+2-2j$ sont \'equivalentes;

(8) les conditions $(j+h(j),1)<_{I}(i+k(i),0)$ et $\beta'_{j}+B+2-2j>\alpha'_{i}+A-2i$ sont \'equivalentes.

De nouveau, il suffit de d\'emontrer le sens "implique" des deux assertions. Supposons $(i+k(i),0)<_{I}(j+h(j),1)$. D'apr\`es (5), on a $k(i)\leq t$ et $(i+k(i),0)<_{I}(b_{k(i)},1)$. Par d\'efinition de l'ordre $<_{I}$, c'est-\`a-dire de $<_{I,A,B,\alpha,\beta}$, on a 
$$(9) \qquad \alpha_{i+k(i)}+A+2-2i-2k(i)> \beta_{b_{k(i)}}+B+2-2b_{k(i)}.$$
 On a aussi $k(i)\leq h(j)$ donc $\beta_{b_{k(i)}}\geq \beta_{b_{h(j)}}$. D'o\`u
 $$\beta_{b_{k(i)}}+B+2-2b_{k(i)}\geq (\beta_{b_{h(j)}}+B+2-2b_{h(j)})+2b_{h(j)}-2b_{k(i)}.$$
 Puisque $b_{h(j)}<j+h(j)$, on a 
 $$\beta_{b_{h(j)}}+B+2-2b_{h(j)}\geq 2+\beta_{j+h(j)}+B+2-2j-2h(j).$$
 Puisque la suite $b_{l}$ est strictement croissante, on a aussi $b_{h(j)}-b_{k(i)}\geq h(j)-k(i)$. On obtient alors
$$\beta_{b_{k(i)}}+B+2-2b_{k(i)}\geq 2+\beta_{j+h(j)}+B+2-2j-2h(j)+2h(j)-2k(i)=2+\beta_{j+h(j)}+B+2-2j-2k(i).$$
Avec (9), cela entra\^{\i}ne 
$$\alpha_{i+k(i)}+A-2i> \beta_{j+h(j)}+B+2-2j.$$
Le membre de gauche est $\alpha'_{i}+A-2i$, celui de droite est $\beta'_{j}+B+2-2j$, ce qui d\'emontre le sens voulu de (7).  Supposons maintenant $(j+h(j),1)<_{I}(i+k(i),0)$. D'apr\`es (6), on a $h(j)+1\leq T$ et $(j+h(j),1)<_{I}(a_{h(j)+1},0)$. Par d\'efinition de l'ordre $<_{I}$, on a donc
$$(10)\qquad \beta_{j+h(j)}+B+2-2j-2h(j)> \alpha_{a_{h(j)+1}}+A+2-2a_{h(j)+1}.$$
On a aussi $h(j)+1\leq k(i)$ donc $\alpha_{a_{h(j)+1}}\geq \alpha_{a_{k(i)}}$. D'o\`u
$$ \alpha_{a_{h(j)+1}}+A+2-2a_{h(j)+1}\geq \alpha_{a_{k(i)}}+A+2-2a_{k(i)}+2a_{k(i)}-2a_{h(j)+1}.$$
Puisque $a_{k(i)}< i+k(i)$, on a
$$\alpha_{a_{k(i)}}+A+2-2a_{k(i)}\geq 2+\alpha_{i+k(i)}+A+2-2i-2k(i).$$
Puisque la suite $a_{l}$ est strictement croissante, on a aussi $a_{k(i)}-a_{h(j)+1}\geq k(i)-h(j)-1$. D'o\`u
 $$\alpha_{a_{h(j)+1}}+A+2-2a_{h(j)+1}\geq  2+\alpha_{i+k(i)}+A+2-2i-2k(i)+2k(i)-2h(j)-2=2+\alpha_{i+k(i)}+A-2i-2h(j).$$
 Avec (10), cela entra\^{\i}ne
 $$\beta_{j+h(j)}+B+2-2j>  2+\alpha_{i+k(i)}+A-2i.$$
 Le membre de gauche est \'egal \`a $\beta'_{j}+B+2-2j$. Celui de droite est \'egal \`a $2+\alpha'_{i}+A-2i$, ce qui d\'emontre une in\'egalit\'e plus forte que la conclusion de (8). Cela ach\`eve la preuve de (7) et (8). 
 
 {\bf Remarque.} Comme on vient de le dire, nous avons en fait d\'emontr\'e une assertion plus forte que (8), \`a savoir
 
 (11) les trois conditions $(j+h(j),1)<_{I}(i+k(i),0)$,  $\beta'_{j}+B+2-2j>\alpha'_{i}+A-2i$ et $\beta'_{j}+B+2-2j>2+\alpha'_{i}+A-2i$ sont \'equivalentes.
 
 Nous l'utiliserons dans le paragraphe suivant.
 
 \bigskip

L'une ou l'autre des conditions $(i+k(i),0)<_{I}(j+h(j),1)$ ou $(j+h(j),1)<_{I}(i+k(i),0)$ est v\'erifi\'ee. Alors (7) et (8) entra\^{\i}nent que $\alpha'_{i}+A-2i\not=\beta'_{j}+B+2-2j$. Donc la suite $p\Lambda_{A-2,B;2}^{n',m'}(\alpha',\beta')$ est sans multiplicit\'es. Par d\'efinition de l'ordre $<_{I'}$, on a $i<_{I'}j$ si et seulement si $(i+k(i),0)<_{I}(j+h(j),1)$. Par d\'efinition de l'ordre $<_{I',A-2,B,\alpha',\beta'}$, on a $i<_{I',A-2,B,\alpha',\beta'}j$ si et seulement si $\alpha'_{i}+A-2i>\beta'_{j}+B+2-2j$. L'assertion (7) implique que les deux ordres sont \'egaux. Cela v\'erifie les hypoth\`eses de r\'ecurrence, donc $P_{A-2,B;2}(<_{I'},\alpha',\beta')$ n'a qu'un \'el\'ement, ce qui ach\`eve la d\'emonstration. $\square$

\bigskip
\subsection{Tensorisation par le caract\`ere signature}
Soit $(\lambda,\epsilon)\in \boldsymbol{{\cal P}^{symp}}(2N)$. Posons $k=k(\lambda,\epsilon)$. Il existe un unique \'el\'ement de $\boldsymbol{{\cal P}^{symp}}(2N)$, que nous notons $({^s\lambda},{^s\epsilon})$ tel que $k({^s\lambda},{^s\epsilon})=k$ et $\rho({^s\lambda},{^s\epsilon})=sgn\otimes \rho(\lambda,\epsilon)$. 

{\bf Remarque.} Nous adoptons cette notation faute de mieux. Elle n'est pas tr\`es bonne car $^s\lambda$ ne d\'epend pas seulement de $\lambda$ et $^s\epsilon$ ne d\'epend pas seulement de $\epsilon$: ils d\'ependent tous deux du couple $(\lambda,\epsilon)$. 
\bigskip

 On fixe $r$ assez grand  et on pose $A=2r+k$ et $B=2r-1-k$ comme en 5.2. Posons $\underline{n}=A+1$, $\underline{m}=B+1$.  Notons $(\alpha,\beta)\in {\cal P}_{2}(N-k(k+1)/2)$ le couple param\'etrant $\rho_{\lambda,\epsilon}$. On peut consid\'erer que $\alpha$ a $\underline{n}$ termes et que $\beta$ en a $\underline{m}$. On peut alors d\'efinir les suites de r\'eels $U_{\lambda,\epsilon}=\alpha+[A/2, 0]_{1/2}$, $V_{\lambda,\epsilon}=\beta+[B/2, 0]_{1/2}$ (on rappelle que $[A/2,0]_{1/2}=\{A/2,A/2-1/2,..., 1/2,0\}$). Pour toute suite de r\'eels $\mu=(\mu_{1},...,\mu_{t})$, posons $2\mu=(2\mu_{1},...,2\mu_{t})$ et $\mu^2=(\mu_{1},\mu_{1},\mu_{2},\mu_{2},...,\mu_{t-1},\mu_{t-1},\mu_{t},\mu_{t})$. 
 \ass{Lemme}{On a l'\'egalit\'e
 $$2U_{\lambda,\epsilon}\sqcup 2V_{\lambda,\epsilon}={^t(^s\lambda)}+([2r,0]_{1}\sqcup[2r-1,0]_{1}).$$}
 
 Preuve. On aura besoin de la propri\'et\'e suivante. Soient $\mu$ une partition et $u,v\in {\mathbb N}$. En supposant que $u$ et $v$ sont assez grands, on consid\`ere que $\mu$ a $u+1$ termes et que la partition transpos\'ee $^t\mu$ en a $v+1$. On d\'efinit $X=\mu+[u,0]_{1}$ et $Y={^t\mu}+[v,0]_{1}$. Posons $\tilde{Y}=\{u+v+1-y; y\in Y\}$. Alors

(1) $X\sqcup\tilde{Y}=[u+v+1,0]_{1}$. 

Preuve. On v\'erifie que $\tilde{Y}$ est contenu dans $[u+v+1,0]_{1}$. Il en est de m\^eme de $X$. Les partitions des deux membres de (1) ont le m\^eme nombre d'\'el\'ements. Il suffit donc de prouver que $X$ et $\tilde{Y}$ n'ont pas d'\'el\'ements communs. Les \'el\'ements de $X$ sont les $\mu_{i}+u+1-i$ pour $i=1,...,u+1$. Ceux de $Y$ sont les $^t\mu_{j}+v+1-j$ pour $j=1,...,v+1$, ceux de $\tilde{Y}$ sont les $u+j-{^t\mu}_{j}$. Il faut donc prouver que $\mu_{i}+{^t\mu}_{j}\not=i+j-1$. On voit que $^t\mu_{j}$ est le plus grand $k\in \{0,..,u+1\}$ tel que $\mu_{k}\geq j$, avec la convention $\mu_{0}=\infty$. Si $\mu_{i}\geq j$, on a donc $^t\mu_{j}\geq i$ et $\mu_{i}+{^t\mu}_{j}\not=i+j-1$. Si $\mu_{i}<j$, on a $^t\mu_{j}<i$ et encore $\mu_{i}+{^t\mu}_{j}\not=i+j-1$. Cela prouve (1).

Posons $X =(^t\beta)^2+[A,0]_{1}$, $Y=(^t\alpha)^2+[B,0]_{1}$, $\tilde{X} =\{A+B+1-x; x\in  X \}$, $\tilde{Y} =\{A+B+1-y; y\in Y \}$. On sait que pour toute partition $\mu$, on a $^t(2\mu)=(^t\mu)^2$. En appliquant (1), on voit que
$$2U_{\lambda,\epsilon}\sqcup \tilde{Y} =[A+B+1,0]_{1},$$
$$2V_{\lambda,\epsilon}\sqcup \tilde{X} =[A+B+1,0]_{1}.$$
On a aussi $A+B+1=4r$. D'o\`u
$$(2) \qquad 2U_{\lambda,\epsilon}\sqcup 2V_{\lambda,\epsilon}\sqcup \tilde{X}\sqcup \tilde{Y}=[4r,0]_{1}^2.$$
Posons $\eta_{A}=0$, $\eta_{B}=1$ si $k$ est pair  et $\eta_{A}=1$, $\eta_{B}=0$  si $k$ est impair (rappelons que $A$ et $k$ sont de m\^eme parit\'e tandis que $B$ et $k$ sont de parit\'e oppos\'ee).  Par construction, $X =({^t\beta}+[A,\eta_{A}]_{2})\sqcup ({^t\beta}+[A-1,1-\eta_{A}]_{2})$, $Y =({^t\alpha}+[B,\eta_{B}]_{2})\sqcup ({^t\alpha}+[B-1,1-\eta_{B}]_{2})$. En posant
$$\Lambda_{{^t\beta},{^t\alpha}}=A_{{^t\beta},{^t\alpha}}\sqcup B_{{^t\beta},{^t\alpha}}=({^t\beta}+[A,\eta_{A}]_{2})\sqcup ({^t\alpha}+[B,\eta_{B}]_{2})$$
et en notant 
$\Lambda_{{^t\beta},{^t\alpha}}^-$ la partition form\'ee des $x-1$ pour $x\in \Lambda_{{^t\beta},{^t\alpha}}$, $x\not=0$, on obtient
$$X \sqcup Y =\Lambda_{{^t\beta},{^t\alpha}}\sqcup \Lambda_{{^t\beta},{^t\alpha}}^-.$$
Le couple $ (\alpha,\beta)$ param\`etre la repr\'esentation $\rho_{\lambda,\epsilon}$ et $({^t\beta},{^t\alpha})$ param\`etre la repr\'esentation $sgn\otimes \rho_{\lambda,\epsilon}$. Autrement dit, $({^t\beta},{^t\alpha})$ est le couple associ\'e \`a $({^s\lambda},{^s\epsilon})$. Pour simplifier, posons $(\mu,\tau)=({^s\lambda},{^s\epsilon})$. On a donc $\Lambda_{{^t\beta},{^t\alpha}}=\Lambda_{\mu,\tau}$, o\`u $\Lambda_{\mu,\tau}=A_{\mu,\tau}\sqcup B_{\mu,\tau}$. D'o\`u, en d\'efinissant $\Lambda_{\mu,\tau}^-$ comme ci-dessus,

(3) $X \sqcup Y =\Lambda_{\mu,\tau}\sqcup \Lambda_{\mu,\tau}^-$.

Introduisons les partitions  $\mu'$ et $\mu''$ d\'efinies par $\mu'_{j}=[\mu_{j}/2]$ pour $j=1,...,2r+1$, $\mu''_{j}=[(\mu_{j}+1)/2]$ pour $j=1,...,2r$. Montrons que
$$(4) \qquad \Lambda_{\mu,\tau}\sqcup \Lambda_{\mu,\tau}^-=(\mu'+[2r,0]_{1})\sqcup(\mu''+[2r-1,0]_{1})
.$$
On  a 
$\Lambda_{\mu,\tau}=A_{\mu,\tau}\sqcup B_{\mu,\tau} =A^{\sharp}_{\mu}\sqcup B^{\sharp}_{\mu}$ et ce dernier terme est calcul\'e par les   formules (1)(a) et (1)(b) de 5.1. Pour un indice $j=1,...,2r$ tel que $\mu_{j}$ est pair, $\mu_{j}$ contribue \`a $ \Lambda_{\mu,\tau}$ par un terme $\mu_{j}/2+2r+1-j$, donc \`a $\Lambda_{\mu,\tau}\sqcup \Lambda_{\mu,\tau}^-$ par deux termes $\mu_{j}/2+2r+1-j$ et $\mu_{j}/2+2r-j$. Mais le premier, resp. le second, est la contribution de l'indice $j$ \`a $\mu'+[2r,0]_{1}$, resp. $\mu''+[2r-1,0]_{1}$. Consid\'erons la contribution des termes impairs de $\mu$. Puisque $\mu$ est symplectique, on peut les regrouper par paires de termes \'egaux. Consid\'erons donc un couple de termes impairs de la forme $\mu_{2j-1}=\mu_{2j}$. Ils contribuent \`a $\Lambda_{\mu,\tau}$ par les termes $A_{\mu,j}^{\sharp}$ et $B_{\mu,j}^{\sharp}$, tous deux \'egaux \`a $(\mu_{2j}+1)/2+2r+1-2j$. Ils contribuent \`a $\Lambda_{\mu,\tau}\sqcup \Lambda_{\mu,\tau}^-$ par quatre termes:  deux fois $(\mu_{2j}+1)/2+2r+1-2j$ et deux fois  $(\mu_{2j}-1)/2+2r+1-2j$. Les indices $2j-1,2j$  contribuent \`a $\mu'+[2r,0]_{1}$ par $(\mu_{2j}+1)/2+2r+1-2j$ et $(\mu_{2j}-1)/2+2r+1-2j$. Ils contribuent \`a $\mu''+[2r-1,0]_{1}$ par les m\^emes termes. Les contributions  aux deux membres de (4) sont donc les m\^emes. Consid\'erons maintenant  un couple de termes impairs de la forme $\mu_{2j}=\mu_{2j+1}$. Ils contribuent \`a $\Lambda_{\mu,\tau}$ par les termes $A_{\mu,j+1}^{\sharp}$ et $B_{\mu,j}^{\sharp}$, tous deux \'egaux \`a $(\mu_{2j}+1)/2+2r-2j$. Ils contribuent \`a  $\Lambda_{\mu,\tau}\sqcup \Lambda_{\mu,\tau}^-$ par quatre termes: deux fois $(\mu_{2j}+1)/2+2r-2j$ et deux fois $(\mu_{2j}-1)/2+2r-2j$.  Les indices $2j,2j+1$  contribuent \`a $\mu'+[2r,0]_{1}$ par $(\mu_{2j}+1)/2+2r-2j$ et $(\mu_{2j}-1)/2+2r-2j$ et ils contribuent \`a $\mu''+[2r-1,0]_{1}$ par les m\^emes termes. 
Les contributions  aux deux membres de (4) sont encore une fois les m\^emes. Enfin, l'indice $2r+1$ contribue \`a $\Lambda_{\mu,\tau}$ par $A_{\mu,r+1}^{\sharp}=0$, qui dispara\^{\i}t dans $\Lambda_{\mu,\tau}^-$. De m\^eme, cet indice contribue \`a $\mu'+[2r,0]_{1}$ par $0$ et ne contribue pas \`a $\mu''+[2r-1,0]_{1}$. Cela prouve (4). 

Par d\'efinition,  
$\tilde{X}\sqcup\tilde{Y}$ est la famille des $4r-x$ pour $x\in X\sqcup Y$. D'apr\`es (3) et (4), c'est donc la r\'eunion de $\tilde{F}'=\{4r-x; x\in \mu'+[2r,0]_{1}\}$ et de $\tilde{F}''=\{4r-x; x\in \mu''+[2r-1,0]_{1}\}$. En posant $E'={^t\mu'}+[2r-1,0]_{1}$ et $E''={^t\mu''}+[2r,0]_{1}$, (1) implique que $E'\sqcup \tilde{F}'=[4r,0]_{1} =E''\sqcup \tilde{F}'$. D'o\`u
$$E'\sqcup E''\sqcup \tilde{X}\sqcup \tilde{Y}=[4r,0]_{1}^2.$$
Avec (2), on obtient
$$(5) \qquad 2U_{\lambda,\epsilon}\sqcup 2V_{\lambda,\epsilon}=E'\sqcup E''.$$
 Montrons que

(6) pour tout $j=1,...,2r$, on a ${^t\mu''}_{j}\geq {^t\mu'}_{j}\geq {^t\mu''}_{j+1}$. 

On rappelle que, pour toute partition $\nu$, ${^t\nu}_{j}$ est le nombre d'indices $h$ tels que $\nu_{h}\geq j$. Si $\mu'_{h}=[\mu_{h}/2]\geq j$, on a aussi $\mu''_{h}=[(\mu_{h}+1)/2]\geq j$, d'o\`u ${^t\mu'}_{j}\leq {^t\mu''}_{j}$. Si $\mu''_{h}=[(\mu_{h}+1)/2]\geq j+1$, on a aussi $\mu'_{h}=[\mu_{h}/2]\geq j$, d'o\`u ${^t\mu''}_{j+1}\leq {^t\mu'}_{j}$. Cela prouve (6). 

Il r\'esulte de (6) que
$${^t\mu'}\sqcup {^t\mu''}=({^t\mu''}_{1},{^t\mu'}_{1},{^t\mu''}_{2},{^t\mu'}_{2},...,{^t\mu''}_{2r},{^t\mu'}_{2r},{^t\mu''}_{2r+1}).$$
Ecrivons $E'=(e'_{1},...,e'_{2r})$, $E''=(e''_{1},...,e''_{r+1})$. Il resulte de (6) que l'on a aussi $e''_{j}>e'_{j}\geq e''_{j+1}$, donc
$$E'\sqcup E''=(e''_{1},e'_{1},e''_{2},e'_{2},...,e''_{2r},e'_{2r},e''_{2r+1})$$
$$=({^t\mu''}_{1}+2r,{^t\mu'}_{1}+2r-1,{^t\mu''}_{2}+2r-1,{^t\mu'}_{2}+2r-2,...,{^t\mu''}_{2r}+1,{^t\mu'}_{2r},{^t\mu''}_{2r+1})$$
$$=({^t\mu'}\sqcup {^t\mu''})+([2r,0]_{1}\sqcup[2r-1,0]_{1}).$$
On sait que la transposition des partitions \'echange la somme et la r\'eunion. Donc $({^t\mu'}\sqcup {^t\mu''})={^t(\mu'+\mu'')}$. Mais, par d\'efinition, $\mu'+\mu''=\mu$. D'o\`u 
$$E'\sqcup E''={^t\mu}+([2r,0]_{1}\sqcup[2r-1,0]_{1}).$$
Avec (5), on obtient l'\'enonc\'e. $\square$

\bigskip
\subsection{Un th\'eor\`eme de minimalit\'e}
 \ass{Th\'eor\`eme}{Soit $(\lambda,\epsilon)\in \boldsymbol{{\cal P}^{symp}}(2N)$. Supposons que $\lambda$ n'a que des termes pairs. Alors il existe un unique \'el\'ement $(\lambda^{min},\epsilon^{min})\in \boldsymbol{{\cal P}^{symp}}(2N)$ v\'erifiant les propri\'et\'es suivantes:
 
 (i) $\mathfrak{mult}(\lambda,\epsilon;{^s\lambda}^{min},{^s\epsilon}^{min})=1$;
 
 (ii) pour tout \'el\'ement $(\lambda',\epsilon')\in \boldsymbol{{\cal P}^{symp}}(2N)$ tel que $\mathfrak{mult}(\lambda,\epsilon;{^s\lambda'},{^s\epsilon'})\not=0$, on a $\lambda^{min}<\lambda'$ ou $(\lambda',\epsilon')=(\lambda^{min},\epsilon^{min})$. 
 
 De plus, on a $({^s\lambda}^{min},{^s\epsilon}^{min})=(\lambda^{max},\epsilon^{max})$.}
 
 Preuve. On pose $k=k(\lambda,\epsilon)$. De nouveau, seuls interviennent les couples $(\lambda',\epsilon')$ tels que $k(\lambda',\epsilon')=k$, c'est-\`a-dire les \'el\'ements de $\boldsymbol{{\cal P}^{symp}}(2N)_{k}$. On fixe $r$ assez grand  et on pose $A=2r+k$ et $B=2r-1-k$ comme en 5.2.    Pour  $(\lambda',\epsilon') \in \boldsymbol{{\cal P}^{symp}}(2N)_{k}$, on  a d\'efini dans le paragraphe pr\'ec\'edent les suites $U_{\lambda',\epsilon'}$ et $V_{\lambda',\epsilon'}$. Montrons que, pour $(\lambda',\epsilon'), (\lambda'',\epsilon'') \in \boldsymbol{{\cal P}^{symp}}(2N)_{k}$, on a
  l'\'equivalence
 
 (1) $^s\lambda'\leq {^s\lambda''}$ si et seulement si $U_{\lambda'',\epsilon''}\sqcup V_{\lambda'',\epsilon''}\leq U_{\lambda',\epsilon'}\sqcup V_{\lambda',\epsilon'}$.
 
 Evidemment $U_{\lambda'',\epsilon''}\sqcup V_{\lambda'',\epsilon''}\leq U_{\lambda',\epsilon'}\sqcup V_{\lambda',\epsilon'}$ \'equivaut \`a $2U_{\lambda'',\epsilon''}\sqcup 2V_{\lambda'',\epsilon''}\leq 2U_{\lambda',\epsilon'}\sqcup 2V_{\lambda',\epsilon'}$. D'apr\`es le lemme 5.6, cela \'equivaut \`a $^t(^s\lambda'')\leq{^t(^s\lambda')}$. On sait bien que la transposition dans ${\cal P}(2N)$ est d\'ecroissante. La condition pr\'ec\'edente \'equivaut donc \`a $^s\lambda'\leq{^s\lambda''}$. Cela prouve (1).

On peut modifier l\'eg\`erement l'assertion (1). On note $(\alpha',\beta')$ et $(\alpha'',\beta'')$ les \'el\'ements de ${\cal P}_{2}(N-k(k+1)/2)$ correspondant \`a $(\lambda',\epsilon')$ et $(\lambda'',\epsilon'')$. On a d\'efini en 5.2 les entiers $n=r+[k/2]+1$ et $m=r-[k/2]$. En 5.6, on a pos\'e  $\underline{n}=2r+k+1$ et $\underline{m}=2r-k$. Puisqu'on suppose $r$ grand, on a $\underline{n}\geq n$ et $\underline{m}\geq m$. On a consid\'er\'e en 5.6  que $\alpha'$ et $\beta'$ avaient respectivement $\underline{n}$ et $\underline{m}$ termes. Mais, d'apr\`es  la construction de  5.1, $\alpha'$, resp. $\beta'$, a au plus $n$, resp. $m$, termes non nuls. En consid\'erant maintenant que $\alpha'\in {\cal P}_{n}$ et $\beta'\in {\cal P}_{m}$, on a
$$U_{\lambda',\epsilon'}=(\alpha'+[A/2,(A+1-n)/2]_{1/2})\sqcup [(A-n)/2,0]_{1/2},$$
$$ V_{\lambda',\epsilon'}=(\beta'+[B/2,(B+1-m)/2]_{1/2})\sqcup [(B-m)/2,0]_{1/2}.$$
Evidemment, l'in\'egalit\'e $U_{\lambda'',\epsilon''}\sqcup V_{\lambda'',\epsilon''}\leq U_{\lambda',\epsilon'}\sqcup V_{\lambda',\epsilon'}$ \'equivaut \`a
$$(\alpha''+[A/2,(A+1-n)/2]_{1/2})\sqcup (\beta''+[B/2,(B+1-m)/2]_{1/2})\leq $$
$$(\alpha'+[A/2,(A+1-n)/2]_{1/2})\sqcup(\beta'+[B/2,(B+1-m)/2]_{1/2}).$$
Remarquons que 
$$  (\alpha'+[A/2,(A+1-n)/2]_{1/2})\sqcup(\beta'+[B/2,(B+1-m)/2]_{1/2})= p\Lambda^{n,m}_{A/2,B/2;1/2}(\alpha',\beta').$$
On a donc

(2) ${s\lambda'}\leq {^s\lambda''}$ si et seulement si $p\Lambda^{n,m}_{A/2,B/2;1/2}(\alpha'',\beta'')\leq p\Lambda^{n,m}_{A/2,B/2;1/2}(\alpha',\beta')$.

On peut reformuler la premi\`ere assertion de l'\'enonc\'e sous la forme

il existe un unique \'el\'ement $(\underline{\lambda},\underline{\epsilon})\in \boldsymbol{{\cal P}^{symp}}(2N)_{k}$ v\'erifiant les propri\'et\'es

(i) $\mathfrak{mult}(\lambda,\epsilon;\underline{\lambda},\underline{\epsilon})=1$;

(ii) pour tout $(\lambda',\epsilon')\in \boldsymbol{{\cal P}^{symp}}(2N)_{k}$ tel que $\mathfrak{mult}(\lambda,\epsilon;\lambda',\epsilon')\not=0$, on a ${^s\underline{\lambda}}<{^s\lambda'}$ ou $(\lambda',\epsilon')=(\underline{\lambda},\underline{\epsilon})$.

Le couple $(\lambda^{min},\epsilon^{min})$ de l'\'enonc\'e est alors $(^s\underline{\lambda},{^s\underline{\epsilon}})$. La derni\`ere assertion de l'\'enonc\'e affirme que $(\underline{\lambda},\underline{\epsilon})=(\lambda^{max},\epsilon^{max})$.
 
  L'assertion (2) et la proposition 5.2 permet de reformuler les assertions ci-dessus sous la forme suivante:

il existe un unique couple $(\underline{\alpha},\underline{\beta})\in {\cal P}_{n}\times {\cal P}_{m}$ v\'erifiant les conditions

(i) $mult(<_{I,\alpha,\beta};\alpha,\beta;\underline{\alpha},\underline{\beta})=1$;

(ii) pour tout couple $(\alpha',\beta')\in {\cal P}_{n}\times {\cal P}_{m}$ tel que $mult(<_{I,\alpha,\beta};\alpha,\beta;\alpha',\beta')\not=0$, on a $p\Lambda_{A/2,B/2;1/2}^{n,m}(\alpha',\beta')< p\Lambda_{A/2,B/2;1/2}^{n,m}(\underline{\alpha},\underline{\beta})$ ou $(\alpha',\beta')=(\underline{\alpha},\underline{\beta})$.

La proposition 4.1 entra\^{\i}ne que ceci est v\'erifi\'e si et seulement si l'ensemble $P_{A/2,B/2;1/2}(<_{I,\alpha,\beta};\alpha,\beta)$ a un unique \'el\'ement et, dans ce cas, $(\underline{\alpha},\underline{\beta})$ est cet unique \'el\'ement. On veut de plus prouver que cet \'el\'ement est \'egal \`a l'\'el\'ement $(\alpha^{max},\beta^{max})$ param\'etrant $(\lambda^{max},\epsilon^{max})$. Dans la  preuve du th\'eor\`eme 5.5, on a vu que $(\alpha^{max},\beta^{max})$ \'etait l'unique \'el\'ement de $P_{A,B;2}(<_{I,\alpha,\beta};\alpha,\beta)$. Il nous suffit donc de prouver que 
$$(3)\qquad P_{A/2,B/2;1/2}(<_{I,\alpha,\beta};\alpha,\beta)=P_{A,B;2}(<_{I,\alpha,\beta};\alpha,\beta).$$

Pour prouver cette assertion, g\'en\'eralisons un peu les hypoth\`eses. On consid\`ere des entiers $n,m\in {\mathbb N}$ et des partitions $\alpha\in {\cal P}_{n}$, $\beta\in {\cal P}_{m}$. On fixe des entiers $A,B$ tels que $A\geq 2(n-1)$ et $B\geq 2(m-1)$, on d\'efinit la  partition $p\Lambda_{A,B;2}^{n,m}(\alpha,\beta)$ et on suppose que cette suite est sans multiplicit\'es. On d\'efinit l'ordre $<_{I,A,B,\alpha,\beta}$ sur l'ensemble d'indices $I=(\{1,...,n\}\times\{0\})\cup (\{1,...,m\}\times \{1\})$ de la fa\c{c}on suivante: pour $i\in \{1,...,n\}$ et $j\in \{1,...,m\}$, on a $(i,0)<_{I,A,B,\alpha,\beta}(j,1)$ si et seulement si $\alpha_{i}+A+2-2i> \beta_{j}+B+2-2j$. On fixe des entiers $e,f\in {\mathbb N}$ et on pose les hypoth\`eses suivantes, pour tout $i\in\{1,...,n\}$ et $j\in \{1,...,m\}$:

(4) si $(i,0)<_{I,A,B,\alpha,\beta}(j,1)$, alors $\alpha_{i}+A+2-2i\geq 2f+1+ \beta_{j}+B+2-2j$;

(5) si $(j,1)<_{I,A,B,\alpha,\beta}(i,0)$, alors $\beta_{j}+B+2-2j\geq 2e+1+\alpha_{i}+A+2-2i$.

Sous ces hypoth\`eses, on va prouver que
$$(6)\qquad P_{(A+e)/2,(B+f)/2;1/2}(<_{I,A,B,\alpha,\beta};\alpha,\beta)=P_{A,B;2}(<_{I,A,B\alpha,\beta};\alpha,\beta).$$

{\bf Remarque.} Dans la situation qui nous int\'eresse, les hypoth\`eses   ci-dessus sont v\'erifi\'ees pour $e=f=0$ (les conditions (4) et (5) sont alors tautologiques par d\'efinition de l'ordre  $<_{I,A,B,\alpha,\beta}$. La conclusion (6) co\"{\i}ncide avec (3).
\bigskip

Si $n=0$ ou $m=0$, la conclusion est triviale puisque $P(\alpha,\beta)$ n'a qu'un \'el\'ement. On suppose $n\geq1$ et $m\geq1$ et on note simplement $<_{I}$ l'ordre  $<_{I,A,B,\alpha,\beta}$. On ne perd rien \`a supposer que $(1,0)<_{I}(1,1)$. On reprend les notations de la preuve de 5.5. Le proc\'ed\'e (a) est le seul autoris\'e pour construire les \'el\'ements de $P_{A,B;2}(<_{I};\alpha,\beta)$. Ce proc\'ed\'e construit un  \'el\'ement $\nu_{1}$ et des partitions $(\alpha',\beta')\in {\cal P}_{n'}\times {\cal P}_{m'}$. L'ordre $<_{I}$ induit un ordre $<_{I'}$ sur l'ensemble d'indices $I'=(\{1,...,n'\}\times \{0\})\cup(\{1,...,m'\}\times \{1\})$. Les \'el\'ements de   $P_{A,B;2}(<_{I};\alpha,\beta)$ sont les $(\{\nu_{1}\}\sqcup\nu'\sqcup\{0^{n-n'-1}\},\mu'\sqcup\{0^{m-m'}\})$, pour $(\nu',\mu')\in P_{A-2,B;2}(<_{I'},\alpha',\beta')$. Montrons que le proc\'ed\'e (a) est aussi le seul autoris\'e pour construire les \'el\'ements de $P_{(A+e)/2,(B+f)/2;1/2}(<_{I},\alpha,\beta)$, autrement dit que l'on a $\alpha_{1}+(A+e)/2> (B+f)/2$. Il suffit de prouver que $\alpha_{1}> (B-A+f)/2$. C'est automatique si $B-A+f<0$. Supposons $B-A+f\geq0$. On a alors
$$(B-A+f)/2\leq B-A+f<B-A+2f+1+\beta_{1}\leq \alpha_{1}$$
d'apr\`es (4) pour $i=j=1$. Cela d\'emontre l'in\'egalit\'e cherch\'ee. Donc les \'el\'ements de   $P_{(A+e)/2,(B+f)/2;1/2}(<_{I};\alpha,\beta)$ sont les $(\{\nu_{1}\}\sqcup\nu'\sqcup\{0^{n-n'-1}\},\mu'\sqcup\{0^{m-m'}\})$, pour $(\nu',\mu')\in P_{(A+e-1)/2,(B+f)/2;1/2}(<_{I'},\alpha',\beta')$. Il nous suffit de prouver que
$$ P_{(A+e-1)/2,(B+f)/2;1/2}(<_{I'},\alpha',\beta')=P_{A-2,B;2}(<_{I'},\alpha',\beta').$$
De nouveau, c'est trivial si $n'=0$ ou $m'=0$. On suppose $n'\geq1$, $m'\geq1$. 
On pose $A'=A-2$, $B'=B$, $e'=e+1$, $f'=f$. On a vu dans la preuve du th\'eor\`eme 5.5 que la partition $p\Lambda_{A',B';2}^{n',m'}(\alpha',\beta')$ \'etait sans multiplicit\'es et que l'ordre $<_{I'}$ co\"{\i}ncidait avec $<_{I',A',B'\alpha',\beta'}$.
L'\'egalit\'e ci-dessus se r\'ecrit donc
$$P_{(A'+e')/2,(B'+f')/2;1/2}(<_{I',A',B'\alpha',\beta'};\alpha',\beta')=P_{A',B';2}(<_{I',A',B',\alpha',\beta'};\alpha',\beta').$$
C'est la relation (6) pour les donn\'ees affect\'ees d'un $'$. En raisonnant par r\'ecurrence, il suffit de prouver que les analogues des conditions (4) et (5) sont v\'erifi\'ees par ces donn\'ees. Soient donc $i\in \{1,...,n'\}$ et $j\in \{1,...,m'\}$. Supposons $(i,0)<_{I'}(j,1)$. Avec les notations de 5.5(7), on a $(i+k(i),0)<_{I}(j+h(j),1)$. On a  alors prouv\'e en 5.5(9) que $\alpha'_{i}+A-2i> \beta'_{j}+B+2-2j$. La preuve utilisait l'in\'egalit\'e 
$$\alpha_{i+k(i)}+A+2-2i-2k(i)>\beta_{b_{k(i)}}+B+2-2b_{k(i)}.$$
D'apr\`es (4), on a l'in\'egalit\'e plus forte
 $$\alpha_{i+k(i)}+A+2-2i-2k(i)>2f+\beta_{b_{k(i)}}+B+2-2b_{k(i)}.$$
Alors, la  m\^eme preuve qu'en 5.5(9) conduit \`a l'in\'egalit\'e $\alpha'_{i}+A-2i> 2f+\beta'_{j}+B+2-2j$. Puisqu'il s'agit d'entiers, cela \'equivaut \`a $\alpha'_{i}+A'+2-2i\geq 1+2f+\beta'_{j}+B'+2-2j$. C'est l'analogue de (4). Supposons maintenant $(j,1)<_{I'}(i,0)$. Avec les notations de 5.5(11), on a $(j+h(j),1)<_{I}(i+k(i),0)$. On a alors prouv\'e que $\beta'_{j}+B+2-2j>2+\alpha'_{i}+A-2i$. La preuve utilisait l'in\'egalit\'e
$$\beta_{j+h(j)}+B+2-2j-2h(j)>\alpha_{a_{h(j)+1}}+A+2-2a_{h(j)+1}.$$
 D'apr\`es (5), on a l'in\'egalit\'e plus forte
$$\beta_{j+h(j)}+B+2-2j-2h(j)>2e+\alpha_{a_{h(j)+1}}+A+2-2a_{h(j)+1}.$$
Alors, la suite  la m\^eme preuve qu'en 5.5(11) conduit \`a l'in\'egalit\'e $\beta'_{j}+B+2-2j>2+2e+\alpha'_{i}+A-2i$. Cela \'equivaut \`a $\beta'_{j}+B'+2-2j\geq 1+2e'+\alpha'_{i}+A'+2-2i$. C'est l'analogue de (5). Cela ach\`eve la preuve de (6) et du th\'eor\`eme. $\square$

\bigskip
   
 \section{Calcul explicite du couple $(\lambda^{max},\epsilon^{max})$}

\bigskip
\subsection{D\'efinition d'un couple $(\bar{\lambda},\bar{\epsilon})$}

Soit $(\lambda,\epsilon)\in \boldsymbol{{\cal P}^{symp}}(2N)$. On suppose que tous les termes de $\lambda$ sont pairs.  On va associer par r\'ecurrence \`a $(\lambda,\epsilon)$ un autre  couple $(\bar{\lambda},\bar{\epsilon})\in \boldsymbol{{\cal P}^{symp}}(2N)$.

Si $N=0$, $(\bar{\lambda},\bar{\epsilon})=(\lambda,\epsilon)$ est l'unique \'el\'ement de $\boldsymbol{{\cal P}^{symp}}(0)$.

On suppose $N>0$. On repr\'esente $\lambda$ sous la forme $\lambda=(\lambda_{1},...,\lambda_{2r+1})$, o\`u $\lambda_{2r+1}=0$. Comme en 5.1, on consid\`ere que $\epsilon$ est d\'efini  soit sur $Jord^{bp}(\lambda)\cup\{0\}=Jord(\lambda)\cup\{0\}$ (avec $\epsilon_{0}=1$), soit sur l'ensemble d'indices $\{1,...,2r+1\}$: $\epsilon(j)=\epsilon_{\lambda_{j}}$. Notons $J^{a}=\{j=1,...,2r+1; (-1)^{j+1}\epsilon(j)=1\}$ et $J^{b}=\{j=1,...,2r+1; (-1)^j\epsilon(j)=1\}$.  
Notons $\mathfrak{S}$ la r\'eunion de $\{1\}$ et de l'ensemble des $j\in \{2,...,2r+1\}$ tels que $\epsilon(j)(-1)^j\not=\epsilon(j-1)(-1)^{j-1}$. Ecrivons $\mathfrak{S}$ comme une suite croissante $1=s_{1}<s_{2}<...<s_{S}$. Posons
$$\bar{\lambda}_{1}=\left\lbrace\begin{array}{cc}(\sum_{h=1,...,S}\lambda_{s_{h}})+S-1-2\vert J^{b}\vert ,& \text{ si }\epsilon(1)=1;\\ (\sum_{h=1,...,S}\lambda_{s_{h}})+S-2\vert J^{a}\vert ,& \text{ si }\epsilon(1)=-1;\\ \end{array}\right.$$
Nous montrerons que

(1) $\bar{\lambda}_{1}$ est pair et on a $2\leq \bar{\lambda}_{1}\leq 2N$.

Posons $N'=N-\bar{\lambda}_{1}/2$. Notons $\lambda'$ la r\'eunion des $\lambda_{j}$ pour $j\in \{1,...,2r+1\}$, $j\not\in \mathfrak{S}$, $\epsilon(j)(-1)^j=-\epsilon(1)$ et des $\lambda_{j}+2$ pour $j\in \{1,...,2r+1\}$, $j\not\in \mathfrak{S}$, $\epsilon(j)(-1)^j=\epsilon(1)$. Les termes de cette partition sont tous pairs.  Ils sont en nombre $2r+1-S$.  Soit $i\in Jord^{bp}(\lambda')$. Par d\'efinition, on peut fixer $j\in \{1,...,2r+1\}$, $j\not\in \mathfrak{S}$, tel que, soit $\epsilon(j)(-1)^j=-\epsilon(1)$ et $i=\lambda_{j}$, soit $\epsilon(j)(-1)^j=\epsilon(1)$ et $i=\lambda_{j}+2$. On note $h[j]$ le plus grand entier $h\in \{1,...,S\}$ tel que $s_{h}<j$. On pose $\epsilon'_{i}=(-1)^{h[j]+1}\epsilon(j)$. Nous montrerons que

(2) cette d\'efinition ne d\'epend pas du choix de $j$;

(3) le couple $(\lambda',\epsilon')$ ainsi d\'efini appartient \`a $\boldsymbol{{\cal P}^{symp}}(2N')$.

Puisque $N'<N$ d'apr\`es (1), on sait par r\'ecurrence associer \`a $(\lambda',\epsilon')$ un couple $(\bar{\lambda}',\bar{\epsilon}')$. Nous montrerons que

(4) $\bar{\lambda}_{1}\geq\bar{\lambda}'_{1}$.

On note $\bar{\lambda}$ la partition $\{\bar{\lambda}_{1}\}\sqcup \bar{\lambda}'$. Elle appartient \`a ${\cal P}^{symp}(2N)$ d'apr\`es (1) et (3). Nous montrerons que

(5) il existe un unique $\bar{\epsilon}\in \{\pm 1\}^{Jord^{bp}(\bar{\lambda})}$ tel que, pour $i\in Jord^{bp}(\bar{\lambda}')$, on ait $\bar{\epsilon}_{i}=\bar{\epsilon}'_{i}$ et que $\bar{\epsilon}_{\bar{\lambda}_{1}}=\epsilon(1)$.

On a ainsi d\'efini le couple $(\bar{\lambda},\bar{\epsilon})\in \boldsymbol{{\cal P}^{symp}}(2N)$. A priori, ce couple d\'epend de l'entier $r$ que l'on a choisi plus haut. En fait

(6) ce couple est ind\'ependant de $r$.

 On reporte \`a 6.3 les d\'emonstrations des assertions (1) \`a (6).

\bigskip

\subsection{Comparaison avec  le couple du th\'eor\`eme 5.5}
Soit $(\lambda,\epsilon)\in \boldsymbol{{\cal P}^{symp}}(2N)$. On suppose que tous les termes de $\lambda$ sont pairs.  

\ass{Proposition}{Le couple  $(\lambda^{max},\epsilon^{max})$ intervenant dans le th\'eor\`eme 5.5 est \'egal au couple $(\bar{\lambda},\bar{\epsilon})$ d\'efini dans le paragraphe pr\'ec\'edent.}

La d\'emonstration sera donn\'ee en 6.5.

\bigskip

\subsection{D\'emonstration des assertions de 6.1}
On reprend les hypoth\`eses et notations de 6.1 en supposant $N>0$. Remarquons que, si on d\'ecompose $J^{a}$ et $J^{b}$ en r\'eunions d'intervalles non cons\'ecutifs, les \'el\'ements de $\mathfrak{S}$ sont exactement les plus petits termes de ces intervalles. On a $1\in J^{a}$ si $\epsilon(1)=1$ et $1\in J^{b}$ si $\epsilon(1)=-1$. On a toujours $2r+1\in J^{a}$. En posant formellement $s_{S+1}=2r+2$, on en d\'eduit

(1) si $\epsilon(1)=1$, $S$ est impair, $J^{a}=\bigcup_{h=1,...,(S+1)/2}\{s_{2h-1},...,s_{2h}-1\}$, 

\noindent $J^{b}=\bigcup_{h=1,...,(S-1)/2}\{s_{2h},...,s_{2h+1}-1\}$; si $\epsilon(1)=-1$, $S$ est pair, 

\noindent $J^{a}=\bigcup_{h=1,...,S/2}\{s_{2h},...,s_{2h+1}-1\}$, $J^{b}=\bigcup_{h=1,...,S/2}\{s_{2h-1},...,s_{2h}-1\}$.
\bigskip

On a

(2) si $j,j+1$ sont deux indices cons\'ecutifs appartenant tous deux \`a $J^{a}$ ou tous deux \`a $J^{b}$, on a $\lambda_{j}\geq \lambda_{j+1}+2$.

En effet, la fonction $j'\mapsto \epsilon(j')(-1)^{j'+1}$ est constante sur $J^{a}$ comme sur $J^{b}$. On en d\'eduit $\epsilon(j)=-\epsilon(j+1)$. Cela interdit \`a $\lambda_{j}$ et $\lambda_{j+1}$ d'\^etre \'egaux. L'assertion s'ensuit puisque ces termes sont pairs par hypoth\`ese.

On en d\'eduit

(3) pour $h=1,...,S$, on a $\lambda_{s_{h}}-\lambda_{s_{h+1}}\geq 2(s_{h+1}-s_{h}-1)$.

C'est \'evident si $s_{h+1}=s_{h}+1$. Sinon, les termes $s_{h},s_{h}+1,...,s_{h+1}-1$ sont tous soit dans $J^{a}$, soit dans $J^{b}$. D'apr\`es (2), on a donc $\lambda_{s_{h}}-\lambda_{s_{h+1}-1}\geq 2(s_{h+1}-s_{h}-1)$. L'assertion s'ensuit puisque $\lambda_{s_{h+1}-1}\geq \lambda_{s_{h+1}}$.

En vertu de (1), on a

si $\epsilon(1)=1$, $2\vert J^b\vert+1-S =2\sum_{h=1,...(S-1)/2}(s_{2h+1}-s_{2h}-1)$,

si $\epsilon(1)=-1$, $2\vert J^{a}\vert -S=2\sum_{h=1,...,S/2}(s_{2h+1}-s_{2h}-1)$.

La d\'efinition de $\bar{\lambda}_{1}$ donn\'ee en 6.1 peut se r\'ecrire

$$(4) \qquad \bar{\lambda}_{1}=(\sum_{h=1,...,S}\lambda_{s_{h}})-2\sum_{h=1,...,[S/2]}(s_{2h+1}-s_{2h}-1).$$

D\'emontrons l'assertion (1) de 6.1, que nous r\'ep\'etons:

(5) $\bar{\lambda}_{1}$ est pair et on a $2\leq \bar{\lambda}_{1}\leq 2N$.

\noindent  Tous les termes de (4) sont pairs donc $\bar{\lambda}_{1}$ l'est. La relation (4) entra\^{\i}ne aussi $\bar{\lambda}_{1}\leq S(\lambda)=2N$. 
Cette formule (4) se r\'ecrit
$$\bar{\lambda}_{1}=(\sum_{h=1,...,[(S+1)/2]}\lambda_{s_{2h-1}})+(\sum_{h=1,...,[S/2]}(\lambda_{s_{2h}}-2(s_{2h+1}-s_{2h}-1)).$$
En vertu de (3), tous les termes des sommes ci-dessus sont positifs ou nuls. La premi\`ere somme contient $\lambda_{1}\geq 2$. Donc $\bar{\lambda}_{1} \geq 2$.  Cela prouve (5).

Montrons que

(6) $S(\lambda')=2N'$.

En effet, 
$$S(\lambda')=(\sum_{j=1,...,2r+1; j\not\in \mathfrak{S}}\lambda_{j})+2X$$
$$=2N-(\sum_{h=1,...,S}\lambda_{s_{h}})+2X,$$
o\`u $X$ est le nombre de $j\in \{1,...,2r+1\}$ tels que $j\not\in \mathfrak{S}$ et $\epsilon(j)(-1)^j=\epsilon(1) $.
D'apr\`es (1), on calcule $2X=2\vert J^{b}\vert -S+1$ si $\epsilon(1)=1$, $2X=2\vert J^{a}\vert -S$ si $\epsilon(1)=-1$. On obtient  alors $S(\lambda')=2N-\bar{\lambda}_{1}=2N'$. 

Les termes de $\lambda'$ sont des termes de $\lambda$ auxquels on ajoute \'eventuellement $2$. Ces  additions de $2$  ne perturbent pas l'ordre des termes. Pr\'ecis\'ement, 

(7) soient $j,k\in \{1,...,2r+1\}$, $j,k\not\in \mathfrak{S}$. Supposons $\epsilon(j)(-1)^j=-\epsilon(1)$ et $\epsilon(k)(-1)^k=\epsilon(1)$; alors $j<k$ si et seulement si $\lambda_{j}\geq \lambda_{k}+2$.

Si $\lambda_{j}\geq \lambda_{k}+2$, on a a fortiori $\lambda_{j}> \lambda_{k}$ donc $j<k$. Inversement, supposons $j<k$. Puisque $\epsilon(j)(-1)^j\not=\epsilon(k)(-1)^k$, il y a au moins un terme  de la suite $\mathfrak{S}$ dans l'intervalle $\{j+1,...,k-1\}$. Soit $s_{h}$ le plus grand terme de la suite $\mathfrak{S}$ dans cet intervalle.  D'apr\`es (2), on  a $\lambda_{s_{h}}\geq \lambda_{k}+2$. Or $\lambda_{j}\geq \lambda_{s_{h}}$. D'o\`u la conclusion $\lambda_{j}\geq \lambda_{k}+2$. 

Ecrivons $\lambda'=(\lambda'_{1},...,\lambda'_{2r+1-S})$. En vertu de (7), pour $j=1,...,2r+1-S$, on calcule $\lambda'_{j}$ de la fa\c{c}on suivante. Notons $h(j)$ le plus grand entier $h=1,...,S$ tel que
 $s_{h}-h<j$. Alors $\lambda'_{j}=\lambda_{j+h(j)}$ si $\epsilon(j+h(j))(-1)^{j+h(j)}=-\epsilon(1)$, $\lambda'_{j}=\lambda_{j+h(j)}+2$ si $\epsilon(j+h(j))(-1)^{j+h(j)}=\epsilon(1)$. On pose
 $$\tau'(j)=\epsilon(j+h(j))(-1)^{h(j)+1}.$$
 On va montrer que cette fonction $\tau'$  sur l'ensemble d'indices de $\lambda'$ provient d'un \'el\'ement  de $  \{\pm 1\}^{Jord^{bp}(\lambda')}$.  Autrement dit
 
 (8) pour $j,k\in \{1,...,2r+1-S\}$ tels que $\lambda'_{j}=\lambda'_{k}$, on a $\tau'(j)=\tau'(k)$; pour $j\in \{1,...,2r+1-S\}$ tel que $\lambda'_{j}=0$, on a $\tau'(j)=1$. 
 
 Pour la premi\`ere assertion, il suffit de traiter le cas o\`u $k=j+1$. Alors $j+h(j)$ et $j+1+h(j+1)$ sont deux termes cons\'ecutifs de l'ensemble $\{1,...,2r+1\}-\mathfrak{S}$. Autrement dit tous les \'el\'ements de $\{j+h(j)+1,...,j+h(j+1)\}$ appartiennent \`a $\mathfrak{S}$. On peut pr\'eciser qu'ils sont diff\'erents de $s_{1}=1$.  Par d\'efinition de la suite $\mathfrak{S}$, on a $\epsilon(s_{h})=\epsilon(s_{h}-1)$ pour tout $h\geq 2$. Donc $\epsilon$ est constant sur $\{j+h(j),..., j+h(j+1)\}$. Par contre, $\epsilon(j+1+h(j+1))=-\epsilon( j+h(j+1))$, sinon $j+1+h(j+1)$ serait dans la suite $\mathfrak{S}$. D'o\`u $\epsilon(j+h(j))=-\epsilon(j+1+h(j+1))$. A fortiori $\lambda_{j+h(j)}\not=\lambda_{j+1+h(j+1)}$, d'o\`u $\lambda_{j+h(j)}\geq \lambda_{j+1+h(j+1)}+2$. Par d\'efinition de la partition $\lambda'$, l'\'egalit\'e $\lambda'_{j}=\lambda'_{j+1}$ implique alors que  $\lambda'_{j}=\lambda_{j+h(j)}$ et $\lambda'_{j+1}=\lambda_{j+1+h(j+1)}+2$. Cela ne se produit que si $\epsilon(j+h(j))(-1)^{j+h(j)}=-\epsilon(1)$ et $\epsilon(j+1+h(j+1))(-1)^{j+1+h(j+1)}=\epsilon(1)$. Alors $\epsilon(j+h(j))(-1)^{h(j)+1}=\epsilon(j+1+h(j+1))(-1)^{h(j+1)+1}$, c'est-\`a-dire $\tau'(j)=\tau'(j+1)$. Cela prouve la premi\`ere assertion de (8). Supposons maintenant $\lambda'_{j}=0$. On a n\'ecessairement $\lambda'_{j}=\lambda_{j+h(j)}=0$. L'\'egalit\'e $\lambda'_{j}=\lambda_{j+h(j)}$ implique $\epsilon(j+h(j))(-1)^{j+h(j)}=-\epsilon(1)$. Puisque $j+h(j)$ appartient \`a l'intervalle $\{s_{h(j)},...,s_{h(j)+1}-1\}$, l'\'egalit\'e pr\'ec\'edente et (1) impliquent que $h(j)$ est impair. De plus, l'\'egalit\'e $\lambda_{j+h(j)}=0$ implique $\epsilon(j+h(j))=1$. Mais alors $\tau'(j)=\epsilon(j+h(j))(-1)^{h(j)+1}=1$. Cela d\'emontre (8). 
 
 Soit $i\in Jord^{bp}(\lambda')$. Fixons 
   $j\in \{1,...,2r+1\}$, $j\not\in \mathfrak{S}$, tel que, soit $\epsilon(j)(-1)^j=-\epsilon(1)$ et $i=\lambda_{j}$, soit $\epsilon(j)(-1)^j=\epsilon(1)$ et $i=\lambda_{j}+2$. Alors il existe un unique $k\in \{1,...,2r+1-S\}$ tel que $j=k+h(k)$. On a $i=\lambda'_{k}$ et $h[j]=h(k)$. Alors
   $$\epsilon'_{i}=(-1)^{h[j]+1}\epsilon(j)=(-1)^{h(k)+1}\epsilon(k+h(k))=\tau'(k).$$
   La relation (8) montre que ceci ne d\'epend pas du choix de $j$ et que $\epsilon'$ co\"{\i}ncide avec $\tau'$ vu comme un \'el\'ement de  $  \{\pm 1\}^{Jord^{bp}(\lambda')}$. Cela d\'emontre l'assertion (2) de 6.1. Avec (6), on a aussi l'assertion (3) de ce paragraphe.  Puisque $\epsilon'=\tau'$, nous abandonnons d\'esormais  la notation $\tau'$.

Montrons que

(9) on a $\bar{\lambda}_{1}\geq \bar{\lambda}'_{1}$; si $\bar{\lambda}_{1}=\bar{\lambda}'_{1}$, on a $\bar{\epsilon}'_{\bar{\lambda}'_{1}}= \epsilon(1)$.

L'assertion est triviale si $N'=0$ puisqu'alors $\bar{\lambda}'_{1}=0<\bar{\lambda}_{1}$ d'apr\`es (5). On suppose $N'>0$. Remarquons que, puisque $\bar{\epsilon}'_{\bar{\lambda}'_{1}}=\epsilon'(1)$ par d\'efinition, la deuxi\`eme assertion de (9) peut se reformuler: si $\bar{\lambda}_{1}=\bar{\lambda}'_{1}$, on a $ \epsilon'_{1}= \epsilon(1)$.

Pour calculer le terme $\bar{\lambda}'_{1}$, il faut repr\'esenter la partition $\lambda'$ sous la m\^eme forme que l'on a repr\'esent\'e $\lambda$ (remarquons que l'on sait par r\'ecurrence que $\bar{\lambda}'_{1}$ ne d\'epend pas de la forme choisie). Si $\epsilon(1)=1$, (1) implique que le nombre de termes $2r+1-S$ de $\lambda'$ est pair. On lui adjoint donc un terme $\lambda'_{2r+2-S}=0$ et on pose $r'=r+(1-S)/2$. Si $\epsilon(1)=-1$, (1) implique que $2r+1-S$ est impair. Mais on ne sait pas si le dernier terme $\lambda'_{2r+1-S}$ est nul. On adjoint donc deux termes $\lambda'_{2r+2-S}=\lambda'_{2r+3-S}=0$ et on pose $r'=r+1-S/2$. Cela \'etant, on  d\'efinit la suite $\mathfrak{S}'=\{s'_{1},...,s'_{S'}\}$ et les ensembles d'indices ${J'}^{a}$ et ${J'}^b$ comme on a d\'efini $\mathfrak{S}$, $J^{a}$ et $J^{b}$. On note  $\tilde{\mathfrak{S}}'=\{1,...,2r+1-S\}\cap \mathfrak{S}'$, $\tilde{J'}^{a}=\{1,...,2r+1-S\}\cap {J'}^{a}$ et $\tilde{J'}^{b}=\{1,...,2r+1-S\}\cap {J'}^{b}$. On v\'erifie que

 (10)(a) si $\epsilon(1)=1$, ${J'}^{a}=\tilde{J'}^{a}\sqcup\{2r+2-S\}$ et ${J'}^{b}=\tilde{J'}^{b}$;

(10)(b) si $\epsilon(1)=-1$, ${J'}^{a}=\tilde{J'}^{a}\sqcup\{2r+3-S\}$ et ${J'}^{b}=\tilde{J'}^{b}\sqcup\{2r+2-S\}$.

Notons aussi $\tilde{J}^{a}$, resp. $\tilde{J}^b$, l'ensemble des $j\in J^{a}$, resp. $j\in J^b$, qui n'appartiennent pas \`a la suite $\mathfrak{S}$. On a

(11) pour $j\in \{1,...,2r+1-S\}$, les propri\'et\'es $j\in \tilde{J'}^{a}$ et $j+h(j)\in \tilde{J}^{b}$ sont \'equivalentes et les propri\'et\'es $j\in \tilde{J'}^{b}$ et $j+h(j)\in \tilde{J}^{a}$ sont \'equivalentes.

Par exemple, $j\in \tilde{J'}^{a}$ si et seulement si $\epsilon'(j)(-1)^{j+1}=1$ et $j+h(j)\in \tilde{J}^{b}$ si et seulement si $\epsilon(j+h(j))(-1)^{j+h(j)+1}=-1$. La d\'efinition $\epsilon'(j)=\epsilon(j+h(j))(-1)^{h(j)+1}$ entra\^{\i}ne l'\'equivalence de ces relations.

Soit $\tilde{S}'$ le plus grand indice tel que $s'_{\tilde{S}'}$ appartienne \`a $\tilde{\mathfrak{S}}'$. Montrons que

(12) pour $k=1,...,\tilde{S}'$, on a  $s'_{k}+h(s'_{k})-1=s_{h(s'_{k})}$; cet indice appartient \`a $\mathfrak{S}\cap J^{b}$ si $s'_{k}\in \tilde{J'}^{a}$ et \`a $\mathfrak{S}\cap J^{a}$ si $s'_{k}\in \tilde{J'}^b$.

Supposons par exemple $s'_{k}\in \tilde{J'}^{a}$. D'apr\`es (11), $s'_{k}+h(s'_{k})$ appartient \`a $J^{b}$. Par d\'efinition de $\lambda'$, $s'_{k}+h(s'_{k})$ n'appartient pas \`a $\mathfrak{S}$. Donc $s'_{k}+h(s'_{k})-1\in J^{b}$. Si cet indice n'appartenait pas \`a $\mathfrak{S}$, $\lambda_{s'_{k}+h(s'_{k})-1}$ appara\^{\i}trait dans la suite $\lambda'$, forc\'ement \`a l'indice $s'_{k}-1$. On aurait alors $h(s'_{k}-1)=h(s'_{k})$ et, d'apr\`es (11), $s'_{k}-1$ appartiendrait \`a $\tilde{J'}^{a}$. Mais, par d\'efinition de la suite $\mathfrak{S}'$, les conditions $s'_{k}-1\in {J'}^{a}$ et $s'_{k}\in \mathfrak{S}'$ entra\^{\i}nent $s'_{k}\in {J'}^{b}$ contrairement \`a hypoth\`ese. Cette contradiction prouve que $s'_{k}+h(s'_{k})-1$ appartient \`a $\mathfrak{S}$.  Par d\'efinition de $h(s'_{k})$, c'est alors forc\'ement $s_{h(s'_{k})}$. Cela prouve (12). 

Il r\'esulte de (12) que l'application $k\mapsto h(s'_{k})$ est une injection de $\{1,...,\tilde{S}'\}$ dans $\{1,...,S\}$. Elle est \'evidemment croissante. 
On a

(13) si 
 $\epsilon(1)=\epsilon'(1)$, alors $1$ n'appartient pas \`a l'image de cette injection.
 
Puisque $s'_{1}=1$, il s'agit de prouver que $h(1)>1$. Supposons par exemple $\epsilon(1)=\epsilon'(1)=1$. La relation (1) appliqu\'ee \`a $\lambda$ et \`a $\lambda'$ implique que $1\in J^{a}$ et $1\in {J'}^{a}$. D'apr\`es (12), cette derni\`ere relation implique $s_{h(1)}\in J^{b}$ (puisque $s'_{1}+h(s'_{1})-1=h(1)$). Alors la relation $1\in J^{a}$ interdit \`a $h(1)$ d'\^etre \'egal \`a $1$. Cela prouve (13).

Supposons $\epsilon(1)=\epsilon'(1)$.
 On veut prouver que $\bar{\lambda}_{1}\geq \bar{\lambda}'_{1}$. D'apr\`es (4) appliqu\'e \`a $\lambda'$, on a $\bar{\lambda}'_{1}\leq \sum_{k=1,...,S'}\lambda'_{s'_{k}}$. De plus, $\lambda'_{s'_{k}}=0$ si $k> 2r+1-S$. Il suffit donc de prouver
$$\bar{ \lambda}_{1}\geq \sum_{k=1,...,\tilde{S}'}\lambda'_{s'_{k}}.$$
D'apr\`es (11) et (12), pour $k=1,...,\tilde{S}'$, les indices $s_{h(s'_{k})}$ et $s_{h(s'_{k})}+1=s'_{k}+h(s'_{k})$ appartiennent tous deux \`a $J^{a}$ ou tous deux \`a $J^{b}$.  D'apr\`es  (2), cela implique $\lambda_{s_{h(s'_{k})}}\geq \lambda_{s'_{k}+h(s'_{k})}+2$. Puisque $\lambda'_{s'_{k}}$ est \'egal \`a $\lambda_{s'_{k}+h(s'_{k})}$ ou \`a $\lambda_{s'_{k}+h(s'_{k})}+2$, on a en tout cas $\lambda'_{s'_{k}}\leq \lambda_{s_{h(s'_{k})}}$. Puisque l'application $k\mapsto h(s'_{k})$ est injective et que son image ne contient pas $1$, on obtient
$$\sum_{k=1,...,\tilde{S}'}\lambda'_{s'_{k}}\leq \sum_{ h=2,...,S}\lambda_{s_{h}}.$$
Il suffit alors de prouver que
$$ \bar{\lambda}_{1}\geq \sum_{ h=2,...,S}\lambda_{s_{h}}.$$
D'apr\`es (4), il suffit encore de prouver que
$$\lambda_{1}\geq 2\sum_{h=1,...,[S/2]}(s_{2h+1}-s_{2h}-1).$$
Mais la somme sur $h$ des relations (3) conduit \`a l'in\'egalit\'e plus forte
$$\lambda_{1}\geq 2\sum_{h=1,...,S}(s_{h+1}-s_{h}-1).$$
Cela d\'emontre   (9) dans le cas $\epsilon(1)=\epsilon'(1)$. 

On suppose maintenant $\epsilon(1)=1$ et $\epsilon'(1)=-1$. Pour $k=1,...,\tilde{S'}$, on a vu ci-dessus que $\lambda_{s_{h(s'_{k})}}\geq \lambda_{s'_{k}+h(s'_{k})}+2$. Par d\'efinition de la partition $\lambda'$, on a $\lambda'_{s_{k}}= \lambda_{s'_{k}+h(s'_{k})}$ si $s'_{k}+h(s'_{k})\in J^{a}$, ou encore si $s'_{k}\in {J'}^{b}$ d'apr\`es (11),  et  $\lambda'_{s_{k}}= \lambda_{s'_{k}+h(s'_{k})}+2$ si $s'_{k}+h(s'_{k})\in J^{b}$, ou encore $s'_{k}\in {J'}^{a}$ d'apr\`es (11). On a donc toujours $\lambda'_{s'_{k}}\leq \lambda_{s_{h(s'_{k})}}$ et m\^eme $\lambda'_{s'_{k}}\leq \lambda_{s_{h(s'_{k})}}-2$ si $s'_{k}\in {J'}^{b}$. En utilisant l'injectivit\'e de l'application $k\mapsto h(s'_{k})$, on obtient
$$\sum_{k=1,...,S'}\lambda'_{s'_{k}}=\sum_{k=1,...,\tilde{S}'}\lambda'_{s'_{k}}\leq (\sum_{h=1,...,S}\lambda_{s_{h}})-2\vert {J'}^{b}\cap \tilde{\mathfrak{S}}'\vert.$$ 
D'apr\`es (10)(a), ${J'}^{b}\cap \tilde{\mathfrak{S}}'={J'}^{b}\cap \mathfrak{S}'$. D'apr\`es (1) appliqu\'e \`a $\lambda'$ (en se souvenant de l'hypoth\`ese  $\epsilon'(1)=-1$), le nombre d'\'el\'ements de cet ensemble est $S'/2$. D'o\`u 
$$\sum_{k=1,...,S'}\lambda'_{s'_{k}}\leq (\sum_{h=1,...,S}\lambda_{s_{h}})-S'.$$
On a aussi
$$\bar{\lambda}'_{1}=(\sum_{k=1,...,S'}\lambda'_{s'_{k}})+S'-2\vert {J'}^{a}\vert ,$$
 d'o\`u
 $$\bar{\lambda}'_{1}\leq (\sum_{h=1,...,S}\lambda_{s_{h}})-2\vert {J'}^{a}\vert.$$
D'apr\`es (10)(a), $\vert {J'}^{a}\vert =\vert \tilde{J'}^{a}\vert +1$. Mais (11) \'etablit une bijection de $\tilde{J'}^{a}$ sur le compl\'ementaire $\tilde{J}^b$ de $J^b\cap \mathfrak{S}$ dans $J^b$. D'apr\`es (1), on a $\vert J^b\cap \mathfrak{S}\vert =(S-1)/2$. On obtient donc $\vert {J'}^{a}\vert =1+\vert J^b\vert -(S-1)/2=\vert J^b\vert +(3-S)/2$, puis
$$\bar{\lambda}'_{1}\leq (\sum_{h=1,...,S}\lambda_{s_{h}})-3+S-2\vert J^{b}\vert .$$
 Mais
 $$\bar{\lambda}_{1}=(\sum_{h=1,...,S}\lambda_{s_{h}})+S-1-2\vert J^{b}\vert$$
 d'o\`u
 $$\bar{\lambda}'_{1}\leq \bar{\lambda}_{1}-2,$$
 ce qui d\'emontre (9) dans le cas $\epsilon(1)=1$, $\epsilon'(1)=-1$.
 
 La preuve du cas $\epsilon(1)=-1$, $\epsilon'(1)=1$ est similaire. Pour $k=1,...,\tilde{S}'$, on a cette fois $\lambda'_{s'_{k}}\leq \lambda_{s_{h(s'_{k})}}$ et m\^eme $\lambda'_{s'_{k}}\leq \lambda_{s_{h(s'_{k})}}-2$ si $s'_{k}\in {J'}^{a}$. D'o\`u
 $$\sum_{k=1,...,S'}\lambda'_{s'_{k}} \leq (\sum_{h=1,...,S}\lambda_{s_{h}})-2\vert {J'}^{a}\cap \tilde{\mathfrak{S}}'\vert.$$
 On utilise (10)(b). 
  Puisque $2r+3-S\in {J'}^{a}$ et $2r+2-S\in {J'}^b$, on a $2r+3-S\in \mathfrak{S}'$ par d\'efinition de cette suite. Donc $ {J'}^{a}\cap \mathfrak{S}'=({J'}^{a}\cap \tilde{\mathfrak{S}}')\sqcup\{2r+3-S\}$ et $\vert  {J'}^{a}\cap \tilde{\mathfrak{S}}'\vert =\vert  {J'}^{a}\cap \mathfrak{S}'\vert -1$. D'apr\`es (1) appliqu\'e \`a $\lambda'$ (en se  rappelant l'hypoth\`ese  $\epsilon'(1)=1)$, on a $\vert  {J'}^{a}\cap \mathfrak{S}'\vert=(S'+1)/2$. D'o\`u
 $$\sum_{k=1,...,S'}\lambda'_{s'_{k}} \leq (\sum_{h=1,...,S}\lambda_{s_{h}})-S'+1.$$
On a aussi
$$\bar{\lambda}'_{1}=(\sum_{k=1,...,S'}\lambda'_{s'_{k}})+S'-1-2\vert {J'}^{b}\vert ,$$
d'o\`u
 $$\bar{\lambda}'_{1}\leq (\sum_{h=1,...,S}\lambda_{s_{h}})-2\vert {J'}^{b}\vert.$$
D'apr\`es (10)(b), $\vert {J'}^{b}\vert =\vert \tilde{J'}^{b}\vert +1$. Mais (11) \'etablit une bijection de $\tilde{J'}^{b}$ sur le compl\'ementaire  $\tilde{J}^{a}$ de $J^a\cap \mathfrak{S}$ dans $J^a$. D'apr\`es (1), on a $\vert J^a\cap \mathfrak{S}\vert =S/2$. On obtient donc $\vert {J'}^{b}\vert =1+\vert J^a\vert -S/2$, puis
$$\bar{\lambda}'_{1}\leq (\sum_{h=1,...,S}\lambda_{s_{h}})-2+S-2\vert J^{a}\vert .$$
 Mais
 $$\bar{\lambda}_{1}=(\sum_{h=1,...,S}\lambda_{s_{h}})+S-2\vert J^{a}\vert$$
 d'o\`u
 $$\bar{\lambda}'_{1}\leq \bar{\lambda}_{1}-2,$$
 ce qui d\'emontre (9) dans le cas $\epsilon(1)=1$, $\epsilon'(1)=-1$. Cela ach\`eve la preuve de (9). 

L'assertion (9) implique les deux assertions (4) et (5) de 6.1. Montrons enfin que les constructions ne d\'ependent pas de l'entier $r$ choisi. Il suffit de montrer qu'elles ne changent pas si l'on remplace $r$ par $r+1$. Par ce remplacement,  on ajoute \`a $\lambda$ deux termes $\lambda_{2r+2}=\lambda_{2r+3}=0$. On voit que les deux indices $2r+2$ et $2r+3$  s'ajoutent \`a la suite $\mathfrak{S}$, le premier s'ajoutant \`a $J^{b}$ et le second \`a $J^{a}$. Le nombre $S$ s'accro\^{\i}t de $2$ et on voit que ces modifications ne changent pas $\bar{\lambda}_{1}$. Puisque les indices ajout\'es appartiennent \`a $\mathfrak{S}$, la partition $\lambda'$ ne change pas, ni bien s\^ur $\epsilon'$.   Cela prouve l'assertion (6) de 6.1.

\bigskip

\subsection{Conservation de l'entier $k$}
Les hypoth\`eses sont les m\^emes qu'en 6.1.

\ass{Lemme}{On a l'\'egalit\'e $k(\lambda,\epsilon)=k(\bar{\lambda},\bar{\epsilon})$.}

Preuve. C'est trivial si $N=0$. On suppose $N>0$. Posons $k=k(\lambda,\epsilon)$, $\bar{k}=k(\bar{\lambda},\bar{\epsilon})$, $M=M(\lambda,\epsilon)$, $\bar{M}=M(\bar{\lambda},\bar{\epsilon})$, cf. 5.1. Rappelons que $k=sup(2M,-2M-1)$, $\bar{k}=sup(2\bar{M},-2\bar{M}-1)$. Montrons d'abord que 

(1) $2M=\vert J^{a}\vert -\vert J^b\vert -1$.  

  D'apr\`es  5.1(2), on a
 $$2M=(\sum_{j=1,...,2r+1}\epsilon(j)(-1)^{j+1})-\sum_{j=1,...,2r+1}(-1)^{j+1}.$$
 La premi\`ere somme vaut $\vert J^{a}\vert -\vert J^b\vert $, la seconde vaut $1$ et l'assertion (1) s'en d\'eduit. 

Supposons d'abord $N'=0$. Alors $\bar{\lambda}_{1}=2N$. Puisque 
$$\bar{\lambda}_{1}\leq \sum_{h=1,...,S}\lambda_{s_{h}}\leq \sum_{j=1,...,2r+1}\lambda_{j},$$
 ces in\'egalit\'es sont des \'egalit\'es. En comparant avec la d\'efinition de $\bar{\lambda}_{1}$, la premi\`ere \'egalit\'e entra\^{\i}ne que  $\vert J^b\vert =(S-1)/2$ si $\epsilon(1)=1$ et $\vert J^{a}\vert =S/2$ si $\epsilon(1)=-1$.  La deuxi\`eme \'egalit\'e entra\^{\i}ne que
  $\lambda_{j}=0$ pour tout \'el\'ement $j\in \{1,...,2r+1\}$ qui n'appartient pas \`a la suite $\mathfrak{S}$.    Si un tel \'el\'ement existe, consid\'erons le plus grand de ces \'el\'ements, notons-le $j$. On  a $j\geq2$ puisque $j$ n'appartient pas \`a $\mathfrak{S}$. Supposons par exemple $j\in J^{a}$. Puisque $j\not\in \mathfrak{S}$, on a aussi $j-1\in J^{a}$ donc $\lambda_{j-1}\geq \lambda_{j}+2$ d'apr\`es 6.3(2). Les $\lambda_{j'}$ pour $j'<j$ sont donc non nuls et appartiennent donc forc\'ement \`a $\mathfrak{S}$. Cela prouve qu'il y a au plus un \'el\'ement de $\{1,...,2r+1\}$ qui n'appartient pas \`a $\mathfrak{S}$. D'o\`u $S=2r+1$ ou $S=2r$. On  a d\'etermin\'e la parit\'e de $S$ en 6.3(1) et on en d\'eduit $S=2r+1$ si $\epsilon(1)=1$,  $S=2r$ si $\epsilon(1)=-1$. On a vu ci-dessus que $\vert J^{b}\vert=(S-1)/2 $ dans le premier cas et que $\vert J^{a}\vert =S/2$ dans le second. Puisque $\vert J^{a}\vert +\vert J^{b}\vert =2r+1$,  on en d\'eduit $\vert J^{a}\vert =(S+1)/2$ dans le premier cas, $\vert J^{b}\vert =S/2+1$ dans le second. Avec (1) et l'\'egalit\'e $k=sup(2M,-2M-1)$, on calcule $k=0$ si $\epsilon(1)=1$, $k=1$ si $\epsilon(-1)=-1$. D'autre part, on a $\bar{\lambda}=(\bar{\lambda}_{1})$ et $\bar{\epsilon}_{\bar{\lambda}_{1}}=\epsilon(1)$. Avec la d\'efinition de 5.1, on calcule $\bar{k}=0$ si $\epsilon(1)=1$ et $\bar{k}=1$ si $\epsilon(1)=-1$. L'\'egalit\'e de l'\'enonc\'e s'ensuit.

 Supposons maintenant $N'>0$. On pose $k'=k(\lambda',\epsilon')$, $M'=M(\lambda',\epsilon')$.  Montrons que
  
  $$(2) \qquad M'=\left\lbrace\begin{array}{cc}-M,&\text{ si }\epsilon(1)=1,\\ -M-1,&\text{ si }\epsilon(1)=-1.\\ \end{array}\right.$$
  
  D'apr\`es (1), on a $2M=\vert J^{a}\vert -\vert J^b\vert -1$ et  $2M'=\vert {J'}^{a}\vert -\vert {J'}^b\vert -1$. Supposons par exemple $\epsilon(1)=1$. Alors, d'apr\`es les relations (10)(a) et (11) de 6.3, 
  $$2M'=\vert \tilde{J'}^{a}\vert -\vert \tilde{J'}^b\vert =\vert \tilde{J}^{b}\vert -\vert \tilde{J}^a\vert .$$ D'apr\`es 6.3(1), $\vert \tilde{J}^{a}\vert =\vert J^{a}\vert -(S+1)/2$, $\vert \tilde{J}^{b}\vert =\vert J^b\vert -(S-1)/2$. D'o\`u $2M'=\vert J^b\vert -\vert J^{a}\vert +1=-2M$. La preuve est similaire dans le cas $\epsilon(1)=-1$. Cela prouve (2). 
  
  Posons  $\bar{M}'=M(\bar{\lambda}',\bar{\epsilon}')$. L'\'egalit\'e du lemme \'equivaut \`a $M=\bar{M}$. Par r\'ecurrence, on peut supposer que $M'= \bar{M}'$. En vertu de (2), pour d\'emontrer le lemme, il suffit de prouver
  
   $$(3) \qquad \bar{M}=\left\lbrace\begin{array}{cc}-\bar{M}',&\text{ si }\epsilon(1)=1,\\ -\bar{M}'-1,&\text{ si }\epsilon(1)=-1.\\ \end{array}\right.$$  
   
   Notons $i_{1}>...,>i_{m}$ les \'el\'ements $i\in Jord^{bp}(\bar{\lambda})$ tels que $mult_{\bar{\lambda}}(i)$ soit impaire et $i'_{1}>....>i'_{m'}$  la suite analogue pour $\bar{\lambda}'$.  Puisque $\bar{\lambda}=\{\bar{\lambda}_{1}\}\sqcup \bar{\lambda}'$, on a deux possibilit\'es:
   
 (4)  si  $\bar{\lambda}_{1}>\bar{\lambda}'_{1}$, ou si $\bar{\lambda}_{1}=\bar{\lambda}'_{1}$ et $mult_{\bar{\lambda}'}(\bar{\lambda}'_{1})$ est paire, alors $i_{1}=\bar{\lambda}_{1}$, $m=m'+1$ et $i_{l}=i'_{l-1}$ pour $l=2,...,m$;
 
  (5) si $\bar{\lambda}_{1}=\bar{\lambda}'_{1}$ et $mult_{\bar{\lambda}'}(\bar{\lambda}'_{1})$ est impaire,  on a $i'_{1}=\bar{\lambda}'_{1}$, $m=m'-1$ et $i_{l}=i'_{l+1}$ pour $l=1,...,m$.  
 
 Supposons $\epsilon(1)=1$. Dans le cas (4), $i_{1}$ ne contribue pas \`a $\bar{M}$ et 
 $$\bar{M}=\sum_{l=2,...,m; \bar{\epsilon}_{\bar{\lambda}_{l}}=-1}(-1)^l=\sum_{l=1,...,m';\bar{\epsilon}'_{\bar{\lambda}'_{l}}=-1}(-1)^{l+1}=-\bar{M}'.$$
 Dans le  cas (5), on a forc\'ement $\epsilon'(1)=1$ et $i'_{1}$ ne contribue pas \`a $\bar{M}'$. Alors
 $$\bar{M}=\sum_{l=1,...,m;\bar{\epsilon}_{\bar{\lambda}_{l}}=-1}(-1)^l=\sum_{l=2,...,m';\bar{\epsilon}'_{\bar{\lambda}'_{l}}=-1}(-1)^{l+1}=-\bar{M}'.$$
  
 Supposons $\epsilon(1)=-1$. Dans le cas (4), $i_{1}$  contribue  \`a $\bar{M}$ par $-1$ et 
 $$\bar{M}=-1+\sum_{l=2,...,m; \bar{\epsilon}_{\bar{\lambda}_{l}}=-1}(-1)^l=-1+\sum_{l=1,...,m';\bar{\epsilon}'_{{\lambda}'_{l}}=-1}(-1)^{l+1}=-1-\bar{M}'.$$
 Dans le cas (5) , on a forc\'ement $\epsilon'(1)=-1$ et $i'_{1}$  contribue  \`a $\bar{M}'$. Alors
 $$\bar{M}=\sum_{l=1,...,m;\bar{\epsilon}_{\bar{\lambda}_{l}}=-1}(-1)^l=\sum_{l=2,...,m';\bar{\epsilon}'_{\bar{\lambda}'_{l}}=-1}(-1)^{l+1}=-(\bar{M}'+1).$$ 
 Cela prouve (3) et ach\`eve la preuve du lemme. $\square$ 
 
\bigskip

\subsection{Preuve de la proposition 6.2}
Les calculs de cette preuve diff\`erent selon les diff\'erents cas $\epsilon(1)=\pm 1$ et $k$ pair ou impair. Pour tenter de les unifier, reprenons les d\'efinitions de 5.1 et 5.2. On fixe un entier  $N>0$,  un entier $k\geq0$ tel que $k(k+1)\leq 2N$ et un entier $r$ assez grand.  On a d\'efini les entiers $n=r+[k/2]+1$, $m=r-[k/2]$, $A=2r+k$, $B=2r-k-1$. Pour un \'el\'ement $(\lambda,\epsilon)\in \boldsymbol{{\cal P}^{symp}}(2N)_{k}$ (c'est-\`a-dire que $k(\lambda,\epsilon)=k$), on suppose que $\lambda$ a $2r+1$ termes avec $\lambda_{2r+1}=0$. On a d\'efini les partitions $A_{\lambda,\epsilon}$ et $B_{\lambda,\epsilon}$. Pour un couple $(\alpha,\beta)\in {\cal P}_{2}(N-k(k+1)/2)$, on suppose que $\alpha$ a $n$ termes et que $\beta$ en a $m$. On a d\'efini les partitions $A_{\alpha,\beta}$ et $B_{\alpha,\beta}$. Fixons un signe $u\in \{\pm 1\}$. On modifie les d\'efinitions ci-dessus en posant

 si $u=1$, $\tilde{A}_{\lambda,\epsilon}=A_{\lambda,\epsilon}$, $\tilde{B}_{\lambda,\epsilon}=B_{\lambda,\epsilon}$;

si $u=-1$, $\tilde{A}_{\lambda,\epsilon}=B_{\lambda,\epsilon}$, $\tilde{B}_{\lambda,\epsilon}=A_{\lambda,\epsilon}$;

si $u(-1)^k=1$, $\tilde{n}=n$, $\tilde{m}=m$, $\tilde{A}=A$, $\tilde{B}=B$, $\tilde{\alpha}=\alpha$, $\tilde{\beta}=\beta$, $\tilde{A}_{\alpha,\beta}=A_{\alpha,\beta}$, $\tilde{B}_{\alpha,\beta}=B_{\alpha,\beta}$;

si $u(-1)^k=-1$, $\tilde{n}=m$, $\tilde{m}=n$, $\tilde{A}=B$, $\tilde{B}=A$, $\tilde{\alpha}=\beta$, $\tilde{\beta}=\alpha$, $\tilde{A}_{\alpha,\beta}=B_{\alpha,\beta}$, $\tilde{B}_{\alpha,\beta}=A_{\alpha,\beta}$.

On a 

$$(1) \qquad \tilde{A}_{\alpha,\beta}=\tilde{\alpha}+[\tilde{A},\tilde{A}+2-2\tilde{n}]_{2},\,\,\tilde{B}_{\alpha,,\beta}=\tilde{\beta}+[\tilde{B},\tilde{B}+2-2\tilde{m}]_{2}.$$

 Dire que $(\alpha,\beta)$ correspond \`a $(\lambda,\epsilon)$ \'equivaut aux \'egalit\'es

$$(2)\qquad \tilde{A}_{\lambda,\epsilon}=\tilde{A}_{\alpha,\beta},\,\, \tilde{B}_{\lambda,\epsilon}=\tilde{B}_{\alpha,\beta}.$$

Consid\'erons les ensembles d'indices $I=(\{1,...,n\}\times\{0\})\cup(\{1,...,m\}\times\{1\})$ et $\tilde{I}=(\{1,...,\tilde{n}\}\times\{0\})\cup(\{1,...,\tilde{m}\}\times\{1\})$. Donnons-nous un ordre $<_{I}$ sur $I$. On en d\'eduit un ordre $<_{\tilde{I}}$ sur $\tilde{I}$: si $u(-1)^k=1$, on a $\tilde{I}=I$ et l'ordre  $<_{\tilde{I}}$ co\"{\i}ncide avec $<_{I}$; si $u(-1)^k=-1$, on impose que $(i,0)<_{\tilde{I}}(j,1)$ si et seulement si $(i,1)<_{I}(j,0)$. On voit alors que l'application $(\nu,\mu)\mapsto (\tilde{\nu},\tilde{\mu})$ est une bijection de $P_{A,B;2}(\alpha,\beta)$ sur $P_{\tilde{A},\tilde{B};2}(\tilde{\alpha},\tilde{\beta})$. 
 
 Fixons maintenant $(\lambda,\epsilon)\in \boldsymbol{{\cal P}^{symp}}(2N)_{k}$ et supposons que tous les termes de  $\lambda$ sont pairs. Soit $(\alpha,\beta)\in {\cal P}_{2}(N-k(k+1)/2)$ le couple correspondant \`a $(\lambda,\epsilon)$. On munit $I$ de l'ordre $<_{I,\alpha,\beta,A,B}$ d\'efini par $(i,0)<_{I}(j,1)$ si et seulement si $\alpha_{i}+A+2-2i>\beta_{j}+B+2-2j$. On note simplement $<_{I}$ cet ordre. On voit que l'ordre $<_{\tilde{I}}$ qui s'en d\'eduit n'est autre que $<_{\tilde{I},\tilde{\alpha},\tilde{\beta},\tilde{A},\tilde{B}}$. D'apr\`es la preuve du th\'eor\`eme 5.5, il y a un unique \'el\'ement $(\alpha^{max},\beta^{max})$ dans $P_{A,B;2}(\alpha,\beta)$ et $(\lambda^{max},\epsilon^{max})$ est l'\'el\'ement de $\boldsymbol{{\cal P}^{symp}}(2N)_{k}$ qui correspond \`a $(\alpha^{max},\beta^{max})$. Alors $(\tilde{\alpha}^{max},\tilde{\beta}^{max})$ est l'unique \'el\'ement de $P_{\tilde{A},\tilde{B};2}(\tilde{\alpha},\tilde{\beta})$ et on a les \'egalit\'es
 $$(3) \qquad \tilde{A}_{\lambda^{max},\epsilon^{max}}=\tilde{A}_{\alpha^{max},\beta^{max}},\,\, \tilde{B}_{\lambda^{max},\epsilon^{max}}=\tilde{B}_{\alpha^{max},\beta^{max}}.$$
 
 Les constructions ci-dessus d\'ependent du signe $u$. On choisit d\'esormais $u=\epsilon(1)$. On a calcul\'e les termes $A_{\lambda,\epsilon}$ et $B_{\lambda,\epsilon}$ en 5.1.  On en d\'eduit que  $\tilde{A}_{\lambda,\epsilon}=\{\lambda_{j}/2+2r+1-j; j=1,...,2r+1, \epsilon(j)(-1)^{j+1}=\epsilon(1)\}$ et $\tilde{B}_{\lambda,\epsilon}=\{\lambda_{j}/2+2r+1-j; j=1,...,2r+1, \epsilon(j)(-1)^j=\epsilon(1)\}$. En particulier, le plus grand terme de $\tilde{A}_{\lambda,\epsilon}\sqcup \tilde{B}_{\lambda,\epsilon}$ appara\^{\i}t dans $\tilde{A}_{\lambda,\epsilon}$. Gr\^ace \`a 6.3(1), on obtient que $\tilde{A}_{\lambda,\epsilon}$ est la r\'eunion sur $h=1,...,[(S+1)/2]$ des $ \{\lambda_{j}/2+2r+1-j; j=s_{2h-1},...,s_{2h}-1\}$ et que $\tilde{B}_{\lambda,\epsilon}$ est la r\'eunion sur $h=1,...,[S/2]$ des $\{\lambda_{j}/2+2r+1-j; j=s_{2h},...,s_{2h+1}-1\}$. 
Divisons l'ensemble $\{1,...,\tilde{n}\}$ en $[(S+1)/2]$ intervalles cons\'ecutifs $X_{1},...,X_{[(S+1)/2]}$, l'intervalle $X_{h}$ ayant $s_{2h}-s_{2h-1}$ termes. Notons $x_{h}$ le plus petit \'el\'ement de $X_{h}$.   Divisons l'ensemble $\{1,...,\tilde{m}\}$ en $[S/2]$ intervalles cons\'ecutifs $Y_{1},...,Y_{[S/2]}$, l'intervalle $Y_{h}$ ayant $s_{2h+1}-s_{2h}$ termes.  Notons $y_{h} $ le plus petit \'el\'ement de $Y_{h}$.  Pour $j\in X_{h}$, le $j$-i\`eme terme de $\tilde{A}_{\lambda,\epsilon}$ est $\lambda_{j+s_{2h-1}-x_{h}}/2+2r+1-j-s_{2h-1}+x_{h}$. D'apr\`es  (1) et (2), c'est aussi $\tilde{\alpha}_{j}+\tilde{A}+2-2j$.  De m\^eme,  pour $j\in Y_{h}$, le $j$-i\`eme terme de $\tilde{B}_{\lambda,\epsilon}$ est $\lambda_{j+s_{2h}-y_{h}}/2+2r+1-j-s_{2h}+y_{h}$.  D'apr\`es (1) et (2), c'est aussi $\tilde{\beta}_{j}+\tilde{B}+2-2j$. L'ordre $<_{\tilde{I}}$ est tel que, pour $i\in X_{h}$ et $j\in Y_{k}$, on a $(i,0)<_{\tilde{I}}(j,1)$ si et seulement si $h\leq k$. En particulier, $(1,0)<_{\tilde{I}}(1,1)$ et l'\'el\'ement $(\tilde{\alpha}^{max},\tilde{\beta}^{max})$ est construit \`a l'aide du proc\'ed\'e (a). Les suites d'indices associ\'ees \`a ce proc\'ed\'e sont $x_{1}=1,...,x_{[(S+1)/2]}$ et $y_{1}=1,...,y_{[S/2]}$ et on obtient
$$(4) \qquad \tilde{\alpha}^{max}_{1}=(\sum_{h=1,...,[(S+1)/2]}\tilde{\alpha}_{x_{h}})+(\sum_{h=1,...,[S/2]}\tilde{\beta}_{y_{h}}).$$

On sait d'apr\`es 1.3(1) que $\tilde{\alpha}_{1}^{max}+\tilde{A}$ est le plus grand \'el\'ement de $p\Lambda^{\tilde{n},\tilde{m}}_{\tilde{A},\tilde{B};2}(\tilde{\alpha}^{max},\tilde{\beta}^{max})$. En fait, cet \'el\'ement intervient avec multiplicit\'e $1$ dans cette partition. Pour montrer cela, il suffit de prouver que $\tilde{\alpha}_{1}^{max}+\tilde{A}> \tilde{\beta}_{1}^{max}+\tilde{B}$. Or, d'apr\`es 1.2(3), $\tilde{\beta}_{1}^{max}$ est major\'e par le terme $\mu_{1}$ issu du proc\'ed\'e (b) appliqu\'e \`a $(\tilde{\alpha},\tilde{\beta})$. D'apr\`es 1.2(4), on a $\mu_{1}=\tilde{\alpha}_{1}^{max}-\tilde{\alpha}_{1}$. Donc 
$$\tilde{\beta}_{1}^{max}+\tilde{B}\leq \tilde{\alpha}_{1}^{max}-\tilde{\alpha}_{1}+\tilde{B}\leq (\tilde{\alpha}_{1}^{max}+\tilde{A})-(\tilde{\alpha}_{1}+\tilde{A})+(\tilde{\beta}_{1}+\tilde{B})<\tilde{\alpha}_{1}^{max}+\tilde{A}$$
puisque $\tilde{\alpha}_{1}+\tilde{A}>\tilde{\beta}_{1}+\tilde{B}$. Cela d\'emontre l'assertion.  D'apr\`es (2), $\tilde{\alpha}_{1}^{max}+\tilde{A}$ est aussi l'unique \'el\'ement maximal de $A_{\lambda^{max},\epsilon^{max}}\sqcup B_{\lambda^{max},\epsilon^{max}}$. On ne sait pas encore que $\lambda^{max}$ a tous ses termes pairs. Mais \`a l'aide de 5.1(a) et 5.1(b), on v\'erifie les propri\'et\'es suivantes, valables pour tout $(\mu,\tau)\in \boldsymbol{{\cal P}^{symp}}(2N)$:

si $\mu_{1}$ est pair, $A_{\mu,\tau}\sqcup B_{\mu,\tau}$ a un unique \'el\'ement maximal, qui est $\mu_{1}/2+2r$; il intervient dans $A_{\mu,\tau}$ si $\tau_{\mu_{1}}=1$ et dans $B_{\mu,\tau}$ si $\tau_{\mu_{1}}=-1$;

si $\mu_{1}$ est impair, l'\'el\'ement maximal de $A_{\mu,\tau}\sqcup B_{\mu,\tau}$ intervient avec multiplicit\'e $2$.

Ceci entra\^{\i}ne que $\lambda^{max}_{1}$ est pair et \'egal \`a $2(\tilde{\alpha}_{1}^{max}+\tilde{A}-2r)$. Puisqu'il intervient dans $\tilde{A}_{\alpha^{max},\beta^{max}}=\tilde{A}_{\lambda^{max},\epsilon^{max}}$, on voit en se reportant \`a la d\'efinition de cet ensemble que $\epsilon^{max}_{\lambda^{max}_{1}}=\epsilon(1)$. Calculons $\lambda_{1}^{max}$ en utilisant l'\'egalit\'e (4). 
Avec les formules ci-dessus, on a $\tilde{\alpha}_{x_{h}}=\lambda_{s_{2h-1}}/2+2r+1-s_{2h-1}-\tilde{A}-2+2x_{h}$, $\tilde{\beta}_{y_{h}}=\lambda_{s_{2h}}/2+2r+1-s_{2h}-\tilde{B}-2+2y_{h}$. D'o\`u
$$\tilde{\alpha}^{max}_{1}=(\sum_{h=1,...,S}\lambda_{s_{h}}/2)-(\sum_{h=1,...,S}s_{h})+2(\sum_{h=1,...,[(S+1)/2]}x_{h})+2(\sum_{h=1,...,[S/2]}y_{h})$$
$$+S(2r+1)-[(S+1)/2](\tilde{A}+2)-[S/2](\tilde{B}+2).$$
Puisque $X_{h}$ a $s_{2h}-s_{2h-1}$ termes, on a $x_{h+1}-x_{h}=s_{2h}-s_{2h-1}$, en posant par convention $x_{[(S+1)/2]+1}=\tilde{n}+1$. De m\^eme, $y_{h+1}-y_{h}=s_{2h+1}-s_{2h}$, en posant par convention $y_{[S/2]+1}=\tilde{m}+1$.  On en d\'eduit ais\'ement $s_{2h-1}=x_{h}+y_{h}-1$ et $s_{2h}=x_{h+1}+y_{h}-1$. Puis
$$\sum_{h=1,...,S}s_{h}=2(\sum_{h=1,...,[(S+1)/2]}x_{h})+2(\sum_{h=1,...,[S/2]}y_{h})+\left\lbrace\begin{array}{cc}\tilde{n}-S,&\text{ si }S{\rm\,\, est\,\, pair;}\\ \tilde{m}-S,&\text{ si }S{\rm \,\,est\,\, impair.}\\ \end{array}\right.$$
D'o\`u
$$\tilde{\alpha}^{max}_{1}=(\sum_{h=1,...,S}\lambda_{s_{h}}/2)+S(2r+1)-[(S+1)/2](\tilde{A}+2)-[S/2](\tilde{B}+2)+\left\lbrace\begin{array}{cc}S-\tilde{n},&\text{ si }S{\rm\,\, est\,\, pair;}\\S- \tilde{m},&\text{ si }S{\rm \,\,est\,\, impair.}\\ \end{array}\right.$$
D'apr\`es 6.3(1), $S$ est impair, resp. pair, si $\epsilon(1)=1$, resp. $\epsilon(1)=-1$.  En distinguant les diff\'erents cas $\epsilon(1)=\pm 1$, $k$ est pair ou impair, on d\'eduit par calcul de la formule ci-dessus l'\'egalit\'e
$$2(\tilde{\alpha}^{max}_{1}+\tilde{A}-2r)=(\sum_{h=1,...,S}\lambda_{s_{h}})+\left\lbrace\begin{array}{cc}S-1-2\tilde{m},&\text{ si } \epsilon(1)=1\\S- 2\tilde{m},&\text{ si } \epsilon(1)=-1.\\ \end{array}\right.$$
L'entier $\tilde{m}$ est le nombre d'\'el\'ements de $\tilde{B}_{\lambda,\epsilon}$, c'est-\`a-dire celui de $B_{\lambda,\epsilon}$ si $\epsilon(1)=1$, de $A_{\lambda,\epsilon}$ si $\epsilon(1)=-1$. Ces derniers sont $\vert J^{b}\vert $, resp. $\vert J^{a}\vert $. Le terme de droite de l'\'egalit\'e ci-dessus est le nombre not\'e $\bar{\lambda}_{1}$ en 6.1.  D'o\`u $\lambda_{1}^{max}=\bar{\lambda}_{1}$.   On a aussi montr\'e plus haut que $\epsilon^{max}_{\lambda^{max}_{1}}=\epsilon(1)=\bar{\epsilon}_{\bar{\lambda}_{1}}$.

Si $\lambda^{max}_{1}=2N$, $\lambda^{max}$ est enti\`erement d\'etermin\'e et co\"{\i}ncide avec l'\'el\'ement  $\bar{\lambda}$  de 6.1. Supposons d\'esormais $\lambda^{max}_{1}<2N$. Le proc\'ed\'e (a) que l'on a appliqu\'e \`a $(\tilde{\alpha},\tilde{\beta})$ a cr\'e\'e le terme $\tilde{\alpha}_{1}^{max}$. Il cr\'ee aussi un couple de partitions, notons-le $(\boldsymbol{\alpha},\boldsymbol{\beta})\in {\cal P}_{{\bf n}}\times {\cal P}_{{\bf m}}$. D'apr\`es la preuve du th\'eor\`eme 5.5, l'ensemble $P_{\tilde{A}-2,\tilde{B};2}(\boldsymbol{\alpha},\boldsymbol{\beta})$ a un unique \'el\'ement, notons-le $(\nu,\mu)$. On a
$$\tilde{\alpha}^{max}=\{\tilde{\alpha}_{1}^{max}\}\sqcup \nu\sqcup\{0^{\tilde{n}-{\bf n}-1}\},\,\,\tilde{\beta}^{max}=\mu\sqcup\{0^{\tilde{m}-{\bf m}}\}.$$
D'o\`u
$$(5) \qquad \tilde{A}_{\lambda^{max},\epsilon^{max}}=\{\lambda_{1}^{max}/2+2r\}\sqcup((\nu\sqcup\{0^{\tilde{n}-{\bf n}-1}\})+[\tilde{A}-2,\tilde{A}+2-2\tilde{n}]_{2}),$$
$$\tilde{B}_{\lambda^{max},\epsilon^{max}}=(\mu\sqcup\{0^{\tilde{m}-{\bf m}}\})+[\tilde{B},\tilde{B}+2-2\tilde{m}]_{2}.$$

En 6.1, on a introduit l'entier $N'=N-\bar{\lambda}_{1}/2=N-\lambda_{1}^{max}/2$ et un \'el\'ement $(\lambda',\epsilon')\in \boldsymbol{{\cal P}^{symp}}(2N')$. On peut repr\'esenter $\lambda'$ comme une partition \`a $2r+1-S$ termes mais, pour faire bonne mesure, on adjoint des termes nuls et on le repr\'esente comme une partition \`a $2r'+1$ termes, avec $2r'+1> 2r+1-S$. On pose $k'=k(\lambda',\epsilon')$. Avec l'\'egalit\'e 6.4(2) et les relations $k=sup(2M,-2M-1)$, $k'=sup(2M',-2M'-1)$, on calcule facilement
$$(6) \qquad k'=\left\lbrace\begin{array}{cc}k-1,&\text{si } \epsilon(1)(-1)^k=1\text{ et }k\not=0,\\ k+1,&\text{ si }\epsilon(1)(-1)^k=-1,\\ 0,&\text{ si }\epsilon(1)=1\text{ et }k=0.\\ \end{array}\right.$$
On introduit les entiers $n'=r'+[k'/2]+1$, $m'=r'-[k'/2]$, $A'=2r'+k'$, $B'=2r'-k'-1$ et le couple de partitions $(\alpha',\beta')\in {\cal P}_{2}(N'-k'(k'+1)/2)$ correspondant \`a $ (\lambda',\epsilon')$. Fixons  un signe $u'\in \{\pm 1\}$. En utilisant ce signe et l'entier $k'$, on construit comme au d\'ebut du paragraphe des entiers $\tilde{n}'$, $\tilde{m}'$, $\tilde{A}'$, $\tilde{B}'$ et des partitions $\tilde{A}_{\lambda',\epsilon'}$, $\tilde{B}_{\lambda',\epsilon'}$, $\tilde{\alpha}'$, $\tilde{\beta}'$. On choisit $u'=-u=-\epsilon(1)$. 

A l'aide de (6) et d'un calcul cas par cas, on voit que

(7) $\tilde{A}'=\tilde{A}+2r'-2r-1,\,\,\tilde{B}'=\tilde{B}+2r'-2r+1$.

Montrons que

(8) $\tilde{\alpha}'=\boldsymbol{\alpha}\sqcup\{0^{\tilde{n}'-{\bf n}}\}$, $\tilde{\beta}'=\boldsymbol{\beta}\sqcup\{0^{\tilde{m}'-{\bf m}}\}$,

\noindent ce qui, d'apr\`es l'analogue de (2), \'equivaut \`a

(9) $\tilde{A}_{\lambda',\epsilon'}=(\boldsymbol{\alpha}\sqcup\{0^{\tilde{n}'-{\bf n}}\})+[\tilde{A}',\tilde{A}'+2-2\tilde{n}']_{2}$, $\tilde{B}_{\lambda',\epsilon'}=(\boldsymbol{\beta}\sqcup\{0^{\tilde{m}'-{\bf m}}\})+[\tilde{B}',\tilde{B}'+2-2\tilde{m}']_{2}$.

On a
$$(10)\qquad \tilde{A}_{\lambda',\epsilon'}=\{\lambda'_{j}/2+2r'+1-j; j=1,...,2r'+1, \epsilon'(j)(-1)^{j}=\epsilon(1)\}$$
$$=\{\lambda'_{j}/2+2r'+1-j; j=1,...,2r+1-S, \epsilon'(j)(-1)^{j}=\epsilon(1)\}\cup\{j\in [2r'-2r+S-1,0]_{1}; (-1)^j=\epsilon(1)\};$$
$$\tilde{B}_{\lambda',\epsilon'}=\{\lambda'_{j}/2+2r'+1-j; j=1,...,2r'+1, \epsilon'(j)(-1)^{j}=-\epsilon(1)\}$$
$$=\{\lambda'_{j}/2+2r'+1-j; j=1,...,2r+1-S, \epsilon'(j)(-1)^{j}=-\epsilon(1)\}\cup\{j\in [2r'-2r+S-1,0]_{1}; (-1)^j=-\epsilon(1)\}.$$
Un calcul cas par cas montre que les contributions des termes compl\'ementaires correspondant aux $0$ que l'on ajoute \`a nos partitions co\"{\i}ncident dans les deux membres des \'egalit\'es (9). C'est-\`a-dire que
$$(11)\qquad \{j\in [2r'-2r+S-1,0]_{1}; (-1)^j=\epsilon(1)\}=[\tilde{A}'-{\bf n},\tilde{A}'+2-2\tilde{n}']_{2},$$
$$\{j\in [2r'-2r+S-1,0]_{1}; (-1)^j=-\epsilon(1)\}=[\tilde{B}'-{\bf m},\tilde{B}'+2-\tilde{m}']_{2}.$$
Pour $j\in \{1,...,2r+1-S\}$, on a $j+h(j)\in \{s_{h(j)}+1,...,s_{h(j)+1}-1\}$ et $\epsilon'(j)=\epsilon(j+h(j))(-1)^{h(j)+1}$ par d\'efinition. La condition $\epsilon'(j)(-1)^j= \epsilon(1)$ \'equivaut \`a $\epsilon(j+h(j))(-1)^{j+h(j)+1}=\epsilon(1)$, dont on sait d'apr\`es 6.3(1) qu'elle \'equivaut \`a $h(j)$ impair. En reprenant la d\'efinition de $\lambda'_{j}$, on voit qu'alors $\lambda'_{j}=\lambda_{j+h(j)}$. On en d\'eduit
$$\{\lambda'_{j}/2+2r'+1-j; j=1,...,2r+1-S, \epsilon'(j)(-1)^{j}=\epsilon(1)\}=$$
$$\cup_{h=1,...,[(S+1)/2]}\{\lambda_{j}/2+2r'-j+2h; j\in \{s_{2h-1}+1,...,s_{2h}-1\}\}.$$
 La condition $\epsilon'(j)(-1)^j= -\epsilon(1)$ \'equivaut \`a $\epsilon(j+h(j))(-1)^{j+h(j)+1}=-\epsilon(1)$, dont on sait d'apr\`es 6.3(1) qu'elle \'equivaut \`a $h(j)$ pair. En reprenant la d\'efinition de $\lambda'_{j}$, on voit qu'alors $\lambda'_{j}=\lambda_{j+h(j)}+2$. On en d\'eduit
$$\{\lambda'_{j}/2+2r'+1-j; j=1,...,2r+1-S, \epsilon'(j)(-1)^{j}=-\epsilon(1)\}=$$
$$\cup_{h=1,...,[S/2]}\{\lambda_{j}/2+2r'-j+2h+2; j\in \{s_{2h}+1,...,s_{2h+1}-1\}\}.$$
 On a donn\'e plus haut les formules exprimant les termes $\tilde{\alpha}_{j}+\tilde{A}+2-2j$ et $\tilde{\beta}_{j}+\tilde{B}+2-2j$ en fonction des $\lambda_{j}/2+2r+1-j$. On traduit alors:
$$\{\lambda'_{j}/2+2r'+1-j; j=1,...,2r+1-S, \epsilon'(j)(-1)^{j}=\epsilon(1)\}=$$
$$\cup_{h=1,...,[(S+1)/2]}\{\tilde{\alpha}_{j}+\tilde{A}-2j+2r'-2r+2h+1; j\in \{x_{h}+1,...,x_{h+1}-1\}\},$$
$$\{\lambda'_{j}/2+2r'+1-j; j=1,...,2r+1-S, \epsilon'(j)(-1)^{j}=-\epsilon(1)\}=$$
$$\cup_{h=1,...,[S/2]}\{\tilde{\beta}_{j}+\tilde{B}-2j+2r'-2r+2h+3; j\in \{y_{h}+1,...,y_{h+1}-1\}\}.$$
Par d\'efinition de $\boldsymbol{\alpha}$ et $\boldsymbol{\beta}$, on a $\tilde{\alpha}_{j}=\boldsymbol{\alpha}_{j-h}$ pour $j\in \{x_{h}+1,...,x_{h+1}-1\}$ et $\tilde{\beta}_{j}=\boldsymbol{\beta}_{j-h}$ pour $j\in \{y_{h}+1,...,y_{h+1}-1\}$. On obtient 
$$\{\lambda'_{j}/2+2r'+1-j; j=1,...,2r+1-S, \epsilon'(j)(-1)^{j}=\epsilon(1)\}=\boldsymbol{\alpha}+[\tilde{A}+2r'-2r-1,\tilde{A}+2r'-2r+1-2{\bf n}]_{2},$$
$$\{\lambda'_{j}/2+2r'+1-j; j=1,...,2r+1-S, \epsilon'(j)(-1)^{j}=-\epsilon(1)\}=\boldsymbol{\beta}+[\tilde{B}+2r'-2r+1,\tilde{B}+2r'-2r+3-2{\bf m}]_{2}.$$
D'o\`u, en utilisant (7), 
$$(12)\qquad\{\lambda'_{j}/2+2r'+1-j; j=1,...,2r+1-S, \epsilon'(j)(-1)^{j}=\epsilon(1)\}=\boldsymbol{\alpha}+[\tilde{A}',\tilde{A}'+2-2{\bf n}]_{2},$$
$$\{\lambda'_{j}/2+2r'+1-j; j=1,...,2r+1-S, \epsilon'(j)(-1)^{j}=-\epsilon(1)\}=\boldsymbol{\beta}+[\tilde{B}',\tilde{B}'+2-2{\bf m}]_{2}.$$
En rassemblant (10), (11) et (12), on obtient (9), d'o\`u (8). 

 Introduisons l'\'el\'ement $({\lambda'}^{max},{\epsilon'}^{max})\in \boldsymbol{{\cal P}^{symp}}(2N')_{k'}$. 
D'apr\`es la preuve du th\'eor\`eme 5.5, il y a unique \'el\'ement $({\alpha'}^{max},{\beta'}^{max})$ dans $P_{A',B';2}(\alpha',\beta')$ et c'est le couple correspondant \`a $({\lambda'}^{max},{\epsilon'}^{max})$. On d\'efinit comme ci-dessus les partitions ${\tilde{\alpha'}}^{max}$, ${\tilde{\beta'}}^{max}$, $\tilde{A}_{{\alpha'}^{max},{\beta'}^{max}}$, $\tilde{B}_{{\alpha'}^{max},{\beta'}^{max}}$, $\tilde{A}_{{\lambda'}^{max},{\epsilon'}^{max}}$, $\tilde{B}_{{\lambda'}^{max},{\epsilon'}^{max}}$, relatives \`a l'entier $k'$ et au signe $u'=-\epsilon(1)$. On a l'analogue de la relation (3), c'est-\`a-dire

(13) $\tilde{A}_{{\lambda'}^{max},{\epsilon'}^{max}}=\tilde{A}_{{\alpha'}^{max},{\beta'}^{max}}$,  $\tilde{B}_{{\lambda'}^{max},{\epsilon'}^{max}}=\tilde{B}_{{\alpha'}^{max},{\beta'}^{max}}$.

A ce point, on applique par r\'ecurrence la proposition 6.2 au couple $(\lambda',\epsilon')$: on a l'\'egalit\'e $({\lambda'}^{max},{\epsilon'}^{max})=(\bar{\lambda}',\bar{\epsilon}')$.  En se rappelant que $u'=-\epsilon(1)$, on calcule
$$\tilde{A}_{\bar{\lambda}',\bar{\epsilon}'}=\{\bar{\lambda}'_{j}/2+2r'+1-j; j=1,...,2r'+1,\bar{\epsilon}'(j)(-1)^j=\epsilon(1)\},$$
$$\tilde{B}_{\bar{\lambda}',\bar{\epsilon}'}=\{\bar{\lambda}'_{j}/2+2r'+1-j; j=1,...,2r'+1,\bar{\epsilon}'(j)(-1)^j=-\epsilon(1)\}.$$
Par construction, $\bar{\lambda}'_{j}=\bar{\lambda}_{j+1}$ et $\bar{\epsilon}'(j)=\bar{\epsilon}(j+1)$. A l'aide de ces formules, (13) se r\'ecrit
$$(14) \qquad \tilde{A}_{{\alpha'}^{max},{\beta'}^{max}}=\{\bar{\lambda}_{j}/2+2r'+2-j; j=2,...,2r'+2,\bar{\epsilon}(j)(-1)^{j+1}=\epsilon(1)\},$$
$$\tilde{B}_{{\alpha'}^{max},{\beta'}^{max}}=\{\bar{\lambda}_{j}/2+2r'+2-j; j=2,...,2r'+2,\bar{\epsilon}(j)(-1)^{j+1}=-\epsilon(1)\}.$$

Montrons que

(15) on a les \'egalit\'es ${\tilde{\alpha'}}^{max}=\nu\sqcup \{0^{\tilde{n}'-{\bf n}}\}$, ${\tilde{\beta'}}^{max}=\mu\sqcup\{0^{\tilde{m}'-{\bf m}}\}$.

Le couple $( {\tilde{\alpha'}}^{max},{\tilde{\beta'}}^{max})$ est l'unique \'el\'ement de $P_{\tilde{A}',\tilde{B}';2}(\tilde{\alpha}',\tilde{\beta}')$ tandis que $(\nu,\mu)$ est l'unique \'el\'ement de $P_{\tilde{A}-2,\tilde{B};2}(\boldsymbol{\alpha},\boldsymbol{\beta})$.  Pour tout $c\in {\mathbb N}$, on a l'\'egalit\'e $P_{\tilde{A}-2,\tilde{B};2}(\boldsymbol{\alpha},\boldsymbol{\beta})=P_{\tilde{A}-2+c,\tilde{B}+c;2}(\boldsymbol{\alpha},\boldsymbol{\beta})$: cela r\'esulte simplement de la construction de ces ensembles. En prenant $c=2r'-2r+1$ et en utilisant (7), on obtient que  $(\nu,\mu)$ est aussi l'unique \'el\'ement de $P_{\tilde{A}',\tilde{B}';2}(\boldsymbol{\alpha},\boldsymbol{\beta})$. En tenant compte de (8), l'assertion (15) r\'esulte du lemme 2.5, pourvu que les hypoth\`eses de ce lemme soient v\'erifi\'ees pour l'ensemble d'indice $\tilde{I}'=(\{1,...,\tilde{n}'\}\times \{0\})\cup(\{1,...,\tilde{m}'\}\times \{1\})$ et les entiers ${\bf n}\leq \tilde{n}'$ et ${\bf m}\leq \tilde{m}'$. Il s'agit de voir que
$$\tilde{\alpha}'_{{\bf n}}+\tilde{A}'+2-2{\bf n}> \tilde{B}'-2{\bf m}$$
et
$$\tilde{\beta}'_{{\bf m}}+\tilde{B}'+2-2{\bf m}>\tilde{A}'-2{\bf n}.$$
Or il r\'esulte de la preuve de (9) ci-dessus que les termes de gauche sont de la forme $\lambda'_{j}+2r'+1-j$ pour un $j\leq 2r+1-S$ tandis que les termes de droite sont de la forme $\lambda'_{j}+2r'+1-j$ pour un $j> 2r+1-S$. Les in\'egalit\'es cherch\'ees s'ensuivent, d'o\`u (15). 

Il r\'esulte de (7) et  (15) que $(\nu\sqcup\{0^{\tilde{n}-{\bf n}-1}\})+[\tilde{A}-2,\tilde{A}+2-2\tilde{n}]_{2}$ est l'ensemble des $\tilde{n}-1$ plus grands termes de $\tilde{A}_{{\alpha'}^{max},{\beta'}^{max}}$, auxquels on ajoute $2r-1-2r'$ et que $(\mu\sqcup\{0^{\tilde{m}-{\bf m}}\})+[\tilde{B},\tilde{B}+2-2\tilde{m}]_{2}$ est l'ensemble des $\tilde{m}$ plus grands termes de $\tilde{B}_{{\alpha'}^{max},{\beta'}^{max}}$, auxquels on ajoute $2r-1-2r'$. Gr\^ace \`a (14), on obtient que $(\nu\sqcup\{0^{\tilde{n}-{\bf n}-1}\})+[\tilde{A}-2,\tilde{A}+2-2\tilde{n}]_{2}$ est l'ensemble des $\tilde{n}-1$ plus grands termes de $\{\bar{\lambda}_{j}/2+2r+1-j; j=2,...,2r'+2,\bar{\epsilon}(j)(-1)^{j+1}=\epsilon(1)\}$ tandis que $(\mu\sqcup\{0^{\tilde{m}-{\bf m}}\})+[\tilde{B},\tilde{B}+2-2\tilde{m}]_{2}$ est l'ensemble des $\tilde{m}$ plus grands termes de $\{\bar{\lambda}_{j}/2+2r+1-j; j=2,...,2r'+2,\bar{\epsilon}(j)(-1)^{j+1}=-\epsilon(1)\}$. Remarquons que dans ce dernier ensemble, on peut aussi bien fair varier $j$ de $1$ \`a $2r'+2$: le terme $j=1$ ne v\'erifie pas la condition $\bar{\epsilon}(j)(-1)^{j+1}=-\epsilon(1)$ puisque $\bar{\epsilon}(1)=\epsilon(1)$ par d\'efinition. On utilise (5) en se rappelant que l'on a d\'ej\`a d\'emontr\'e l'\'egalit\'e $\lambda_{1}^{max}=\bar{\lambda}_{1}$. Alors $\tilde{A}_{\lambda^{max},\epsilon^{max}}$ est l'ensemble des $\tilde{n}$ plus grands termes de $\{\bar{\lambda}_{j}/2+2r+1-j; j=1,...,2r'+2,\bar{\epsilon}(j)(-1)^{j+1}=\epsilon(1)\}$ et $\tilde{B}_{\lambda^{max},\epsilon^{max}}$ est 
  l'ensemble des $\tilde{m}$ plus grands termes de $\{\bar{\lambda}_{j}/2+2r+1-j; j=1,...,2r'+2,\bar{\epsilon}(j)(-1)^{j+1}=-\epsilon(1)\}$.

Repr\'esentons $\bar{\lambda}$ comme une partition \`a $2r+1$ termes. On a prouv\'e en 6.4 que $k(\bar{\lambda},\bar{\epsilon})=k$. On peut donc d\'efinir les partitions $\tilde{A}_{\bar{\lambda},\bar{\epsilon}}$ et $\tilde{B}_{\bar{\lambda},\bar{\epsilon}}$ relatives \`a $r$, $k$ et au signe $u=\epsilon(1)$. On a $\tilde{A}_{\bar{\lambda},\bar{\epsilon}}=\{\bar{\lambda}_{j}/2+2r+1-j; j=1,...,2r+1,\bar{\epsilon}(j)(-1)^{j+1}=\epsilon(1)\}$ et, par construction, cette partition a $\tilde{n}$ termes. Donc $\tilde{A}_{\bar{\lambda},\bar{\epsilon}}$ est l'ensemble des $\tilde{n}$ plus grands termes de
$\{\bar{\lambda}_{j}/2+2r+1-j; j=1,...,2r'+2,\bar{\epsilon}(j)(-1)^{j+1}=\epsilon(1)\}$, autrement dit
$$\tilde{A}_{\bar{\lambda},\bar{\epsilon}}=\tilde{A}_{\lambda^{max},\epsilon^{max}}.$$
De la m\^eme fa\c{c}on,
$$\tilde{B}_{\bar{\lambda},\bar{\epsilon}}=\tilde{B}_{\lambda^{max},\epsilon^{max}}.$$
Mais ces \'egalit\'es \'equivalent \`a l'\'egalit\'e  $(\lambda^{max},\epsilon^{max})=(\bar{\lambda},\bar{\epsilon})$. Cela ach\`eve la d\'emonstration de la proposition 6.2. $\square$

CNRS-Institut de Math\'ematiques de Jussieu-PRG

4 place Jussieu

750005 Paris

jean-loup.waldspurger@imj-prg.fr

\end{document}